%% file: JMB_arxiv-11-jun-2021.tex
\documentclass[smallextended]{svjour3}
\usepackage{preamble_arxiv}

\begin{document}

\titlerunning{Operon Dynamics with State Dependent Delays}
\authorrunning{T.~Gedeon, A.R.~Humphries, M.C.~Mackey, H.-O.~Walther and Z.~Wang}

\title{Operon Dynamics with State Dependent Transcription and/or Translation Delays\\
\thanks{This
        work was supported by the Natural Sciences and Engineering Research Council (NSERC) of Canada (ARH and MCM) and the Alexander von Humboldt Stiftung of Germany (H-OW).
        TG was partially supported by NSF grant DMS-1951510. }}
\author{Tom\'{a}\v{s} Gedeon
\and Antony R. Humphries
\and Michael C. Mackey
\and Hans-Otto Walther
\and Zhao (Wendy) Wang}

\institute{Tom\'{a}\v{s} Gedeon
\at Department of Mathematics, Montana State University, Bozeman, MT 59717\\ \email{gedeon@math.montana.edu}
\and
Antony R. Humphries
\at
Department of Mathematics and Statistics and Department of Physiology, McGill University, Montreal, QC, Canada H3A 0B9 \\
\email{tony.humphries@mcgill.ca}
\and
Michael C. Mackey
\at
Departments of Physiology, Physics and Mathematics \& Statistics, McGill University, 3655 Promenade Sir William Osler, Montreal, Quebec H3G 1Y6, Canada \\
\email{michael.mackey@mcgill.ca}
\and
Hans-Otto Walther
\at
Mathematisches Institut, Universit{\"{a}}t Giessen, Arndtstrasse 2, 35392 Giessen, Germany\\
\email{Hans-Otto.Walther@math.uni-giessen.de}
\and
Zhao (Wendy) Wang
\at
Department of Mathematics and Statistics, McGill University, Montreal, QC, Canada H3A 0B9 \\
\email{zhao.wang3@mail.mcgill.ca}
}

\date{\today,\currenttime}
\maketitle

\input{abstract.tex}

\tableofcontents

\input{intro.tex}

\input{operonequations}

\input{equilibria.tex}

\input{quasilinearization.tex}

\input{numerics.tex}

\input{dynamics.tex}

\input{discussion.tex}

\section*{Acknowledgments}

MCM would like to thank the Universities of Giessen and Bremen for their hospitality during the time that this work was initiated.

\appendix
\setcounter{tocdepth}{2}

\section*{Appendices}
\addcontentsline{toc}{section}{Appendices}
\addtocontents{toc}{\protect\setcounter{tocdepth}{0}}

\input{app_semiflow.tex}

\input{app_attractor.tex}

\input{app_linearization.tex}


\end{document}

%% file: abstract.tex
\begin{abstract}
Transcription and  translation retrieve and operationalize gene encoded information in cells. These processes are not instantaneous and incur significant delays.
In this paper we study  Goodwin models of both inducible and repressible operons with state-dependent delays. The paper provides justification and derivation of
the model, detailed analysis of the appropriate setting of the corresponding dynamical system, and extensive numerical analysis of its dynamics. Comparison with constant delay models shows
significant differences in dynamics that include existence of stable periodic orbits in inducible systems and multistability in repressible systems. A combination of parameter space exploration, numerics, analysis of  steady state linearization and bifurcation theory indicates the likely presence of Shilnikov-type homoclinic bifurcations in the repressible operon model.
 \end{abstract}

%

%% file: intro.tex
\section{Introduction}

\begin{figure}[thp]
\centering
\includegraphics[width=1.0\textwidth]{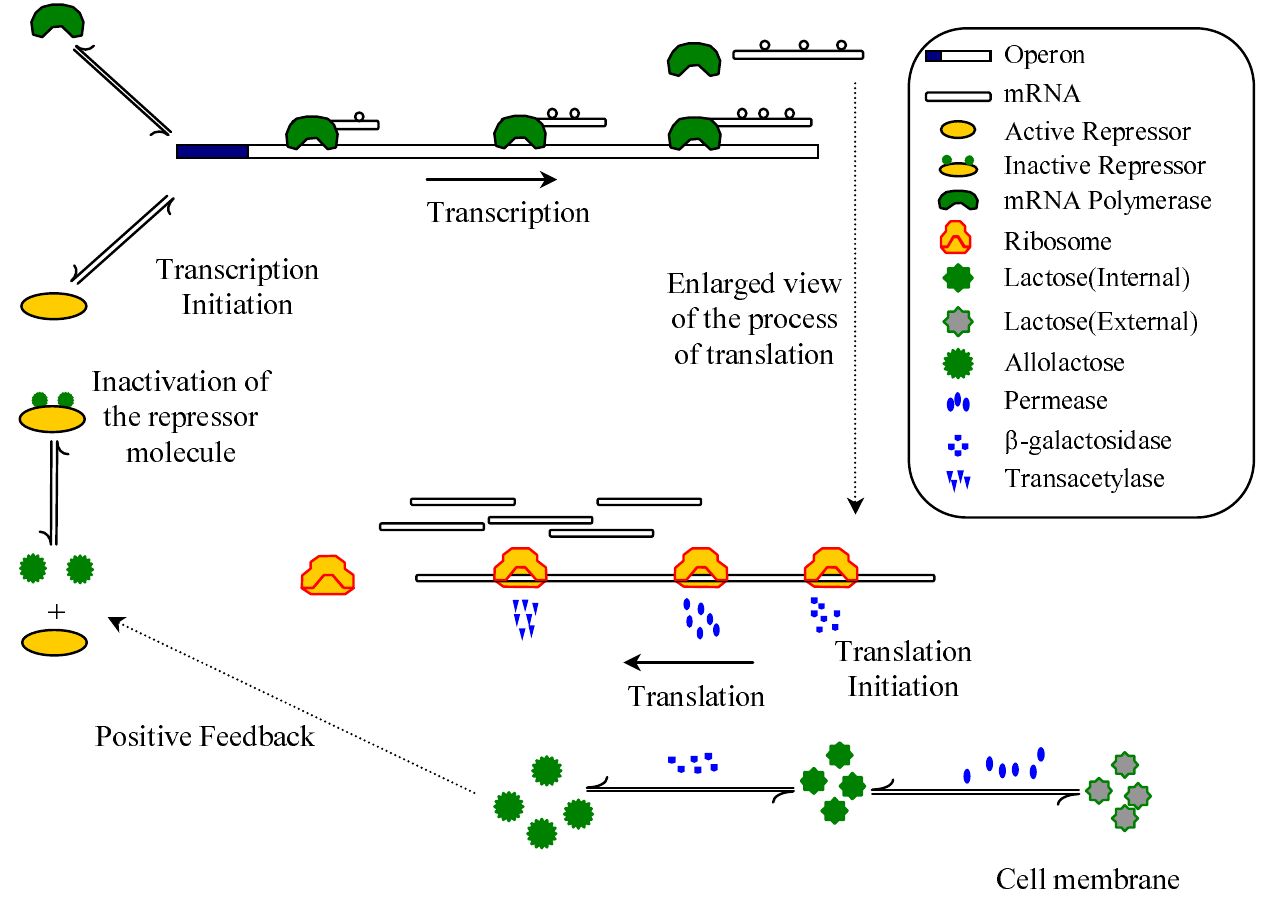}
\caption{A cartoon representation of the operation of the {\it lac} operon enabling bacteria to utilize lactose as an energy source in the absence of glucose, as elucidated by \citet{Jacob1960operon, jacob-1961}.  The process starts (upper left) when the operator region (dark-blue) is free of active repressor molecules so mRNA polymerase can attach to the DNA and start moving along the structural genes to produce mRNA.  Once the mRNA is fully formed, ribosomes start the translation process which, for the {\it lac} operon, produces, in sequence, $\beta$-galactosidase, permease, and transacetylase.  The permease facilitates the transport of extracellular lactose to the intracellular space (bottom right), while the permease is essential for the conversion of the internalized lactose into allolactose.  Allolactose, in turn, is able to bind to active repressor molecules thereby inactivating them and giving rise to the positive feedback nature of the {\it lac} operon.  Modified from \cite{yildirim2003feedback}.}
\label{fig:lac-operon}
\end{figure}

%
The regulation of information encoding and transmission in biological systems has intrigued and occupied mathematicians and physicists for decades.  One of the earliest published papers along these lines is \citet{timofeeff1935natur} which played a major role in the motivation of Erwin
Schr\"{o}dinger to give the 1943 Dublin lectures that are immortalized in \citet{schrodinger1943}.  The regulation of information retrieval started to become understood very quickly after the seminal work of Jacob and Monod \citep{Jacob1960operon,jacob-1961} elucidating the nature of the regulation of lactose production in bacteria.  The molecular apparatus carrying out this procedure in bacteria, involving transcription of DNA to produce mRNA and the translation of the mRNA to ultimately produce an effector protein, was named an 'operon' by them.
In Figure \ref{fig:lac-operon} we have illustrated the operon concept using the lactose ({\it lac}) operon as an example \citep{Jacob1960operon,jacob-1961}.

Rather astonishingly, mathematical models of the process of transcription and translation rapidly appeared~\citep{goodwin1963temporal,goodwin1965}. These first attempts were swiftly followed by an analysis of a simple repressible operon \citep{Griffith68a} and an inducible operon \citep{Griffith68b}. These and other results were summarized in the \citet{tyson-othmer-1978} review which is still relevant today.

 Though Goodwin clearly noted the existence of significant delays in both transcription and translation in \cite{goodwin1963temporal}, and thought that the delays might have significant dynamic influences\footnote{Personal conversation with MCM, November, 1994}, he did not examine their potential effects.  Apparently the first to incorporate constant transcriptional and translational delays into the Goodwin model was \cite{banks77} and then \cite{macdonald1977time} followed in rapid succession by \cite{banks1978global,banks1978stability}, \citet{adh79,adh83} and \cite{mahaffy1984models}.  These were followed by a number of subsequent investigations.

Since the processes of transcription and translation are rather complicated, the assumption of constant delay may limit our ability to appreciate the richness of dynamics that the process of protein production can impose on the cell. The goal of this paper is to derive a Goodwin-like   delay-differential equation (DDE) model  with {\it state dependent delays} that we feel may more closely correspond to biological reality, explore the potential dynamics both in repressible and inducible cases and contrast these dynamics with that of a system with constant delays.

This paper is rather long and detailed, and a summary of the contents may be of help to the reader.  Section~\ref{sec:basic} outlines the basic operon equations starting with a summary of the Goodwin model in Section~\ref{ssec:Goodwin} and a summary of the equations we derive in this paper in
Section~\ref{ssec:delay}.  Section~\ref{ssec:evolution} details the full derivation of the model equations we study here while Section~\ref{ssec:trans-init} summarizes the functional forms for the transcription initiation rates that we use for inducible and repressible operons.

Section~\ref{ssec:trans-trans-vel} contains our arguments for the nature and form of the transcriptional and translational velocities in the two types of operons that lead naturally to the state dependent delays that are central to our study.  Section~\ref{ssec:equilibria} deals with the quantitative and qualitative nature of the possible equilibrium states of our model equations developed in Section~\ref{ssec:evolution}, and the following Section~\ref{ssec:quasilinear} gives a linearization procedure in the neighborhood of these steady states that is not fully justified mathematically but easily understood by most readers.  (The analytically exact linearization is to be found in Appendix~\ref{app:linearization} and leads to precisely the same result).  These linearizations are needed for stability determinations based on the eigenvalues evaluated at the equilibria.

Section~\ref{sec:method} contains the details of the numerical methods we have used in our numerical studies of this  model, while the following Section~\ref{sec:examples} contains the extensive details of our numerical studies for both the repressible (Section \ref{sec:repexamples}) and inducible (Section \ref{sec:indexamples}) operon models. In both cases we have found that significant new types of dynamics are introduced by the state dependency of the delays.
 In the inducible operon model we found a stable periodic orbit as well as tristability between a periodic orbit and two steady states. In the repressible operon model we found bistability between two steady states as well as between a periodic orbit and a steady state. All of these results are obtained with one state dependent delay and are not present in the corresponding DDE model with a constant delay.  In addition, in the repressible operon model we found evidence of a homoclinic bifurcation of Shilnikov type~\citep{Shilnikov65,Kuznetsov2004}, indicating the potential for complex dynamics. Finally, in both types of operons there are stable  periodic orbits, where a short burst of  transcription and translation is interspersed with longer periods of quiescence. These orbits represent a  pulse-generating mechanism on a sub-cellular level and may be connected to the phenomena of transcriptional bursting~\citep{Chubb2020}.    The main body of the paper concludes with a discussion and summary in Section~\ref{sec:disc} and is followed by three mathematical appendices.  Appendix~\ref{app:semiflow} treats the semiflows arising from our basic model with state dependent delays, Appendix \ref{app:dissipativity} considers some aspects of the nature of the model solutions including positivity and the global attractor, while Appendix \ref{app:linearization} treats the linearization mentioned above.

%% file: operonequations.tex
\section{Basic operon equations} \label{sec:basic}

\subsection{The Goodwin model}\label{ssec:Goodwin}

The \citet{goodwin1965}  model for operon dynamics    considers a large population of cells, each of which contains
one copy of a particular operon,  and we use that as a basis for discussion.  We let $(M,I,E)$ respectively
denote the mRNA, intermediate protein, and effector protein {\it concentrations}. For a generic operon  
the dynamics are assumed to be  given by \citep{goodwin1963temporal,goodwin1965,Griffith68a,Griffith68b,othmer76,selgrade79}
\begin{align}
    \dfrac{dM\!}{dt}(t) &= {\cal F}(E(t))
    -\gamma_M M(t),\label{eq:mrna}\\
    \dfrac{dI}{dt}(t) &= \beta_I M(t) -\gamma_I I(t) ,\label{eq:intermed}\\
    \dfrac{dE}{dt}(t) &= \beta_E I(t) - \gamma_E E(t).\label{eq:effector}
\end{align}
It is assumed here that the flux ${\cal F}$ (in units of $[\mbox{concentration}\cdot \mbox{time}^{-1}]$) of initiation of mRNA production is a function of the effector level $E$.  Furthermore, the model assumes that the flux of protein and metabolite production  are proportional (at rates $\beta_I,\beta_E$ respectively) to the amount of mRNA and
intermediate protein respectively. All three of the components $(M,I,E)$ are subject to degradation at rates $\gamma_M, \gamma_I, \gamma_E$.  The parmeters $\beta_I,\beta_M,\gamma_I$ and
 $\gamma_M$ have dimensions [time$^{-1}$].

\subsection{The effects of cell growth and state dependent transcription and translation delays}\label{ssec:delay}

We will study an extended Goodwin model taking into account the effects of cell growth and delays which are introduced by state dependent transcription and translation processes.
 The cell growth affects  the volume and hence  the concentrations of all the molecules in the cell.

The following sections are devoted to the derivation of the generalization of the Goodwin model:
\begin{align}
\dfrac{dM\!}{dt}(t) & = \label{eq:mrna-delay}
\beta_M e^{-\mu \tau_M}\dfrac{v_M(E(t))}{v_M(E(t-\tau_M))}f(E(t-\tau_M))-\bgamma_M M(t), \\
\dfrac{dI}{dt}(t) & = \label{eq:intermed-delay}
    \beta_I  e^{-\mu\tau_I}\dfrac{v_I(M(t))}{v_I(M(t-\tau_I))} M(t-\tau_I) -\bgamma_I I(t),\\
\dfrac{dE}{dt}(t) & = \beta_E I(t) -\bgamma_E E(t).\label{eq:effector-delay}
\end{align}

In equations \eqref{eq:mrna-delay}-\eqref{eq:effector-delay} there are several changes to be noted relative to the original Goodwin model \eqref{eq:mrna}-\eqref{eq:effector}. The first is
the introduction of the two delay terms $E(t-\tau_M)$ and $M(t-\tau_I)$ indicating that   $E$ and $M$ are now to be evaluated at a time in the past due to the non-zero times required for transcription and translation.  From a dynamic point of
view, the presence of these delays can have a dramatic effect.

The
second change\footnote{We  note that previous inclusions of cell growth and transcription/translation delays \citep{yildirim2003feedback,yildirim2004dynamics} contained an error which we correct here. In the original \cite{yildirim2003feedback} paper, the term $e^{-\mu \tau_M}$ in \eqref{eq:mrna-delay} was mistakenly inserted within the argument of $f$ instead of being in front of $f$.} is the appearance of the terms $e^{-\mu \tau_M}$ and $e^{-\mu \tau_I}$
which respectively account for an effective dilution of the maximal mRNA production and intermediate protein
fluxes because the cell is growing at a rate $\mu$
(in units of [time$^{-1}$]).

The third change from \eqref{eq:mrna}-\eqref{eq:effector} to
\eqref{eq:mrna-delay}-\eqref{eq:effector-delay} is
the alteration of the decay rates  $\gamma_i$ to $\bgamma _i \equiv \gamma_i + \mu$ because the dilution
due to cell growth
leads to an effective increase in the rate of destruction.

A fourth change is the replacement
of ${\cal F}$ in \eqref{eq:mrna} by $\beta_M f$ in \eqref{eq:mrna-delay}.
Here $\beta_M$ is the {\it maximal} production rate
(in units of $[\mbox{concentration}\cdot \mbox{time}^{-1}]$)
possible and $f$ is the fraction of free operator sites on the operon, a function that will vary between a maximal value of $1$ and a minimal value in $[0,1)$. We remark that $\beta_M$ thus has different units
than the linear rate constants $\beta_I$ and $\beta_E$.  

Finally, velocity ratio terms of the form $\frac{v_j(w)}{v_j(w(t-\tau))}$ appear in
\eqref{eq:mrna-delay} and \eqref{eq:intermed-delay} as a consequence of the delays $\tau_M$ and $\tau_I$ being non-constant and depending on the state variables,
as explained in Section~\ref{ssec:trans-trans-vel}.

\subsection{Evolution equations incorporating state-dependent transcription and translation rates}\label{ssec:evolution}

Transcription is initiated when RNA polymerase (RNAP) is recruited to the promoter region by one or more transcription factors, partially unwinds the promoter DNA to form the transcription bubble, and subsequently leaves the promoter region, moving along the DNA. If multiple initiations take place in rapid succession, then transcribing RNAPs  start to interfere with each other, and as a result the average velocity of individual RNAPs transcription events will decrease. This, in turn, leads to an increased time of transcription.

Translation is initiated by assembly of the ribosome on the initiation region of the mRNA.
Ribosomes catalyze subsequent  binding of codon specific transfer RNAs (tRNA) to the mRNA and  transfer of the attached amino acid to the nascent polypeptide. Subsequent translocation of  the ribosome completes the cycle.

The result is a bio-polymerization process whose velocity depends on current demand for both ribosomes and  tRNAs, which is affected both  by the number of actively translated mRNAs and the  growth rate of the cell.
Both transcription and translation share key characteristics that lead to a common model of these processes.
The most basic model, from which other models are derived, is a stochastic Totally Asymmetric Simple Exclusion Process (TASEP) model for particles hopping on a strand with a finite number of discrete sites, that represent nucleotides ~\citep{derrida93,schutz93,kolomeisky98,Shaw03,Zia11}.

It should be noted  that in eukaryotes (as opposed to our consideration in this paper of prokaryotic gene regulation) transcription takes place in the nucleus  while translation takes place in the cytoplasm and the consequent transport of intermediate from cytoplasm into the nucleus gives rise to a {\it transport delay} that may, on occasion, be considered as state dependent \citep{verriest17,wang21}.

\subsubsection*{mRNA dynamics}

For the mRNA molecules we start with the mRNA transcripts and consider their density  $r(t,a)$ at time $t$ and location $a$ along the DNA,
so
$$
\int_{a_1}^{a_2}r(t,a)da
$$
is the number of mRNA molecules with positions between $a_1$  and $a_2$, $0\le a_1<a_2\le a_M$, where   $a_M$ is the end of the transcription region.
The velocity of transcription along the DNA is given by a function $v_M$, and we assume that the actual velocity of the process depends on the value $w(t)$ of a function $w$, to be determined later.
If the transcription process takes place without any loss of mRNA transcripts, then the evolution equation for the density $r(t,a)$ is given by
\begin{equation}
\dfrac{\partial r}{\partial t}(t,a) + v_M(w(t)) \dfrac{\partial r}{\partial a}(t,a) = 0.
\label{e:m(t,a)}
\end{equation}
We look for a differential equation
$$
\frac{dm}{dt}(t)=p(t)-\gamma_Mm(t)
$$
for the number $m(t)$ of complete mRNA molecules at time $t$, with a constant rate $\gamma_M>0$ of degradation
and a \emph{production function} $p$ which describes the contribution of the release of completed mRNA molecules at time $t$ to the rate of change $\frac{dm}{dt}(t)$.
In order to determine $p(t)$ consider the number of mRNA molecules undergoing transcription at time $t$, which is
$$
J(t)=\int_0^{a_M}r(t,a)da.
$$
Using Eq. \eqref{e:m(t,a)}  we have a balance equation for $J$,
$$
\frac{dJ}{dt}(t)=v_M(w(t))r(t,0)-v_M(w(t))r(t,a_M),
$$
 where the term $v_M(w(t)r(t,0)$ represents  the initiation rate of transcription of mRNA molecules to $\frac{dJ}{dt}(t)$, and $-v_M(w(t))r(t,a_M)$ is the
release rate of completed mRNA molecules. Therefore the term $v_M(w(t))r(t,a_M)$ is the desired contribution $p(t)$ to $\frac{dm}{dt}(t)$.
Using characteristics we obtain
$$
p(t)=v_M(w(t))r(t,a_M)=v_M(w(t))r(t-\tau,0),
$$
with the time $\tau=\tau_M(t)$ needed for production of mRNA molecules which reach the final length $a_M$ at time $t$,
\be
a_M  =  \int^t_{t-\tau_M(t)}v(w(s))ds
 =  \int^0_{-\tau_M(t)}v(w(t+s))ds.\label{eq:delay-tau-M}
\ee
We arrive at
\begin{align*}
	p(t) & =	v_M(w(t))r(t-\tau_M(t),0)\\
	& = \frac{v_M(w(t))}{v_M(w(t-\tau_M(t)))}[v_M(w(t-\tau_M(t)))r(t-\tau_M(t),0)],	
\end{align*}
where the term $[v_M(w(t-\tau_M(t)))r(t-\tau_M(t),0)]=F(t-\tau_M(t))$ stands for the onset of transcription of mRNA molecules at time $t-\tau_M(t)$. The differential equation for $m$ thus becomes
\begin{align}
\frac{dm}{dt}(t) & = \frac{v_M(w(t))}{v_M(w(t-\tau_M(t)))}[v_M(w(t-\tau_M(t)))r(t-\tau_M(t),0)]-\gamma_Mm(t)\nonumber\\
& =  \frac{v_M(w(t))}{v_M(w(t-\tau_M(t)))}F(t-\tau_M(t))-\gamma_Mm(t).\label{eq:m}	
\end{align}

We now switch to a description of transcription in terms of molecule concentration, rather than numbers of molecules. Since the concentration $M$ is related to the number of molecules $m$ by $M=m/V$, we have
$$
\dfrac{dm}{dt}(t) = \dfrac{dM\!}{dt}(t) V(t) + M(t)\dfrac{dV\!}{dt}(t)= V(t)\dfrac{dM\!}{dt}(t) + \mu V(t)M(t),
$$
under the assumption that the cells are growing exponentially with $\tfrac{dV}{dt}(t) = \mu V(t)$.
Consequently, noting that $V(t)=e^{\mu\tau_M(t)}V(t - \tau_M(t))$, we can rewrite \eqref{eq:m} as
\begin{align}
\dfrac{dM\!}{dt}(t) \hspace*{-2em}\mbox{} & \mbox{}\hspace*{2em} = \frac{1}{V(t)}\dfrac{dm(t)}{dt}-\mu M(t)
\nonumber \\
&=\dfrac{v_M(w(t))}{v_M(w(t-\tau_M(t)))} e^{-\mu\tau_M(t)}
\dfrac {   F (t-\tau_M(t))}{V(t - \tau_M(t))}
- (\gamma_M + \mu)  M(t) \nonumber
\end{align}

We express  the {\it  initiation flux}
$\frac {F(\ldots)}{V(\ldots)}$ in concentration units as
$$
\dfrac {F(t-\tau_M(t))}{V(t - \tau_M(t))}=:\beta_M f(w(t-\tau_M(t)))
$$
where $\beta_M$ is the {\it maximal} initiation flux (units of [$\mbox{concentration}\cdot\mbox{time}^{-1}$]) and $f$ stands for the fraction of free operator sites on the operon, a function that will vary between  a minimal value in $(0,1)$ and a maximal value of $1$.

As derived in \citet[Chapter 1]{mackey2016simple}, the initiation flux is a function of concentration of
the  effector molecule $E$, and the velocity $v_M$ also depends on $E$ (Section~\ref{ssec:trans-trans-vel}). Therefore for the transcription process $w=E$ and we obtain
\be
\dfrac{dM\!}{dt}(t) =\dfrac{v_{M\!}(E(t))}{v_{M\!}(E(t-\tau_{M\!}(t)))} e^{-\mu \tau_{M\!}(t)}
\beta_M f^{\!}(E(t-\tau_{M\!}(t)))
- (\gamma_M + \mu) M(t),
\label{eq:m-concen}
\ee
together with Eq. \eqref{eq:delay-tau-M} for the delay $\tau_M(t)$, which depends on the function $E$.


\subsubsection*{Intermediate dynamics}

We assume that the initiation of the translational production of the mRNA into intermediate protein is a relatively simple process and, unlike the transcription process, not under regulatory control.

For the intermediate molecules,  we use $i(t,b)$  to  describe   their density  at time $t$ and location $b$ along the mRNA. We assume that the translation is proceeding at a velocity $v_I(q)$ along the mRNA. This velocity may depend on $q$, where $q$ is to be determined. Analogous to the transcription process we arrive at
\be
\dfrac{di}{dt}(t) =  \dfrac{v_I(q(t))}{v_I(q(t-\tau_I(t)))} \beta_I m(t-\tau_I(t)) - \gamma_I i(t).
\label{eq:i}
\ee

In Section~\ref{ssec:trans-trans-vel}  we argue that $q=M$, the concentration of mRNA.  Therefore switching to a concentration description  using $I(t)= \frac{i(t)}{V(t)}$,
following the derivation of \eqref{eq:m-concen} we can rewrite \eqref{eq:i} in the form
\be
\dfrac{dI}{dt}(t)
=\dfrac{v_I(M(t))}{v_I(M(t-\tau_I(t)))} e^{-\mu \tau_I(t)} \beta_I M(t-\tau_I(t))) - (\gamma_I + \mu) I(t).
\label{eq:I-concen}
\ee

\subsubsection*{Effector dynamics }

The effector dynamics are the easiest because there is neither transcription nor translation involved.  Rather the production of the effector is assumed to be proportional to the intermediate level $i$ at a rate $\beta_E$, while the effector is destroyed at a rate $\gamma_E$.  Thus
\be
\dfrac{de}{dt}(t) = \beta_E i(t) - \gamma_E e(t),
\label{eq:e}
\ee
and, changing the description from numbers to concentrations, we have simply that
 \be
\dfrac{dE}{dt}(t) = \beta_E I(t) - (\gamma_E + \mu) E(t).
\label{eq:E}
\ee

\subsubsection*{Putting it all together}

Denote the transcriptional velocity  by $v_M(E(t))$ and the translational velocity by $v_I(M(t))$.  Further let   $\bgamma_M = \gamma_M + \mu,\bgamma_I = \gamma_I + \mu,\bgamma_E = \gamma_E + \mu $.  Then  we can write the state dependent forms of
\eqref{eq:mrna-delay}-\eqref{eq:effector-delay} as
\begin{align}
\dfrac{dM\!}{dt}(t) & = \beta_M e^{-\mu \tau_{M}(t)} \dfrac{v_M(E(t))}{v_M(E(t-\tau_{M}(t)))} f(E(t-\tau_M(t))) -\bgamma_M M(t), \label{eq:mrna-delay-var}\\
\dfrac{dI}{dt}(t) & = \beta_I e^{-\mu\tau_I(t)} \dfrac{v_I({M(t)})}{v_I(M(t-\tau_I(t)))} M(t-\tau_I(t)) -\bgamma_I I(t) , \label{eq:intermed-delay-var}\\
\dfrac{dE}{dt}(t) & = \beta_E I(t) -\bgamma_E E(t). \label{eq:effector-delay-var}
\end{align}
These equations are supplemented by the two additional equations  which define the delays  $\tau_M$ and $\tau_I$   by threshold conditions, namely
\begin{align}
a_M &= \int_{t-\tau_M(t)}^{t} v_M(E(s)) ds=\int_{-\tau_M(t)}^{0} v_M(E(t+s)) ds
\label{eq:delay by stateM}\\
a_I &= \int_{t-\tau_I(t)}^{t} v_I(M(s)) ds=\int_{-\tau_I(t)}^{0} v_I(M(t+s)) ds.
\label{eq:delay by stateI}
\end{align}

We write $\tau_M(t)$ and $\tau_I(t)$ for the state-dependent delays, but from
\eqref{eq:delay by stateM} and \eqref{eq:delay by stateI} it is clear that the value of each is determined
by the values of $E(t)$ or $M(t)$ respectively over the whole integration interval. Using the Banach space notation of the Appendices we ought to write $\tau_M(E_t)$ and $\tau_I(M_t)$ for these delays where $E_t$ and $M_t$ are functions defined by $E_t(\theta)=E(t-\theta)$ and $M_t(\theta)=M(t-\theta)$. But, to hopefully make the presentation accessible to readers who are not comfortable with Banach spaces, we will
avoid any Banach space notation in the main body of the text and continue to write $\tau_M(t)$ and $\tau_I(t)$ for the delays at time $t$.

Velocity ratio terms such as those appearing in \eqref{eq:mrna-delay-var} and \eqref{eq:intermed-delay-var}
are ubiquitous in distributed state-dependent DDE problems with either threshold conditions \citep{Craig2016}
or with randomly distributed maturation times \citep{cassidy2019equivalences}. \cite{Bernard2016} explains very clearly why they arise.

\subsection{The control of transcription initiation rates}\label{ssec:trans-init}

The determination of how
effector concentrations modify the fraction of free operator sites, $f$, has been dealt with by a number of authors. Here we merely summarize the nature of $f$ for inducible and repressible systems, see \citet[Chapter 1]{mackey2016simple} for  details.

For a repressible operon, $f$ is a monotone decreasing function
\begin{equation}
f(E) 
= \frac{1+
K_1E^n}{1 + KE^n}, \label{rep-frac2}
\end{equation}
where  $K > K_1$, $n > 1$, so there is maximal repression for large $E$.
For an inducible system
$f$ is a monotone increasing function of the form
\begin{equation}
f(E) 
= \frac{1+
K_1E^n}{K + K_1E^n}, \label{ind-frac2}
\end{equation}
where $K > 1$, $n > 1$.  Maximal induction occurs for very large $E$.

Note that both \eqref{rep-frac2} and  \eqref{ind-frac2} are special cases of
    \begin{equation}
    f(E) =   \dfrac{1+K_1E^n}{A+BE^n} 
      \label{eq:gen-response-fun}
    \end{equation}
The constants $A,B  \geq 0$ are defined in Table \ref{table:delay/velocity behaviour}.


\subsection{Transcriptional and translational velocities}\label{ssec:trans-trans-vel}

In this section we discuss
cellular processes that affect the transcriptional and translational velocities $v_M $ and $v_I $.
Both transcription and translation are polymerization processes where small parts are associated by an enzymatic reaction catalyzed by a large complex into a long polymer chain.

For the transcription process nucleotides A,C,G,T are incorporated by RNA polymerase (RNAP) into an mRNA chain, and for the translation process peptides are incorporated by ribososomes into a polypeptide that, upon folding, becomes a functional protein.
Velocities of both processes depend on a sufficient and timely supply of nucleotides and  peptides, respectively.
The availability, or paucity, of the parts  may result in changes in velocity from position to position along the strand~\citep{Zia11}.
The abundance of the parts reflects the overall growth rate of the cell:
faster growth leads to greater demand on resources and a slower transcription ($v_M$)  and translation ($v_I $) velocity. Therefore  for an inducible operon, the transcription velocity $v_M(E) $ is a decreasing function of the concentration of the effector  $E$ and for a repressible operon  $v_M(E)$ is an increasing function of the concentration of the effector $E$.

The velocity of translation depends on the number of initiations of the translation process, which is directly proportional to the concentration $M$ of mRNA. Since greater demand on peptide availability results in a lower elongation velocity of ribosomes, the translational velocity of ribosomes $v_I(M)$ is a decreasing function of $M$.

There is a second effect  that may affect elongation velocity. This is the effect of  elongation interference by multiple RNAP or multiple ribosomes~\citep{hwa08,klumpp11}.
 The velocity of elongation decreases with the number of   RNAPs and ribosomes that elongate at the same time.  Since this number is proportional to the initiation rate,  the velocity $v_M=v_M(E)$
is a decreasing function of the concentration of $E$ for an inducible operon
and an increasing function of the concentration of $E$
for a repressible operon. The velocity  $v_I(M)$ is  a  decreasing function of $M$.

Both availability of nucleotides and peptides and elongation interference support the following assumptions on the velocities $v_M$ and $v_I$:
\begin{align*}
v_M & =v_M(E) \quad\textrm{is a decreasing function of $E$ for an inducible operon} \\
v_M & =v_M(E) \quad\textrm{is an increasing function of $E$ for a repressible operon} \\
v_I &= v_I(M) \quad\textrm{is always a decreasing function of $M$.}
\end{align*}

There are no analytic expressions for the dependence of $v_M(E)$ and $v_I(M)$ on $E$ and $M$ respectively, but we have made assumptions for modeling purposes and these are detailed in Table \ref{table:delay/velocity behaviour}.  Specifically we have assumed that they can be represented by Hill functions with parameters determining the maximum, minimum and half-maximal values as well as a parameter ($m$ or $m_I$) which controls the slope.  We do not offer any detailed stoichiometric justification for these assumptions, but rather assume that they will capture the essential nature of their dependencies.

\begin{table}
\centering
\begin{tabular}{|c|c|c|} \hline
Quantity & Repressible (e.g. {\it tryp}) & Inducible (e.g. {\it lac})\\
\hline
\parbox{2.2cm}{\centering \vspace*{3pt} Fraction $f(E)$ of free operators\vspace*{1pt}}
 & $ f(E) = \dfrac{1+K_1E^n}{1 + K E^n}$ &  $f(E) = \dfrac{1+K_1E^n}{K + K_1E^n}$
 \\
\hline
\parbox{2.1cm}{\centering Qualitative\\ behaviour of $f$} &
\parbox{3cm}{\vspace*{2pt} \centering Monotone $\downarrow$\\ Max $=1$\\ Min $=K_1/K < 1$ \vspace*{2pt}} &
\parbox{3cm}{\centering Monotone $\uparrow$\\Min $= 1/K< 1$\\ Max $= 1$ \vspace*{2pt}} \\
\hline
\parbox{2.1cm}{\centering transcription\\ velocity\\  $v_M(E)$}
&
\parbox{4.2cm}{\centering $v_M(E)=\dfrac{v_M^{min}E_{50}^m + v_M^{max} E^m}{E_{50}^m+E^m}$\\
Min at $E=0$ is $v_M^{min}$ \\
$v_M(E_{50})=\frac12(v_M^{max}+v_M^{min})$\\
Max as $E \to \infty$ is $v_M^{max} > v_M^{min}$ \vspace*{2pt}}&
\parbox{4.2cm}{\centering $v_M(E)=\dfrac{v_M^{max}E_{50}^m + v_M^{min} E^m}{E_{50}^m+E^m}$\\
Max at $E=0$ is $v_M^{max}$\\
$v_M(E_{50})=\frac12(v_M^{max}+v_M^{min})$\\
Min as $E \to \infty$ is $v_M^{min} < v_M^{max} $}
\\ \hline
\parbox{2.1cm}{\vspace*{3pt} \centering Qualitative $\tau_M$\\ behaviour \vspace*{3pt}}
& Monotone $\downarrow$  with $E$
& Monotone $\uparrow$  with $E$
\\
\hline
\parbox{2.1cm}{\centering translation\\ velocity\\  $v_I(M)$} &
\multicolumn{2}{c|}{\parbox{7cm}{\vspace*{1pt} \centering $v_I(M)=\dfrac{v_I^{max}M_{50}^{m_I} + v_I^{min} M^{m_I}}{M_{50}^{m_I}+M^{m_I}}$\\
Max at $M=0$ is $v_I^{max}$\\
$v_I(M_{50})=\frac12(v_I^{max}+v_I^{min})$\\
Min as $M \to \infty$ is $v_I^{min} < v_I^{max} $ \vspace*{2pt}}} \\ \hline
\parbox{2.1cm}{\vspace*{2pt} \centering Qualitative $\tau_I$\\ behaviour \vspace*{2pt}} &
\multicolumn{2}{c|}{Monotone $\uparrow$  with $M$} \\
\hline
\end{tabular}
\caption{Summary of the form of the fraction $f$ of free operators as determined by their stoichiometry, and the Hill function forms we have assumed for the transcriptional and translational velocities. }
\label{table:delay/velocity behaviour}
\end{table}

The parameters  $v_M^{min}, v_I^{min}$ describe minimal velocity of transcription and translation, respectively. While the individual polymerases and ribosomes may briefly  pause their elongation, in our model where $M,I$ model  concentrations in a large population of cells,  we assume $v_I^{min}>0$ and $v_M^{min}>0$.  Violation of this  assumption would cause significant problems both with our theory and numerical simulations. The minimal velocity being strictly positive ensures that the maximal delay is
bounded since from \eqref{eq:delay by stateM} and \eqref{eq:delay by stateI}
$$\tau_M(t)\leq \frac{a_M}{v_M^{min}}, \qquad  \tau_I(t)\leq \frac{a_I}{v_I^{min}}.$$
Similarly the maximal velocities define the minimal delays. Interesting dynamics can occur when delays become large and in what follows we will often take $v_M^{min}$ as a bifurcation parameter and study what happens as $v_M^{min}\to0$ and consequently $\tau_M(t)$ becomes large.

%% file: equilibria.tex
\subsection{Equilibria}\label{ssec:equilibria}

We next consider the steady states $(M^*,I^*,E^*)$ of \eqref{eq:mrna-delay-var}-\eqref{eq:delay by stateI}.
From \eqref{eq:delay by stateM} and \eqref{eq:delay by stateI}
at steady state the delays satisfy
\begin{equation} \label{eq:eqdels}
\tau_M=\tau_M^*(E^*):=\frac{a_M}{v_M(E^*)}, \qquad \tau_I=\tau_I^*(M^*):=\frac{a_I}{v_I(M^*)}.
\end{equation}
Then, at the steady state, equations \eqref{eq:mrna-delay-var}-\eqref{eq:effector-delay-var}
simplify to
\begin{align}
0 & = \beta_M e^{-\mu \tau_M^*(E^*)}f(E^*) - \bgamma_M M^*, \label{eq:eqM}\\
0 & = \beta_I e^{-\mu \tau_I^*(M^*)} M^* -\bgamma_I I^* , \label{eq:eqI}\\
0 & = \beta_E I^* -\bgamma_E  E^*. \label{eq:eqE}
\end{align}
We rearrange \eqref{eq:eqM} to obtain
\be \label{eq:M-ss}
M^* = \dfrac{\beta_M}{\bgamma_M}e^{-\mu\tau_M^*(E^*)}f(E^*),
\ee
and then substituting this and \eqref{eq:eqE} into \eqref{eq:eqI} we find that
the steady state must satisfy a single equation for $E^*$:
\be \label{eq:gEstar}
0  = g_E(E^*) :=
\frac{\beta_M\beta_I\beta_E}{\bgamma_M\bgamma_I\bgamma_E}
e^{-\mu(\tau_I^*(M^*)+\tau_M^*(E^*))}f(E^*)-E^*
\ee
where the argument of $\tau_I^*$ is given by \eqref{eq:M-ss}.

With the functions $v_M$, $v_I$ and $f$ defined as
in Table~\ref{table:delay/velocity behaviour} then
$v_M(E)\in[v_M^{min},v_M^{max}]$ and $v_I(M)\in[v_I^{min},v_I^{max}]$, so
\begin{equation*}
\tau_M\in[a_M/v_M^{max},a_M/v_M^{min}] \,\,\mbox{and} \,\,
\tau_I\in[a_I/v_I^{max},a_I/v_I^{min}]
\end{equation*}
while $f(E)\in(0,1]$.
Thus $g_E(0)>0$ and
$$g_E(E) \leq
\frac{\beta_M\beta_I\beta_E}{\bgamma_M\bgamma_I\bgamma_E}-E.$$
Therefore,
$g_E(E)<0$ for all $E$ sufficiently large, and by the intermediate value theorem there is at least one solution $E^*>0$ to $g_E(E^*)=0$.  This defines a steady state $(M^*,I^*,E^*)$. It also follows that any steady-state solution must satisfy
$$E^*\leq\frac{\beta_M\beta_I\beta_E}{\bgamma_M\bgamma_I\bgamma_E}, \qquad
I^*\leq\frac{\beta_M\beta_I}{\bgamma_M\bgamma_I}.$$

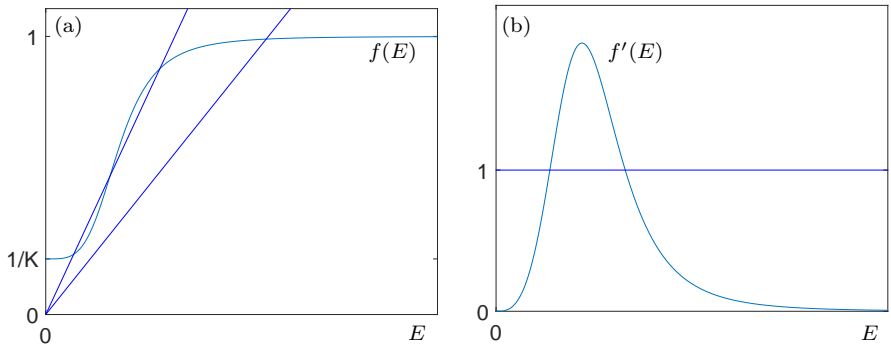
\begin{figure}[tp!]
\centering
\includegraphics[scale=0.45]{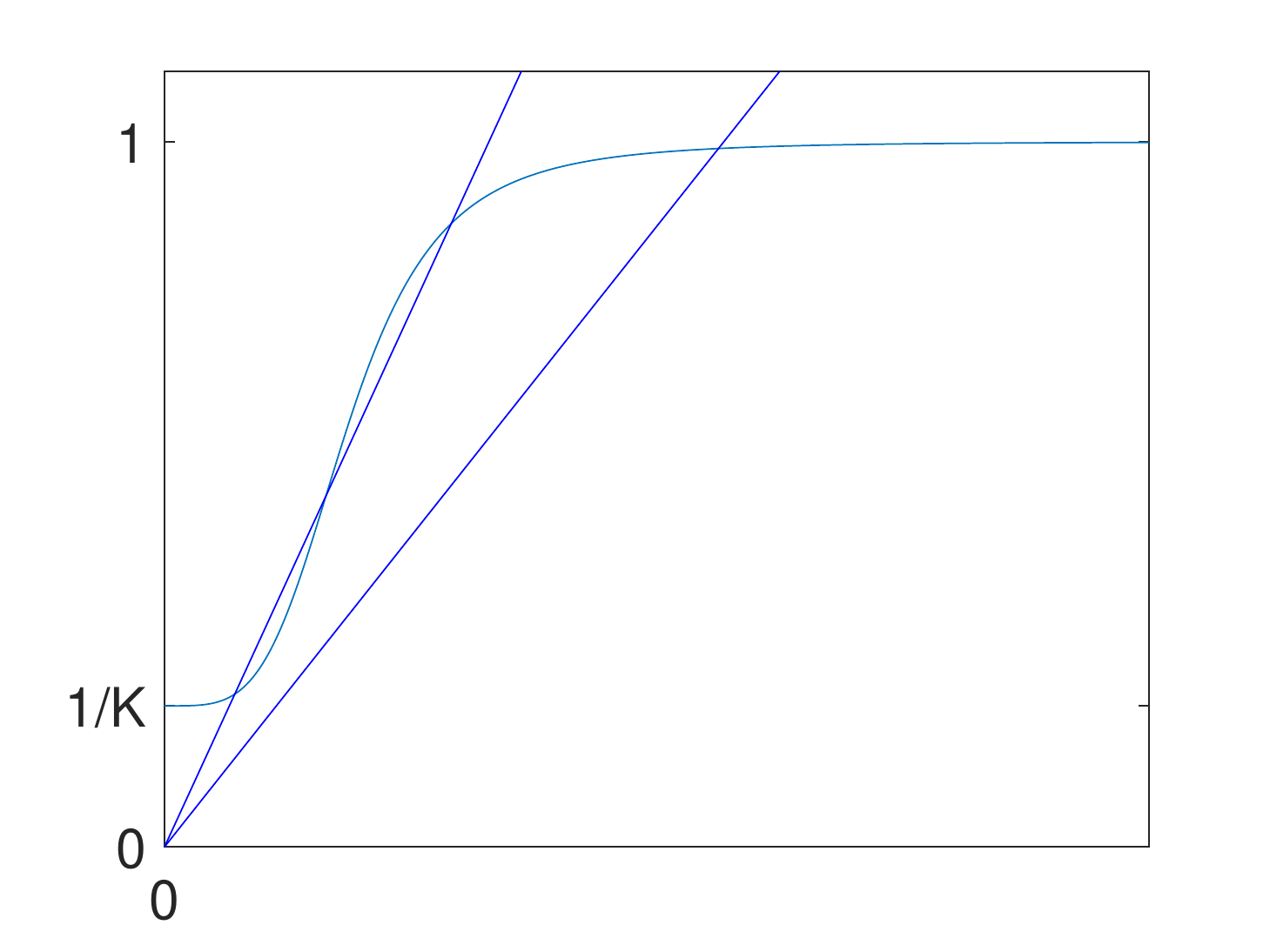}\hspace*{3mm}\includegraphics[scale=0.45]{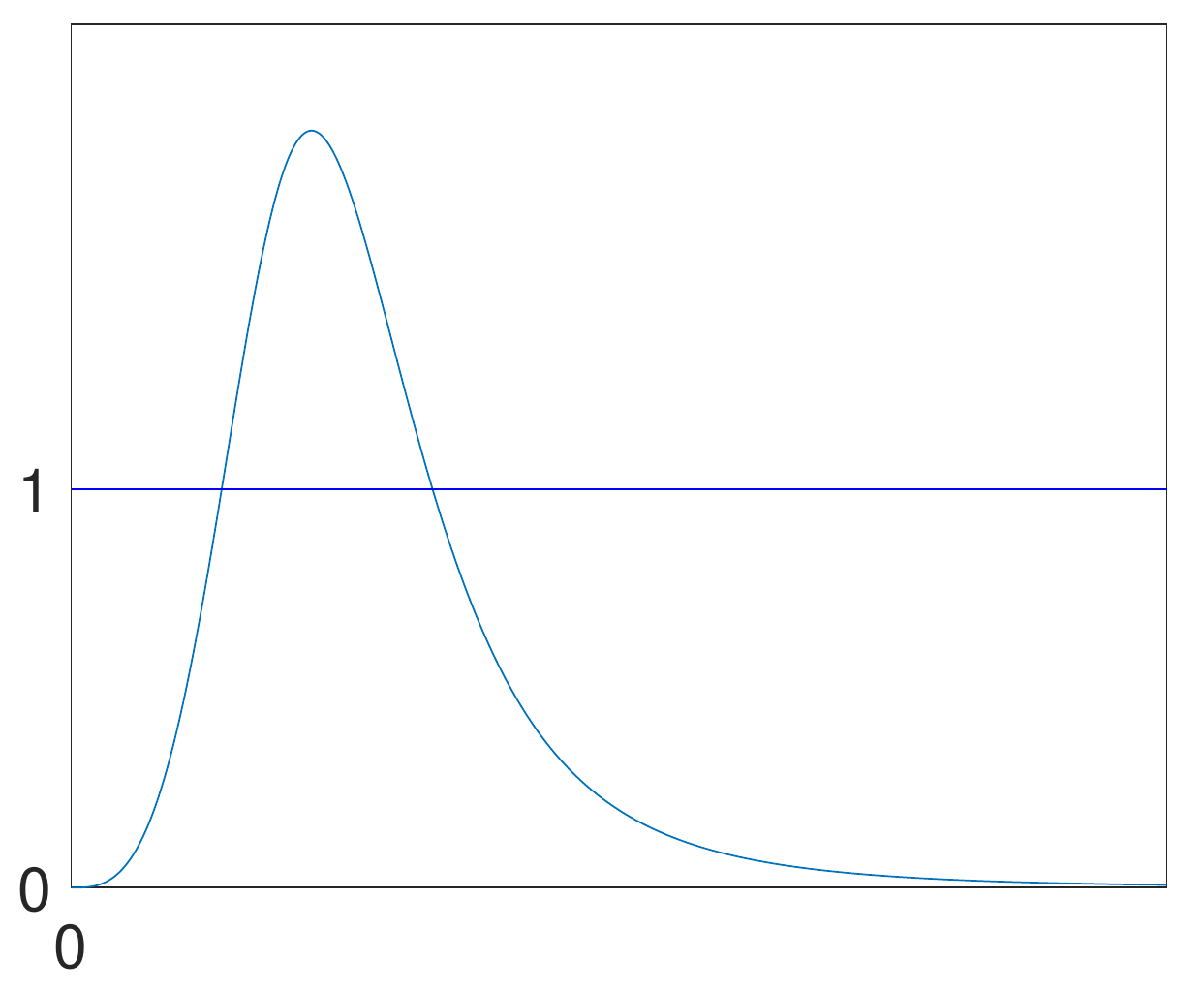}
\put(-185,3){$E$}
\put(-200,110){$f(E)$}
\put(-15,3){$E$}
\put(-110,110){$f'(E)$}
\put(-318,120){(a)}
\put(-150,120){(b)}
\caption{Inducible constant delays with $f$ as defined in Table~\ref{table:delay/velocity behaviour}.
(a) Illustration of one or three solutions to \eqref{eq:gEstarred} for different values of $C_{\beta\gamma}$.
(b) Because $f'$ is unimodal $g'(E)$ has at most two sign changes so \eqref{eq:gEstarred} cannot have more than three solutions.}
\label{fig:induc3sols}
\end{figure}

If there is no cell growth, and thus $\mu=0$, and/or if both delays are constant and independent of the state-variables ($v_M^{max}=v_M^{min}$ and $v_I^{max}=v_I^{min}$) then equation
\eqref{eq:gEstar} reduces to the form
\be \label{eq:gEstarred}
0  = g_E(E) = C_{\beta\gamma}f(E)-E
\ee
for a suitably defined constant $C_{\beta\gamma}>0$. The solutions of \eqref{eq:gEstarred} and similar equations are well-studied
in the context of monotone-cyclic feedback systems both with and without constant delay
(\cite{othmer76,tyson-othmer-1978,yildirim2004dynamics}). For the repressible case $f(E)$ is monotone decreasing,  hence $g_E(E)$ is also monotone decreasing and there is a unique steady state.
For the inducible case $f(E)$ is non-negative and monotone increasing. Here the number of steady states depends on the exact form of $f$. With $f$ defined as in Table~\ref{table:delay/velocity behaviour},
which has
a unique inflection point with $f''(E)=0$ and $E>0$,
there will be at most three steady states
as shown in \cite{yildirim2004dynamics} and illustrated in
Figure~\ref{fig:induc3sols}.

\begin{figure}[tp!]
\mbox{}\hspace*{-3mm}\includegraphics[scale=0.45]{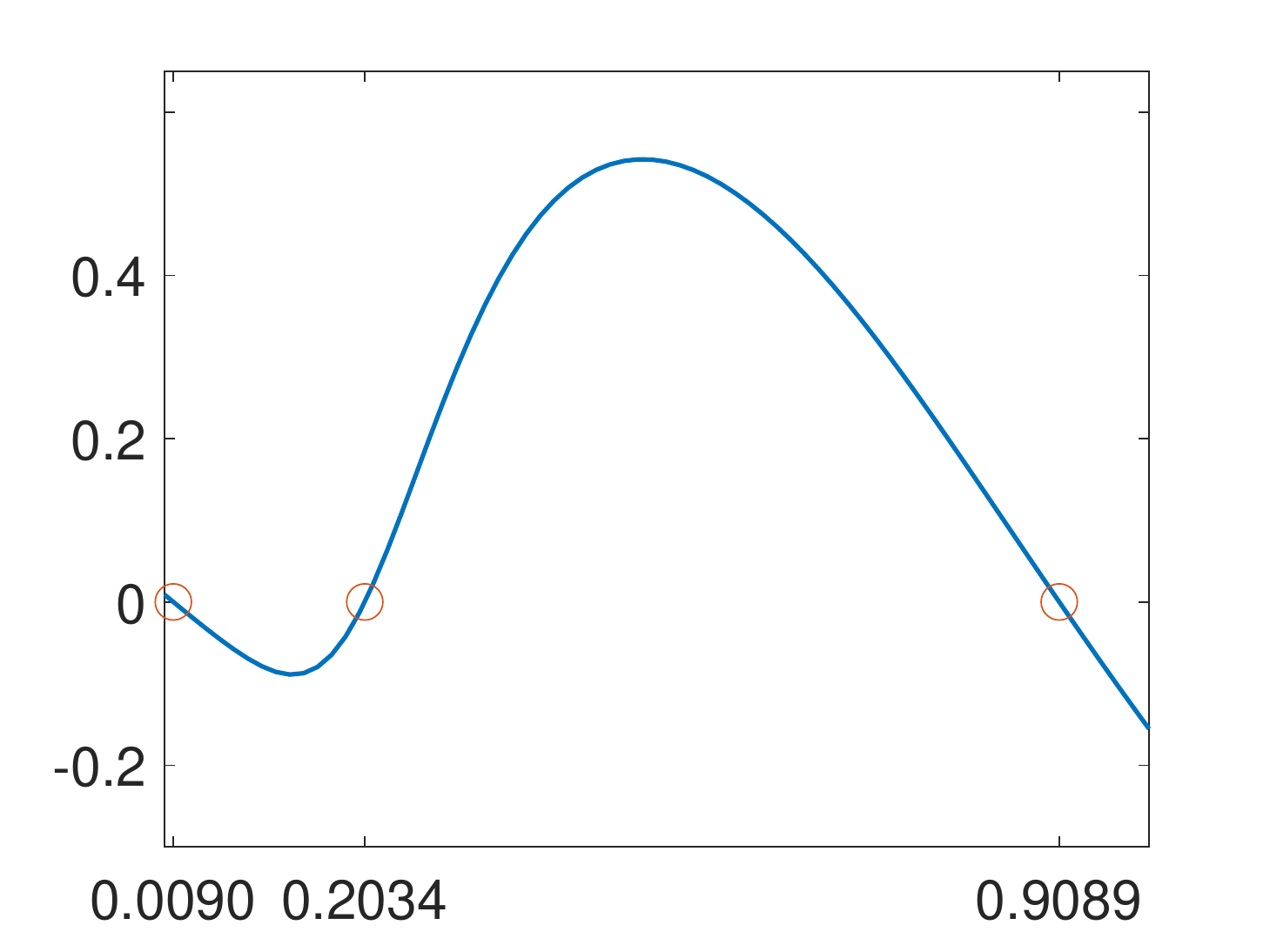}\hspace*{-6mm}\includegraphics[scale=0.45]{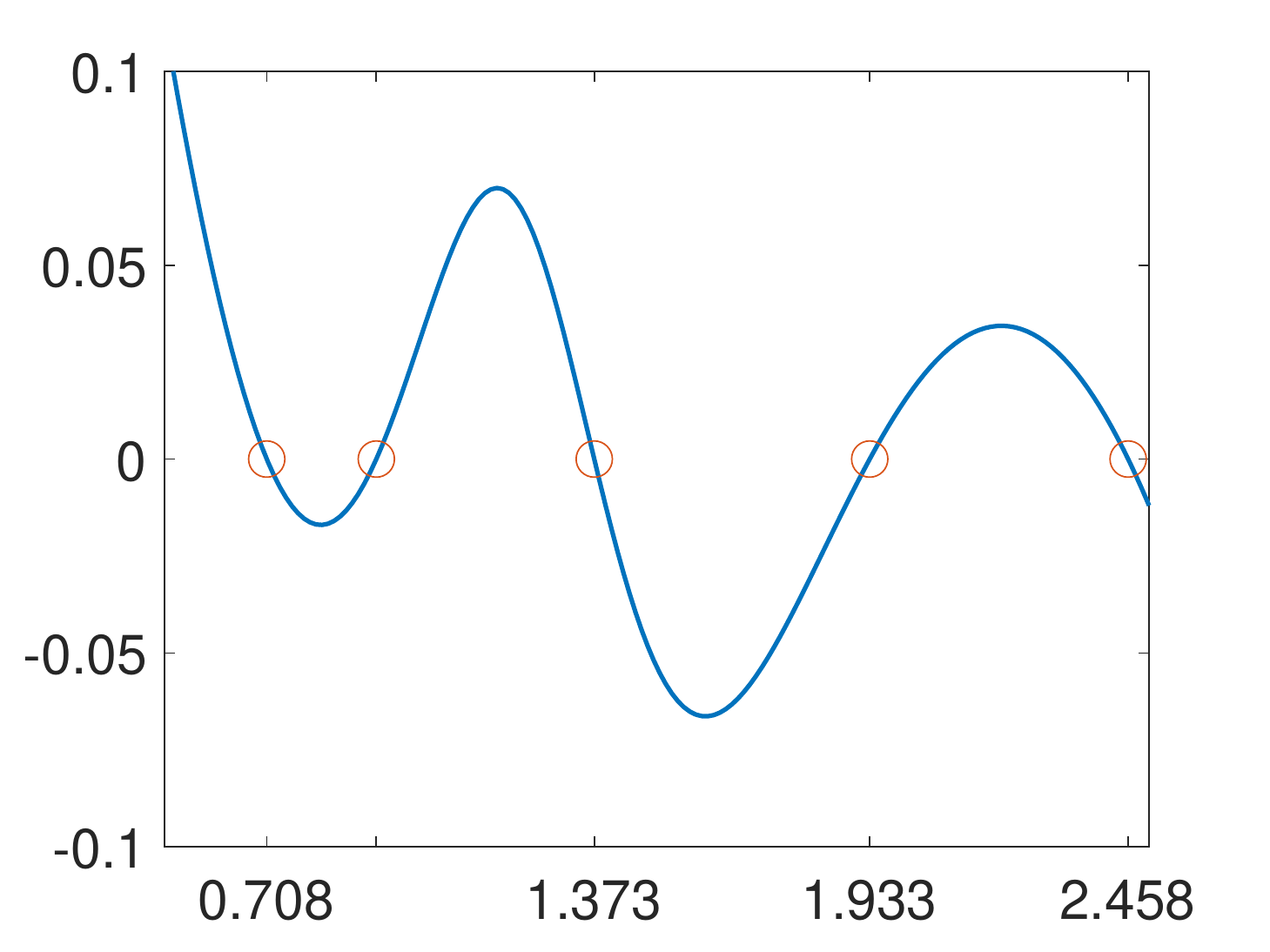}
\put(-332,120){(a) Repressible 3 steady states}
\put(-320,24){$\tau_M$ state-dependent}
\put(-160,120){(b) Inducible 5 steady states}
\put(-150,24){$\tau_M$ state-dependent}
\put(-250,5){$E$}
\put(-44,5){$E$}
\put(-350,108){\rotatebox{90}{$g_E(E)$}}

\mbox{}\hspace*{-3mm}\includegraphics[scale=0.45]{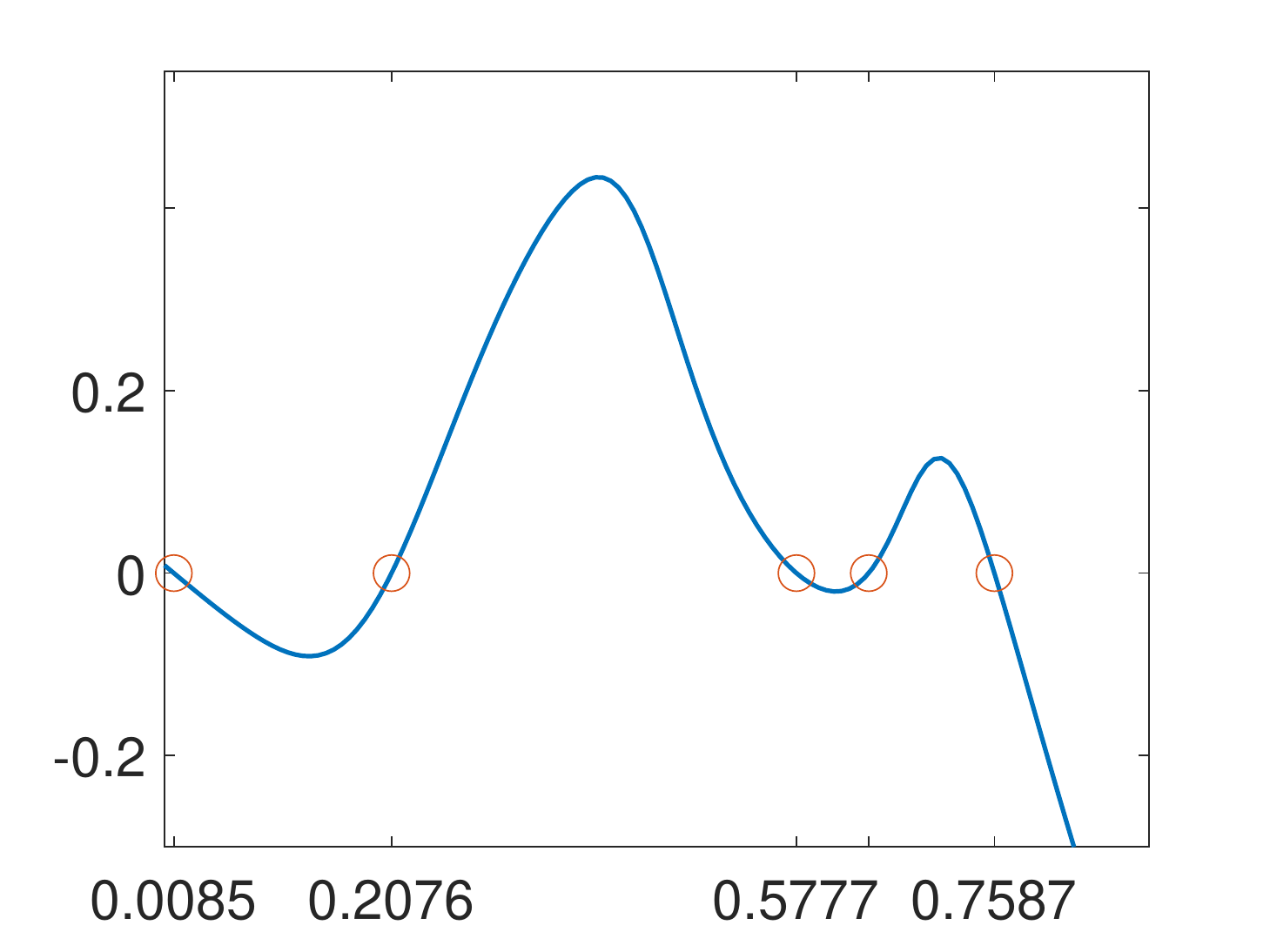}\hspace*{-6mm}\includegraphics[scale=0.45]{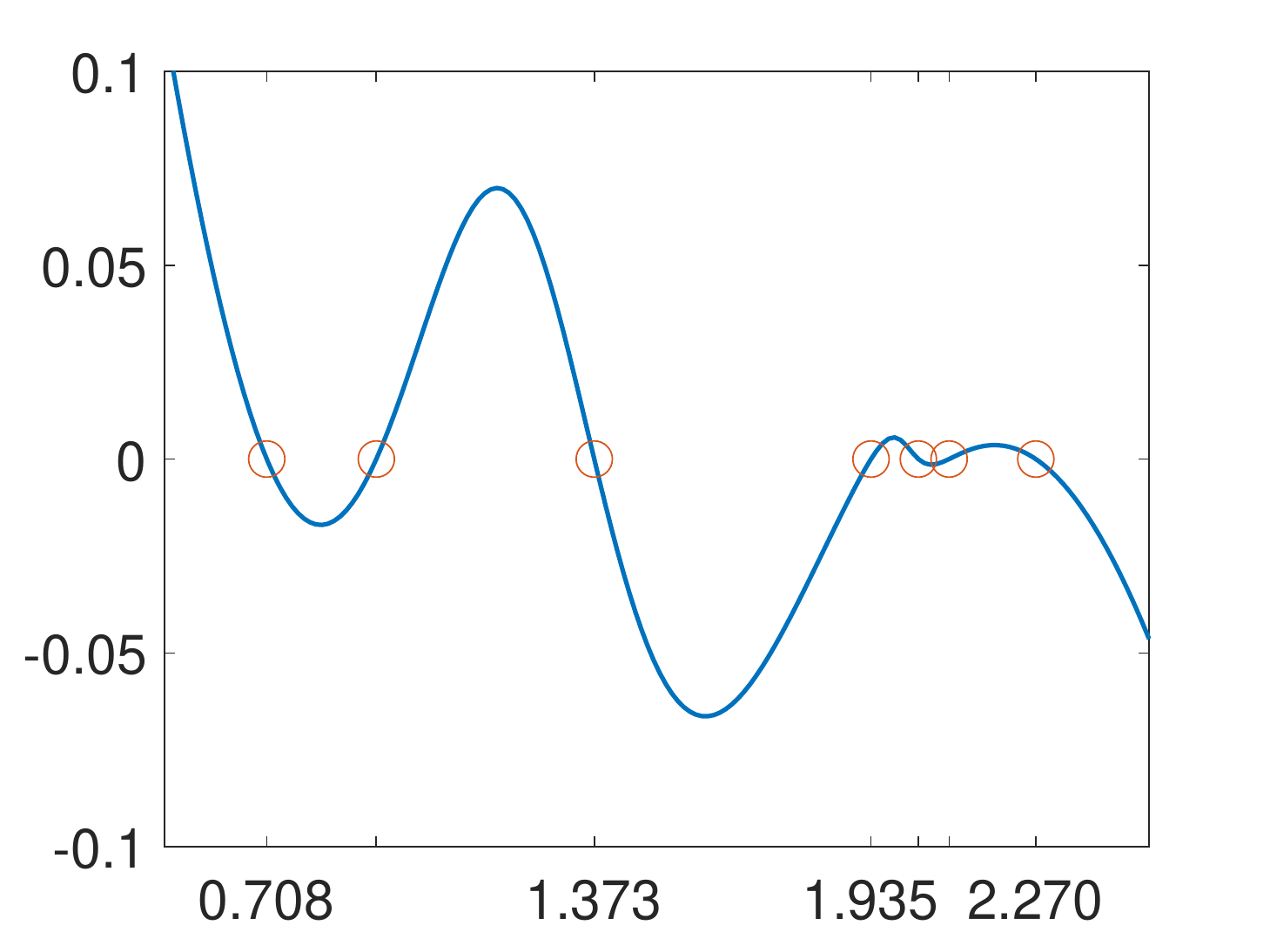}
\put(-332,120){(c) Repressible 5 steady states}
\put(-332,24){\phantom{(c)} $\tau_M$ \& $\tau_I$ state-dependent}
\put(-160,120){(d) Inducible 7 steady states}
\put(-160,24){\phantom{(d)} $\tau_M$ \& $\tau_I$ state-dependent}
\put(-198,5){$E$}
\put(-21,5){$E$}
\put(-350,108){\rotatebox{90}{$g_E(E)$}}
\caption{Examples of the function $g_E(E)$ defined by \eqref{eq:gEstar} with functions defined in
Table~\ref{table:delay/velocity behaviour}. (a) and (b) show a repressible and an inducible example with a single state-dependent delay, and parameter values given in Table~\ref{table:multsspars}. (c) and (d) show examples with two state-dependent delays with the same parameter values, except for those stated in
\eqref{eq:red5pars} and \eqref{eq:ind7pars} respectively.}
\label{fig:gmultsols}
\end{figure}
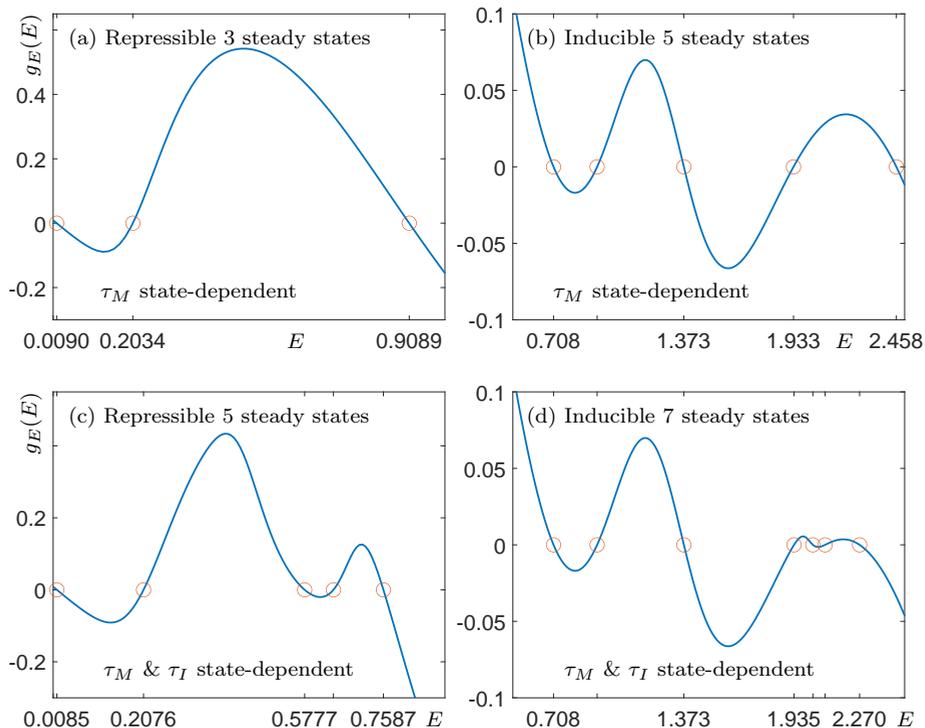

\begin{table}[htp!]
\centering
\begin{tabular}{|c|c|c|}
\hline
Quantity & Repressible & Inducible \\ \hline
		$\mu$ & 0.05 & 0.034  \\
		$\beta_M$ & 1.4 & 1  \\
		$\beta_I$ & 1 & 2  \\
		$\beta_E$ & 1 & 3  \\
		$\bgamma_M$ & 1 & 0.994  \\
		$\bgamma_I$ & 1 & 0.994  \\
		$\bgamma_E$ & 1 & 0.994  \\
		\hline
		$K$ & 2 & 10  \\
		$K_1$ & 1 & 1  \\
		$n$ & 5 & 4  \\
		\hline
		$m$ & 3 & 10  \\
		$a_M$ & 1 & 1  \\
		$v_M^{min}$ & 0.01 & 0.05  \\
		$v_M^{max}$ & 2 & 1  \\
		$E_{50}$ & 1 & 1  \\
		\hline
		$m_I$ & $20$ & $80$  \\
		$a_I$ & 1 & 1 \\
		$v_I^{min}$ & 1 & 2  \\
		$v_I^{max}$ & 1 & 2  \\
		$M_{50}$ & $1$ & $0.3374$ \\
		\hline
\end{tabular}
	\caption{Parameters used for repressible and inducible examples in (a) and (b) of Figure~\ref{fig:gmultsols}(a) and (b)}
     \label{table:multsspars}
\end{table}

In many DDEs the delay(s) only appears in the state variables,
and so do not affect the computation of the steady states. This is also
the case with our model when $\mu=0$, and the computation of the steady states
from \eqref{eq:eqM}-\eqref{eq:eqE} is independent of the delays in that case.

With state-dependent delays and cell growth, and thus $\mu>0$, the behaviour of the model \eqref{eq:mrna-delay-var}-\eqref{eq:delay by stateI} is quite different.
Now, the delays $\tau_I$ and $\tau_M$ enter explicitly into
\eqref{eq:mrna-delay-var}-\eqref{eq:delay by stateI} and hence \eqref{eq:gEstar}.
While the delay will be constant in time on any steady-state solution, with state-dependent delays
the value of the delay will depend on the state variable, which may change the structure of the phase-space of the dynamical system.

As an example consider the repressible case with state-dependent transcription velocity (but constant translation velocity). From
Table~\ref{table:delay/velocity behaviour} the transcription velocity $v_M(E)$ is a monotonic increasing function of $E$, and hence at steady state (from \eqref{eq:eqdels}) the transcription delay is a monotonic decreasing function of $E$. Then $g_E(E^*)$ contains the product of a monotonic increasing function $e^{-\mu \tau_M^*(E^*)}$ and monotonic decreasing function $f(E^*)$.  The product in  $g_E(E)$ need not be monotonic and we can no longer conclude that there is a unique steady state for the repressible case.
This is illustrated in panel (a) of Figure~\ref{fig:gmultsols} which shows an example where $g_E(E)$ has three zeros corresponding to three different steady states of the model for the  repressible case.

For an inducible operon, the situation is reversed.   The velocity $v_M(E)$ is a decreasing function  and so $e^{-\mu\tau_M^*(E^*))}$ a decreasing function of $E^*$, while the function $f(E^*)$ is an increasing function of its argument. This can again lead to additional steady states and
Figure~\ref{fig:gmultsols}(b) shows an example where $g_E(E)$ has five zeros corresponding to five different steady states in the model for the inducible case with state-dependent transcription velocity, but constant translation velocity. The full parameter sets for both of these examples are listed in Table~\ref{table:multsspars}.

In the previous examples we set $v_I^{min}=v_I^{max}=v_I$ so the translation velocity $v_I(M)=v_I$ was constant, as was the translation delay $\tau_I=a_I/v_I$. If we allow $v_I^{min}<v_I^{max}$ then
the translation delay $\tau_I(M)$ becomes a second state-dependent delay, and in
$g_E(E^*)$ the term $e^{-\tau_M^*(E^*))}f(E^*)$ is multiplied by an additional term
$e^{-\mu(\tau_I^*(M^*))}$. With the translation velocity defined as in
Table~\ref{table:delay/velocity behaviour} we see that $e^{-\mu(\tau_I^*(M^*))}$
is a monotonic increasing function of $M^*$. However, $M^*$ itself is defined by \eqref{eq:M-ss}
which again contains the product of $e^{-\tau_M^*(E^*))}$ and $f(E^*)$ that we already discussed above.
Although a full analysis of this case is beyond the scope of this paper, we note that by changing a few parameters from their
values in Table~\ref{table:multsspars}
it is possible to obtain additional steady states.
For the repressible case with
\be \label{eq:red5pars}
K=10, \quad n=10, \quad v_I^{min}=0.05, \quad v_I^{max}=0.5.
\ee
and with both the delays $\tau_M^*(E^*)$ and $\tau_I^*(M^*)$  state-dependent,
we obtain 5 co-existing steady states, as shown in
Figure~\ref{fig:gmultsols}(c). For the inducible case
with
\be \label{eq:ind7pars}
v_I^{min}=1.1.
\ee
we obtain 7 co-existing steady states, where again both delays are state-dependent.

Taken together the examples of Figure~\ref{fig:gmultsols} suggest that there can be
\be \label{eq:numss}
1+2\chi_I+2n_\tau, \qquad
\chi_I=\left\{
\begin{array}{cl}
0, & \text{\ repressible case} \\
1, & \text{\ inducible case}
\end{array}\right.
\ee
steady states where $n_\tau$ is the number of delays which are state-dependent. We cannot
prove that this is the maximum possible number of steady states, but we can construct
examples with this many steady states in a systematic way.

\begin{figure}[thp!]
\mbox{}\hspace*{-2mm}\includegraphics[scale=0.45]{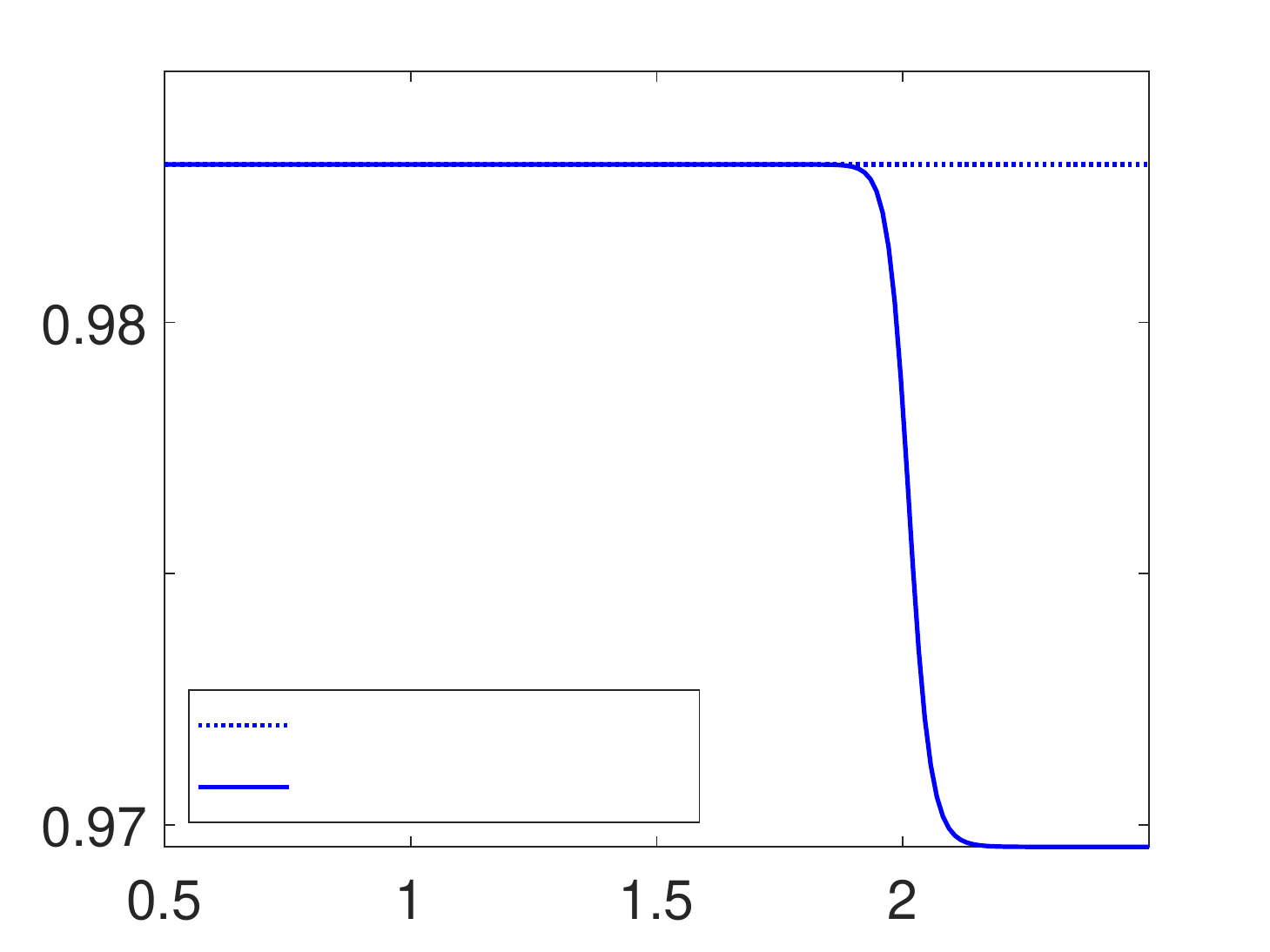}\hspace*{-6mm}\includegraphics[scale=0.45]{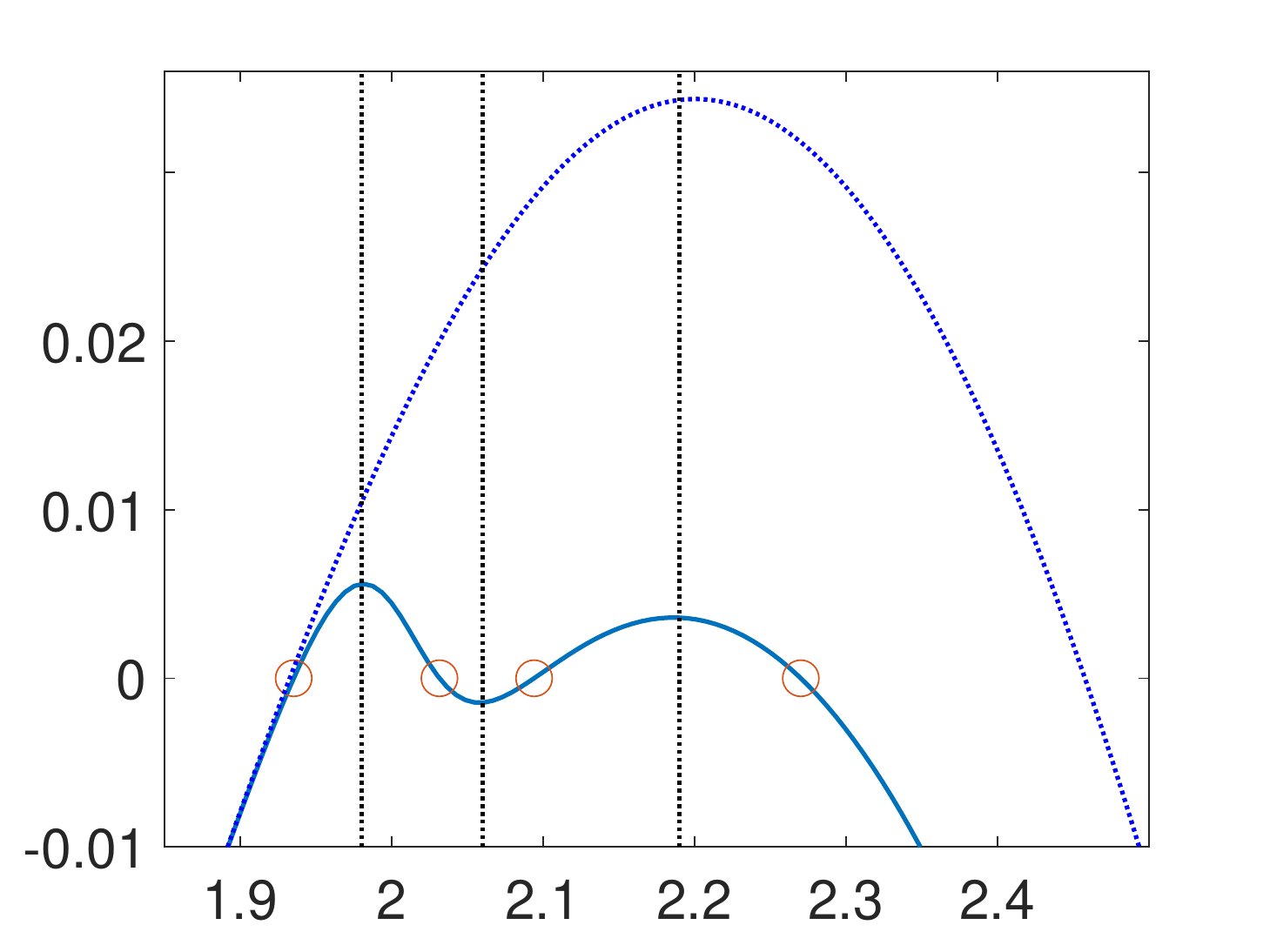}
\put(-334,122){(a)}
\put(-161,122){(b)}
\put(-203,5){$E$}
\put(-25,5){$E$}
\put(-318,30){$e^{-\mu a_I/v_I^{max}}$}
\put(-318,21){$e^{-\mu a_I/v_I(M^*)}$}
\put(-175,105){\rotatebox{90}{$g_E(E)$}}
\caption{(a) For the inducible model, the constant function $e^{-\mu\tau_I(M^*)}$ (broken line)
as in Figure~\ref{fig:gmultsols}(b), and the state-dependent function
$e^{-\mu\tau_I(M^*)}$ (solid line) as in Figure~\ref{fig:gmultsols}(b).
(b) The behaviour of the corresponding functions $g_E(E)$ for $E\approx2$. The black vertical lines separate the intervals on which $g_E(E)$ is increasing or decreasing.}
\label{fig:construct}
\end{figure}
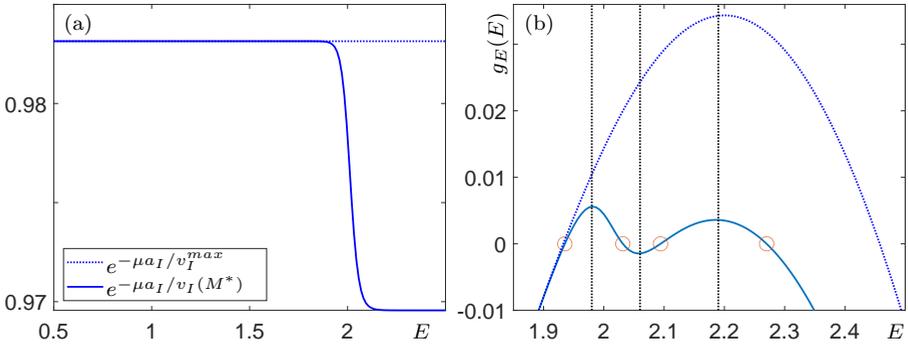

We illustrate this by showing how the example of the inducible operon with two state-dependant delays and 7 steady states in Figure~\ref{fig:gmultsols}(d) is constructed from the example with one state-dependent delay and 5 steady states in Figure~\ref{fig:gmultsols}(b). The only difference between the two examples is that in \eqref{eq:gEstar} the term $e^{-\mu\tau_I^*(M)}$ is constant in the first example, but not in the second. To make $\tau_I$ state-dependent we take $v_I^{min}<v_I^{max}$ and
$m_I\gg0$, so that the translation velocity $v_I(M)$ is
close to a step function, which results in $e^{-\mu\tau_I^*(M^*)}$ also being essentially a switching function. In Figure~\ref{fig:gmultsols} we identify that for $E\approx2$ we have $0<g_E(E)\ll1$ with
$g_E'(E)>0$, and we set the switching function to act at this point by using \eqref{eq:M-ss} to define
$M_{50}$ via
$$
M_{50} = \dfrac{\beta_M}{\bgamma_M}e^{-\mu\tau_M^*(2)}f(2).
$$
Then using \eqref{eq:M-ss} directly we obtain $e^{-\mu\tau_I^*(M^*)}$ as a function of $E^*$ as shown in
Figure~\ref{fig:construct}(a). With this element included in $g_E(E^*)$ the function is modified so that $g_E'(2)<0$, and the function gains an additional maximum and minimum for $E\approx2$, as shown in
Figure~\ref{fig:construct}(b). From there, parameters can be adjusted as needed to ensure $g_E(E)$ has a zero between each sign change of $g_E'(E)$.

\begin{figure}[htp!]
\mbox{}\hspace*{-3mm}\includegraphics[scale=0.45]{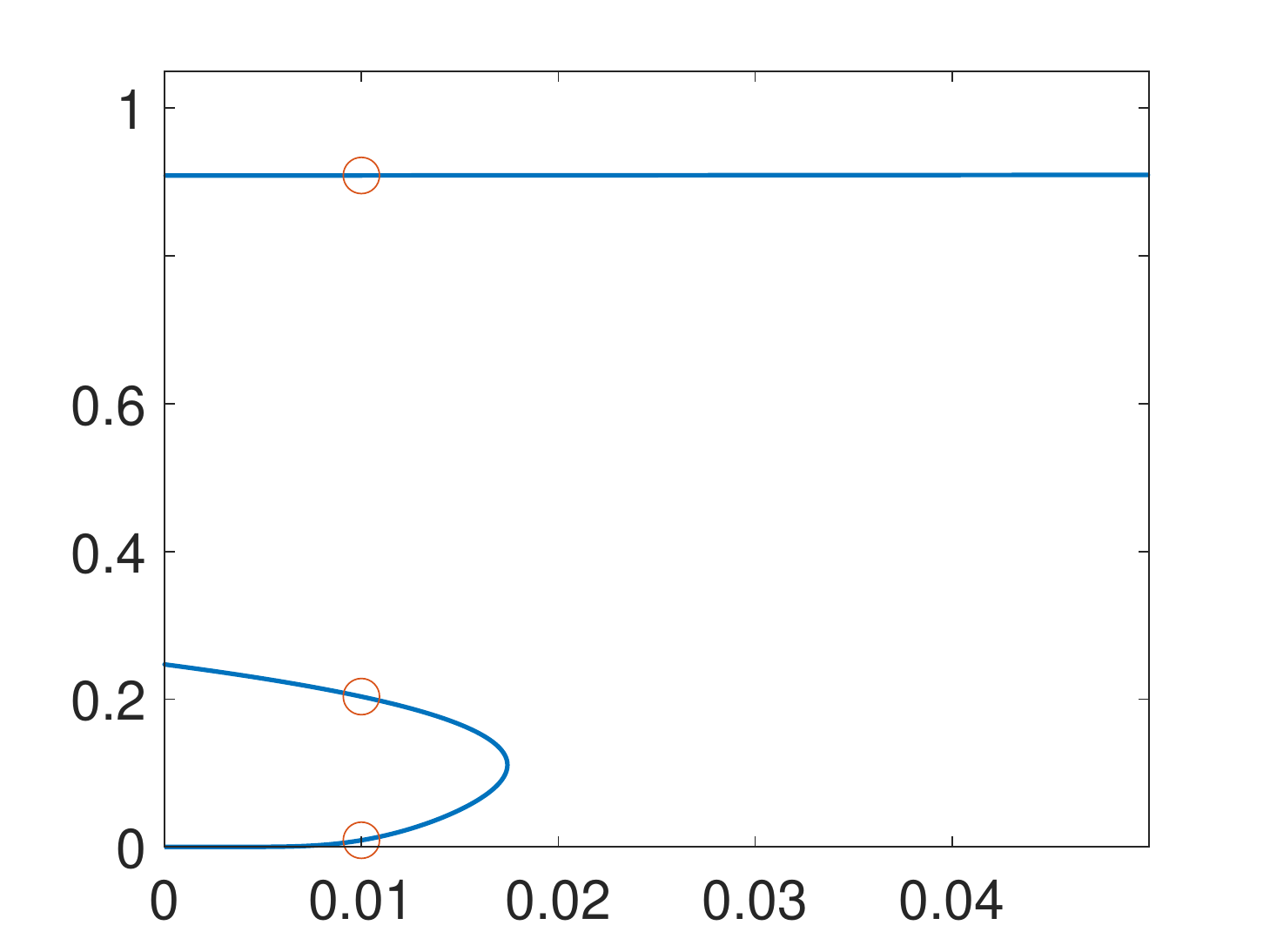}\hspace*{-6mm}\includegraphics[scale=0.45]{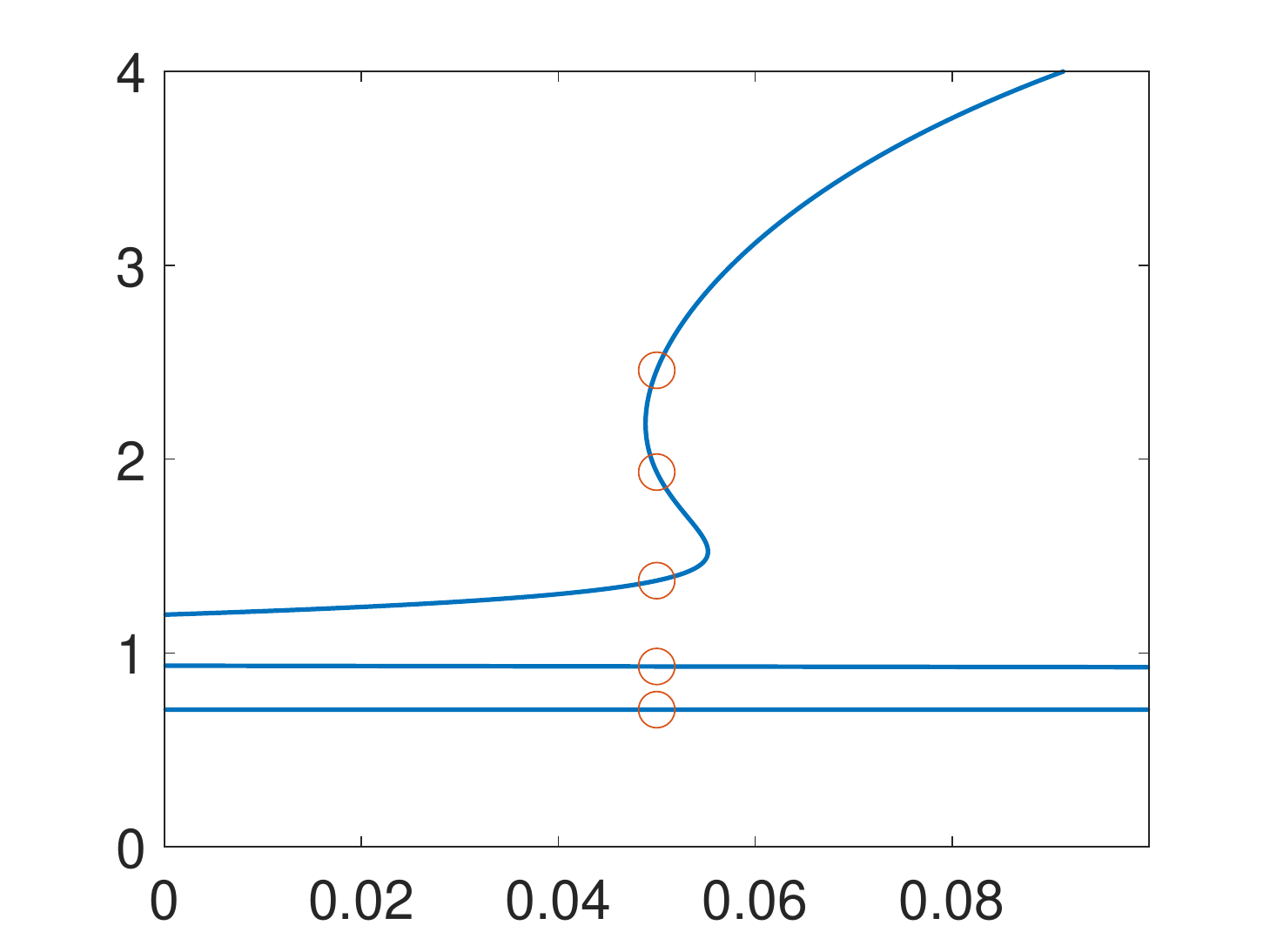}
\put(-332,100){(a) Repressible $\tau_M$ state-dependent}
\put(-160,120){(b) Inducible}
\put(-160,111){$\tau_M$ state-dependent}
\put(-203,5){$v_M^{min}$}
\put(-35,5){$v_M^{min}$}
\put(-348,100){\rotatebox{90}{$E$}}
\put(-173,113){\rotatebox{90}{$E$}}

\mbox{}\hspace*{-3mm}\includegraphics[scale=0.45]{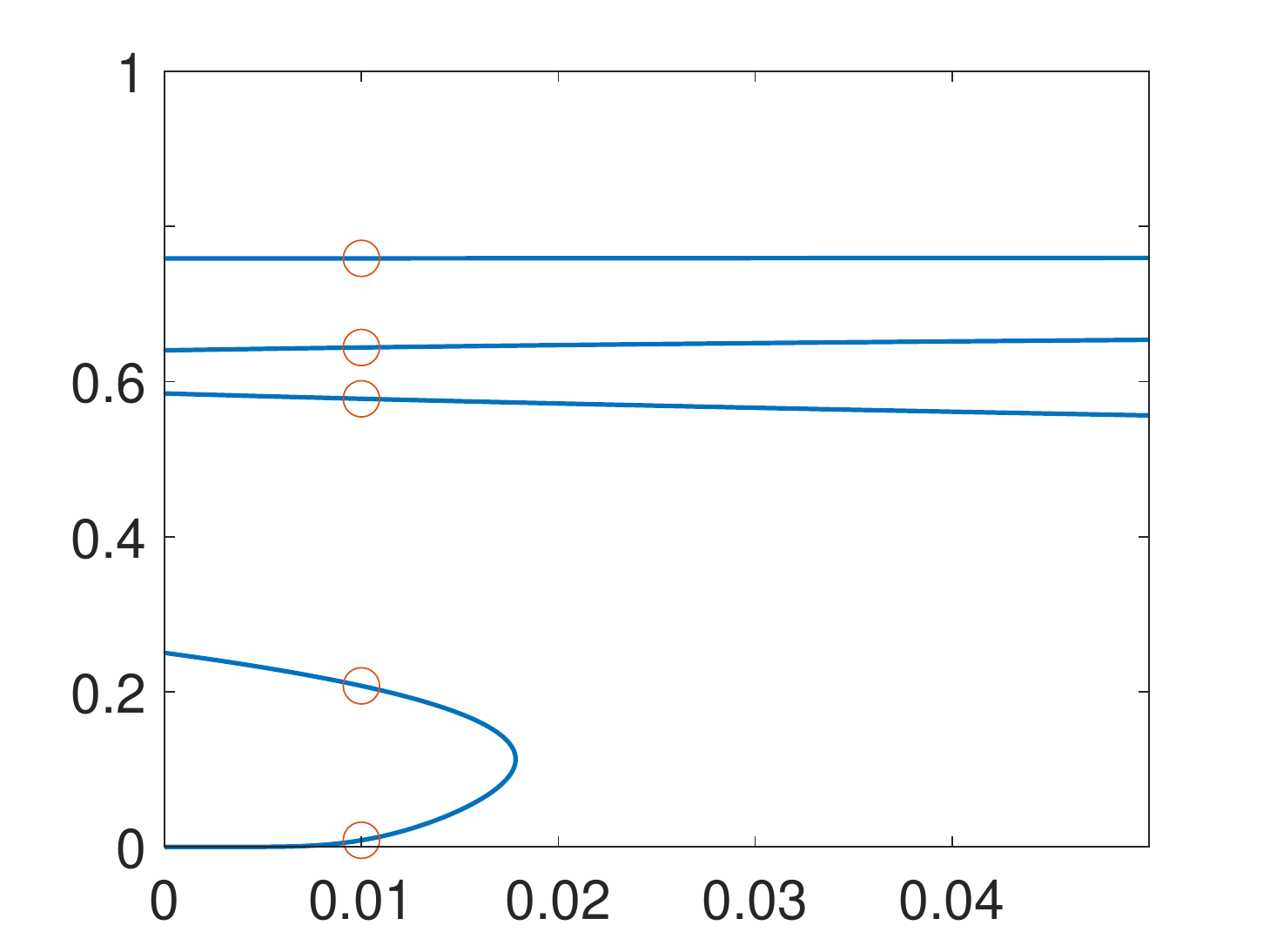}\hspace*{-6mm}\includegraphics[scale=0.45]{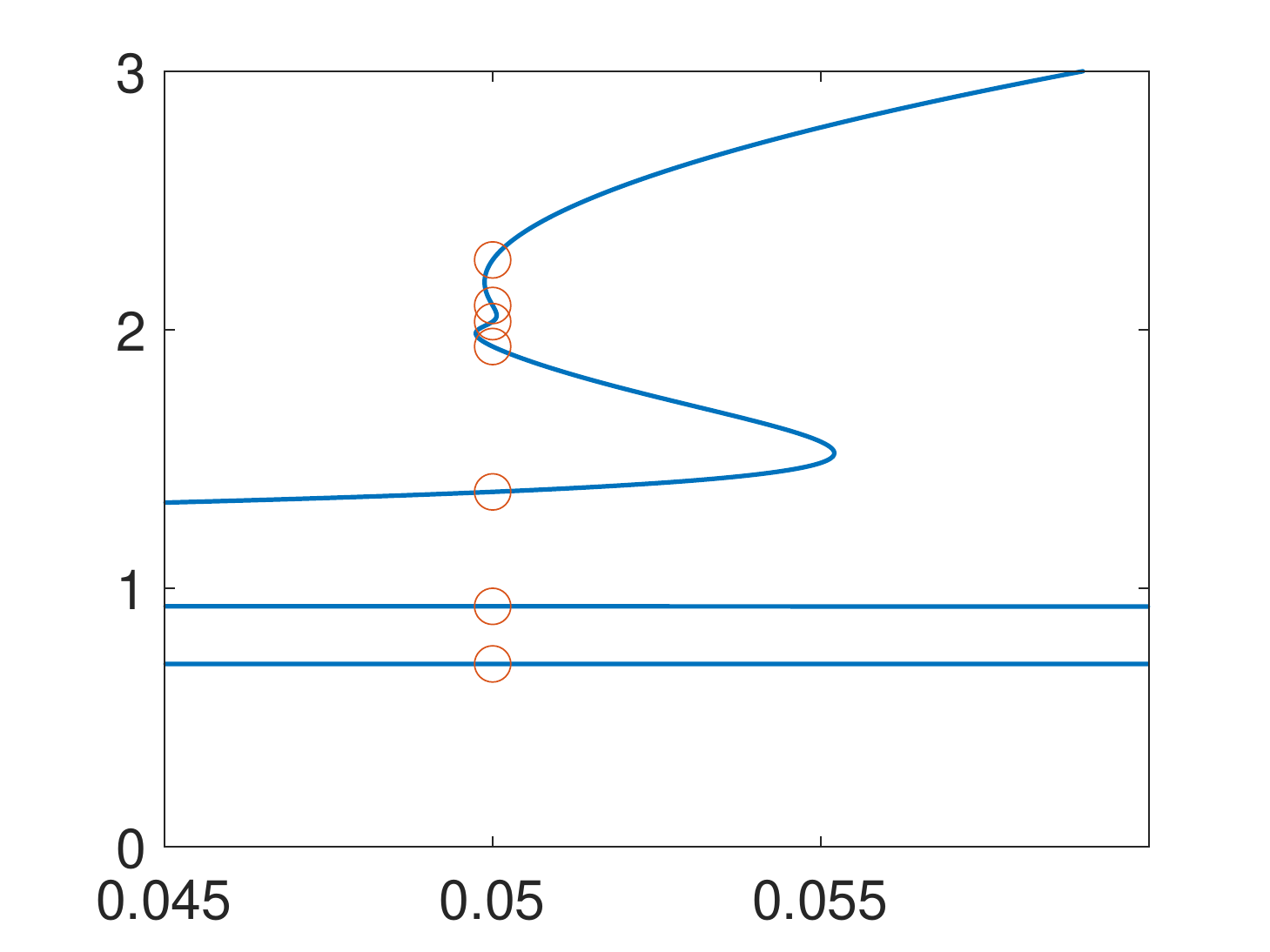}
\put(-332,120){(c) Repressible}
\put(-332,111){$\tau_M$ \& $\tau_I$ state-dependent}
\put(-160,120){(d) Inducible}
\put(-160,111){$\tau_M$ \& $\tau_I$ state-dependent}
\put(-203,5){$v_M^{min}$}
\put(-35,5){$v_M^{min}$}
\put(-348,104){\rotatebox{90}{$E$}}
\put(-174,111){\rotatebox{90}{$E$}}
\caption{One parameter continuations of the steady states as the parameter $v_M^{min}$ is varied, obtained
by plotting the zero contour of the function $g_E$. All the other parameters take the same values as in the corresponding panel of Figure~\ref{fig:gmultsols}. The red circles indicate the co-existing steady states already seen in Figure~\ref{fig:gmultsols}. In the limit as $v_M^{min}\to v_M^{max}$ the delay $\tau_M$ ceases to be state-dependent.}
\label{fig:gcont}
\end{figure}
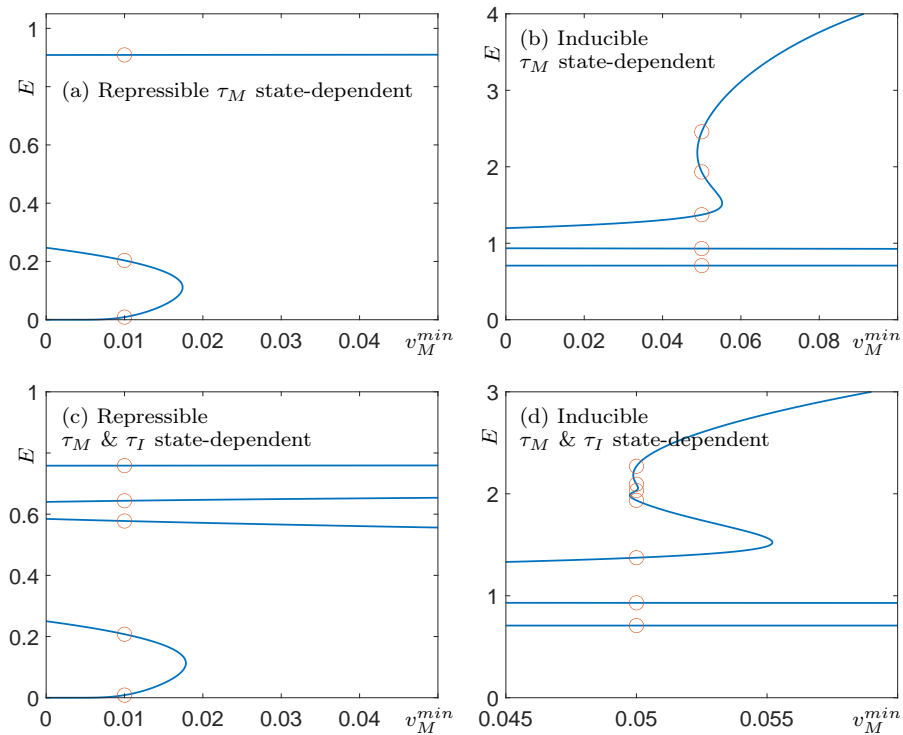

The function $g_E$ can be used to effectively perform a one-parameter continuation of the steady states by varying $E$ and one other parameter, and plotting a single contour of the function corresponding to $g_E=0$. In Figure~\ref{fig:gcont} we demonstrate four examples of one-parameter continuation in the parameter $v_M^{min}$ starting from the cases illustrated in Figure~\ref{fig:gmultsols}. Given that we obtained our examples with several co-existing steady states by constructing the function $g_E(E)$ to have multiple nearby zeros, and hence multiple local extrema close to zero, it should be no surprise to see
that the steady states from Figure~\ref{fig:gmultsols} only co-exist over a small interval of $v_M^{min}$
values and that some of them are destroyed in fold bifurcations. We note also that as $v_M^{min}$ is increased,  in the limit as $v_M^{min}\to v_M^{max}$ the delay $\tau_M$ becomes constant,
and by \eqref{eq:numss} the number of steady states  will be reduced. In particular in the cases of
Figure~\ref{fig:gcont}(a) and (b) when $v_M^{min}=v_M^{max}$ there is no state-dependency in the model and there can only be 1 or 3 steady states in the repressible and inducible cases, respectively.

Figure~\ref{fig:gcont} indicates that the steady states of  \eqref{eq:mrna-delay-var}-\eqref{eq:delay by stateI} undergo fold bifurcations. Hopf bifurcations are also ubiquitous in DDEs, and already known to occur in the repressible case with constant delays. Hence in the following sections we will study the dynamics and bifurcations of the system \eqref{eq:mrna-delay-var}-\eqref{eq:delay by stateI} and will return to the examples of this section.

%% file: quasilinearization.tex
\subsection{Linearization By Expansion}\label{ssec:quasilinear}

To determine the stability of the steady states considered Section~\ref{ssec:equilibria}, we linearize the system \eqref{eq:mrna-delay-var}--\eqref{eq:effector-delay-var} with \eqref{eq:delay by stateM} and \eqref{eq:delay by stateI} in a neighborhood of each steady state and examine the nature of the characteristic values.

This can be done  rigorously using a functional analytic approach in an appropriate Banach space, and this derivation is presented in Appendix~\ref{app:linearization}. However, that approach will not be accessible to many readers, so here we present an alternative heuristic derivation using elementary techniques which arrives at exactly the same characteristic equation as in Appendix~\ref{app:linearization}.

Assuming linear behaviour of the solution for a small perturbation from the steady state $(M^*, I^*, E^*)$,
we begin by setting
\begin{align}
M(t) &= M^*+\mathcal{E}_M e^{\lambda t}, \label{eq:Mpert} \\
I(t) &= I^*+\mathcal{E}_I e^{\lambda t}, \label{eq:Ipert} \\
E(t) &= E^*+\mathcal{E}_E e^{\lambda t}. \label{eq:Epert}
\end{align}

We denote the delays at the steady state by $\tau_M^*(E^*)$ and  $\tau_I^*(M^*)$, as defined in
\eqref{eq:eqdels}, and again write $\tau_M(t)$ and $\tau_I(t)$ for the time varying delays on a solution close to the steady state (even though as noted after \eqref{eq:delay by stateI} these delay terms are properly functions in a Banach space).

From the threshold condition \eqref{eq:delay by stateM}, Taylor expanding the integrand around the steady state we obtain
\begin{align} \notag
v_M(E^*)&\tau_M^*(E^*)  = a_M =  \int_{t-\tau_M(t)}^{t} v_M(E(s))ds
= \int_{-\tau_M(t)}^{0} v_M(E(s+t)) ds \\ \notag
&= \int_{-\tau_M(t)}^{0} v_M(E^*+\cE_Ee^{\lambda(s+t)}) ds  \\ \notag
&= \int_{-\tau_M(t)}^{0} v_M(E^*)+v_M'(E^*)\cE_Ee^{\lambda(s+t)}+\cO(\cE_E^2) ds \\
&= v_M(E^*)\tau_M(t) + \cE_E v_M'(E^*) e^{\lambda t}\int_{-\tau_M(t)}^{0} e^{\lambda s} ds+\cO(\cE_E^2).
\label{eq:tauM-exp-1}
\end{align}
Note that for $\lambda\ne0$,
\begin{equation} \label{eq:int_remov}
\int_{-a}^0 e^{\lambda s}\,ds =
\frac{1}{\lambda}(1-e^{-a\lambda})
\end{equation}
while
$$\lim_{\lambda\to0}\frac{1}{\lambda}(1-e^{-a\lambda})=a=\int_{-a}^0 e^{0 s}\,ds,$$
so $(1-e^{-a\lambda})/\lambda$ has a removable singularity at $\lambda=0$.   Therefore  we can
use \eqref{eq:int_remov} for all $\lambda\in\mathbb{C}$ and
\eqref{eq:tauM-exp-1} becomes
\begin{align}  \notag
v_M(E^*)\tau_M^*(E^*)
& = v_M(E^*)\tau_M(t)\\
& \qquad +\cE_Ee^{\lambda t}
\frac{v_M'(E^*)}{\lambda}(1-e^{-\lambda \tau_M(t)})+\cO(\cE_E^2). \label{eq:tauM-exp}
\end{align}
%
Notice that $v_M(E^*)>0$ is required for $\tau_M(E^*)$, defined by
\eqref{eq:eqdels}, to be finite; this is ensured by  the assumption that
$v_M^{min}>0$. Hence we may rearrange \eqref{eq:tauM-exp} as
$$\tau_M(t) = \tau_M^*(E^*)
-\cE_E e^{\lambda t}\frac{v_M'(E^*)}{\lambda v_M(E^*)}(1-e^{-\lambda \tau_M(t)})+\cO(\cE_E^2).$$
Noting that this implies that $\tau_M(t) = \tau_M^*(E^*)+\cO(\cE_E)$, we obtain
\begin{equation} \label{eq:tauM-epsE}
\tau_M(t) = \tau_M^*(E^*)
-\cE_E e^{\lambda t}\frac{v_M'(E^*)}{\lambda v_M(E^*)}(1-e^{-\lambda \tau_M^*(E^*)})+\cO(\cE_E^2).
\end{equation}
With \eqref{eq:tauM-epsE}, the factor $e^{-\mu\tau_M(t)}$
in \eqref{eq:mrna-delay-var} behaves as
\begin{align}
e^{-\mu\tau_M(t)} \hspace*{-2.5em} & \hspace*{2.5em}
= e^{-\mu\tau_M^*(E^*)}e^{\mathcal{E}_E\frac{\mu v_M'(E^*)}{\lambda v_M(E^*)}e^{\lambda t}(1-e^{-\lambda \tau_M^*(E^*)})}e^{-\mu\mathcal{O}(\mathcal{E}_E^2)} \nonumber \\
&= e^{-\mu\tau_M^*(E^*)}\left[1+\mathcal{E}_E\frac{\mu v_M'(E^*)}{\lambda v_M(E^*)}e^{\lambda t}(1-e^{-\lambda \tau_M^*(E^*)})+\mathcal{O}(\mathcal{E}_E^2)\right]\left[1+\mathcal{O}(\mathcal{E}_E^2)\right] \nonumber \\
&=  e^{-\mu\tau_M^*(E^*)}\left[1+\mathcal{E}_E\frac{\mu v_M'(E^*)}{\lambda v_M(E^*)}e^{\lambda t}(1-e^{-\lambda \tau_M^*(E^*)})\right]+\mathcal{O}(\mathcal{E}_E^2). \label{eq:exp-mutauM}
\end{align}
For the fraction term $v_M(E)/v_M(E(t-\tau_M(t)))$ in the differential equation, we apply Taylor expansion around the steady state to $v_M(E)$ and $1/v_M(E(t-\tau_M(t)))$ separately and then take the product, which gives:
\begin{gather*}
v_M(E) = v_M(E^*)+v_M'(E^*)\mathcal{E}_E e^{\lambda t}+\mathcal{O}(\mathcal{E}_E^2), \\
\frac{1}{v_M(E(t-\tau_M(t)))} = \frac{1}{v_M(E^*)}+ \left(-\frac{v_M'(E^*)}{v_M(E^*)^2}\right)\mathcal{E}_E e^{\lambda(t-\tau_M(t))}+\mathcal{O}(\mathcal{E}_E^2).
\end{gather*}
Thus
\begin{align}
\frac{v_M(E)}{v_M(E(t-\tau_M(t)))} &= 1+\mathcal{E}_E\frac{v_M'(E^*)}{v_M(E^*)}e^{\lambda t}(1-e^{-\lambda\tau_M(t)})+\mathcal{O}(\mathcal{E}_E^2) \nonumber \\
&= 1+\mathcal{E}_E\frac{v_M'(E^*)}{v_M(E^*)}e^{\lambda t}(1-e^{-\lambda\tau_M^*(E^*)})+\mathcal{O}(\mathcal{E}_E^2). \label{eq:vM-vtauM}
\end{align}
Following the derivation of \eqref{eq:exp-mutauM} and \eqref{eq:vM-vtauM}, similarly we have
\begin{gather}
e^{-\mu\tau_I(t)} = e^{-\mu\tau_I^*(M^*)}\left[1+\mathcal{E}_M\frac{\mu v_I'(M^*)}{\lambda v_I(M^*)}e^{\lambda t}(1-e^{-\lambda \tau_I^*(M^*)})\right]+\mathcal{O}(\mathcal{E}_M^2),
\label{eq:exp-mutauI} \\
\frac{v_I(M)}{v_I(M(t-\tau_I(t)))} = 1+\mathcal{E}_M\frac{v_I'(M^*)}{v_I(M^*)}e^{\lambda t}(1-e^{-\lambda\tau_I^*(M^*)})+\mathcal{O}(\mathcal{E}_M^2). \label{eq:vI-vtauI}
\end{gather}

Now we use these expansions to linearize the system \eqref{eq:mrna-delay-var}--\eqref{eq:effector-delay-var} equation by equation.
Substituting the perturbations \eqref{eq:Mpert} and \eqref{eq:Epert}
into \eqref{eq:mrna-delay-var}
and using the expansions \eqref{eq:exp-mutauM} and \eqref{eq:vM-vtauM} we have
\begin{align*}
\cE_M\lambda e^{\lambda t} \hspace*{-2.5em} & \hspace*{2.5em}
= \dfrac{d}{dt}(M^*+\mathcal{E}_Me^{\lambda t})\\
& = \beta_M\left[e^{-\mu\tau_M^*(E^*)}\!\left(\!1+\cE_E\frac{\mu v_M'(E^*)}{\lambda v_M(E^*)}e^{\lambda t}(1-e^{-\lambda \tau_M^*(E^*)})\!\right)\!+\mathcal{O}(\mathcal{E}_E^2)\right] \times  \\
& \qquad \left[1+\mathcal{E}_E\frac{v_M'(E^*)}{v_M(E^*)}e^{\lambda t}(1-e^{-\lambda\tau_M(E^*)})+\mathcal{O}(\mathcal{E}_E^2)\right]
\!f(E^*\!+\mathcal{E}_E e^{\lambda(t-\tau_M(t))}) \\
& \qquad\quad -\bgamma_M(M^*+\mathcal{E}_Me^{\lambda t}) \\
&= \beta_M e^{-\mu\tau_M^*(E^*)}\left[1+\mathcal{E}_E\frac{v_M'(E^*)}{v_M(E^*)}e^{\lambda t}(1-e^{-\lambda \tau_M^*(E^*)})(1+\frac{\mu}{\lambda})\right] \times \\
& \qquad \left[f(E^*)+ \mathcal{E}_E f'(E^*) e^{\lambda(t-\tau_M(t))}\right]-\bgamma_M(M^*+\mathcal{E}_Me^{\lambda t})+\mathcal{O}(\mathcal{E}_E^2)  \\
&= \beta_M e^{-\mu\tau_M^*(E^*)}\left(f(E^*)+\mathcal{E}_E f(E^*)\frac{v_M'(E^*)}{v_M(E^*)}e^{\lambda t}(1-e^{-\lambda \tau_M^*(E^*)})(1+\frac{\mu}{\lambda})\right.  \\
& \qquad + \mathcal{E}_E f'(E^*)e^{\lambda(t-\tau_M^*(E^*))} \bigg) -\bgamma_M(M^*+\mathcal{E}_Me^{\lambda t})+\mathcal{O}(\mathcal{E}_E^2).
\end{align*}
Using the equality \eqref{eq:eqM} and multiplying by $e^{-\lambda t}$, this simplifies to
\begin{align} \notag
\mathcal{E}_M\lambda
& = \mathcal{E}_E\beta_M e^{-\mu\tau_M^*(E^*)}\bigg(f(E^*)\frac{v_M'(E^*)}{v_M(E^*)}(1-e^{-\lambda \tau_M^*(E^*)})(1+\frac{\mu}{\lambda})\\
& \qquad\quad + f'(E^*)e^{-\lambda\tau_M^*(E^*)}\bigg)
-\mathcal{E}_M\bgamma_M+\mathcal{O}(\mathcal{E}_E^2) \label{eq:Mlin}
\end{align}
For the second differential equation \eqref{eq:intermed-delay-var},
substituting the perturbation \eqref{eq:Mpert} and \eqref{eq:Ipert} and the expansion \eqref{eq:exp-mutauI} and \eqref{eq:vI-vtauI} we
similarly find that
\begin{align*}
\mathcal{E}_I\lambda e^{\lambda t} & 
= \dfrac{d}{dt}(I^*+\mathcal{E}_Ie^{\lambda t})\\
&= \beta_I e^{-\mu\tau_I^*(M^*)}\!\left(\!M^*+\mathcal{E}_M e^{\lambda(t-\tau_I^*(M^*))}
+\mathcal{E}_M M^*\frac{v_I'(M^*)}{v_I(M^*)}\times \right.  \\
& \qquad\qquad\qquad
e^{\lambda t}(1-e^{-\lambda \tau_I^*(M^*)})(1+\frac{\mu}{\lambda})\!\Big)-\bgamma_I(I^*\!+\mathcal{E}_Ie^{\lambda t})+\mathcal{O}(\mathcal{E}_M^2).
\end{align*}
Using the equality \eqref{eq:eqI} and multiplying by $e^{-\lambda t}$ we obtain
\begin{align} \notag
\mathcal{E}_I\lambda
&= \mathcal{E}_M\beta_I e^{-\mu\tau_I^*(M^*)}\left(e^{-\lambda \tau_I^*(M^*)}+M^*\frac{v_I'(M^*)}{v_I(M^*)}(1-e^{-\lambda \tau_I^*(M^*)})(1+\frac{\mu}{\lambda})\right)  \\
& \qquad -\bgamma_I\mathcal{E}_I+\mathcal{O}(\mathcal{E}_M^2). \label{eq:Ilin}
\end{align}

Lastly, the case of the differential equation
\eqref{eq:effector-delay-var} is simpler, since it is linear with no delays.
Substituting the perturbations \eqref{eq:Ipert} and \eqref{eq:Epert} we have
$$
\mathcal{E}_E\lambda e^{\lambda t}
= \frac{d}{dt}(E^*+\mathcal{E}_Ee^{\lambda t})
= \beta_E(I^*+\mathcal{E}_I e^{\lambda t})-\bgamma_E(E^*+\mathcal{E}_E e^{\lambda t}).
$$
Using the equality \eqref{eq:eqE}, and multiplying by $e^{-\lambda t}$ this simplifies to
\be \label{eq:Elin}
\mathcal{E}_E\lambda = \mathcal{E}_I\beta_E-\mathcal{E}_E\bgamma_E.
\ee

Combining \eqref{eq:Mlin}, \eqref{eq:Ilin} and \eqref{eq:Elin}, and dropping the higher order terms
gives the linear system
\begin{equation} \label{eq:Alinsys}
A(\lambda)
\begin{pmatrix}
\cE_M \\
\cE_I \\
\cE_E
\end{pmatrix}
=
\left(\begin{matrix}
-\bgamma_M-\lambda & 0 & A_{13}\\
A_{21} & -\bgamma_I-\lambda  & 0\\
0 & \beta_E   &  -\bgamma_E-\lambda
\end{matrix}\right)
\begin{pmatrix}
\cE_M \\
\cE_I \\
\cE_E
\end{pmatrix}
=0
\end{equation}
where the $A_{13}$ and $A_{21}$ entries of the $3\times 3$-matrix $A(\lambda)$ are defined by
\begin{gather*}
A_{13} = \beta_M e^{-\mu\tau_{M\!}^*(E^{*\!})\!}
\!\left(
\!f(E^*)\frac{v_{M\!}'(E^*)}{v_{M\!}(E^*)}(1-e^{-\lambda
\tau_{M\!}^*(E^{*\!})})(1+\frac{\mu}{\lambda})+f'^{\!}(E^{*\!})e^{-\lambda\tau_{M\!}^*(E^{*\!})\!}\!
\right)\hspace*{-0.7mm},
\\
A_{21} = \beta_I e^{-\mu\tau_I^*(M^*)}\left(M^*\frac{v_I'(M^*)}{v_I(M^*)}(1-e^{-\lambda \tau_I^*(M^*)})(1+\frac{\mu}{\lambda}) +e^{-\lambda \tau_I^*(M^*)}\right).
\end{gather*}
%
%
%
The characteristic equation
of \eqref{eq:mrna-delay-var}--\eqref{eq:effector-delay-var} is
\be \label{eq:CharEqA}
\Delta(\lambda)=\det(A(\lambda))=0,
\ee
with $\Delta(\lambda)$ given by
\begin{equation}
\Delta(\lambda) = (\bgamma_M+\lambda)(\bgamma_I+\lambda)(\bgamma_E+\lambda)
+\beta_M\beta_I\beta_E e^{-\mu(\tau_{I\!}^*(M^*)+\tau_{M\!}^*(E^*))}k(\lambda),\label{eq:Delta-lambda-MIE}
\end{equation}
where
\begin{align} \notag
k(\lambda) &=\left(
\frac{v_M'(E^*)}{v_M(E^*)}f(E^*)( 1-e^{-\lambda\tau_M^*(E^*)})\Bigl(1+\frac\mu\lambda\Bigr)
+ f'^{\!}(E^*)e^{-\lambda\tau_M^*(E^*)} \right)\\
& \qquad \times
\left(
\frac{v_I'(M^*)}{v_I(M^*)}M^*(1- e^{-\lambda\tau_I^*(M^*)})\Bigl(1+\frac\mu\lambda\Bigr)
+e^{-\lambda\tau_I^*(M^*)}
\right)  \label{eq:klambda1}
\end{align}

Exactly the same characteristic equation is derived completely rigorously in Appendices
\ref{app:semiflow}-\ref{app:linearization} culminating in equation~\eqref{eq:char-eq-at-MIE}.

In contrast to the rigorous variational approach used in the appendices, 
here we assert without proof that all the quantities of interest can be written as functions of the
perturbation parameters $\mathcal{E}_M$, $\mathcal{E}_I$ and $\mathcal{E}_E$. 
For example in equation \eqref{eq:tauM-epsE} we have Taylor expanded the state-dependent delay $\tau_M$ as a function of $\mathcal{E}_E$. To justify that rigorously requires
functional analysis, and  this is done in Proposition~\ref{prop:1} in Appendix~\ref{app:semiflow}.

Another drawback of the derivation above is that   there is no theory to show that
stability of steady states is determined by the characteristic equation 
\eqref{eq:CharEqA}.
However, since for
our model both approaches lead to the same characteristic equation, the theory relating stability to the characteristic equation applies \citep{HKWW}. Therefore  the stability of equilibria of the system \eqref{eq:mrna-delay-var}-\eqref{eq:effector-delay-var} is determined by characteristic values arising from \eqref{eq:CharEqA}.

%
%

%% file: numerics.tex
\section{Numerical Methods} \label{sec:method}

In this section, we describe numerical methods to study the distributed state-dependent delay  model
\eqref{eq:mrna-delay-var}-\eqref{eq:delay by stateI}.
We would like to conduct one-parameter continuation of steady states and periodic orbits
and compute local stability and bifurcations in Matlab~\citep{matlab}. The standard package for performing numerical bifurcation analysis of DDEs in Matlab is DDE-BIFTOOL~\citep{ddebiftool}.
Unfortunately, although it can handle constant or discrete state-dependent delays, DDE-BIFTOOL cannot be applied directly to problems with distributed state-dependent
delays defined by threshold conditions such as \eqref{eq:delay by stateM} and \eqref{eq:delay by stateI}.
Likewise, the built-in Matlab function \texttt{ddesd} for solving DDE initial value problems is also only implemented for discrete delays.

In equations \eqref{eq:delay by stateM} and \eqref{eq:delay by stateI} the delay $\tau_M$ or $\tau_I$ can be determined by adjusting the lower limit of the integral until the integral has the desired value $a_M$.
Naive numerical implementations would use a bisection or secant iteration to determine the delay, which would necessitate evaluating the integral at each step of the iteration. This would be very slow to compute and would become the main bottleneck slowing down numerical computations. Another problem in evaluating this integral is that the numerical DDE solvers that have been implemented in Matlab are all written for discrete delays and only give access to the value of the solution $u(t-\tau_j)$ at the discrete delays, whereas to evaluate the integral in \eqref{eq:delay by stateM} or \eqref{eq:delay by stateI} we require
the values of the integrand across the whole interval.

We describe below two different implementations of \eqref{eq:mrna-delay-var}-\eqref{eq:delay by stateI}
in DDE-BIFTOOL, neither of which require an iteration to find the delays, and also show how to apply \texttt{ddesd} to solve the initial value problem.

\subsection{Steady State Computations - Linearization Correction}
\label{subsec:correction}

As discussed in Section~\ref{ssec:equilibria}, we can obtain steady states of \eqref{eq:mrna-delay-var}-\eqref{eq:delay by stateI} from the scalar function $g_E(E^*)$ defined in \eqref{eq:gEstar}. Any solution to $g_E(E^*)=0$ gives the $E$ component of a steady state with corresponding $M$ and $I$ components given by \eqref{eq:eqM} and \eqref{eq:eqE}.

We would like to conduct one-parameter continuation of steady states and compute local stability
using DDE-BIFTOOL, but as noted above it cannot be directly applied to solve DDEs when the delay is defined by a threshold condition. Nevertheless, at a steady state the threshold integral conditions
\eqref{eq:delay by stateM} and \eqref{eq:delay by stateI} become integrals of constant functions. Consequently,
the delays are defined by \eqref{eq:eqdels} and can be treated as discrete delays.
Therefore we are able to implement the system \eqref{eq:mrna-delay-var}-\eqref{eq:effector-delay-var} together with \eqref{eq:eqdels} in DDE-BIFTOOL and use it to continue the steady states. This approach also allows us to locate fold bifurcations of steady states. However, although replacing \eqref{eq:delay by stateM} and \eqref{eq:delay by stateI} by \eqref{eq:eqdels} preserves the existence of the steady states, (as detailed in
\cite{wendy-thesis} for a related model) characteristic values and hence stability of the steady state are altered.  The reason is that the integration of exponential perturbations along the solution, see Section~\ref{ssec:quasilinear}, are not included.

To recover the correct stability information when using the modified problem \eqref{eq:mrna-delay-var}-\eqref{eq:effector-delay-var} with discrete delays \eqref{eq:eqdels} we perform
linearization correction using the characteristic equation.
The characteristic roots of the steady state of the modified problem are taken as ``seed values'' which are then corrected by applying the Matlab  nonlinear system solver \texttt{fsolve} to the
exact characteristic equation \eqref{eq:CharEqA} for the original model \eqref{eq:mrna-delay-var}-\eqref{eq:delay by stateI}. This works well at a majority of points along  continuation branch; however, it behaves poorly at some points leading to spurious bifurcations.   In addition,  sometimes the algorithm does not converge,   while sometimes
the solver converges to an already found characteristic value. As well as creating duplicates of characteristic values, this  leads to missing some characteristic values, and so does not reliably classify bifurcations.

To deal with these issues, we remove any characteristic values at which the algorithm fails to converge, as well as
any duplicate values. To resolve the issue of missing characteristic values, we use the corrected characteristic values from the previous point on the branch as a second set of ``seed values'' to compute additional characteristic values, where again we remove duplicates. The removal of duplicate characteristic values is somewhat dangerous, because it could result in missing genuine instances of characteristic values with multiplicity larger than one. However, in practice,
we did not encounter this problem.

Once the corrected characteristic roots are computed at each steady state, we obtain the correct stability information. Hopf bifurcations occur when two complex conjugate characteristic values cross the imaginary axis, or equivalently
when the number of complex characteristic values with positive real part changes by two. Fold bifurcations happen when a real characteristic value changes sign, which we can also detect by the number of complex characteristic values with positive real part changing by one. Because we obtain the stability from a modified problem and we did not alter any of the DDE-BIFTOOL subroutines which use linearization inappropriate for our model,   we are not able to use additional DDE-BIFTOOL subroutines such as those that detect criticality of Hopf bifurcations and perform normal form computations.

\subsection{Steady State and Periodic Orbit Computation - Delay Discretization}
\label{subsec:discretization}

While the approach of Section~\ref{subsec:correction} allows us to compute the stability of steady states and hence to detect fold and Hopf bifurcations, it cannot be used to compute periodic orbits, because the delays
\eqref{eq:delay by stateM}-\eqref{eq:delay by stateI} would not be constant on periodic orbits.

The only way to tackle the full distributed state-dependent delay operon model
\eqref{eq:mrna-delay-var}-\eqref{eq:delay by stateI} is to evaluate the integrals in the threshold conditions
\eqref{eq:delay by stateM}-\eqref{eq:delay by stateI}. While this cannot be done exactly in DDE-BIFTOOL,
 it is enough to evaluate the integral to sufficient accuracy using a numerical quadrature scheme.

As the two delays are of similar form, we will describe the method for approximating $\tau_M$ by discretizing the integral
in \eqref{eq:delay by stateM} using the composite  trapezoidal method and seeking the value of $\tau_M$ that satisfies
\eqref{eq:delay by stateM}. To do this we introduce extra ``dummy'' delays as follows.

\sloppy{With $a_M$ fixed and $v_M(E)  \in [v_M^{min}, v_M^{max}]$, it follows that $\tau_M \in [a_M/v_M^{max}, a_M/v_M^{min}]$. To obtain the state-dependent delay $\tau_M$ that satisfies the threshold condition \eqref{eq:delay by stateM}, we discretize the interval $\left[t-a_M/v_M^{max}, t \right]$} uniformly with a sequence of mesh points $$t=x_0> x_1 > ... > x_N=t-\frac{a_M}{v_M^{max}}$$ and define $N$ constant ``dummy'' delays $$\tau_j=t-x_j=t-\frac{j}{N}\frac{a_M}{v_M^{max}}, \qquad j=1, 2,\ldots, N.$$
In particular, as $\tau_M \geqslant a_M/v_M^{max}$, there is no need for detection of the delay $\tau_M$ over the interval $[t-a_M/v_M^{max}, t]$.

On the interval $[t-a_M/v_M^{min}, t-a_M/v_M^{max}]$ where the delay $\tau_M$ lies, we detect it as follows.
Divide the interval $[t-a_M/v_M^{min}, t-a_M/v_M^{max}]$ into $N$ equal width subintervals,
$$t-\frac{a_M}{v_M^{max}}=x_N > x_{N+1} >... > x_{2N}=t-\frac{a_M}{v_M^{min}},$$
  which   implies that
$$x_{N+j} = t-\frac{a_M}{v_M^{max}}-j \frac{a_M}{N} (\frac{1}{v_M^{min}}-\frac{1}{v_M^{max}}), \qquad j=1, 2, ..., N.$$
We then define another $N$ constant ``dummy'' delays
$$\tau_{N+j}=t-x_{N+j}, \qquad j=1, 2, ..., N.$$

To compute $\tau_M$, we take advantage of the functionality of DDE-BIFTOOL which allows state-dependent delays to be defined as functions of the other delays and the solution values at those delays.   We  let
$$J(j)=\int_{x_j}^{x_0} v_M(E(s))ds = \int_{t- \tau_j}^t v_M(E(s))ds, \qquad j=1, 2, \ldots, 2N$$
and let $J_h(j)$ be the numerical approximation of $J(j)$ using the composite trapezoidal rule,
$$J_h(j) = \sum_{k=1}^{j-1} \frac{1}{2}\left(v_M(E(x_k))+v_M(E(x_{k+1}))\right)(x_k-x_{k+1}).$$ 
We look for the largest $j$ such that $J_h(j) \leqslant a_M$. Since
\begin{align*}
	a_M &> \int_{t-a_M/v_M^{max}}^t v_M(E(s))ds \\
	& \approx \frac{a_M}{N v_M^{max}}\left[\frac{1}{2}v_M(E(t))+\frac{1}{2}v_M(E(t-\frac{a_M}{v_M^{max}}))+\sum_{j=1}^{N-1}v_M(E(t-\tau_j))\right] \\
& = J_h(N),
\end{align*}
we successively add subintervals to the integral until we find $j$ such that
\begin{equation}
a_M \geqslant J_h(j) \quad \text{and} \quad a_M < J_h(j+1).
\label{eq:locatetauM}
\end{equation}
  With such a  $j$, we have   $\tau_M\in[\tau_j,\tau_{j+1})$. To locate $\tau_M$ more precisely, consider
\begin{align*}
	a_M = \int_{t-\tau_M}^t v_M(E(s))ds
	=\int_{t-\tau_M}^{t-\tau_j} v_M(E(s))ds+\int_{t-\tau_j}^t v_M(E(s))ds,
\end{align*}
which implies
\begin{equation}
a_M-J_h(j) \approx \int_{t-\tau_M}^{t-\tau_j} v_M(E(s))ds.
\label{eq:japprox}
\end{equation}
Applying the trapezoidal rule again, we have
\begin{equation}
\int_{t-\tau_M}^{t-\tau_j} v_M(E(s)) ds  \approx \frac{\tau_M-\tau_j}{2}[v_M(E(t-\tau_M))+v_M(E(t-\tau_j))],
\label{eq:trap1}
\end{equation}
and using a linearization in the subinterval (which is consistent with the trapezoidal method) we have
\begin{align} \notag
v_M(E(t-\tau_M)) & = v_M(E(t-\tau_j))\\
& \qquad +(\tau_M-\tau_j)[v_M(E(t-\tau_{j+1}))-v_M(E(t-\tau_j))].
\label{eq:trap2}
\end{align}
Substituting \eqref{eq:trap1} and \eqref{eq:trap2} into \eqref{eq:japprox} gives
\begin{align*} \notag
	a_M- J_h(j) &\approx
\frac{\tau_M-\tau_j}{2}[(2+\tau_j-\tau_M)v_M(E(t-\tau_j)) \notag\\
& \qquad\qquad\qquad\qquad +(\tau_M-\tau_j)v_M(E(t-\tau_{j+1}))]. 	
\end{align*}
Rearranging this we find that $\tau_M$ is given as the solution of
$k(\tau_M)=0$ where
\begin{align}
k(\tau_M) 
& = \frac{(\tau_M-\tau_j)^2}{2}[v_M(E(t-\tau_j))-v_M(E(t-\tau_{j+1}))] \notag\\
& \qquad \qquad\qquad \qquad-(\tau_M-\tau_j)v_M(E(t-\tau_j))+(a_M-J_h(j))
\label{eq:quadtauM}
\end{align}
Note that \eqref{eq:quadtauM} is a quadratic function of $\tau_M$ and the condition \eqref{eq:locatetauM} guarantees that $k(\tau_M)$ has a zero for $\tau_M \in [\tau_j, \tau_{j+1}]$.
Applying the quadratic formula to \eqref{eq:quadtauM}, we obtain the solution
\begin{align}
\tau_M  & = \tau_j + \frac{v_{M\!}(E(x_j))}{v_M(E(x_j))-v_M(E(x_{j+1}))}    \notag \\
& \qquad - \frac{\sqrt{v_{M\!}(E(x_j)^2-2(a_M-J_h(j))(v_{M\!}(E(x_j))-v_{M\!}(E(x_{j+1})))}}{v_M(E(x_j))-v_M(E(x_{j+1}))}
\label{result_tauM}
\end{align}
where the minus sign in the quadratic formula ensures that   the root $\tau_M\in[\tau_j,\tau_{j+1}]$   when   $k(\tau_M)$ is either a concave up or concave down parabola.

With this implementation we are able to apply DDE-BIFTOOL directly to an approximation to
the system \eqref{eq:mrna-delay-var}-\eqref{eq:delay by stateI}. Since the stability computations are carried out within DDE-BIFTOOL (as opposed to the linearization correction technique described in Section~\ref{subsec:correction}) we are able to use the full functionality of DDE-BIFTOOL which allows us to determine criticality of bifurcations and also to compute branches of periodic orbits emanating from
Hopf bifurcations.

The choice of parameters for the numerical discretization is somewhat delicate. If the discretization is too coarse convergence issues arise in the branch continuation, while finer discretizations allow for a smoother continuation of branches with larger continuation steps, but at the cost of each step being very slow. This arises because the numerical linear algebra problems at the heart of the approximate Newton method in each DDE-BIFTOOL continuation step increase in complexity with both the number of delays and the size of the collocation problem. The total number of delays in the discretized problem is $2N+2$, composed of the $2N$ dummy delays, the (assumed constant) delay $\tau_I$, and the computed state-dependent delay $\tau_M$ given by \eqref{result_tauM}.

For the computations in
Section~\ref{sec:examples} we use degree 4 or 5 collocation polynomials and 20 to 40 mesh intervals resulting in 80 to 200 collocation points on the periodic orbit.
For stability computations of steady states we took
$N=32$ which results in 65 constant delays and one state-dependent delay. For computation of periodic orbits we took $N=48$ resulting in close to one hundred delays in the discretized problem.

The computation of each step of the continuation is quite slow compared to the implementation
of Section~\ref{subsec:correction}. The algorithms give consistent results on problems for which both
can be applied (with bifurcation points agreeing to between 3 and 5 significant digits of accuracy), but the algorithm of this section is more widely applicable.
For the results shown in Section~\ref{sec:examples}, we mainly use the discretization method described in this section, with the linearization correction method of Section~\ref{subsec:correction} used to validate the results.

\subsection{Solving Initial Value Problems (IVPs)}
\label{ssec:ddesd}

Simulating IVPs allows us to investigate the dynamics in parameter regimes where none of the steady states are stable.  In Section~\ref{sec:examples} we find stable periodic orbits which do not arise from Hopf bifurcations by   following this procedure.

The Matlab routine \texttt{ddesd}
solves DDE initial value problems with discrete state-dependent delays.  While we would like to use
\texttt{ddesd} to study \eqref{eq:mrna-delay-var}-\eqref{eq:delay by stateI},   we need to address the issue of implicitly defined  delays.

For simplicity, as in the preceding sections, we treat $\tau_I$ as a constant delay which   is defined by \eqref{eq:eqdels}.
We deal with the state-dependent delay $\tau_M$ defined by \eqref{eq:delay by stateM} by
differentiating the integral in \eqref{eq:delay by stateM} with respect to $t$ 
to obtain
$$0=v_M(E(t))-\Big(1-\frac{d\tau_{M}\!}{dt}(t)\Big)v_M(E(t-\tau_M)),$$
which implies that
\begin{equation}
\frac{d\tau_{M}\!}{dt}(t)= 1-\frac{v_M(E(t))}{v_M(E(t-\tau_M))}.
\label{eq:dtauM}
\end{equation}
We can thus solve the system \eqref{eq:mrna-delay-var}-\eqref{eq:delay by stateI} as an initial value problem by considering the system of three equations \eqref{eq:mrna-delay-var}-\eqref{eq:effector-delay-var}
augmented by \eqref{eq:dtauM} to define the evolution of the state-dependent delay $\tau_M$ along with the constant delay $\tau_I=a_I/v_I$ where $v_I=v_I^{min}= v_I^{max}$. The case where $\tau_I$ is state-dependent can be handled similarly.

Although this trick avoids the need to evaluate the integral in \eqref{eq:delay by stateM} during the simulation, care needs to be taken since information is lost when differentiating and while a solution of \eqref{eq:delay by stateM} also solves \eqref{eq:dtauM}, the converse is not necessarily true. To ensure our solution of \eqref{eq:dtauM} also solves \eqref{eq:delay by stateM}, we specify history functions so that \eqref{eq:delay by stateM} is satisfied at time $t=t_0$.
In particular, we require $\tau_M(t_0)$ to satisfy
\begin{equation}
a_M= \int_{t_0-\tau_M(t_0)}^{t_0} v_M(E(s)) ds.
\label{eq:ICtauM}
\end{equation}
This will depend on the choice of the history function $E(t)$ defined for $t \leqslant t_0$. In general we need to evaluate this integral only once. Even this can be avoided if $E(t)=E_0$ is constant for $t \leqslant t_0$ since then \eqref{eq:ICtauM} simplifies to $a_M=\tau_M(t_0) v_M(E_0)$ which implies that
$$\tau_M(t_0)=\frac{a_M}{v_M(E_0)}.$$

Although we do not need to solve the integral threshold condition \eqref{eq:delay by stateM} during the numerical computation, after a numerical solution is computed, it is very easy to evaluate the integral on the right hand side of \eqref{eq:delay by stateM} to check how close it is to $a_M$. In all the examples presented in Section~\ref{sec:examples}, this defect is smaller that $10^{-5}$ at the final time indicating 5 or more digits of accuracy in the computation of the threshold condition across the interval of computation.

To find the period of a stable periodic solution, a simple technique is to take advantage of the idea of the Poincar\'e section, and we implement an event function to detect periodicity. While \texttt{ddesd} has a built-in event detection function which can be used to detect periodicity, it slows down the numerical solution drastically.
Instead, once the simulation is complete, we fit a spline to the numerical solution and use the spline functions within Matlab to obtain the crossings of the Poincar\'e section and maxima and minima of solutions and hence period and amplitude information.

Once we find a stable periodic orbit, the solution may be continued as one parameter is varied either by performing additional numerical IVP solves to find a periodic orbit for a perturbed parameter set, or
by importing the numerically computed periodic solution into DDE-BIFTOOL and use
the discretization of Section~\ref{subsec:discretization} to continue the solution.
The DDE-BIFTOOL discretization has the advantage that it can equally well find stable and unstable periodic orbits, and we will use it in Section~\ref{sec:examples} to detect fold bifurcations of periodic orbits where the stability of the periodic orbit changes.

\texttt{ddesd} can only be used for the continuation of stable periodic orbits, which is useful for validating the DDE-BIFTOOL results. To perform continuation with \texttt{ddesd}
we use the stable periodic solution at each iteration as the history function for the next computation when the continuation parameter is slightly changed. With a small perturbation in the continuation parameter value, we expect to converge to the stable periodic solution as it should still lie in the basin of attraction. Care needs to be taken when doing this, since when making a perturbation of the parameters we need to recompute initial value of the state-dependent delay $\tau_M(t_0)$ so that the integral \eqref{eq:ICtauM} is satisfied at the initial time $t_0$ with the new parameter set and history function given by the numerical solution with the previous parameter set.

%% file: dynamics.tex
\section{Dynamics of Repressible and Inducible Operons with State-Dependent Delays}\label{sec:examples}

In this section, we explore the dynamics of the Goodwin operon model \eqref{eq:mrna-delay-var}-\eqref{eq:delay by stateI} incorporating state-dependent delays. We will mainly focus on the case where the transcription delay $\tau_M$ is state-dependent and the translation delay $\tau_I$ is constant. Then equations
\eqref{eq:mrna-delay-var}-\eqref{eq:delay by stateI} simplify to
\begin{equation}
\begin{aligned}
\dfrac{dM\!}{dt}(t) & = \beta_M e^{-\mu \tau_M(t)} \dfrac{v_M(E(t))}{v_M(E(t-\tau_M(t)))} f(E(t-\tau_M(t))) -\bgamma_M M(t), \\
\dfrac{dI}{dt}(t) & = \beta_I e^{-\mu\tau_I} M(t-\tau_I) -\bgamma_I I(t), \\
\dfrac{dE}{dt}(t) & = \beta_E I(t) -\bgamma_E E(t), \\
a_M &= \int_{t-\tau_M(t)}^{t} v_M(E(s)) ds=\int_{-\tau_M(t)}^{0} v_M(E(t+s)) ds,
\end{aligned}
\label{eq:sysonedel}
\end{equation}
with
$\tau_I=a_I/v_I$ where $v_I=v_I^{min}= v_I^{max}$.
The
respective functions for a repressible or inducible system
are defined in Table~\ref{table:delay/velocity behaviour}.
We will treat the minimum transcription velocity, $v_M^{min}$, as a bifurcation parameter.

\subsection{Repressible Operon with One State-Dependent Delay}
\label{sec:repexamples}

Recall that when there are no state-dependent delays there are only two possibilities for a
repressible system. Namely there is either a globally stable steady state,
or a globally stable limit cycle
which arises through a supercritical Hopf bifurcation from the steady state.
We already showed in Section~\ref{ssec:equilibria} (see Figures~\ref{fig:gmultsols}(a) and~\ref{fig:gcont}(a)) that it is possible for a repressible system with one state-dependent delay to have multiple steady states, as well as fold bifurcations of steady states.
In this section we will explore the dynamics of the repressible system in more depth to reveal the possible dynamics and bifurcations that may arise.





\begin{figure}[thp!]
\centering
\includegraphics[scale=0.6]{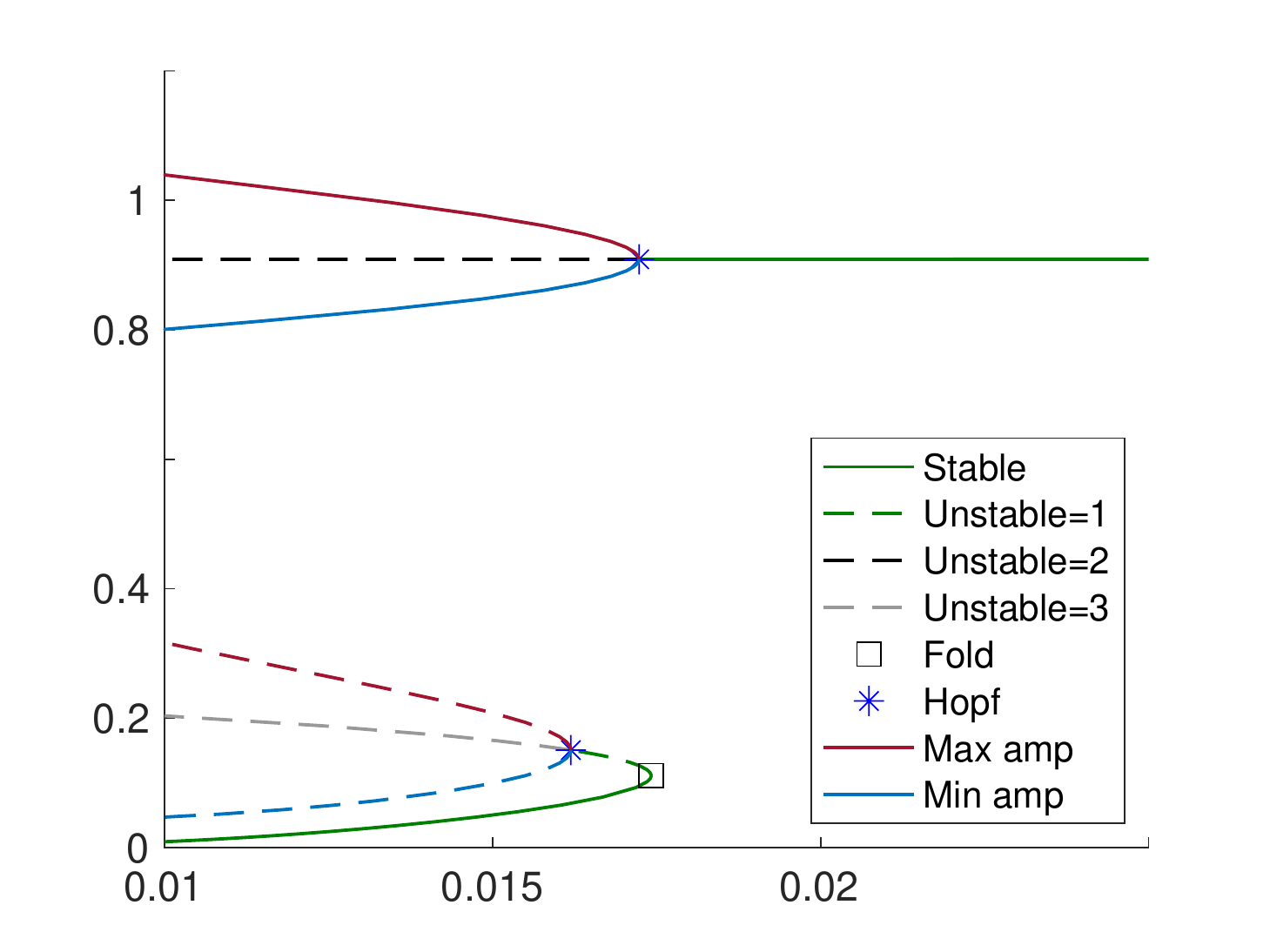}
\put(-230,165){\rotatebox{90}{$E$}}
\put(-40,10){$v_M^{min}$}
\caption{Bifurcation diagram of the model \eqref{eq:sysonedel} for a repressible system,
with parameters defined in Table~\ref{table:multsspars} except $v_M^{min}$ which is taken as the bifurcation parameter.
Solid lines represent stable objects including stable steady state (in green) and stable limit cycle (maximum amplitude in red and minimum amplitude in blue).
Steady states are represented using the $E$-component of the solution, and the amplitude of
periodic solutions is taken from the maximum and minimum of the $E(t)$ on the periodic solution.
Dashed lines represent unstable objects including unstable steady states (depending on the number of eigenvalues with positive real part, green for one, black for two and gray for three and more) and an unstable limit cycle. Bifurcations are listed in Table~\ref{tab:threess_onepara_trp}. }
	\label{fig:threess_onepara_trp}
\end{figure}

\begin{figure}[thp!]
\hspace*{-1em}\includegraphics[scale=0.45]{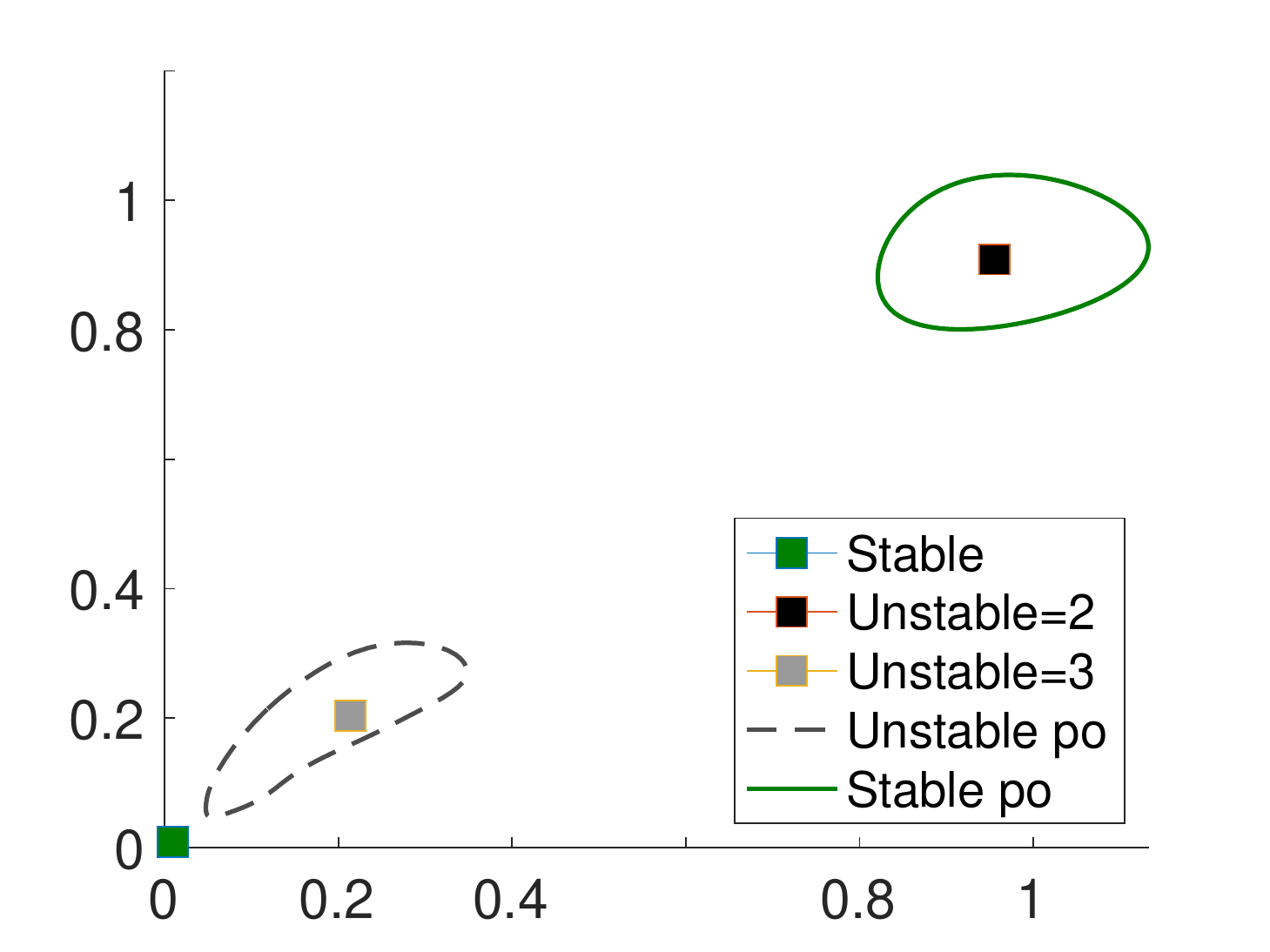}\hspace*{-2em}\includegraphics[scale=0.45]{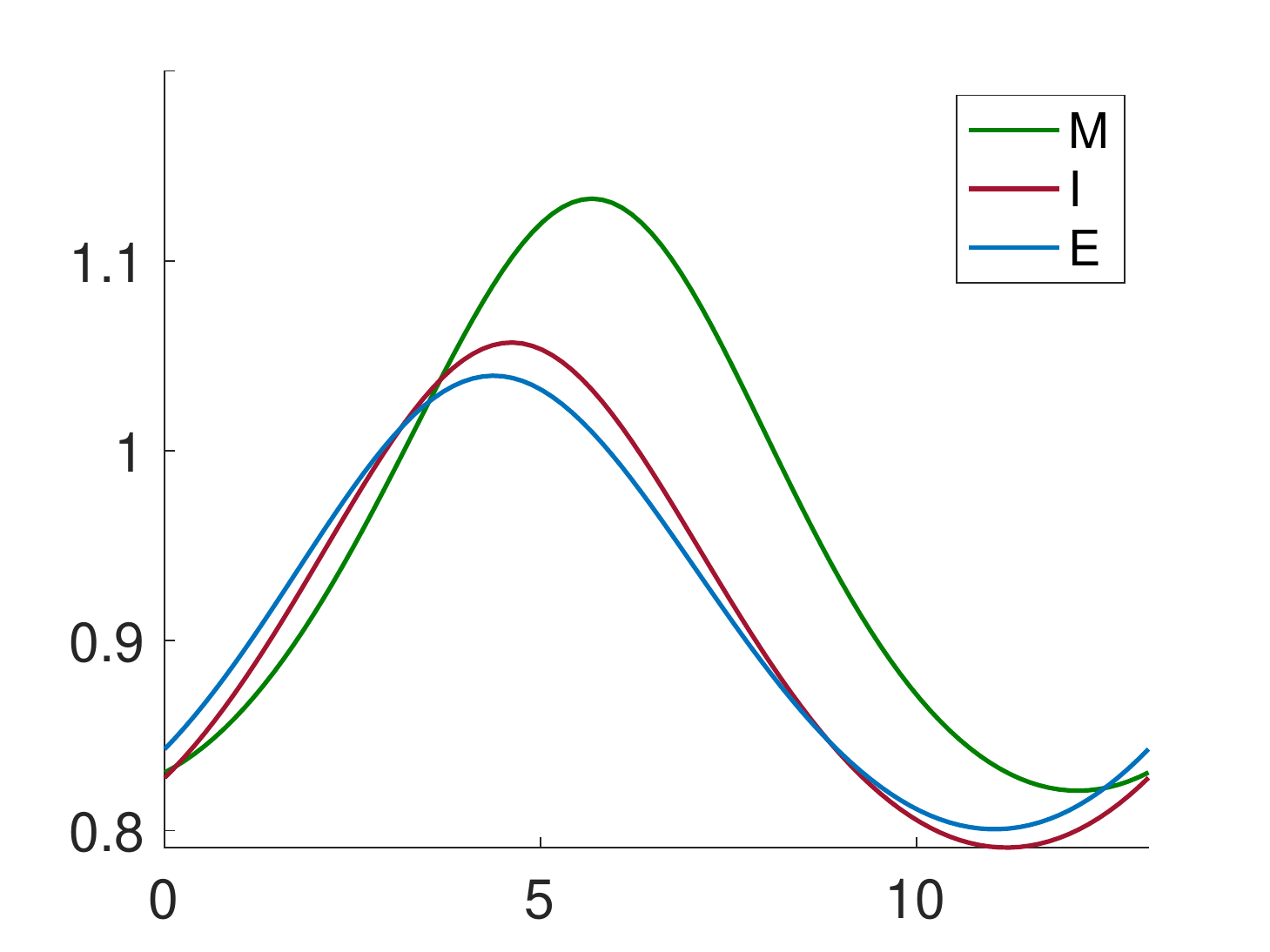}
\put(-330,120){(a) Phase Space}
\put(-155,120){(b) Stable Limit Cycle}
\put(-346,122){\rotatebox{90}{$E$}}
\put(-199,5){$M$}
\put(-20,5){$t$}
\caption{Repressible system \eqref{eq:sysonedel} with parameters
as defined in Table~\ref{table:multsspars} showing the orbits from
Figure~\ref{fig:threess_onepara_trp} at $v_M^{min}=0.01$.
(a) A projection of the phase-space dynamics into the $M$-$E$ plane in $\R^2$  with
curves formed by the points $(M(t),E(t))$, $t \in \R$ along periodic solutions $(M(t),I(t),E(t))$,
with squares denoting steady states (colour-coded according to the dimension of their unstable manifold).
(b) The three components of the stable periodic solution.
}
\label{fig:vMmin01_trp}
\end{figure}

\begin{figure}[thp!]
\centering
\hspace*{0.5em}\includegraphics[scale=0.42]{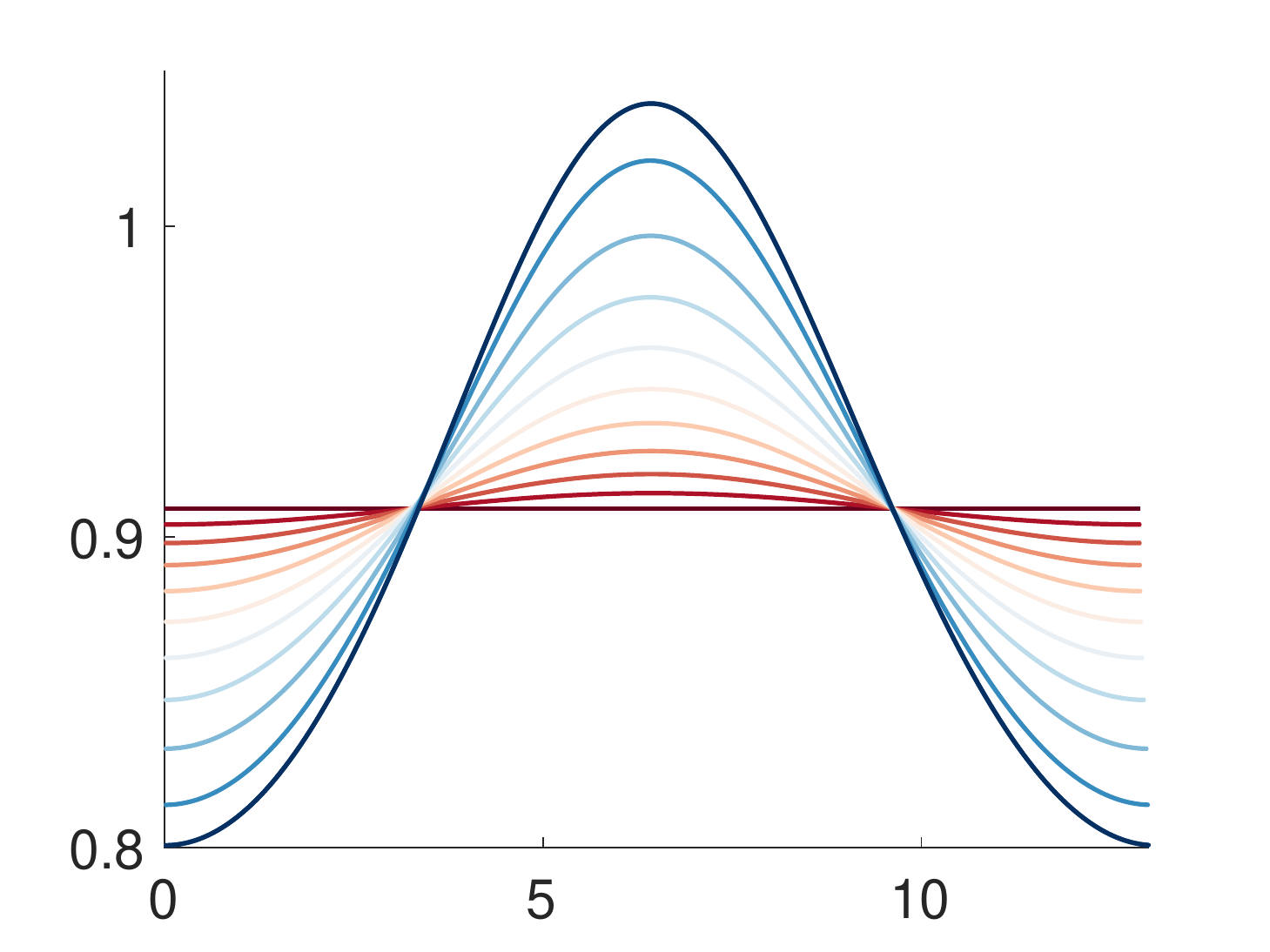}\hspace*{0.5em}\includegraphics[scale=0.42]{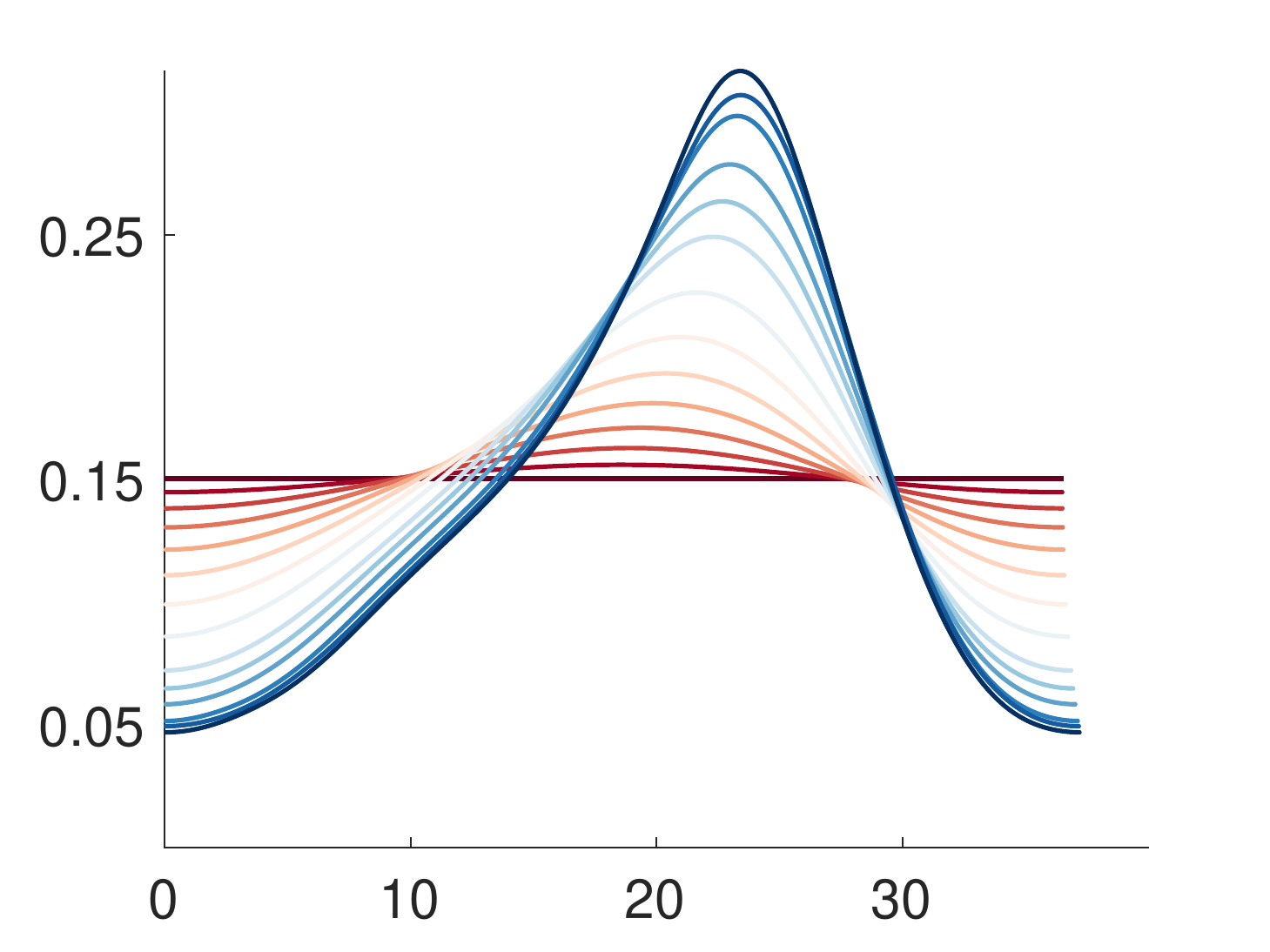}
\put(-12,4){$t$}
\put(-177,3){$t$}
\put(-318,114){\rotatebox{90}{$E$}}
\put(-152,115){\rotatebox{90}{$E$}}

\vspace*{-1ex}

\includegraphics[scale=0.42]{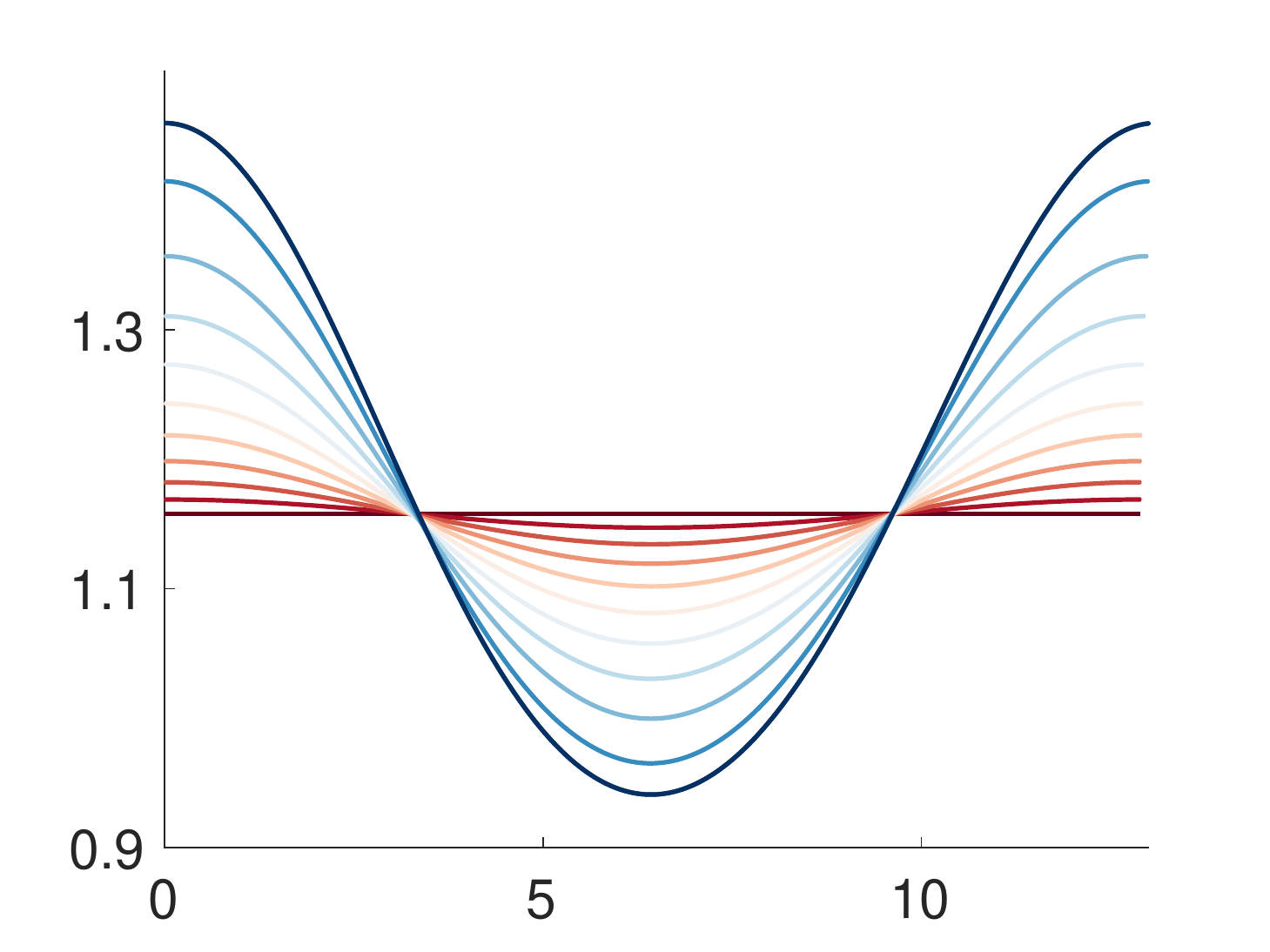}\hspace*{-1em}\includegraphics[scale=0.42]{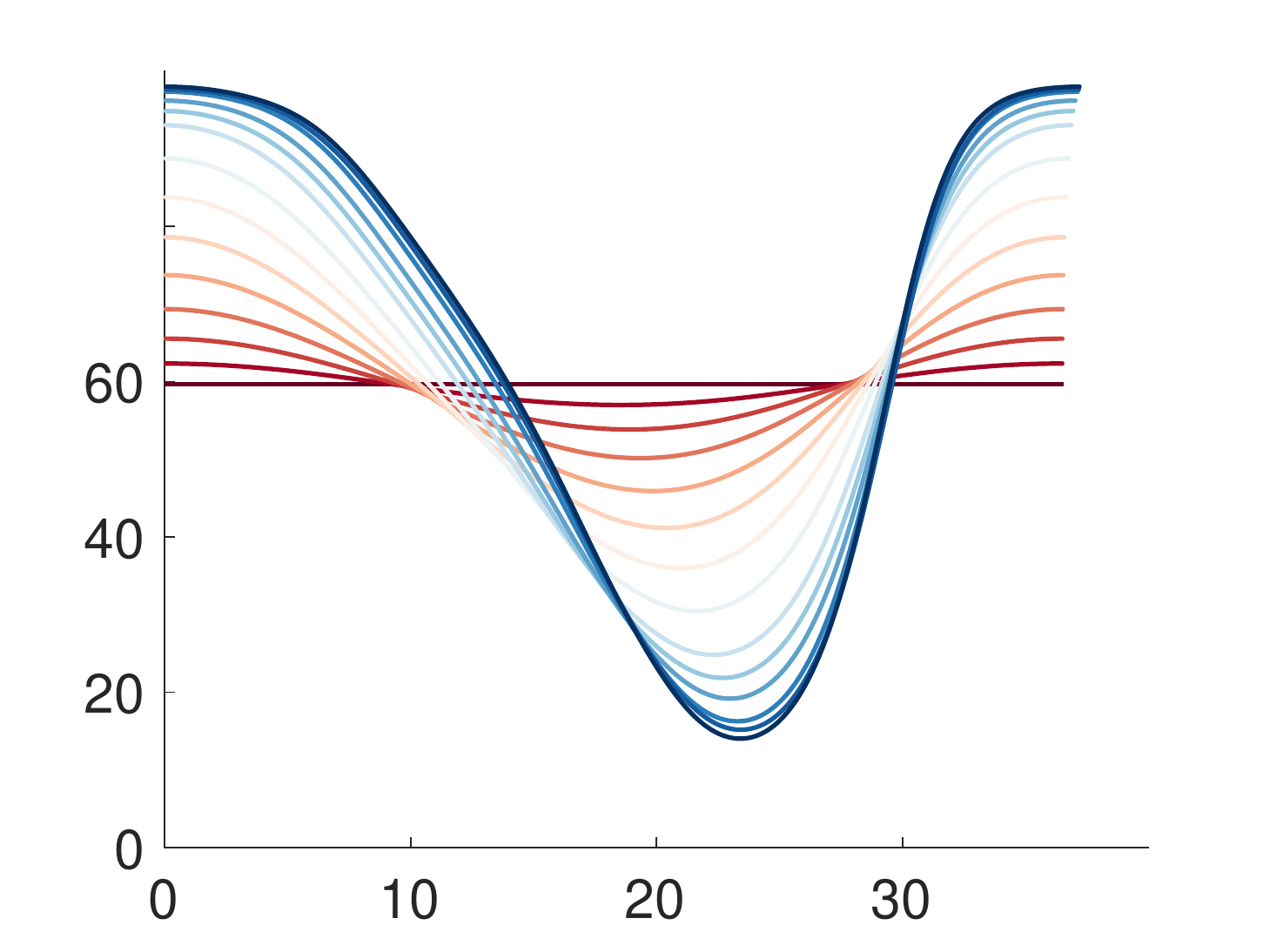}
\put(-22,4){$t$}
\put(-195,4){$t$}
\put(-332,100){\rotatebox{90}{$\tau_M(t)$}}
\put(-164,100){\rotatebox{90}{$\tau_M(t)$}}

\vspace*{-1.2ex}

\includegraphics[scale=0.42]{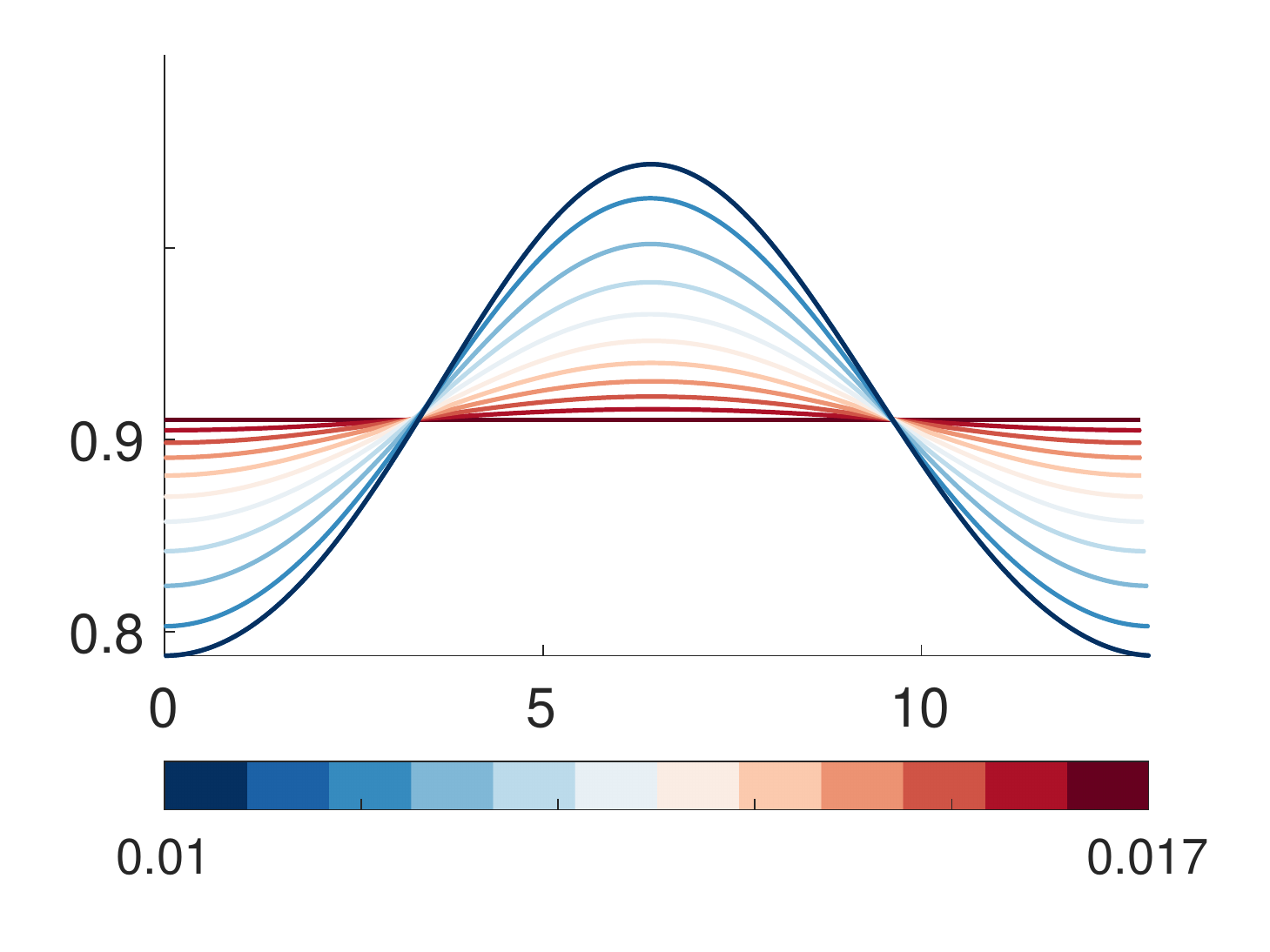}\hspace*{-1em}\includegraphics[scale=0.42]{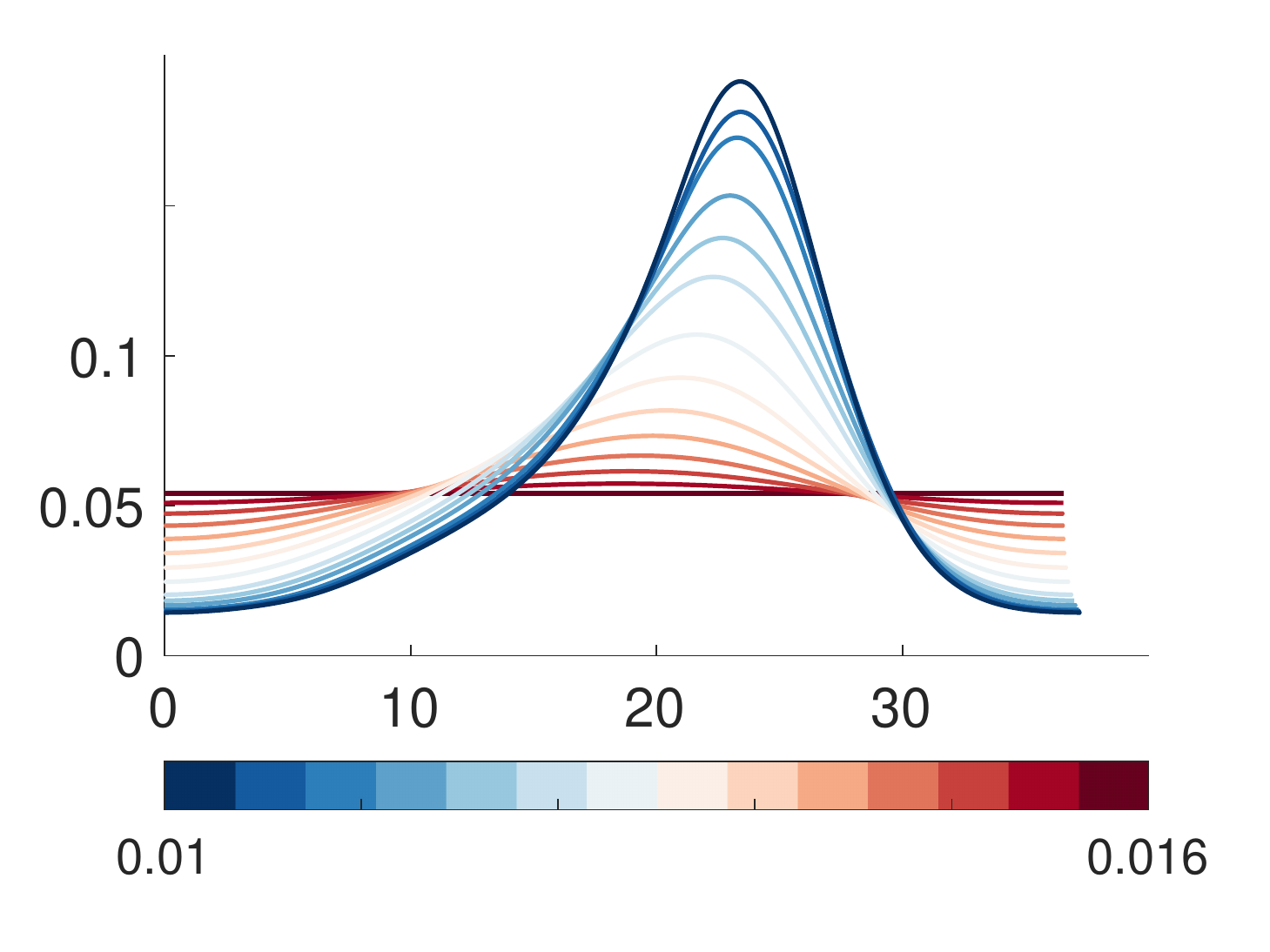}
\put(-24,33){$t$}
\put(-193,30){$t$}
\put(-334,90){\rotatebox{90}{$v_M(E(t))$}}
\put(-165,90){\rotatebox{90}{$v_M(E(t))$}}
\put(-85,11){$v_M^{min}$}
\put(-250,11){$v_M^{min}$}

\vspace*{-2.5ex}

\includegraphics[scale=0.42]{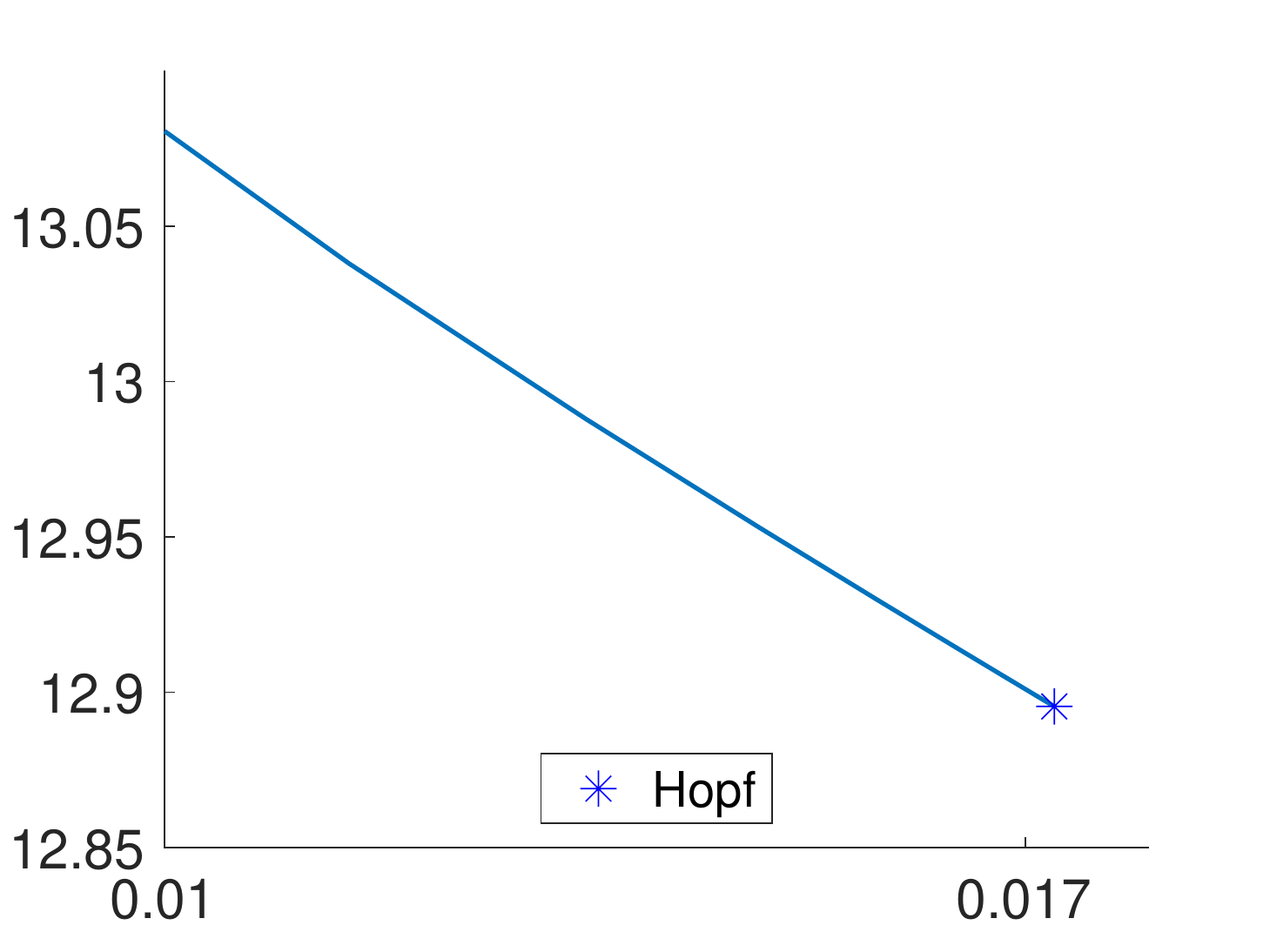}\hspace*{0.5em}\includegraphics[scale=0.42]{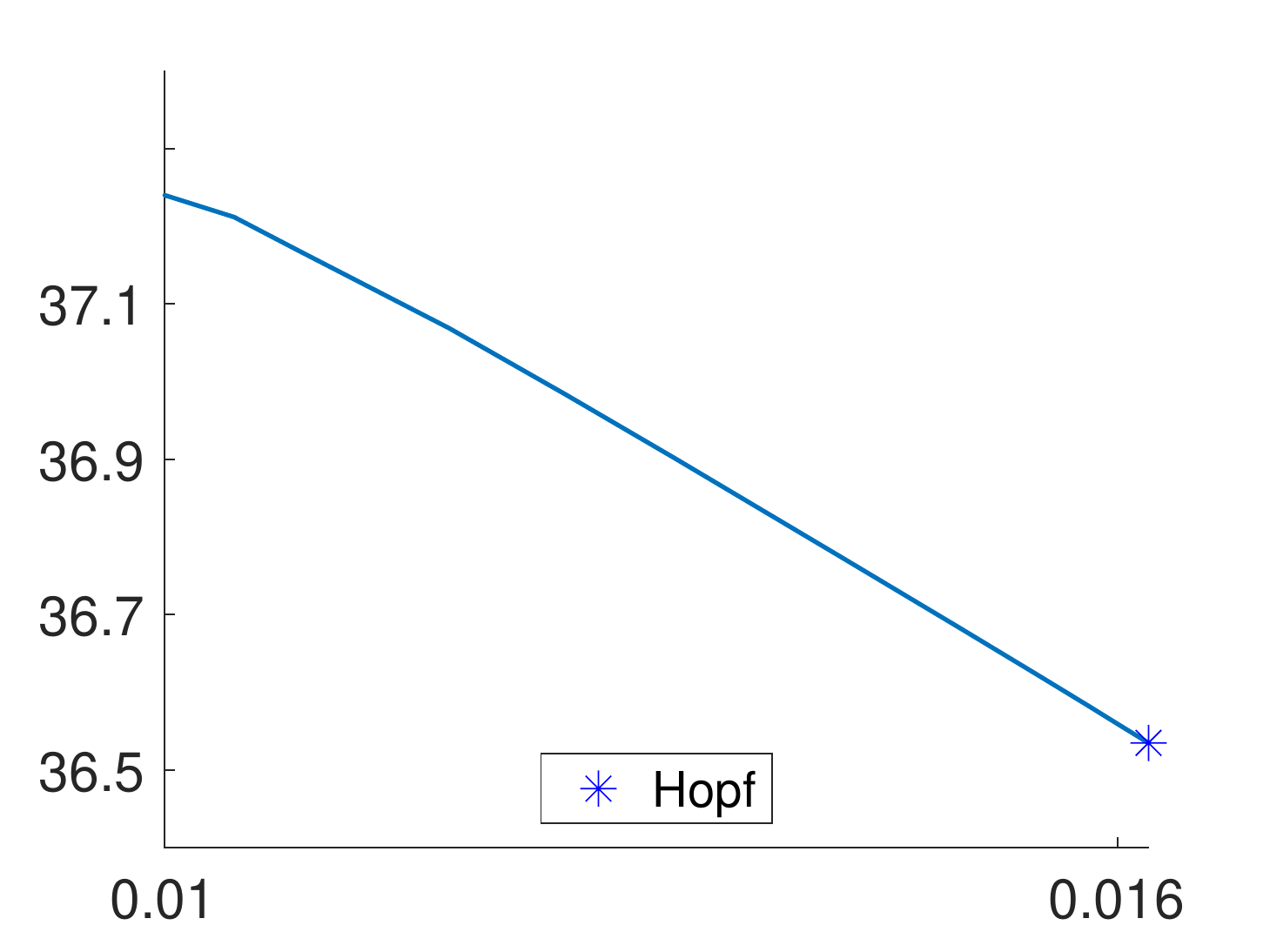}
\put(-85,4){$v_M^{min}$}
\put(-250,4){$v_M^{min}$}
\put(-319,113){\rotatebox{90}{$T$}}
\put(-154,113){\rotatebox{90}{$T$}}
\put(-255,95){Period $=T$}
\put(-90,95){Period $=T$}

\caption{Periodic orbits from the bifurcation diagramsin Figure~\ref{fig:threess_onepara_trp}. Left column: stable periodic orbits. Right column: unstable periodic orbits.
The colormap in each column indicates values of the continuation parameter $v_M^{min}$.}
\label{fig:po_trp}
\end{figure}

We begin by returning to the example from Section~\ref{ssec:equilibria} and consider
the  state-dependent delay system \eqref{eq:sysonedel} with the repressible parameter set defined in Table~\ref{table:multsspars}. The bifurcation diagram in Figure~\ref{fig:threess_onepara_trp} was
computed using DDE-BIFTOOL as detailed in Section~\ref{sec:method}, and extends the diagram previously shown in Figure~\ref{fig:gcont}(a) to
show steady state solutions, periodic orbits along with their stability, as well as Hopf and fold bifurcations. These bifurcations are listed in Table~\ref{tab:threess_onepara_trp}.

\begin{table}[thp!]
	\centering
	\begin{tabular}{| l | l | l | l |}
		\hline
		Bifurcation & Bifurcation parameter value & Unstable eigenvalues& $E^*$ value \\
		\hline
		Fold & $v_M^{min}=0.017416$ & 0 to 1 & 0.1109 \\
		\hline
		Hopf & $v_M^{min}=0.016193,$ period = 36.5348 & 1 to 3 & 0.1505 \\
		\hline\hline
		Hopf & $v_M^{min}=0.017234,$ period = 12.8954 & 0 to 2 & 0.9090 \\
		\hline
	\end{tabular}
	\caption{Steady state bifurcations seen in Figure~\ref{fig:threess_onepara_trp}.}
	\label{tab:threess_onepara_trp}
\end{table}

When $v_M^{min} = v_M^{max}$ both delays $\tau_M$ and $\tau_I$ are constant, and there can only be one steady state.
With the repressible parameter values in Table~\ref{table:multsspars} this steady state is stable. As $v_M^{min}$ is decreased there is a fold bifurcation at $v_M^{min}=0.0174$ giving rise to a pair of additional steady states, one of which is stable. Therefore there is bistability between steady states for the \emph{repressible} model with $\tau_M$ state-dependent. However, the bistability region is very narrow as at
$v_M^{min}=0.0172$ there is a Hopf bifurcation from one of the steady states giving rise to a stable periodic orbit. Consequently, for $v_M^{min}<0.0172$ there is bistability between a steady state and a limit cycle.


There is another Hopf bifurcation at $v_M^{min}=0.0162$ that gives rise to
an unstable limit cycle. Unstable periodic orbits are unlikely to be detected via numerical simulation, but it is possible to compute and follow the unstable periodic orbits in DDE-BIFTOOL
for $v_M^{min} < 0.0162$
as shown in
Figure~\ref{fig:threess_onepara_trp}.
This Hopf bifurcation
results in the coexistence of a stable steady state, two unstable steady states, a stable limit cycle and an unstable limit cycle.

Figure~\ref{fig:vMmin01_trp}(a) shows these coexisting objects at $v_M^{min} = 0.01$
in a projection of phase space onto the $M$-$E$ plane. Since DDEs define infinite dimensional dynamical systems, low dimensional projections of phase space are often used to visualise dynamics, but
the projection will, in general, not be one-to-one. Therefore some orbits may appear to intersect in the projection, even though that is impossible in phase space due to uniqueness of solutions.
As an illustration of the information that is lost in projection consider the stable limit cycle at $v_M^{min}$ which is shown over one period in Figure~\ref{fig:vMmin01_trp}(b), but is represented by the closed green curve in Figure~\ref{fig:vMmin01_trp}(a) and by just two points in
Figure~\ref{fig:threess_onepara_trp}.

Figure~\ref{fig:po_trp} shows the evolution of the stable and unstable limit cycles generated in the Hopf bifurcations as $v_M^{min}$ decreases. Illustrated are the $E$ component of the limit cycle for different
values of $v_M^{min}$ as well as the transcription velocity $v_M(E(t))$
and the delay $\tau_M$ as functions of $t$ on the periodic solution.
Comparing the two columns of Figure~\ref{fig:po_trp} we see that
the stable limit cycles remain fairly sinusoidal over the parameter range, while the unstable limit cycles have larger period than the stable ones and also larger ratios between the maximum and minimum values of the time-dependent components shown.


\subsubsection*{Homoclinic Bifurcation}

Now we change two parameter values from the previous example and consider
the repressible model \eqref{eq:sysonedel}
with parameter values in Table~\ref{table:multsspars} except for $n=15$ and $v_M^{max}=1$.
We again take $v_M^{min}$ as the bifurcation parameter.

When $v_M^{min} = v_M^{max} =1$ both delays $\tau_M$ and $\tau_I$ are constant
with $\tau_M=\tau_I=1$. In this case the constant delay repressible model has an unstable steady state and a globally stable limit cycle. This limit cycle can be found by simulating the DDE system (as described in
Section~\ref{ssec:ddesd}) using \texttt{ddesd} and then continuing the solution using DDE-BIFTOOL
(see Section~\ref{subsec:discretization}).

\begin{figure}[htp!]
\centering
\includegraphics[scale=0.6]{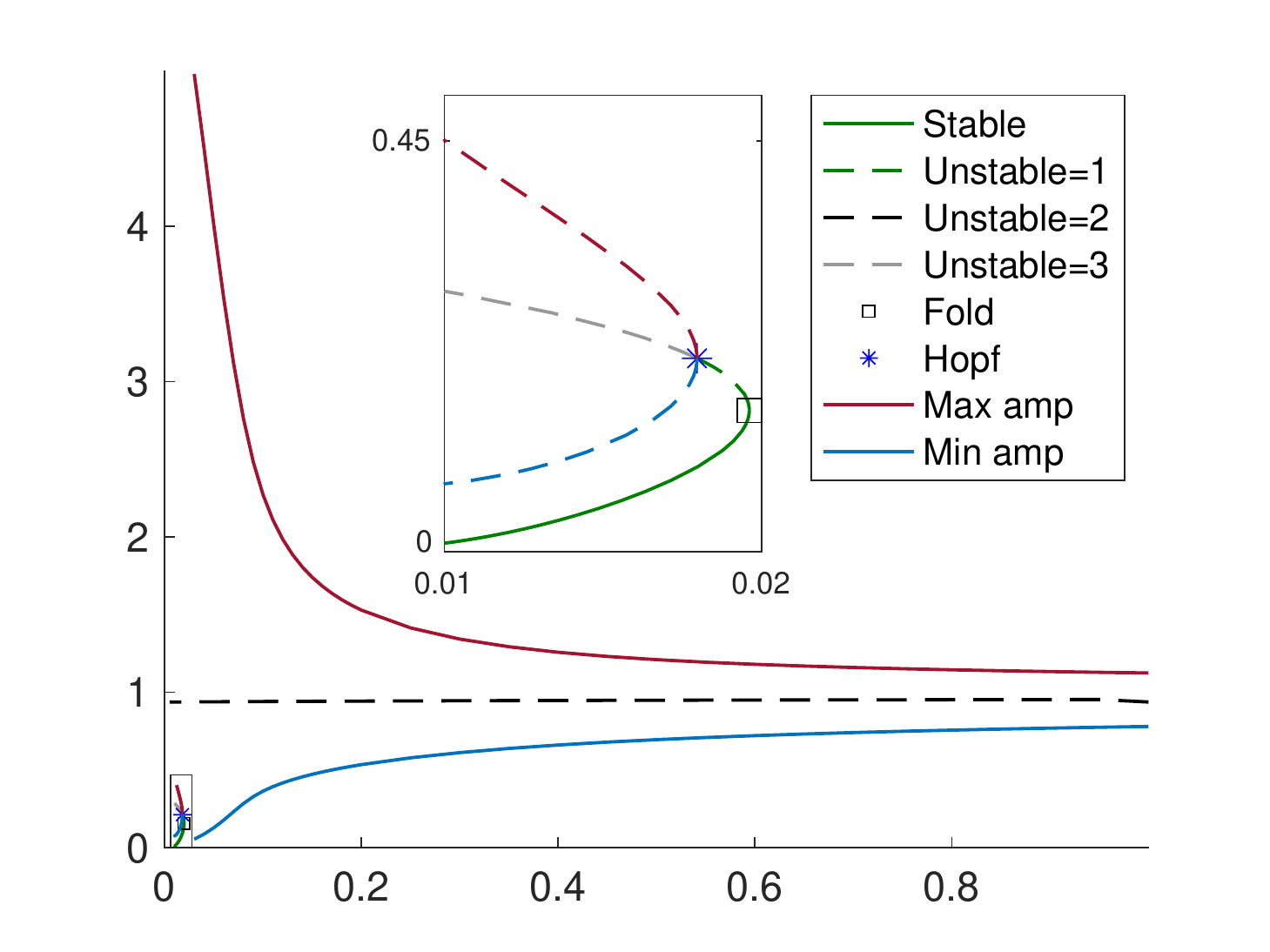}
\put(-230,165){\rotatebox{90}{$E$}}
\put(-40,10){$v_M^{min}$}
\put(-215,30){\vector(3,2){75}}
\caption{Bifurcation diagram of the model \eqref{eq:sysonedel} for a repressible system with constant $\tau_I$. Parameter values are as in Table~\ref{table:multsspars} except $n=15$, $v_M^{max}=1$ and $v_M^{min}$. Line specifications can be found in Figure~\ref{fig:threess_onepara_trp}.}
\label{fig:threess_onepara_trp_II}
\end{figure}

When the parameter value $v_M^{min}$ is decreased the delay $\tau_M$ becomes state-dependent and the amplitude of the stable periodic orbit gradually increases
as shown in the bifurcation diagram in Figure~\ref{fig:threess_onepara_trp_II}.
Bifurcations are listed in Table~\ref{tab:threess_onepara_trp_II}.
Similar to the previous example there is a fold bifurcation when $v_M^{min}$ is very small which leads to two additional steady states, one of which is stable. Thus in this example we obtain two unstable steady states which co-exist with a single stable steady state. There is also an unstable limit cycle
generated by a Hopf bifurcation, also similar to the previous example. We are not able to find stable limit-cycles that co-exist with the stable steady state.

\begin{table}[thp!]
	\centering
	\begin{tabular}{| l | l | l | l |}
		\hline
		Bifurcation & Bifurcation parameter value & Unstable eigenvalues& $E^*$ value \\
		\hline
		Fold & $v_M^{min}=0.019610$ & 0 to 1 & $0.1546$ \\
		\hline
		Hopf & $v_M^{min}=0.017963$, period = 31.4290 & 1 to 3 & $0.2116$ \\
		\hline
	\end{tabular}
	\caption{Steady state bifurcations seen in Figure~\ref{fig:threess_onepara_trp_II}.}
	\label{tab:threess_onepara_trp_II}
\end{table}


\begin{figure}[htp!]
\hspace*{-0.5em}\includegraphics[scale=0.42]{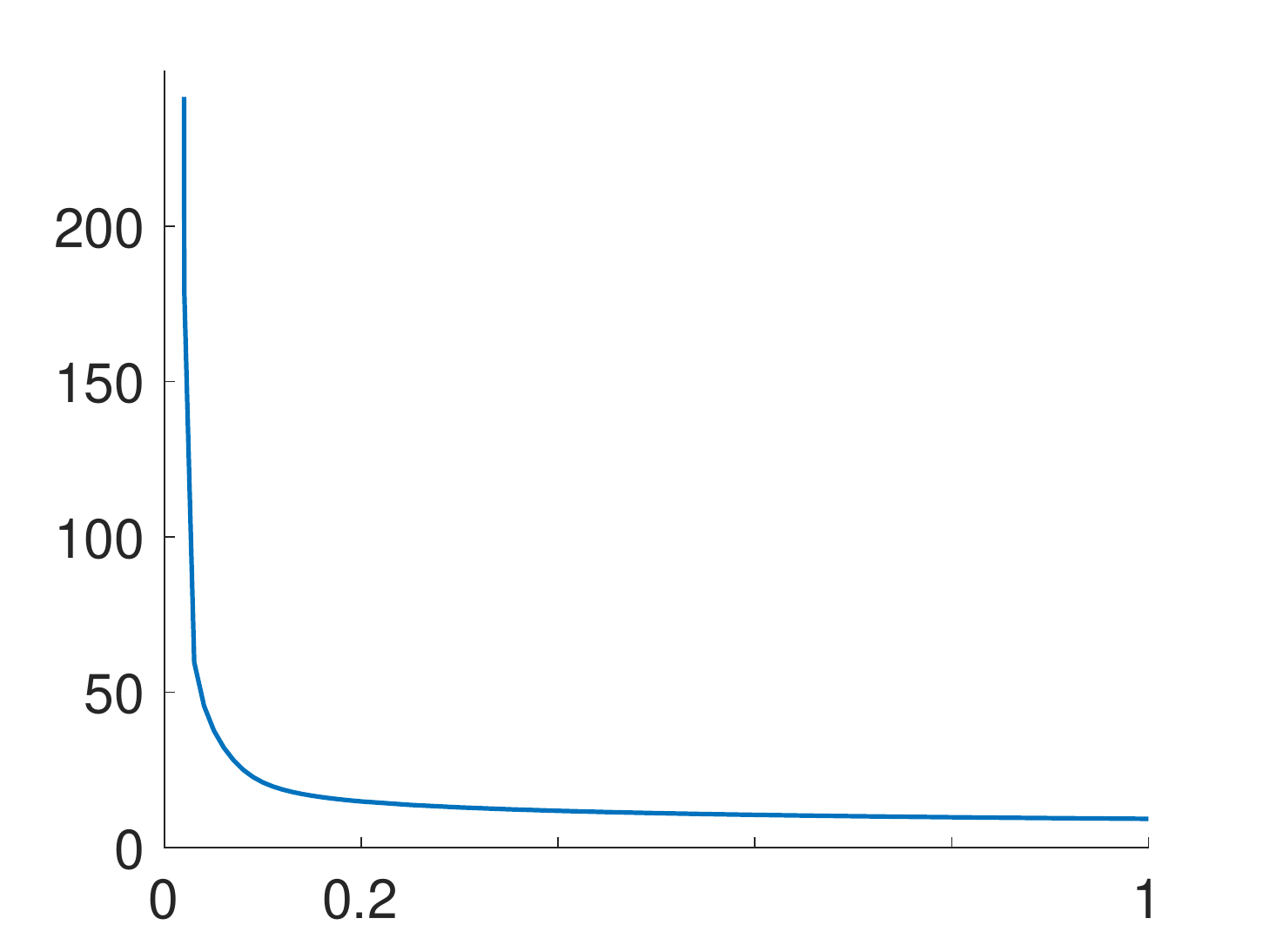}\hspace*{-0.5em}\includegraphics[scale=0.42]{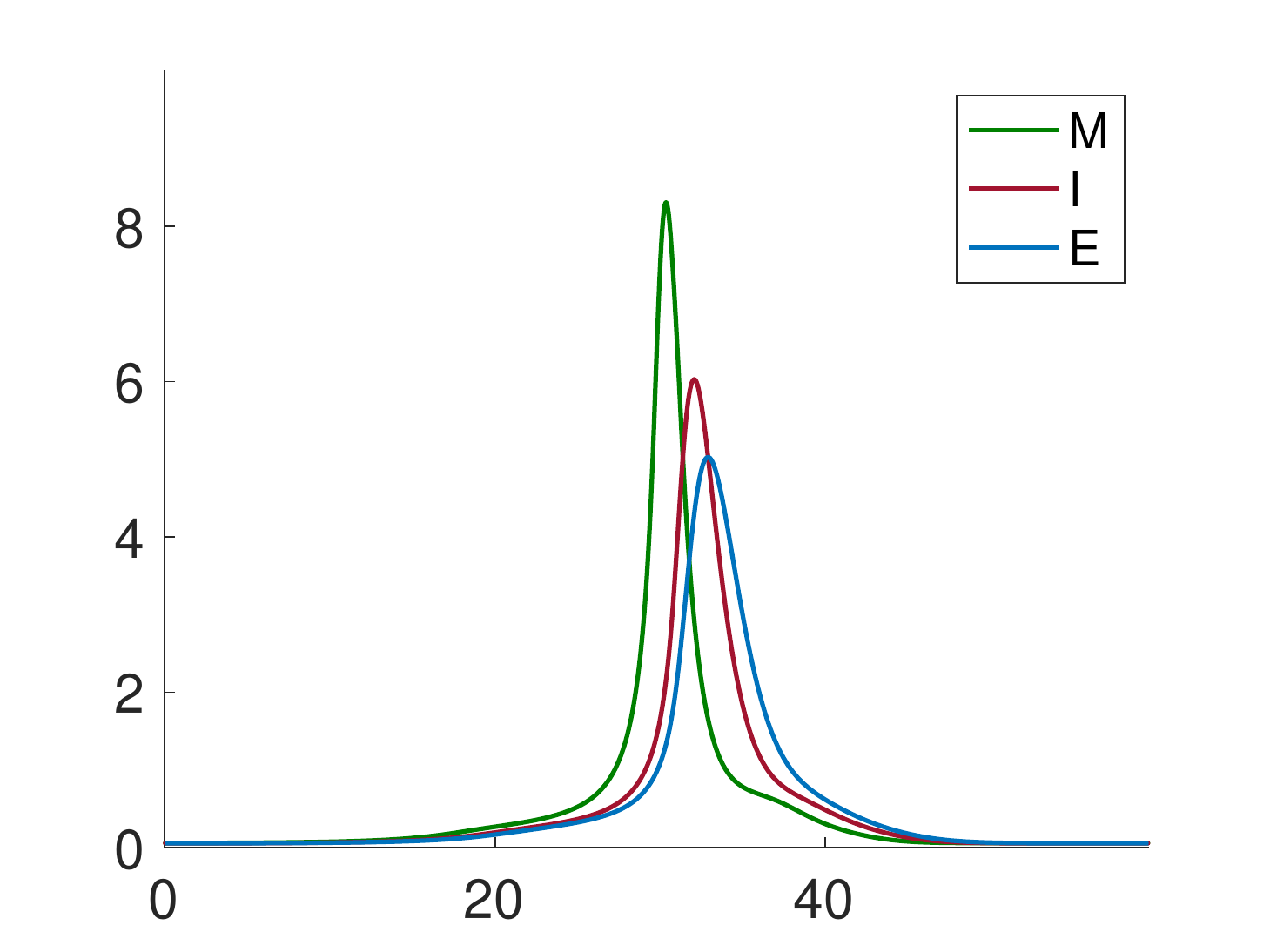}
	\put(-20,5){$t$}
	\put(-230,4){$v_M^{min}$}
	\put(-335,113){\rotatebox{90}{$T$}}
	\put(-295,115){(a) Period $=T$}
	\put(-145,115){(b) Stable Limit Cycle}
	\put(-145,105){$v_M^{min}=0.03$}
	
\hspace*{-0.5em}\includegraphics[scale=0.42]{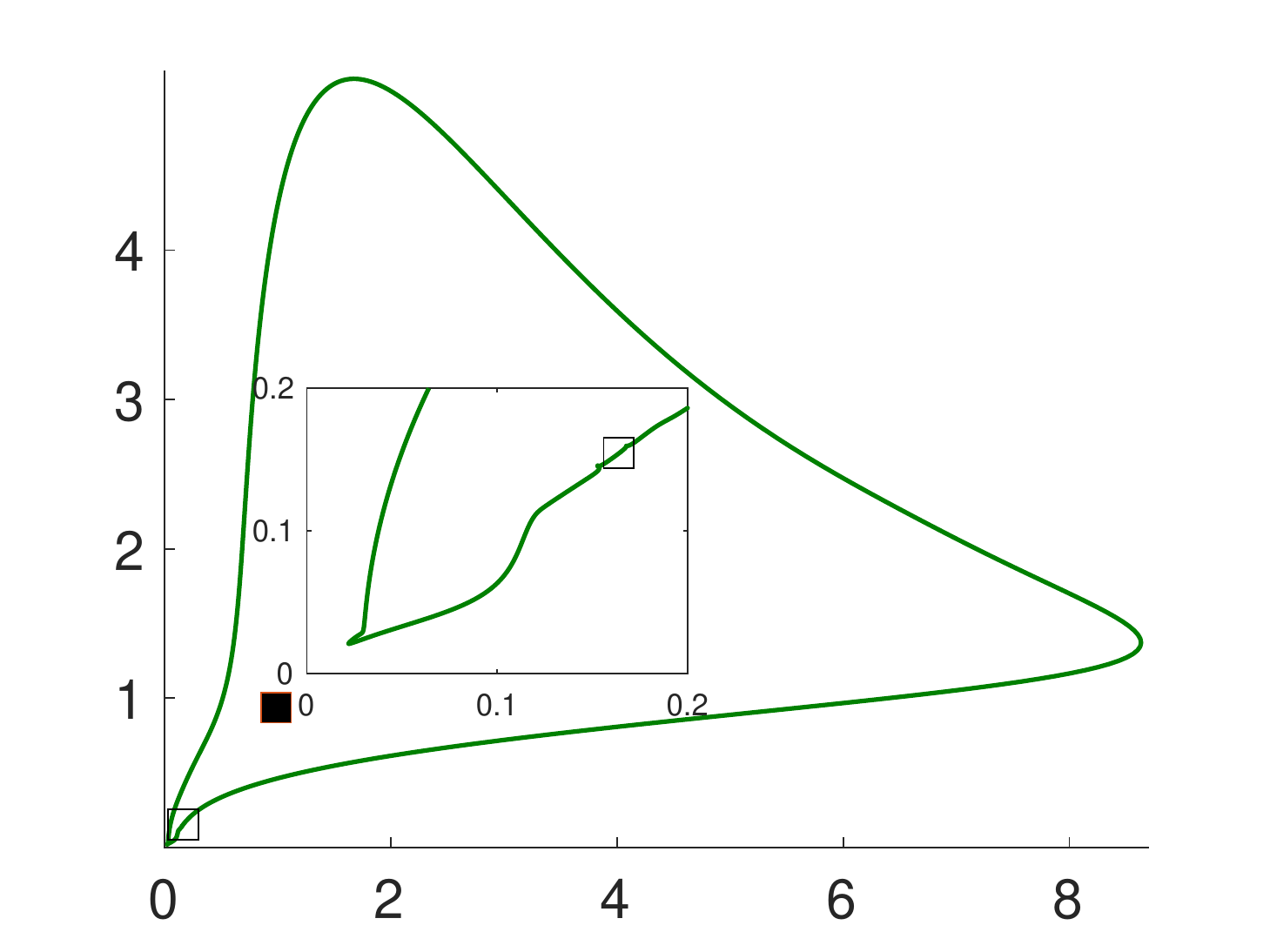}\hspace*{-0.5em}\includegraphics[scale=0.42]{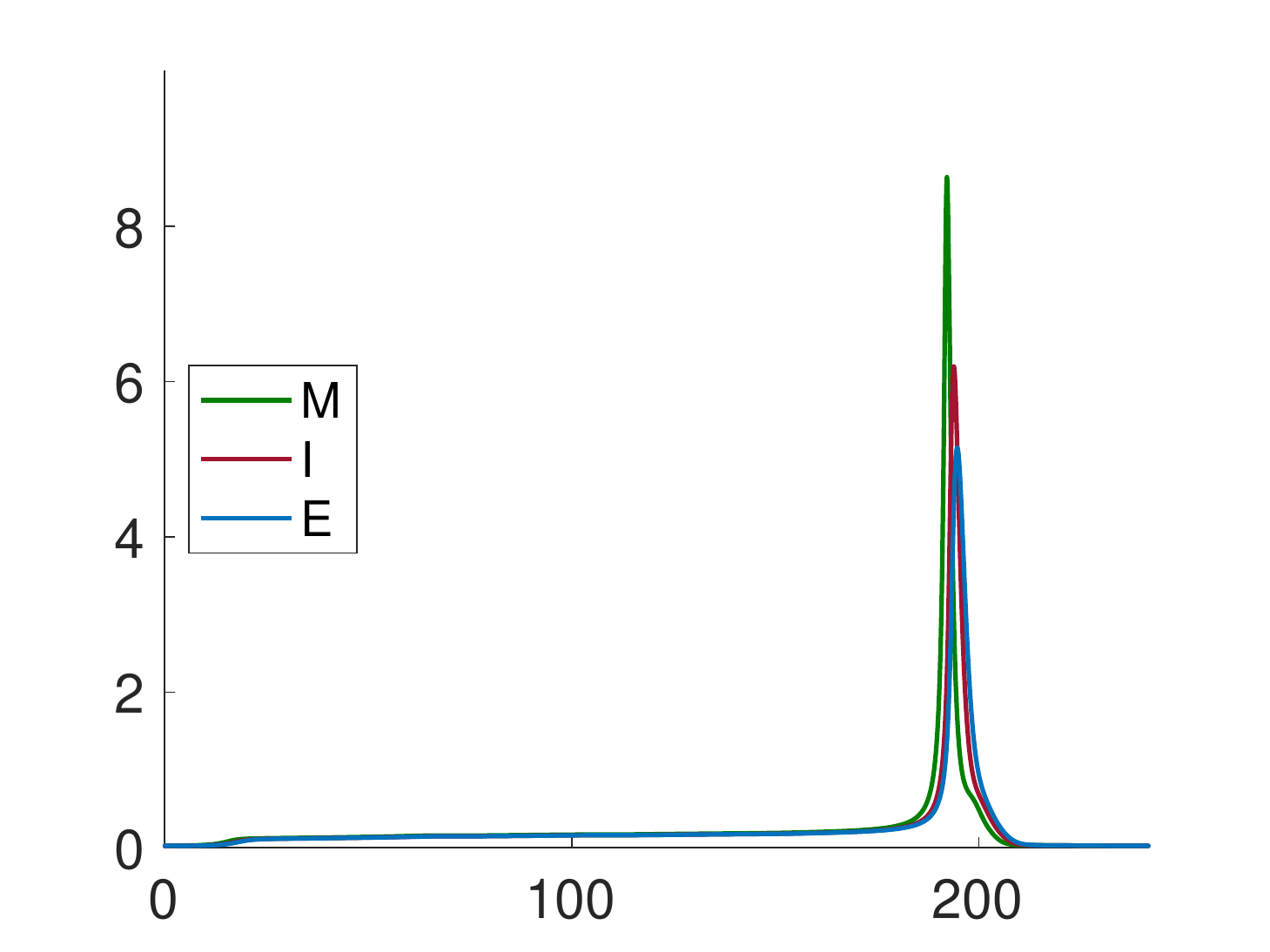}	
	\put(-20,5){$t$}
	\put(-195,5){$M$}
	\put(-145,115){(d) Stable Limit Cycle}
	\put(-275,115){(c) Phase Space}
	\put(-130,105){$v_M^{min}=0.0197$}
	\put(-260,105){$v_M^{min}=0.0197$}
	\put(-335,115){\rotatebox{90}{$E$}}

\caption{(a) The period of the stable periodic orbit (shown in Figure~\ref{fig:threess_onepara_trp_II}) grows dramatically as $v_M^{min}$ decreases. (b) Periodic orbit at $v_M^{min} = 0.03$. (c) Projection of the phase space dynamics into the $M$-$E$ plane at $v_M^{min} = 0.0197$. The open square marks the steady state $(M^*, E^*)$ at the fold bifurcation.
(d) Periodic orbit at $v_M^{min} = 0.0197$. }
	\label{fig:vMmin02_trp_II}
\end{figure}
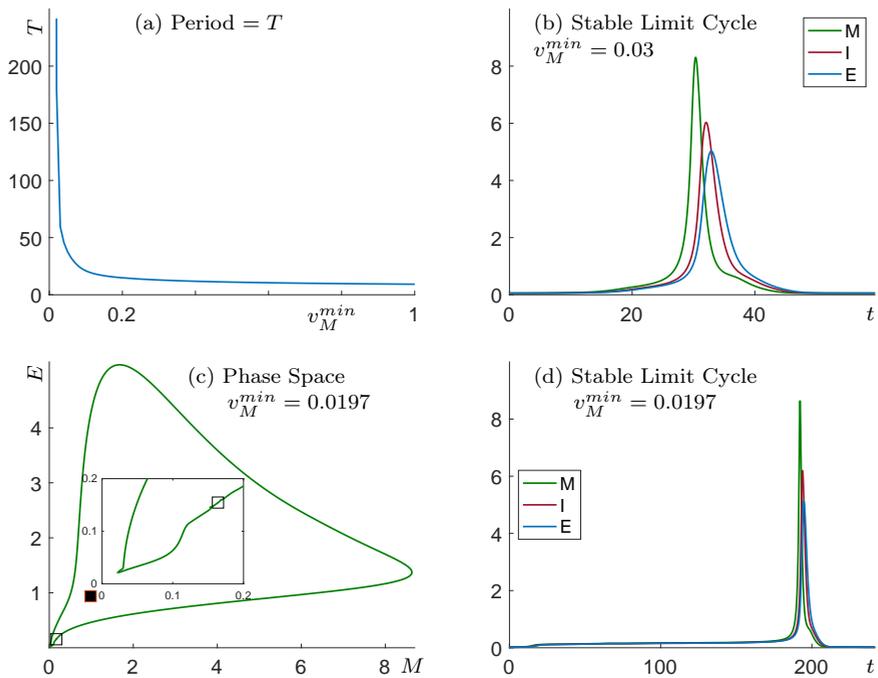

This example differs from the previous example in the behaviour of the stable limit cycle.
We are able to find the limit cycle only for $v_M^{min}\geq0.0197$ with the period increasing dramatically as $v_M^{min}\to0.0197$ as shown in Figure~\ref{fig:vMmin02_trp_II}(a), which suggests that a homoclinic bifurcation may occur.
For $v_M^{min} = 0.03$ the stable limit cycle is shown in Figure~\ref{fig:vMmin02_trp_II}(b), and appears
to behave like a relaxation oscillator with $(M(t),I(t),E(t))\approx(0,0,0)$ for much of the time, with
one burst of production each period.
This periodic solution may have an interesting biological interpretation (see Section~\ref{sec:disc}).
Namely, the burst of transcription is followed in short succession by burst of protein production, and this protein represses the initiation of mRNA transcription for a majority of the period. Only when this repression is released, a burst of transcription follows.  

The last limit cycle that we are able to compute for $v_M^{min}=0.0197$ is shown in Figures~\ref{fig:vMmin02_trp_II}(c) and (d). If there is a homoclinic orbit then the limit cycle  would have to approach a saddle-like steady state. However, in the phase space plot in panel (c), the periodic orbit is always far from the only steady state (denoted by the solid square) that exists for $v_M^{min}=0.0197$. On the other hand, we do observe that the orbit does pass through the region of phase space containing the `ghost' of the saddle steady state destroyed in the fold bifurcation at $v_M^{min}=0.01961$. Panel (d) also shows the solution close to this ghost steady state for $t\in(20,180)$ which is for most of the period.

In this example it seems that a homoclinic bifurcation occurs
very close to the fold bifurcation where the steady state with saddle stability is destroyed.
This suggests that our parameter set is close to a higher co-dimension bifurcation where the homoclinic and saddle bifurcations coincide. We investigate this further in the next example.


\subsubsection*{Zero-Hopf Bifurcation and 3DL Transition}

\begin{figure}[thp!] \includegraphics[scale=0.58]{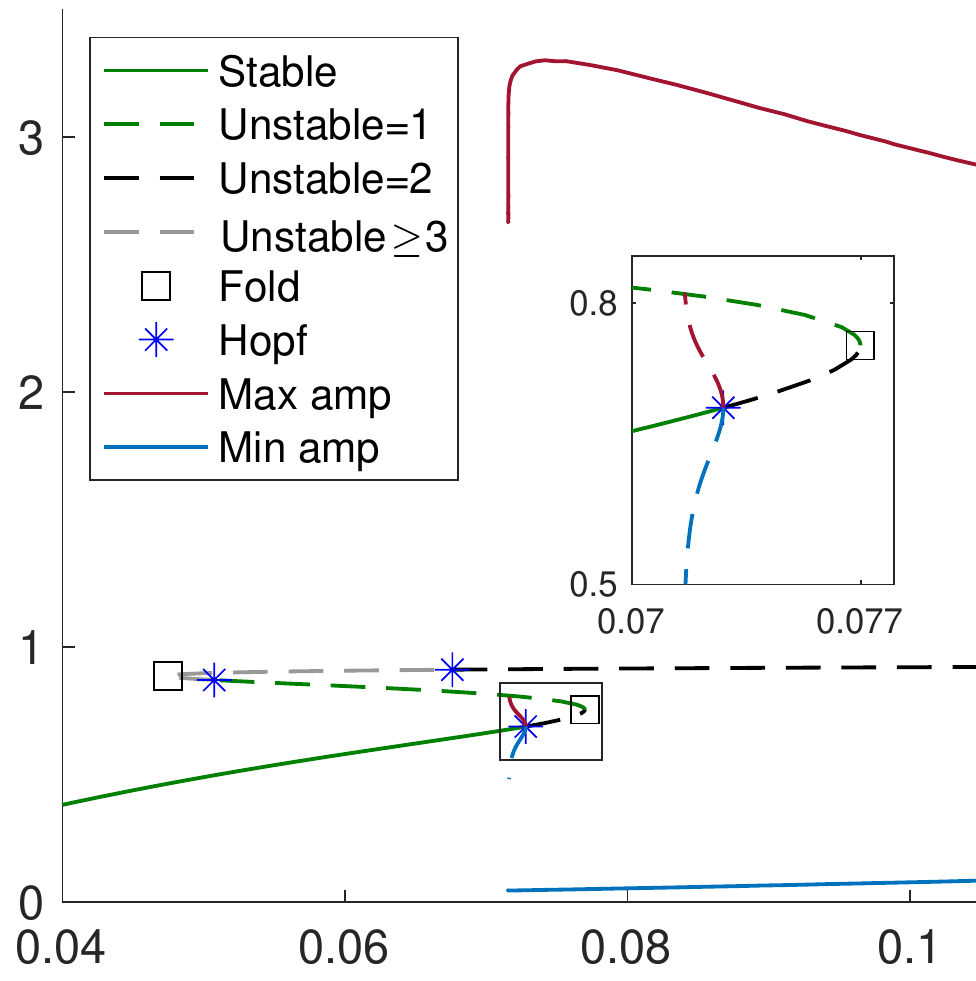}\hspace*{0.5em}\includegraphics[scale=0.58]{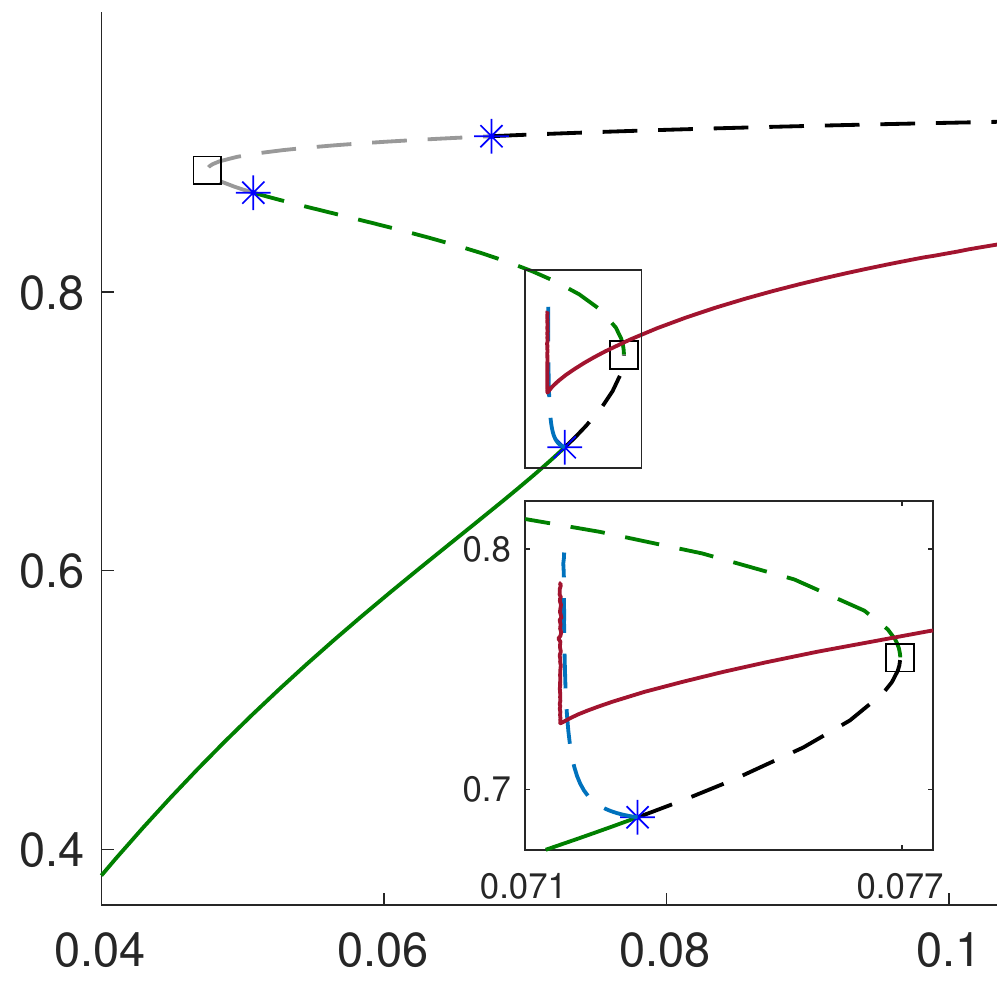}	
   \put(-205,154){(a)}
    \put(-145,154){(b)}
   \put(-159,154){\rotatebox{90}{$E$}}
	\put(-40,5){$v_M^{min}$}
	\put(-332,155){\rotatebox{90}{$E$}}
	\put(-215,4){$v_M^{min}$}
    \put(-60,93){\vector(1,-1){10}}
    \put(-234,42){\vector(1,1){22}}
	\caption{Bifurcation diagram of the model \eqref{eq:sysonedel} for a repressible system with constant $\tau_I$. Parameter values are as in Table~\ref{table:multsspars} except
$m=n=15$, $v_M^{max}=1$ and $v_M^{min}$. (a) Line specifications can be found in Figure~\ref{fig:threess_onepara_trp}. (b) The same as (a) except the stable and unstable periodic orbits
are represented as a solid red and blue dashed curve respectively using the 1-norm \eqref{eq:1norm}
of the periodic solution.}
\label{fig:threess_onepara_trp_III}
\end{figure}
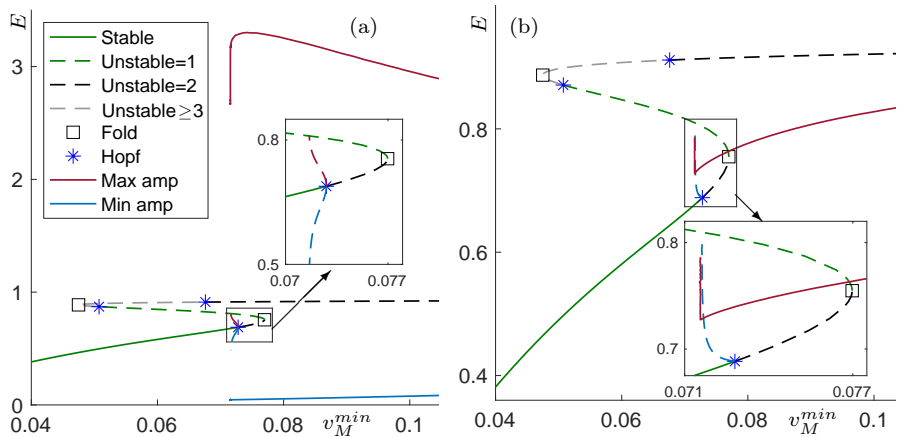

Next we change a single parameter value from the example shown in
Figures~\ref{fig:threess_onepara_trp_II}-\ref{fig:vMmin02_trp_II}
to consider the model \eqref{eq:sysonedel} in the repressible case with the Hill coefficient in the transcription velocity $m=15$ (in both the previous examples we took $m=3$). All the other parameter values remain the same as in the previous example, so $n=15$, $v_M^{max}=1$ and the rest of the parameters as
defined Table~\ref{table:multsspars}.
The resulting bifurcation diagram is shown in Figure~\ref{fig:threess_onepara_trp_III}, and the bifurcations are listed in Table~\ref{table_threess_III_trp}.

\begin{table}[htp!]
	\centering
	\begin{tabular}{| l | l | l | l |}
		\hline
		Bifurcation & Bifurcation parameter value & Unstable eigenvalues& $E^*$ value \\
		\hline
		Hopf & $v_M^{min}=0.072792,$ period = 54.7357  & 0 to 2 & 0.6885 \\
		\hline
		Fold & $v_M^{min}=0.076976$ & 2 to 1 & 0.7548 \\
		\hline
		Hopf & $v_M^{min}=0.050745$ & 1 to 3 & 0.8711 \\
		\hline
		Fold & $v_M^{min}=0.047485$ & 3 to 4 & 0.8871 \\
		\hline
		Hopf & $v_M^{min}=0.067602$ & 4 to 2 & 0.9116 \\
		\hline
	\end{tabular}
	\caption{Bifurcation information associated with Figure~\ref{fig:threess_onepara_trp_III}. }
	\label{table_threess_III_trp}
\end{table}

There are several significant differences between the bifurcation diagram in
Figure~\ref{fig:threess_onepara_trp_III} and the previous case in Figure~\ref{fig:threess_onepara_trp_II}.
Considering first just the steady states, we see that there is an additional fold bifurcation
and that all the steady states now  lie on a single continuous branch of steady states with two fold bifurcations.
As in the previous example there is a single segment of stable steady states, but it loses stability in a subcritical Hopf bifurcation at $v_M^{min}=0.072792$ whereas in the previous example the stable steady state was destroyed in a fold bifurcation. Comparing
the insets in Figures~\ref{fig:threess_onepara_trp_II} and~\ref{fig:threess_onepara_trp_III} we see that the Hopf and the fold bifurcation both occur in each example but, importantly, their order on the branch is reversed.
Therefore there must be an intermediate value of $m\in(3,15)$  where the two bifurcations will coincide in a so-called zero-Hopf or fold-Hopf bifurcation. The codimension-two zero-Hopf bifurcation is known to generate homoclinic orbits and bifurcations \citep{Kuznetsov2004}, which is further evidence for the existence of homoclinic orbits in the state-dependent delay operon model \eqref{eq:sysonedel}.

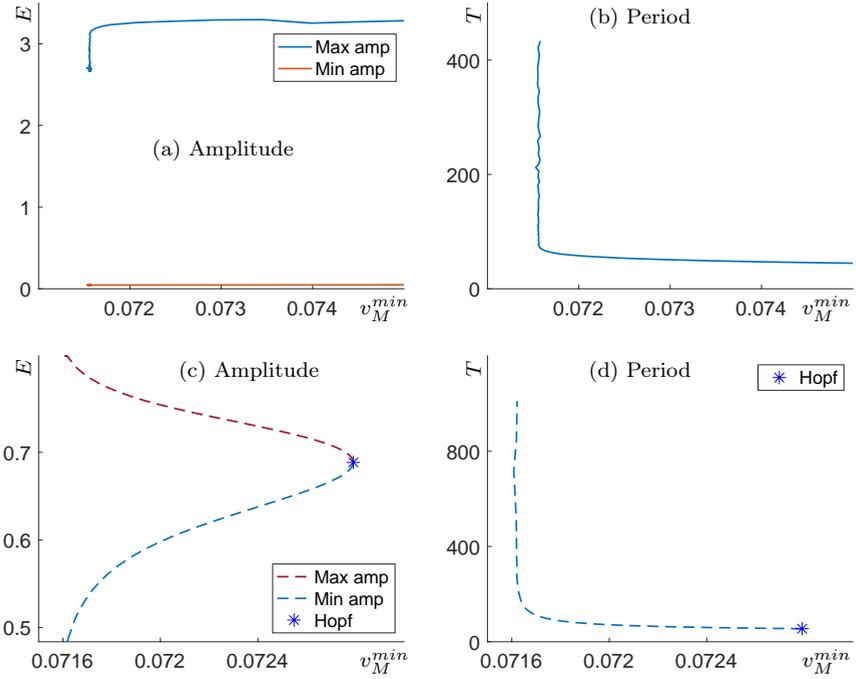
\begin{figure}[thp!] \includegraphics[scale=0.42]{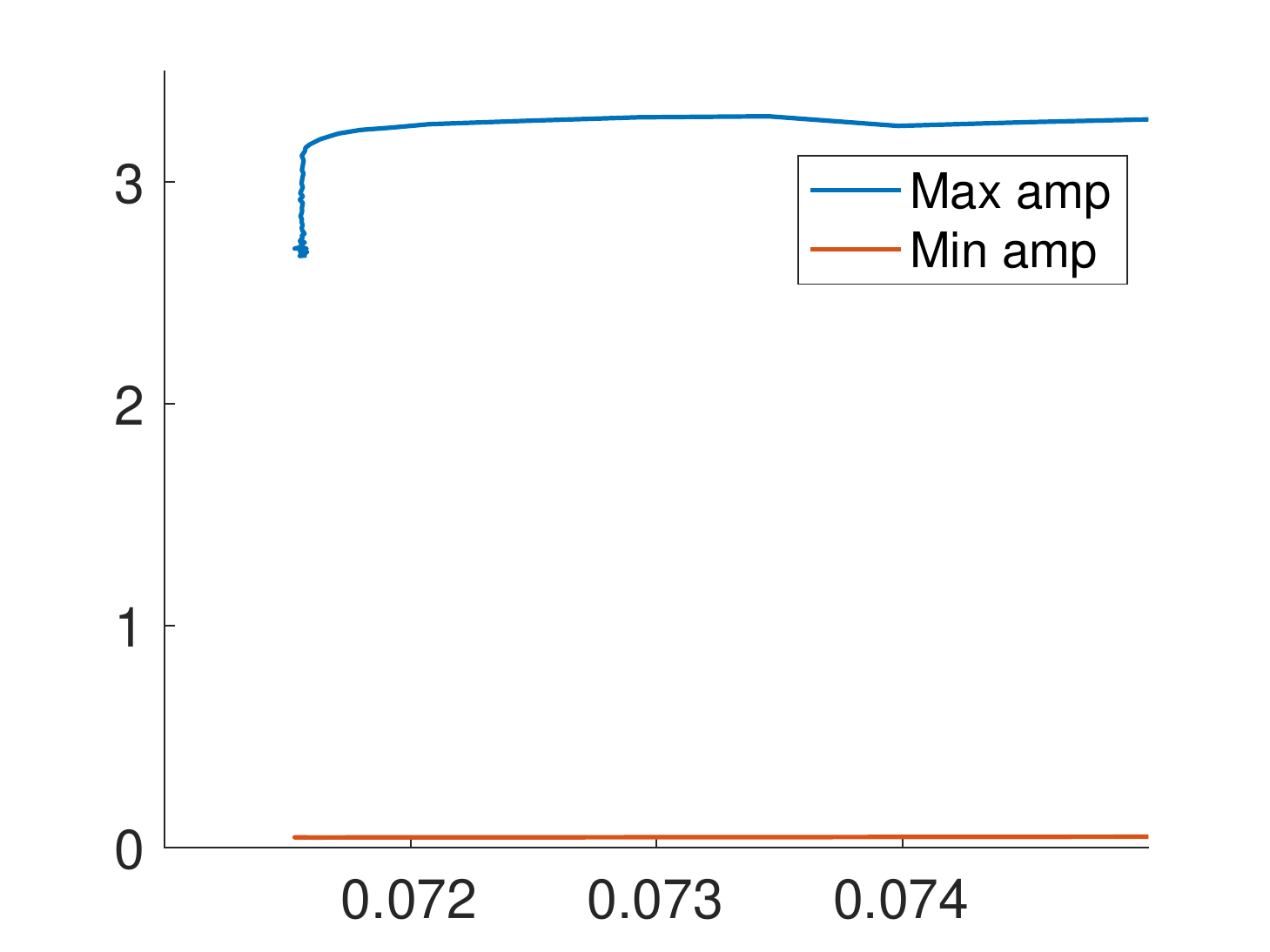}\hspace*{-1em}\includegraphics[scale=0.42]{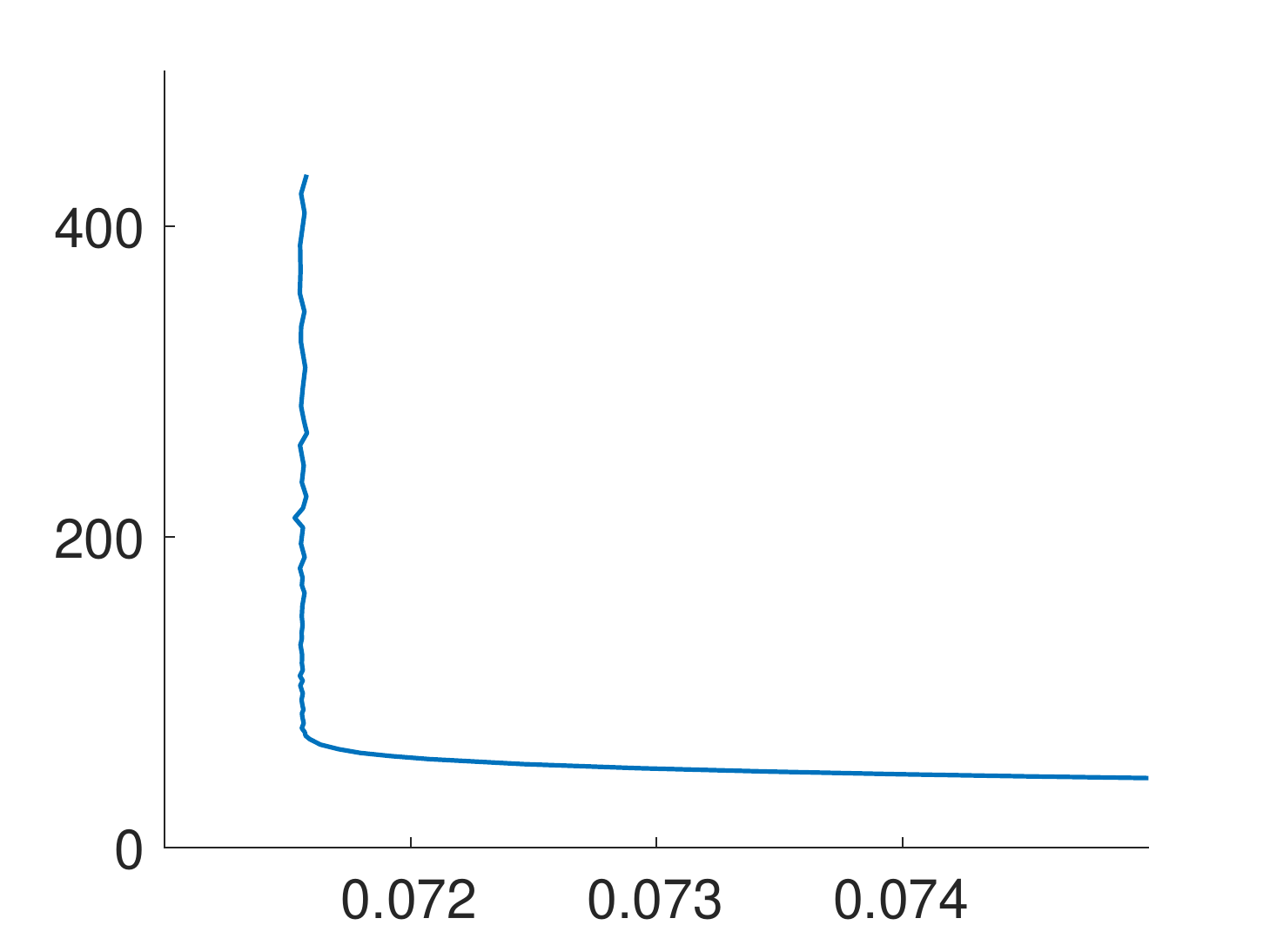}
\put(-280,65){(a) Amplitude}
\put(-116,115){(b) Period}
\put(-36,5){$v_M^{min}$}
\put(-202,5){$v_M^{min}$}
\put(-331,115){\rotatebox{90}{$E$}}
\put(-162,115){\rotatebox{90}{$T$}}

\includegraphics[scale=0.42]{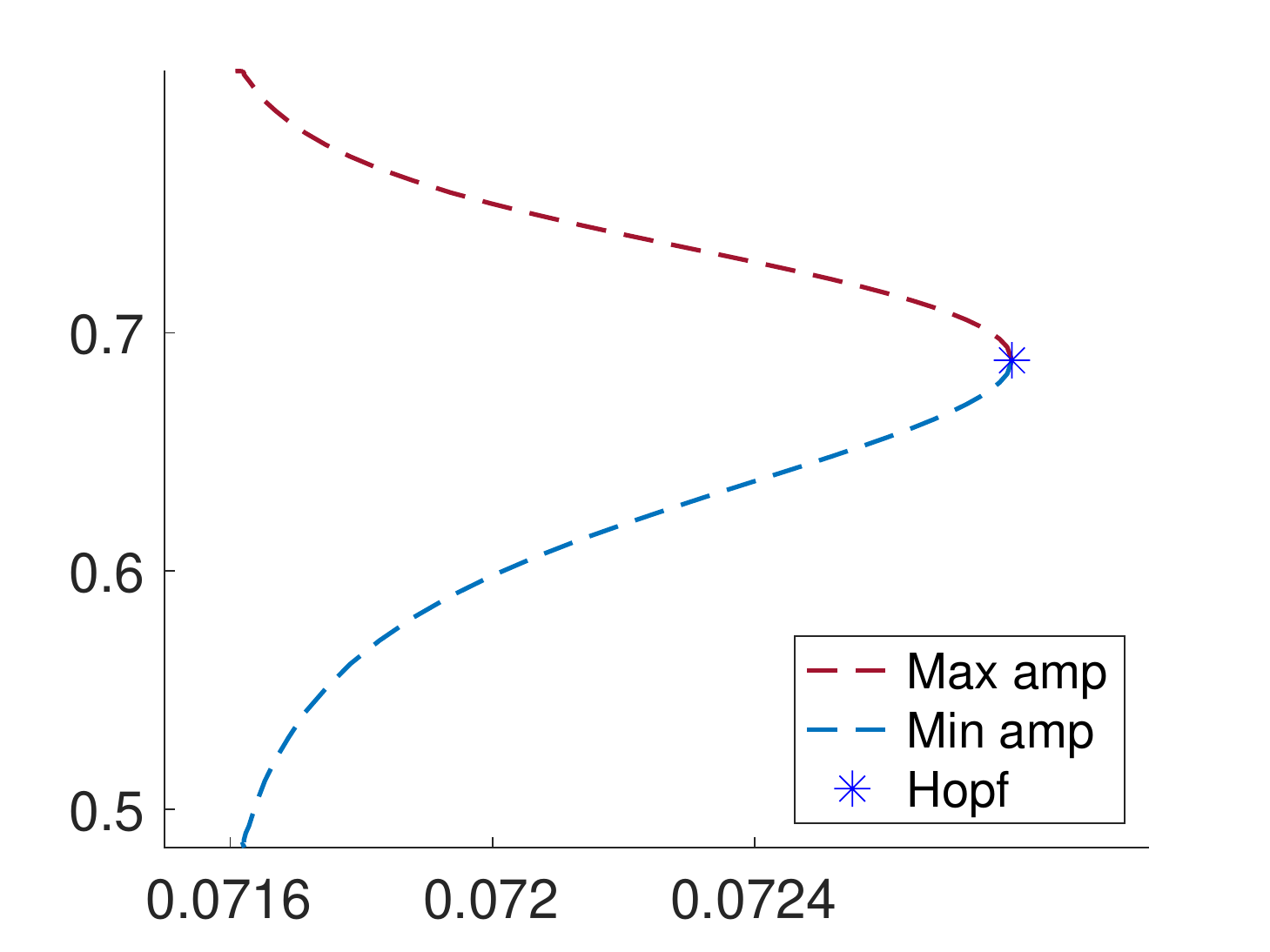}\hspace*{-1em}\includegraphics[scale=0.42]{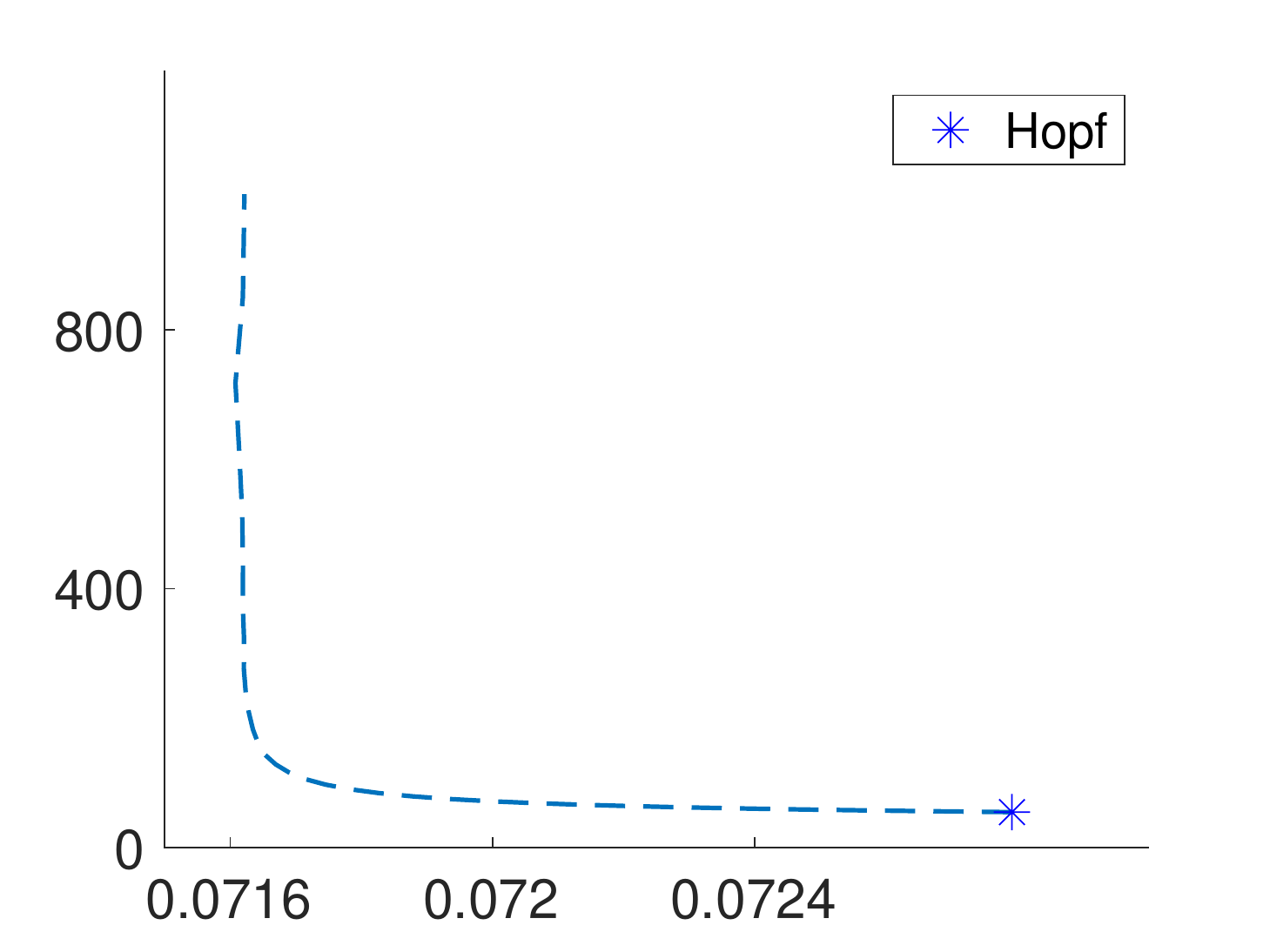}
\put(-270,115){(c) Amplitude}
\put(-116,115){(d) Period}
\put(-36,5){$v_M^{min}$}
\put(-202,5){$v_M^{min}$}
\put(-331,115){\rotatebox{90}{$E$}}
\put(-162,115){\rotatebox{90}{$T$}}
	
\caption{(a) \& (b) Stable, and, (c) \& (d) unstable branches of periodic orbits from the bifurcation diagram in Figure~\ref{fig:threess_onepara_trp_III}.}
\label{fig:perorbs_onepara_trp_III}
\end{figure}

Consideration of the periodic orbits shown in Figure~\ref{fig:threess_onepara_trp_III} provides further evidence supporting existence
of homoclinic orbits. While we could imagine that the two branches of periodic orbits shown in
Figure~\ref{fig:threess_onepara_trp_III}(a) might join up to form one continuous branch that is not what happens.  A different representation of the periodic solutions on the bifurcation diagram is appropriate when considering periodic orbits close to homoclinic. In Figure~\ref{fig:threess_onepara_trp_III}(a) the periodic orbits are represented by two curves, representing their amplitude. Figure~\ref{fig:threess_onepara_trp_III}(b) shows exactly the same bifurcation diagram, except that a periodic orbit of period $T$ is now represented by the 1-norm of its $E(t)$ component:
\begin{equation} \label{eq:1norm}
\|\cdot\|_1=\frac{1}{T}\int_{0}^TE(t)dt
\end{equation}
This representation of periodic orbits is useful because the 1-norm of a periodic orbit approaches  the value of $E^*$ as a periodic orbit approaches either a Hopf bifurcation or a homoclinic bifurcation at the steady state $(M^*,I^*,E^*)$. In Figure~\ref{fig:threess_onepara_trp_III}(b) the stable and unstable periodic orbits are each represented by a single curve using \eqref{eq:1norm}, and in the inset the periodic solution branches can be seen to both be approaching the intermediate steady state.

Figure~\ref{fig:perorbs_onepara_trp_III} shows the evolution of both the amplitude and period of the branch of stable and branch of unstable periodic orbits. The rapidly increasing period at the end of each branch suggests that both terminate in homoclinic bifurcations. This can be seen even more clearly by viewing the periodic orbits in phase space.

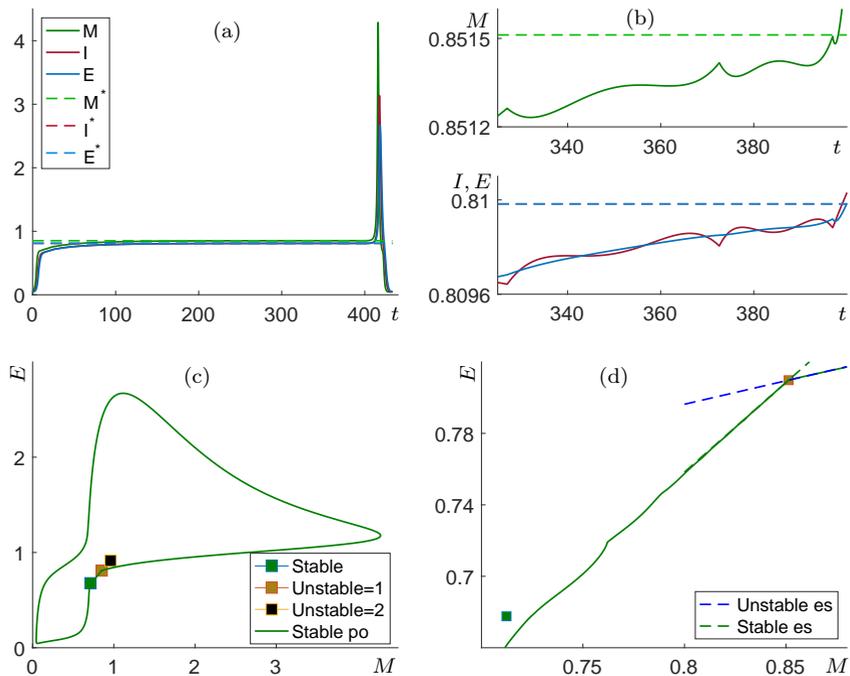
\begin{figure}[thp!]
\includegraphics[scale=0.42]{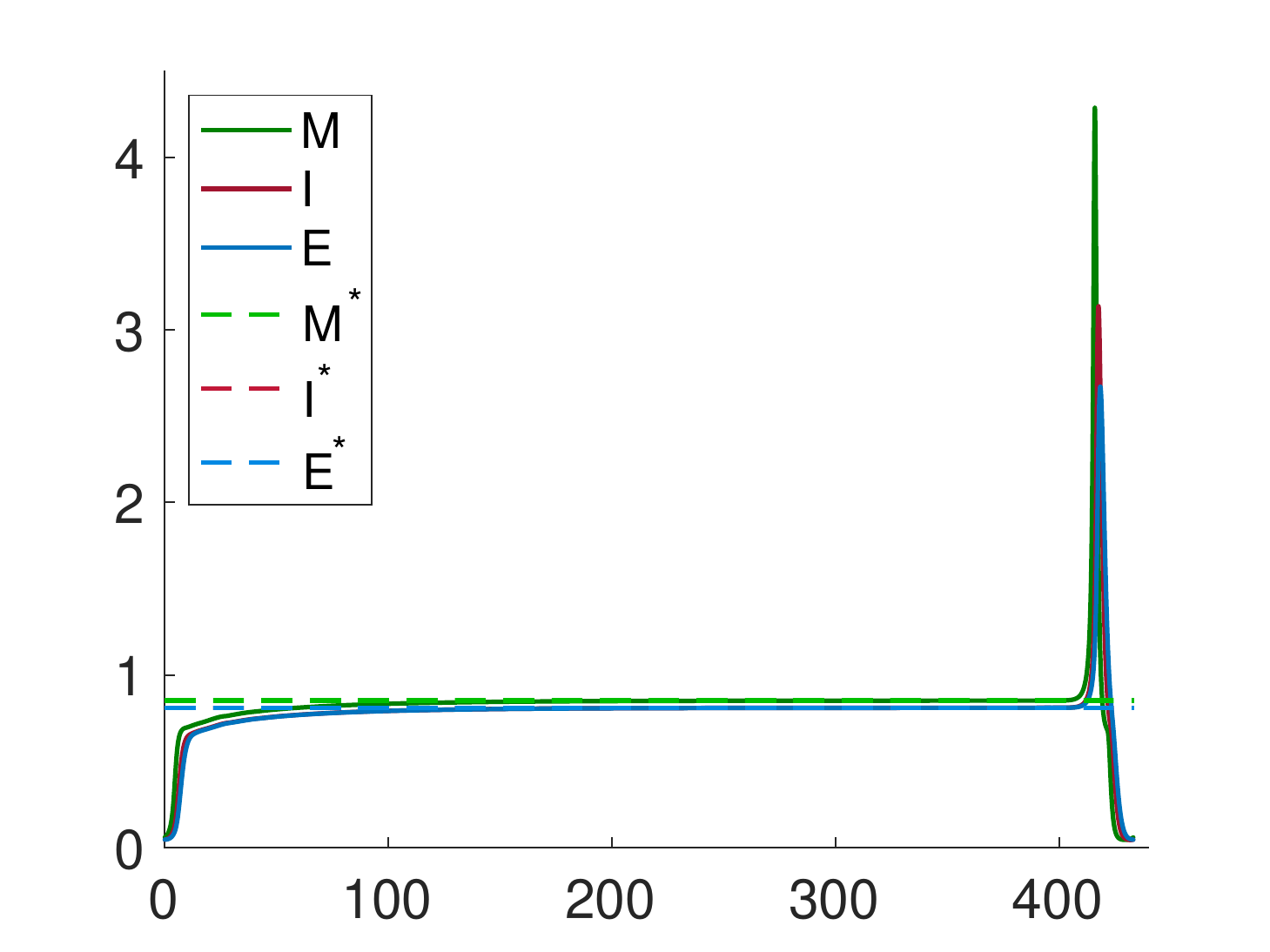}\hspace*{-1em}\includegraphics[scale=0.42]{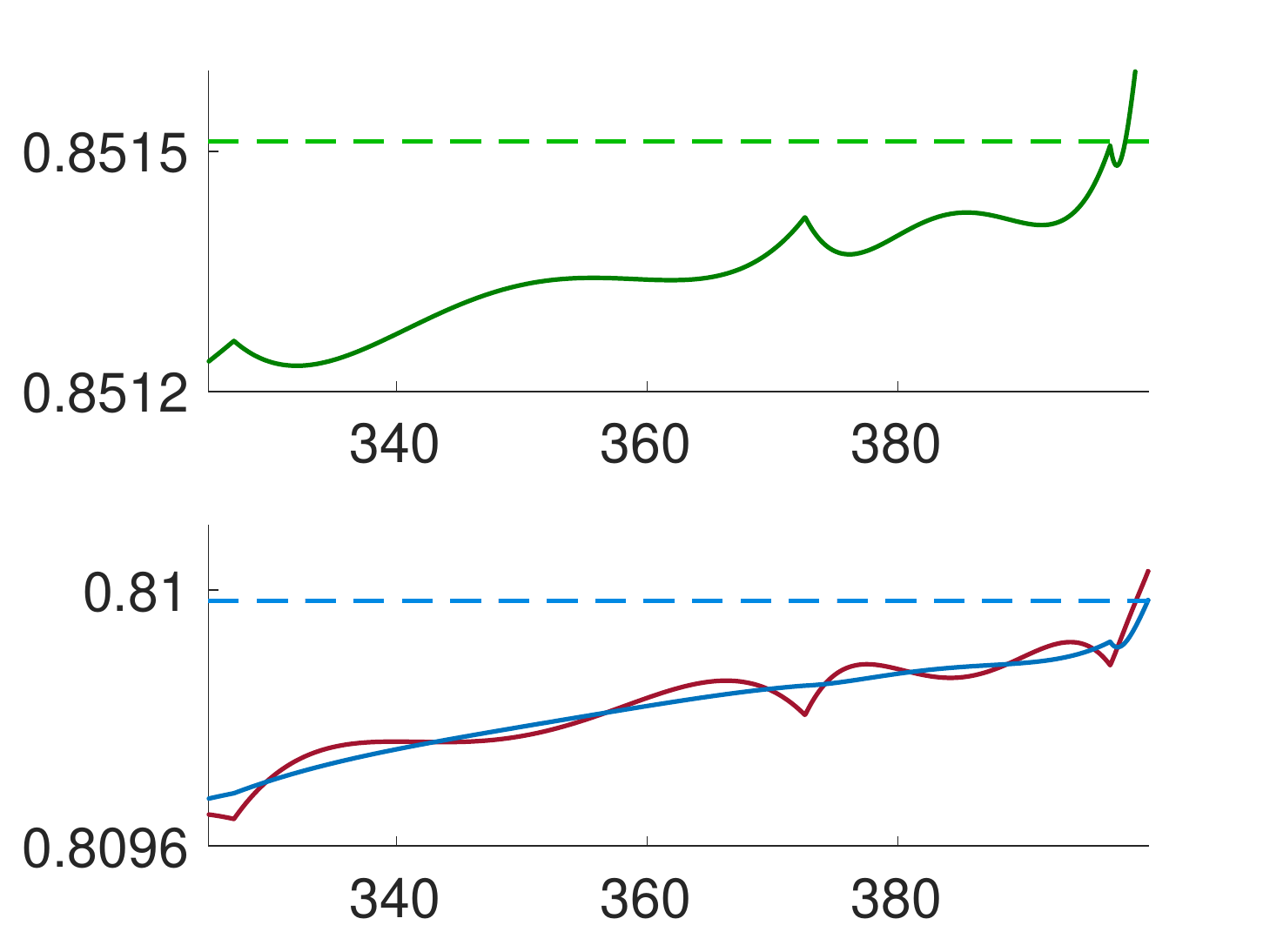}
\put(-255,112){(a)}
\put(-100,117){(b)}
\put(-160,116){$M$}
\put(-165,55){$I,E$}
\put(-188,5){$t$}
\put(-22,68){$t$}
\put(-20,5){$t$}

\includegraphics[scale=0.42]{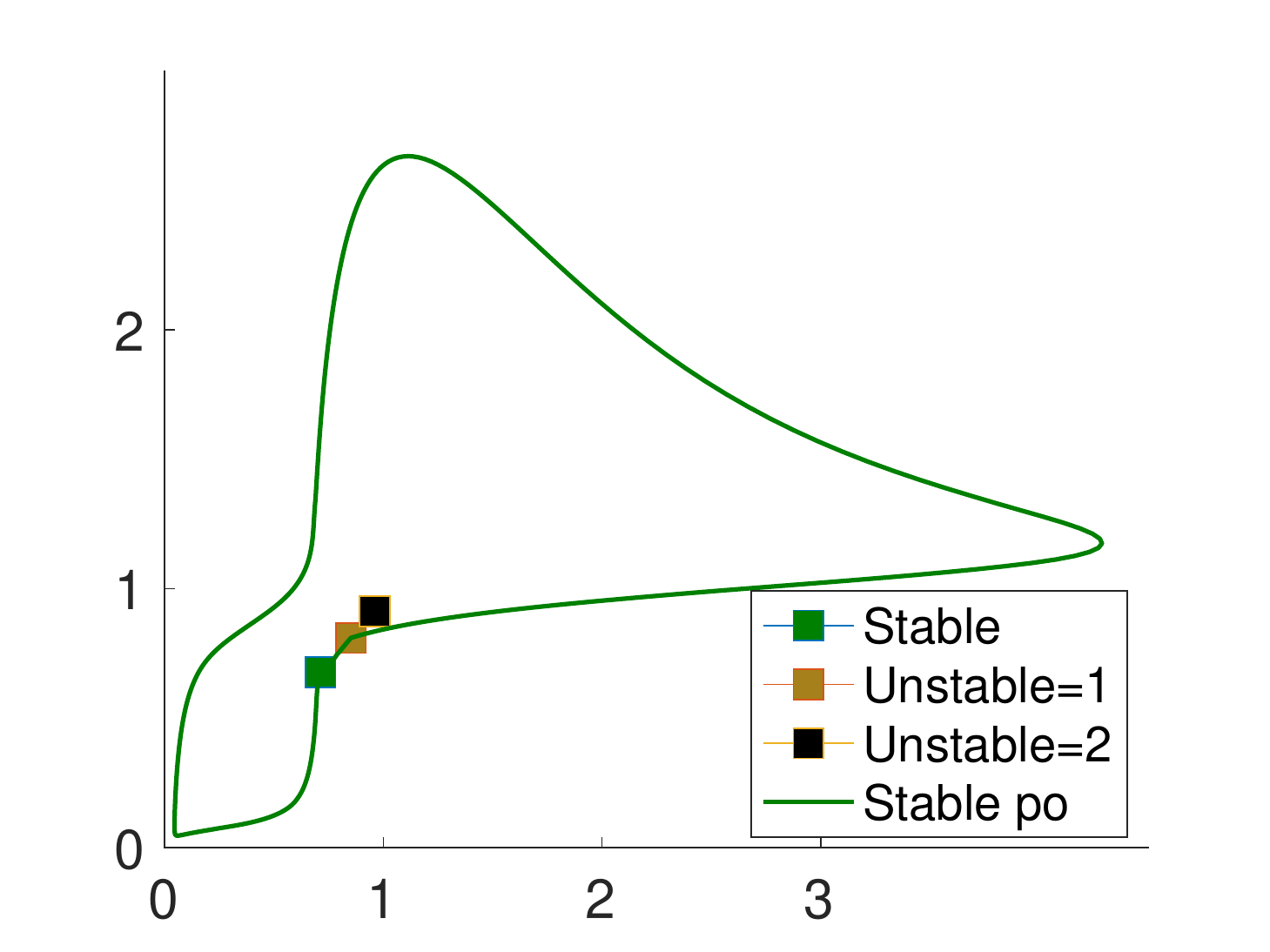}\hspace*{-1em}\includegraphics[scale=0.42]{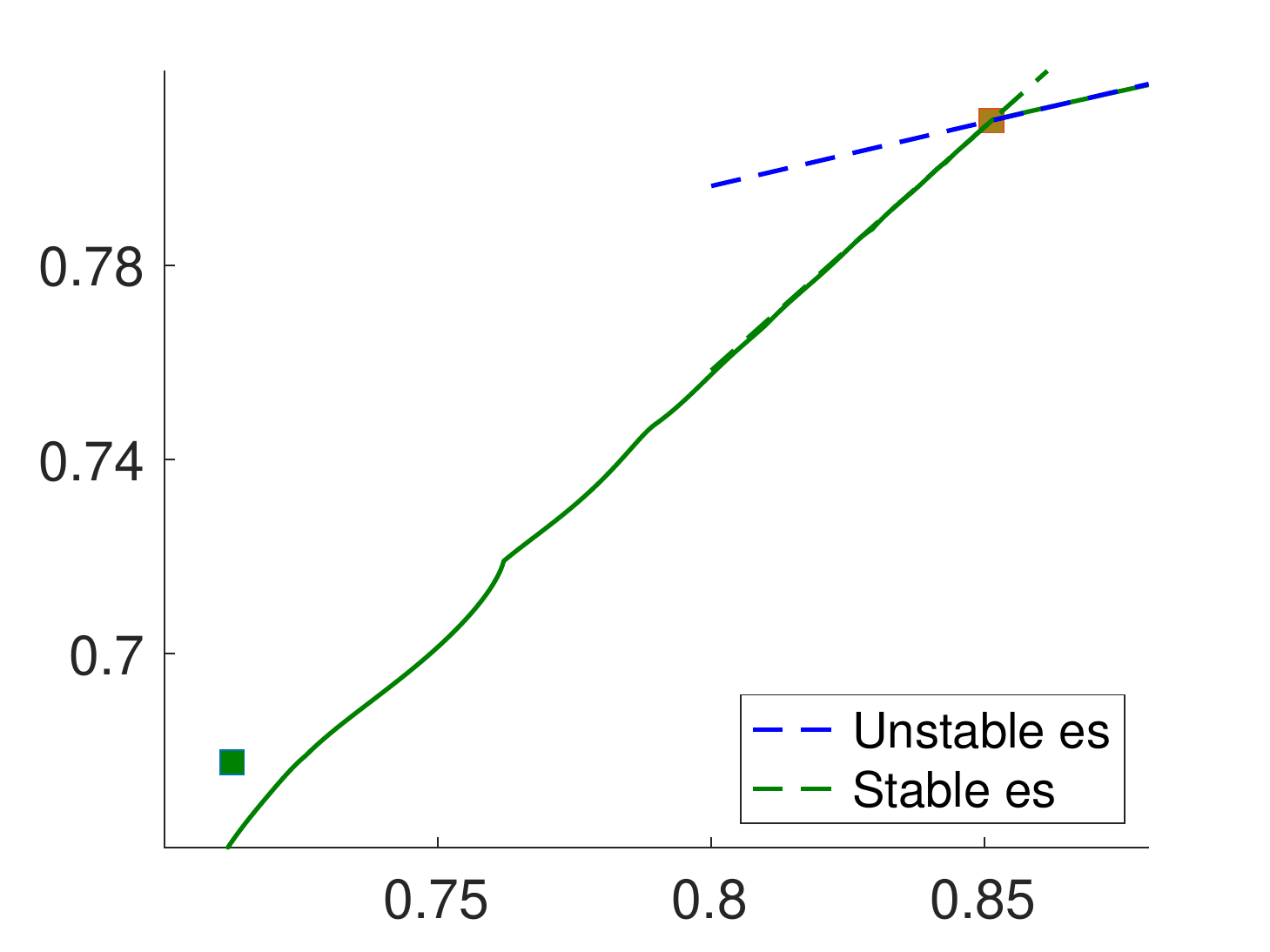}
\put(-266,115){(c)}
\put(-110,115){(d)}
\put(-25,5){$M$}
\put(-195,5){$M$}
\put(-331,115){\rotatebox{90}{$E$}}
\put(-162,115){\rotatebox{90}{$E$}}
\caption{Stable periodic orbit for repressible system \eqref{eq:sysonedel}
with $v_M^{min}=0.071577$,
$v_M^{max}=1$, $m=n=15$, and other parameters as in Table~\ref{table:multsspars}.
(a) The unstable periodic solution.
(b) The part of the periodic orbit very close to the middle steady state.
(c) A projection of the phase-space dynamics into the $M$-$E$ plane
showing the unstable periodic solution and steady states (colour coded according to the dimension of their unstable manifold).
(d) Detail of the phase space showing periodic orbit passing very close to middle steady state.
Also shown is the projection of the linear unstable manifold (in dashed blue) and the leading linear stable manifold (in dashed green) of this steady state.}
\label{stablehm}
\end{figure}

Figure~\ref{stablehm} shows the last limit cycle that we are able to compute on the branch of stable periodic orbits with $v_M^{min}=0.071577$. Panel (a) shows all three components of the periodic orbit as well as the unstable steady state from the middle segment of the branch of steady states. This shows that the system spends most of the time close to this steady state with just a short burst of production once per period.

Recall that DDEs define infinite dimensional dynamical systems whose phase space consists of function segments defined over a time interval equal to the largest delay. It follows that  for two solutions to be close in phase space, it is necessary that they are close in coordinate space for a time interval longer than the largest delay. With the parameters in this example $\tau_I=1$ and at the steady state $E^*=0.81$ and $\tau_M=9.155$. Thus the largest delay is close to $10$. Figure~\ref{stablehm}(b) shows a zoomed view
on the part of the periodic orbit closest to the steady state just before the burst. This shows that all three components of the periodic solution agree with the steady state values to three significant digits over a time interval several times larger than the delay, thus confirming that the periodic orbit passes close to the steady state in phase space.

Figure~\ref{stablehm}(c) shows a projection of the phase space dynamics on to the
$M$-$E$ plane showing the
stable periodic orbit and the three coexisting steady states
at $v_M^{min}=0.071577$. The periodic orbit appears to pass close to all three steady states, but that is an illusion created partly by the projection from infinite dimensions to $\R^2$ and partly by the scale which is very compact to show the large amplitude bursting periodic orbit. Also note that one of the steady states is asymptotically stable and so it is impossible for a periodic orbit to lie in its basin of attraction.

Figure~\ref{stablehm}(d) shows that the periodic orbit passes close to just the middle steady state,
and also shows the behaviour of the solution near to this  steady state.
The leading characteristic values of the intermediate unstable steady
state $(M^*, I^*, E^*)=(0.8515, 0.8100, 0.8100)$
at $v_M^{min}= 0.071577$, where stable periodic orbits cease to exist, are obtained from
\eqref{eq:CharEqA} as
\begin{align*}
\lambda_1 &=0.49051, \\  
\lambda_2 &= -0.017729, \\ 
\lambda_{3,4} 
&=-0.018217 \pm 0.70175i.
\end{align*}
Following the theory of Section~\ref{ssec:quasilinear} this leads to linearized solutions close to the steady state $(M^*,I^*,E^*)$ of the form
$$\left(\begin{array}{c}
M(t) \\ I(t) \\ E(t) \end{array}\right)
=
\left(\begin{array}{c}
M^* \\ I^* \\ E^* \end{array}\right)
+Ce^{\lambda t}
\left(\begin{array}{c}
\mathcal{E}_M \\ \mathcal{E}_I \\ \mathcal{E}_E \end{array}\right)$$
where the constant eigenvector $(\mathcal{E}_M,\mathcal{E}_I,\mathcal{E}_E)$ lies
in the nullspace of the matrix $A(\lambda)$ defined in \eqref{eq:Alinsys}.
The corresponding eigenvectors for $\lambda_1$ and $\lambda_2$ are computed as
\begin{align*}
v_1=\begin{pmatrix}
1\\
0.39077\\  
0.26222   
\end{pmatrix},
\quad
v_2=\begin{pmatrix}
1\\
0.98572\\   
1.0035    
\end{pmatrix}.
\end{align*}
Since $\lambda_1$ is the only characteristic value with positive real part the linear unstable manifold of the steady state is defined by $e^{\lambda_1 t}v_1$. When projecting phase space into the $M$-$E$ plane
this line has slope $\mathcal{E}_E/\mathcal{E}_M$.
The stable manifold of the steady state is infinite-dimensional and so cannot be represented in the $M$-$E$ plane, but the dominant part of the linear stable manifold (with slowest decay) is given by
$e^{\lambda_2 t}v_2$.

The projections of both the dominant part of the linear stable manifold and the linear unstable manifold
are shown in Figure~\ref{stablehm}(d).  The stable periodic orbit is seen to approach the steady state along a direction that is tangential to the dominant stable manifold before leaving   along a direction that is  tangential to the unstable manifold. Since the orbit passes very close to the steady state, the passage through the neighbourhood of the steady state takes a very long time. This results  in the large period of the orbit.

\begin{figure}[thp!]
\includegraphics[scale=0.42]{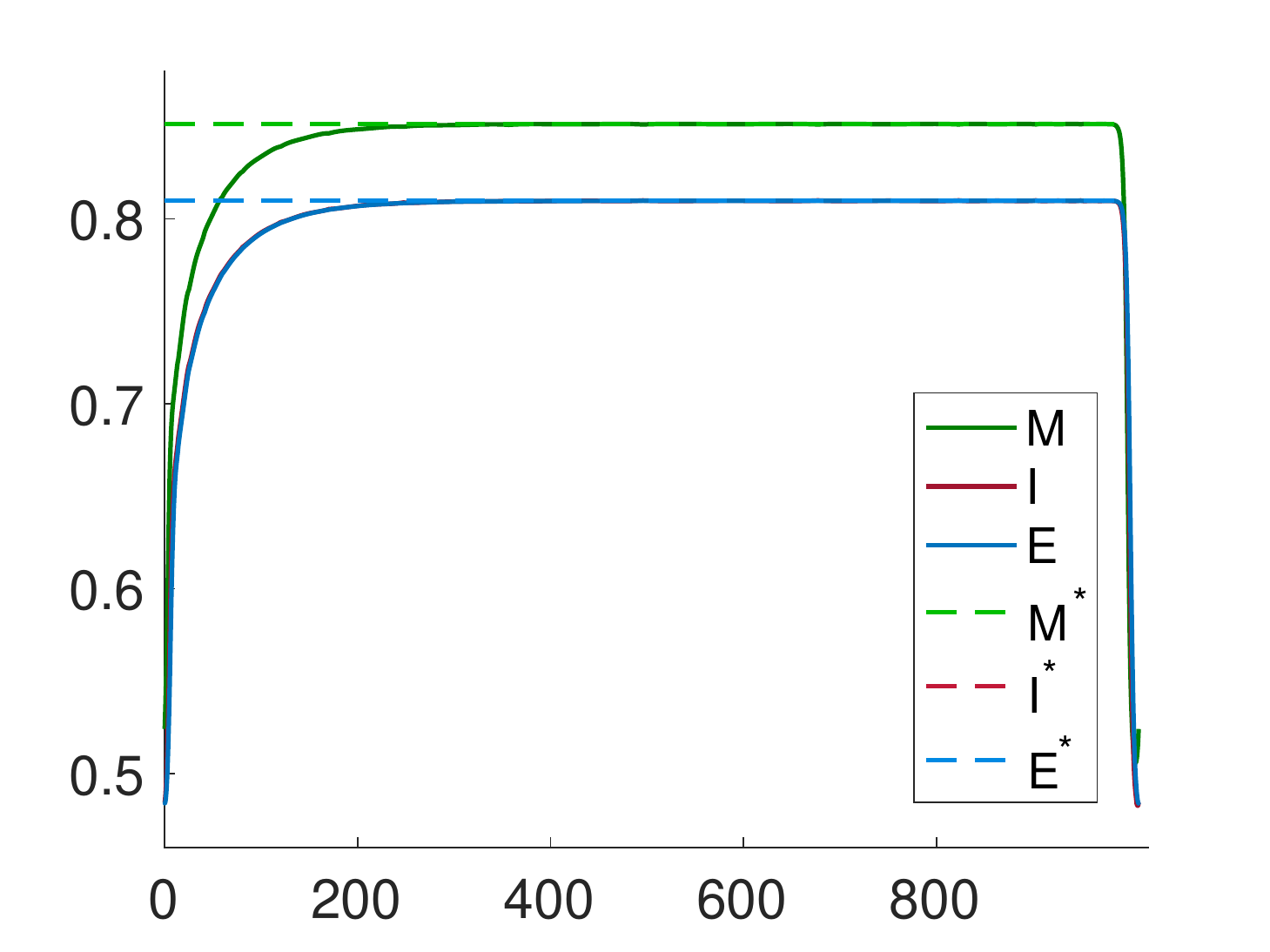}\hspace*{-1em}\includegraphics[scale=0.42]{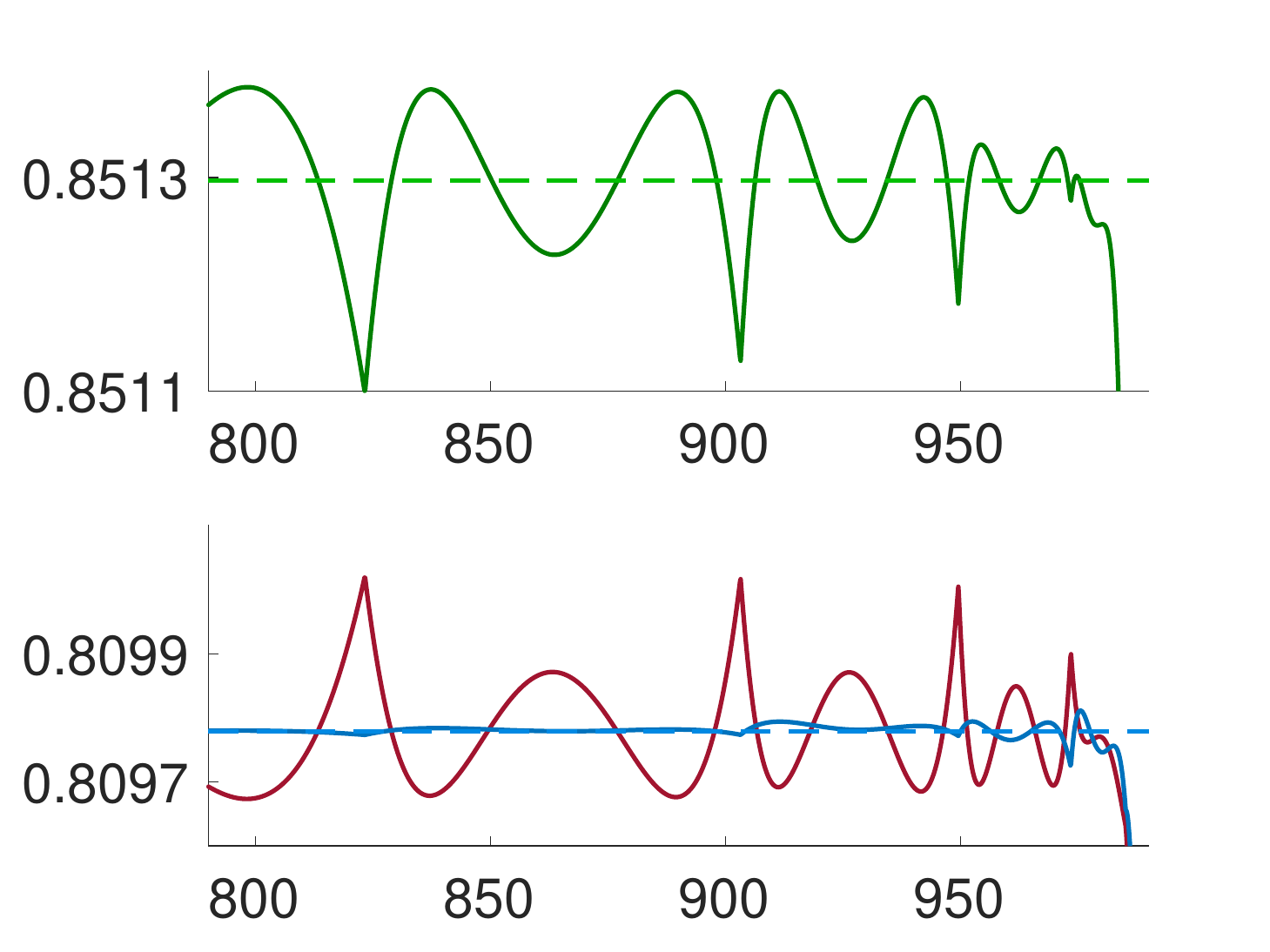}
\put(-270,120){(a)}
\put(-107,114){(b)}
\put(-162,114){$M$}
\put(-168,53){$I,E$}
\put(-192,5){$t$}
\put(-22,68){$t$}
\put(-20,5){$t$}

\includegraphics[scale=0.42]{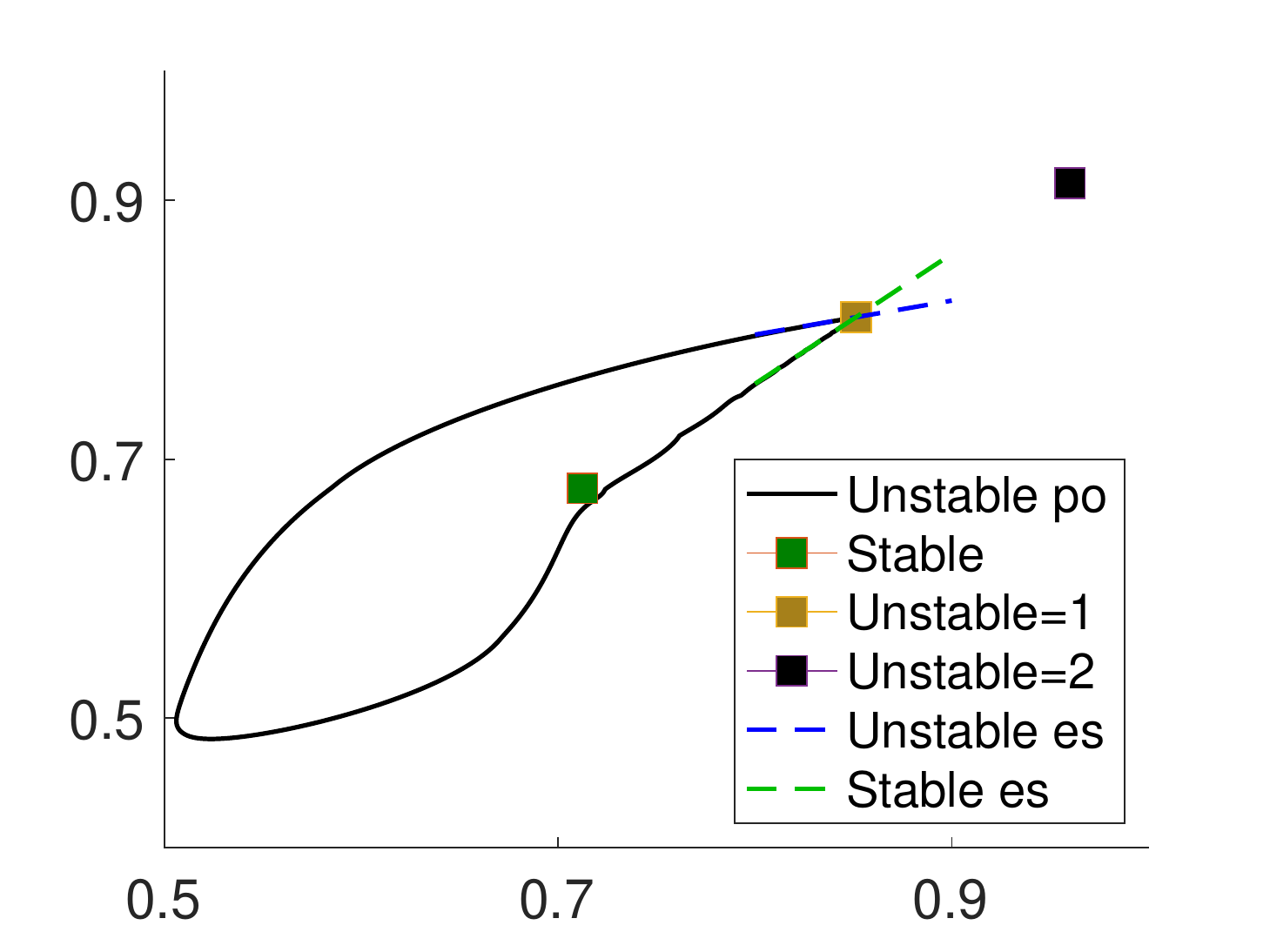}\hspace*{-1em}\includegraphics[scale=0.42]{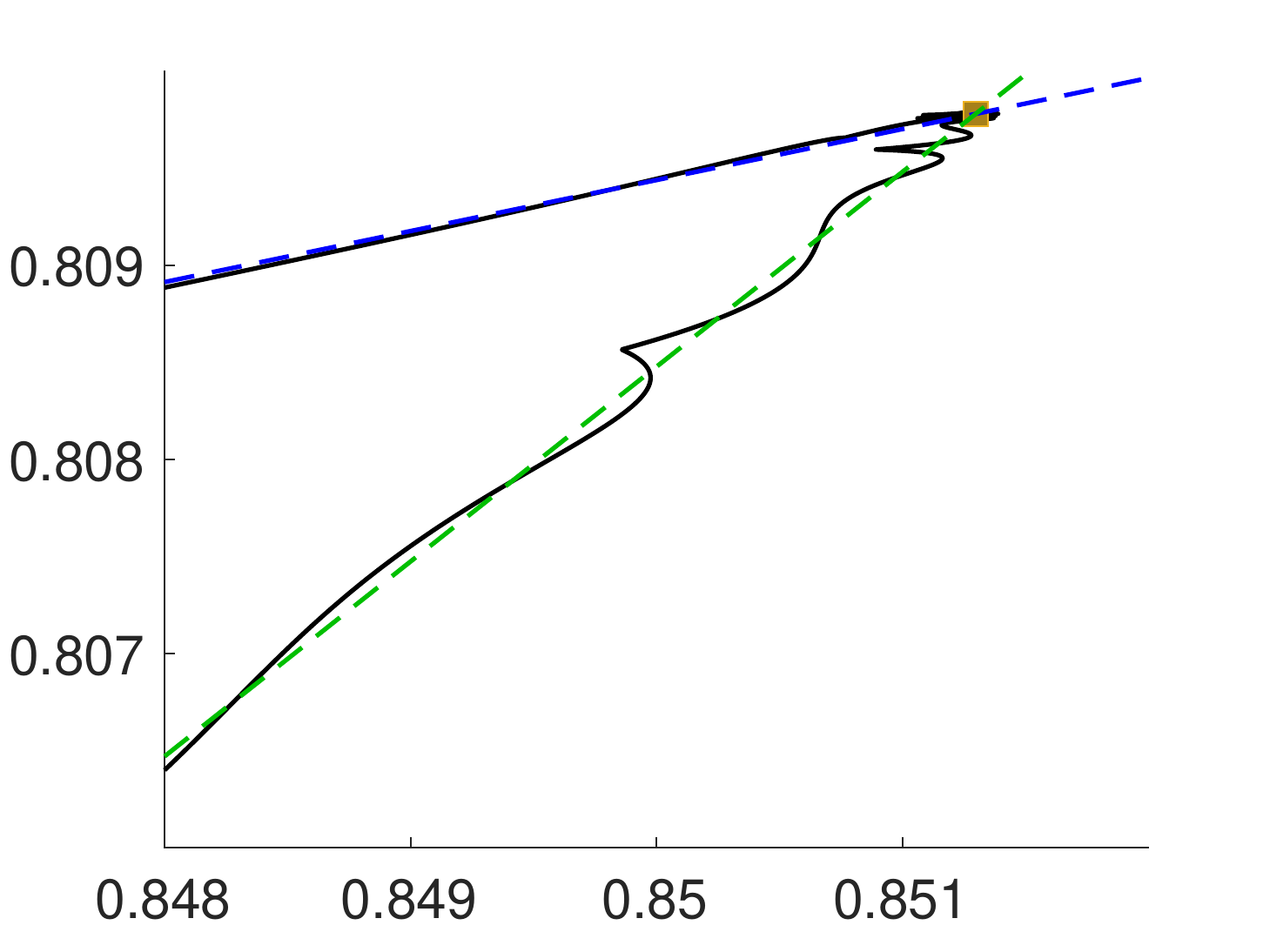}
\put(-270,115){(c)}
\put(-110,115){(d)}
\put(-25,5){$M$}
\put(-195,5){$M$}
\put(-331,115){\rotatebox{90}{$E$}}
\put(-162,115){\rotatebox{90}{$E$}}
		
\caption{Unstable periodic orbit for repressible system \eqref{eq:sysonedel}
with $v_M^{min}=0.071622$,
and the other parameters the same as in Figure~\ref{stablehm}. Panels (a) to (d) as described in caption to Figure~\ref{stablehm}.}
\label{unstablehm}
\end{figure}

\begin{figure}[htp!] \includegraphics[scale=0.28]{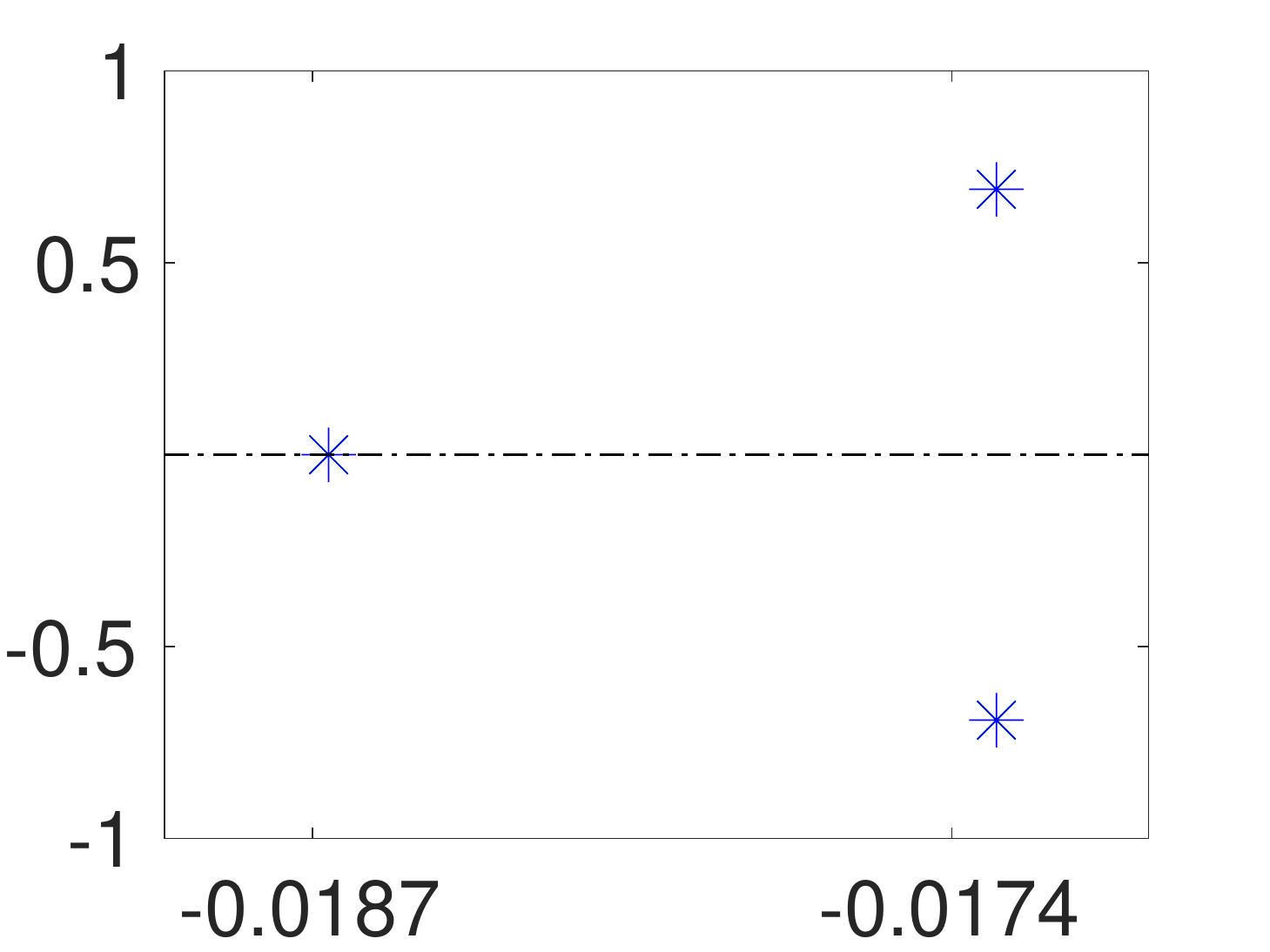}\hspace*{-0.5em}\includegraphics[scale=0.28]{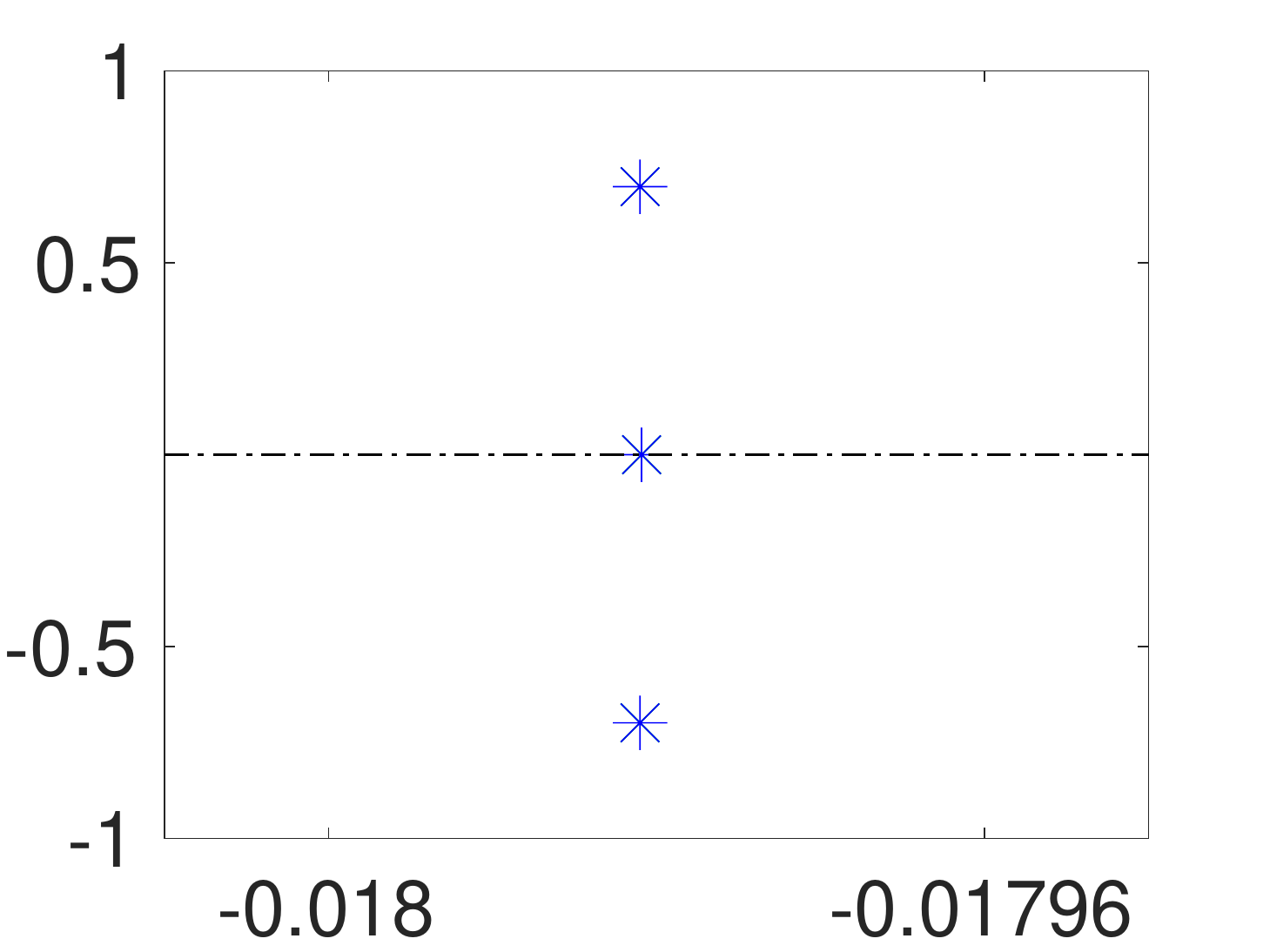}\hspace*{-0.5em}\includegraphics[scale=0.28]{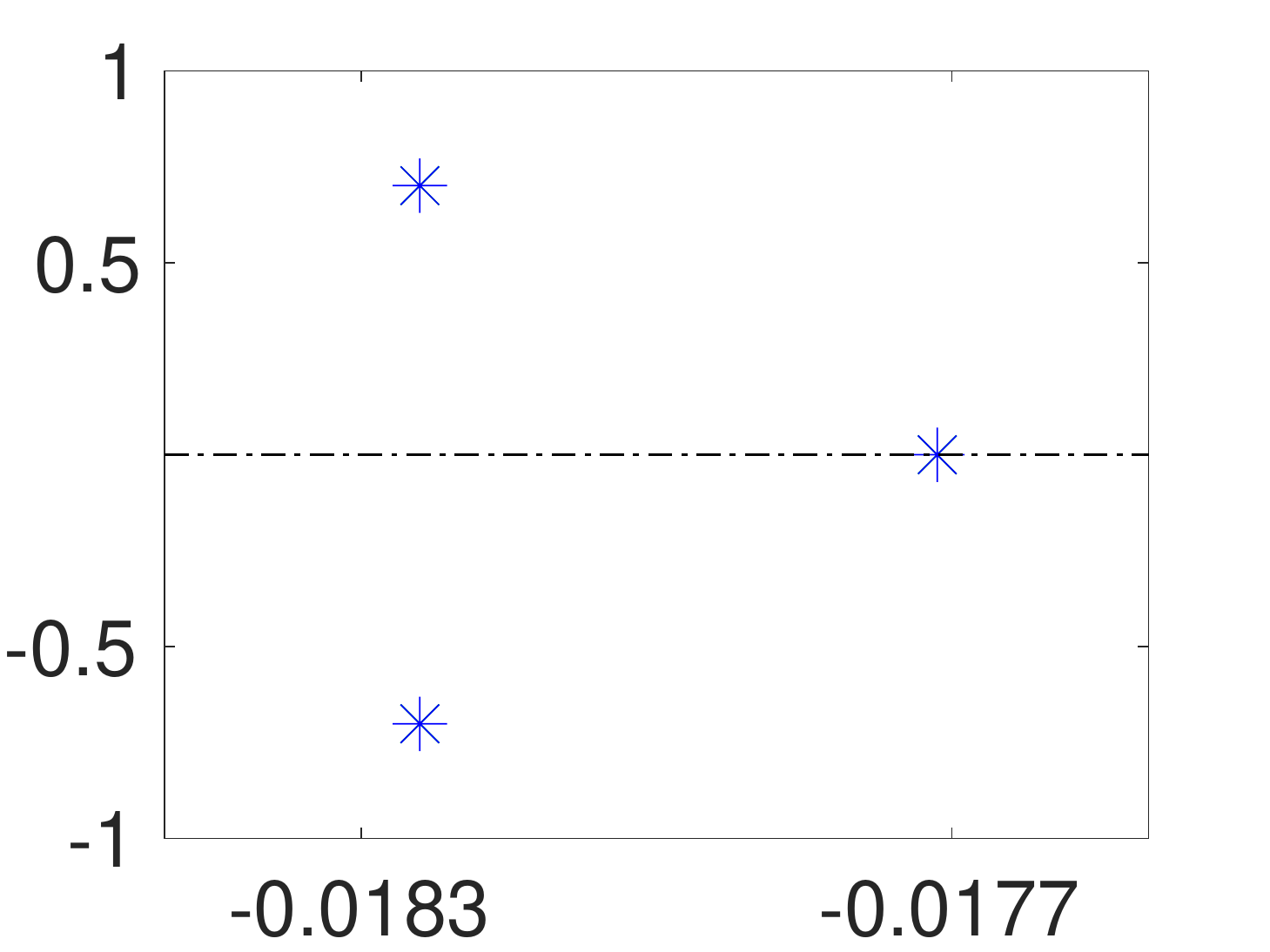}
\put(-70,1){$Im(\lambda)$}
\put(-185,1){$Im(\lambda)$}
\put(-298,1){$Im(\lambda)$}
\put(-112,35){\rotatebox{90}{$Re(\lambda)$}}
\put(-226,35){\rotatebox{90}{$Re(\lambda)$}}
\put(-340,35){\rotatebox{90}{$Re(\lambda)$}}
\put(-315,86){(a) $v_M^{min}=0.07$}
\put(-205,86){(b) $v_M^{min}=0.07103$}
\put(-95,86){(c) $v_M^{min}=0.07162$}
\caption{Configuration of leading negative real and complex-conjugate eigenvalues as $v_M^{min}$ varies. A 3-dimensional transition~\citep{Kalia2019} occurs on the middle branch of unstable steady states as in Figure~\ref{fig:threess_onepara_trp_III}. }
	\label{3DLtrans}
\end{figure}
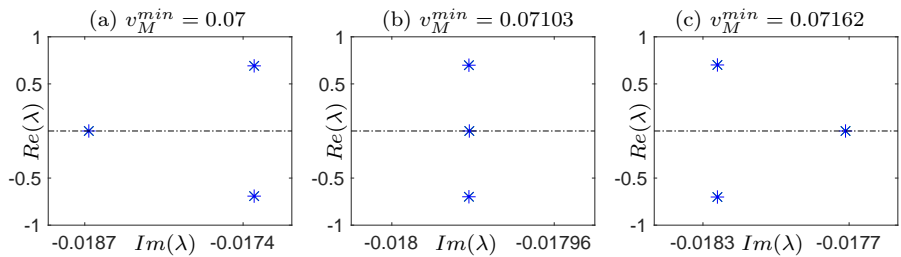

Figure~\ref{unstablehm} is similar to Figure~\ref{stablehm}
but shows the last periodic orbit that we are able to compute
on the branch of unstable periodic orbits with $v_M^{min}=0.071622$.
Comparing the two figures we see that the unstable periodic orbit is quite different to the stable orbit.
Figure~\ref{unstablehm}(a) and (b) show that again the periodic orbit is close to the intermediate steady state for most of the period except for a short burst (or antiburst) of \emph{depressed} production.

The characteristic values and corresponding eigenvectors can also be computed for the unstable periodic orbit, but the value of $v_M^{min}$ only differs in the third significant digit between the two examples. The  characteristic values and eigenvectors agree with those above to the third significant digit.
The phase space plots in Figure~\ref{unstablehm}(c) and (d) again show a periodic orbit close to homoclinic approaching the steady state near the dominant linear stable manifold and leaving tangential to the linear unstable manifold.

There are two significant differences from the stable periodic orbit. Firstly, the unstable periodic orbit leaves the neighbourhood of the steady state in the \emph{opposite} direction to the stable periodic orbit, which results in the production being decreased rather than increased during  the burst. Secondly it is  apparent in Figure~\ref{unstablehm}(b) and (d) that the periodic orbit is not tangential to the dominant part of the linear stable manifiold but rather oscillates about it. This seems to arise because there is only a small difference in the real parts between $\lambda_2$ and the next characteristic values which occur as a complex conjugate pair $\lambda_{3,4}$.
Furthermore, as shown in Figure~\ref{3DLtrans}, for a very nearby value of $v_M^{min}$,
the leading negative real eigenvalue $\lambda_2$ and complex-conjugate eigenvalues $\lambda_{3,4}$ exchange order.
Such a transition in a system also having one real positive eigenvalue is called a 3-dimensional or 3DL transition and is associated with rich Shilnikov homoclinic bifurcation structures
\citep{Kalia2019}.

\subsection{Inducible Operon with One State-Dependent Delay}
\label{sec:indexamples}

We now turn to consider the Goodwin model \eqref{eq:sysonedel} with one-state dependent delay in the case of an inducible operon with functions defined in Table~\ref{table:delay/velocity behaviour}.
Recall that with both delays constant (and also in the absence of delays) an
inducible system with $n>1$ can have either a single globally stable steady state,
or there can be two locally stable steady states and an unstable intermediate steady state. There are no other possibilities when using the functions in Table~\ref{table:delay/velocity behaviour} \citep{yildirim2004dynamics}.

We will show that an inducible operon with state-dependent transcription delay
$\tau_M$ can support stable and unstable periodic orbits and that these can be generated in supercritical or subcritical Hopf bifurcations, or in fold bifurcations of periodic orbits.

\subsubsection*{Inducible Supercritical Hopf Bifurcation}

\begin{figure}[htp!]
	\centering
	\includegraphics[scale=0.6]{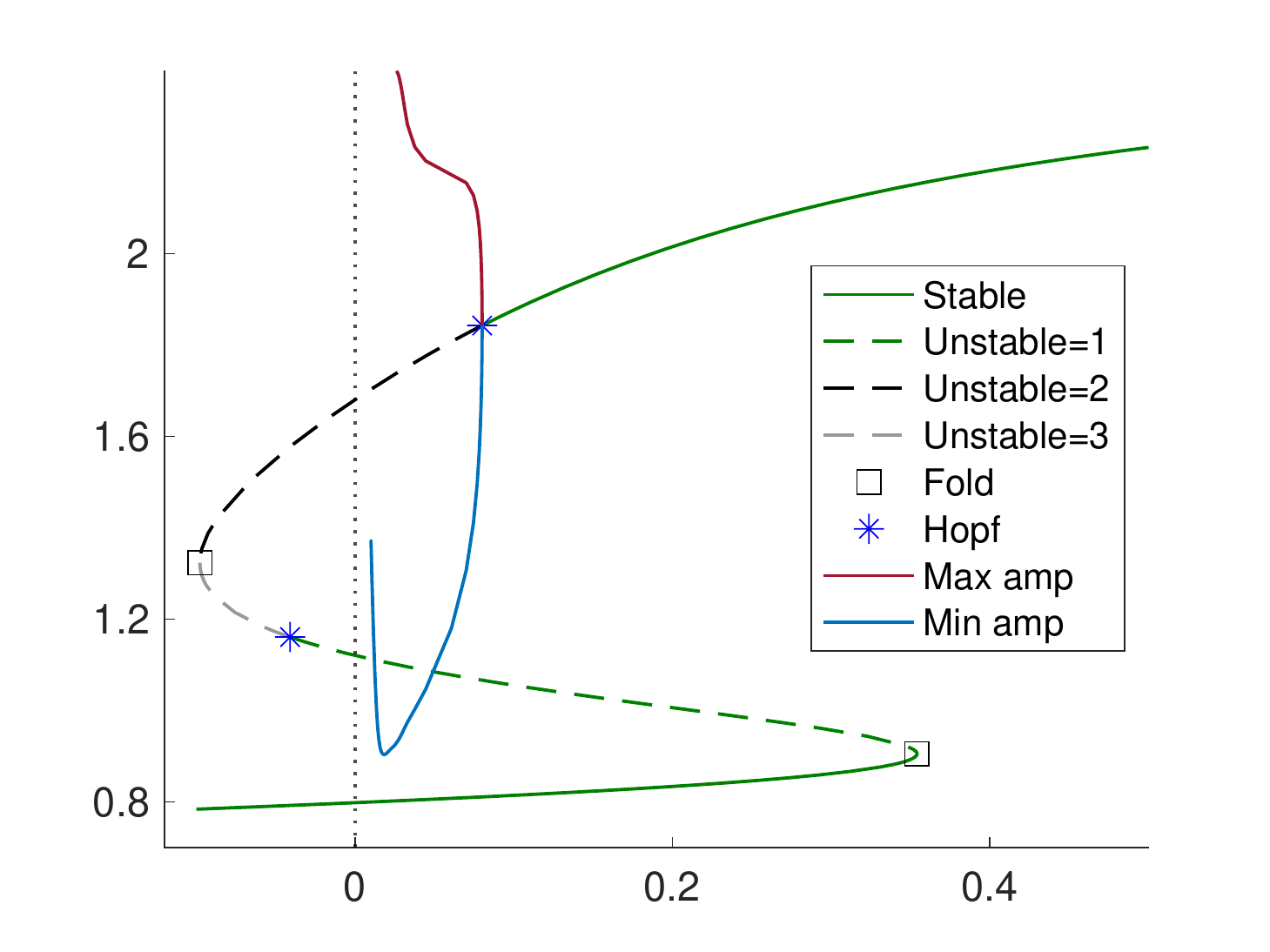}
 \put(-230,165){\rotatebox{90}{$E$}}
\put(-40,10){$v_M^{min}$}
	\caption{Bifurcation diagram of the model \eqref{eq:sysonedel} for an inducible system
             with parameter values as defined in Table \ref{table: parameter_lac}
             and $v_M^{min}$ treated as a continuation/bifurcation parameter.
         Line specifications can be found in Figure \ref{fig:threess_onepara_trp}.
        The amplitudes of $E$-component of periodic solutions $\R \to \R^3$ are shown.
         The vertical dotted line at $v_M^{min} = 0$ separates the biologically realistic case
         $v_M^{min}>0$ from the biologically unrealistic case $v_M^{min}<0$ (see text).
         Bifurcations occurring for $v_M^{min}>0$ are detailed in Table~\ref{table_threess_onepara}.}
	\label{threess_onepara}
\end{figure}

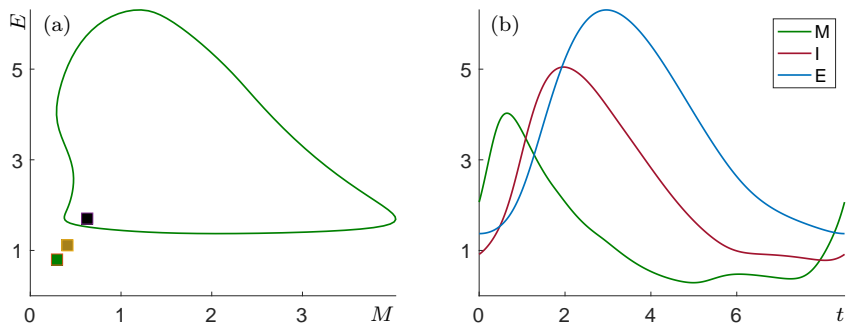
\begin{figure}[htp!]
\hspace*{-0.5em}\includegraphics[scale=0.42]{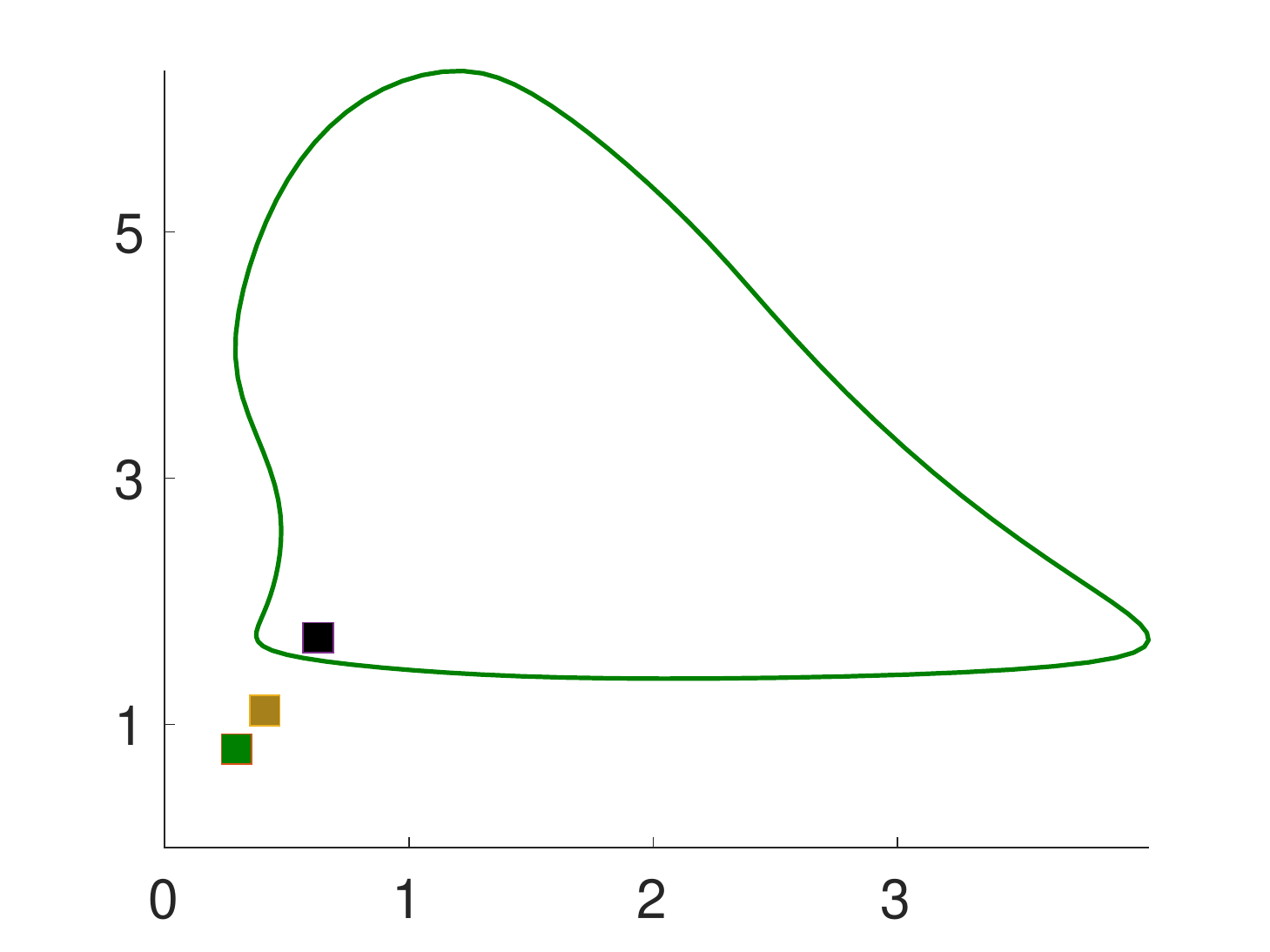}\hspace*{-1em}\includegraphics[scale=0.42]{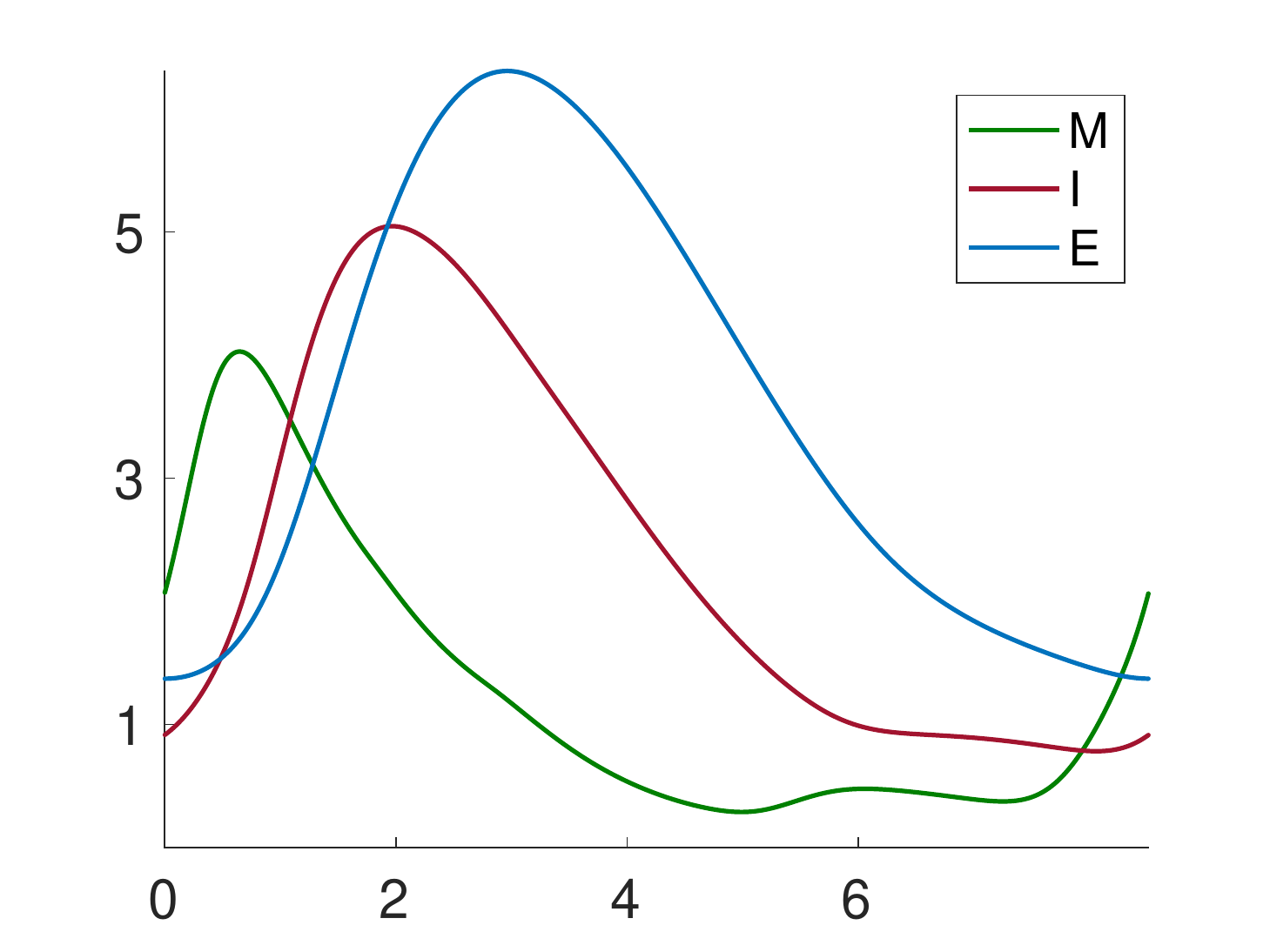}
\put(-318,115){(a)}
\put(-150,115){(b)}
\put(-330,115){\rotatebox{90}{$E$}}
\put(-195,5){$M$}
\put(-20,5){$t$}
\caption{Inducible system \eqref{eq:sysonedel}
with parameters as defined in Table \ref{table: parameter_lac}
showing the orbits from Figure~\ref{threess_onepara}
at $v_M^{min}=0.01$.
(a) A projection of the phase space dynamics into the $M$-$E$ plane showing
periodic orbits represented by closed curves, and steady states by squares (whose colour indicates the number of unstable eigenvalues as in Figure~\ref{stablehm}).
(b) The three components of the stable periodic solution.
}
\label{vMmin001_lac}
\end{figure}

We begin by considering the inducible operon model \eqref{eq:sysonedel} with parameters defined in
Table~\ref{table: parameter_lac}. With this parameter set and $v_M^{min}=v_M^{max}=1$,  both delays are constant  and the model has a single globally stable steady state. For $v_M^{min}<1$ the transcription delay becomes state-dependent and several bifurcations occur, as shown in
Figure~\ref{threess_onepara} and listed in Table~\ref{table_threess_onepara}.

\begin{table}[htp!]
\centering
\begin{tabular}{|c|c|}
\hline
Quantity & Value \\ \hline
		$\mu$ & 0.05   \\
		$\beta_M$ &  1  \\
		$\beta_I$ & 1.8  \\
		$\beta_E$ & 1.5 \\
		$\bgamma_M$ & 1  \\
		$\bgamma_I$ & 0.97 \\
		$\bgamma_E$ & 1  \\
		\hline
		$K$ & 4  \\
		$K_1$ & 1  \\
		$n$ & 4  \\
		\hline
	\end{tabular}
\hspace*{1em}\begin{tabular}{|c|c|}
\hline
Quantity & Value \\ \hline
		$m$ & 2  \\
		$a_M$ & 1  \\
		$v_M^{max}$ & 0.2  \\
		$E_{50}$ & 1   \\
		\hline
		$\tau_I$ & 0.5  \\
		\hline
\end{tabular}
	\caption{Parameters for inducible operon example of Figure~\ref{threess_onepara}.}
	\label{table: parameter_lac}
\end{table}

\begin{table}[htp!]
	\centering
	\begin{tabular}{| l | l | l | l |}
		\hline
		Bifurcation & Bifurcation parameter value & Unstable eigenvalues& $E^*$ value \\
		\hline
		Hopf & $v_M^{min}=0.080031,$ period = 6.2271 & 0 to 2 & 1.8571 \\
		\hline \hline
		Fold & $v_M^{min}=0.35409$ & 1 to 0 & 0.9052 \\
		\hline
	\end{tabular}
	\caption{Bifurcation information associated with Figure \ref{threess_onepara}. Those bifurcations occurring for $v_M^{min}<0$ are not displayed due to lack of physiological meaning. }
	\label{table_threess_onepara}
\end{table}

%
%
%

As $v_M^{min}$ is decreased there is first a fold bifurcation which creates two additional steady states. This results in bistability  between two stable steady states
for $v_M^{min}\in(0.08,0.354)$ separated by an intermediate unstable steady state.
This configuration is well known for inducible operons
with constant delays, but here the bifurcation to three steady states is induced by varying the state-dependency of the delay $\tau_M$.

Reducing $v_M^{min}$ further, an unexpected event occurs; the upper steady state loses stability in a supercritical Hopf bifurcation which creates a stable periodic orbit which exists for $v_M^{min}<0.08$. This stable periodic orbit coexists with one stable and two unstable steady states. Thus we have an interval of bistability between a limit cycle and a steady state for an \emph{inducible} operon.

The stable periodic orbit and a projection of phase space into the $M$-$E$ plane are shown in
Figure~\ref{vMmin001_lac} for $v_M^{min}=0.01$. We suspect that the periodic orbit exists
for all $v_M^{min}>0$, but the numerical discretization of the threshold integral described in
Section~\ref{subsec:discretization} requires $v_M^{min}$ bounded away for zero, and we only compute periodic orbits for $v_M^{min}\geq 0.01$.

The linearization correction method of Section~\ref{subsec:correction} avoids discretizing the
integral and is applicable even when $v_M^{min}<0$. Though $v_M^{min}<0$ leads to negative transcription velocities which is not physiological, this can be computationally useful. This is demonstrated
in Figure~\ref{threess_onepara} where continuation through negative values of $v_M^{min}$ reveals that the
different branches of steady states are joined at a  fold bifurcation with $v_M^{min}<0$. This allows computation  of all the physiological steady states  by continuation of a single branch.


\subsubsection*{Inducible Subcritical Hopf Bifurcation}

We now change just one parameter value from the previous example and consider
the inducible state-dependent transcription delay operon model \eqref{eq:sysonedel}
with parameters as in Table~\ref{table: parameter_lac}, except for the Hill coefficient in the transcription velocity function which we now set to $m=4$.

\begin{figure}[htp!]
	\centering
	\includegraphics[scale=0.6]{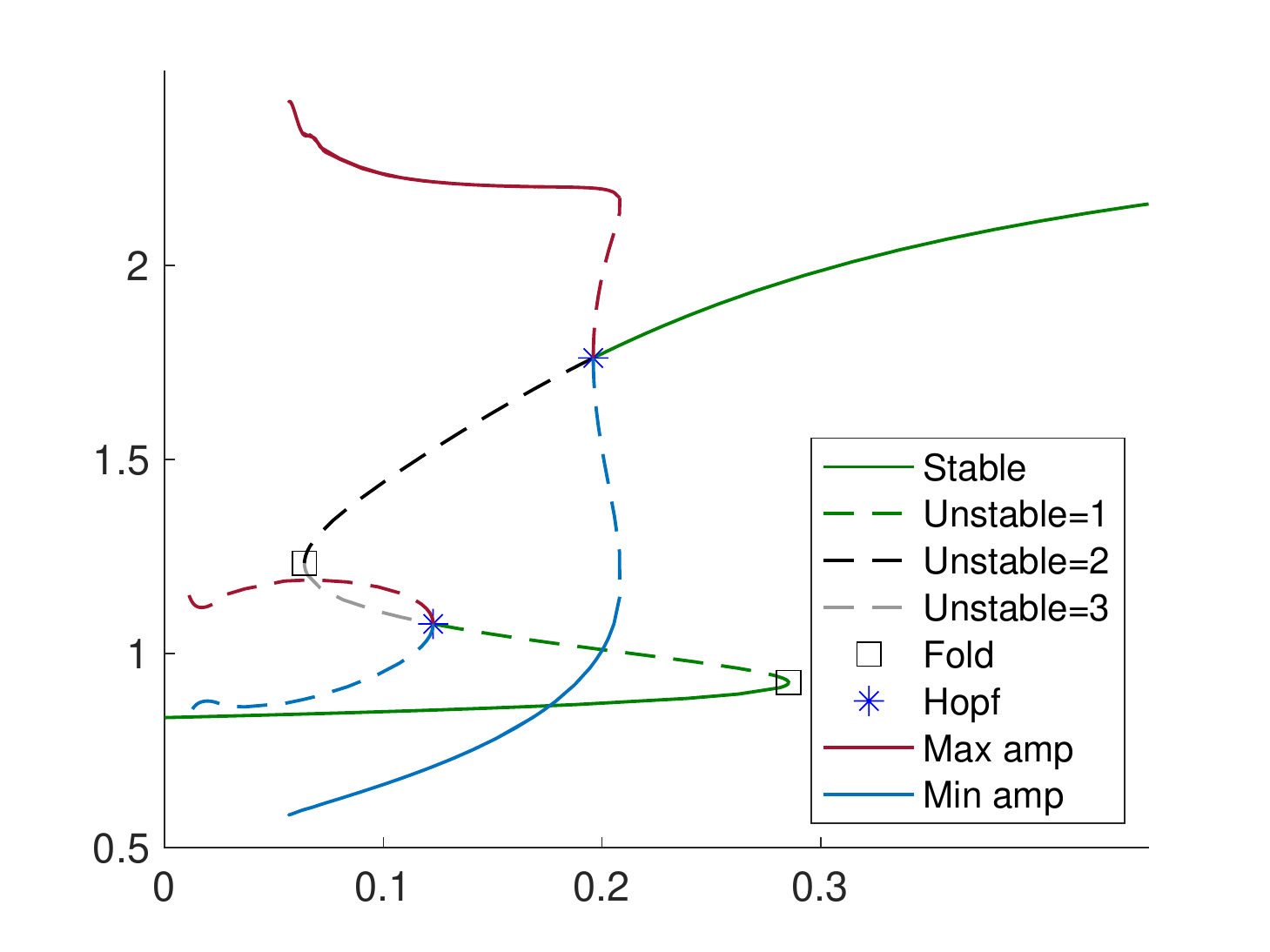}
\put(-230,165){\rotatebox{90}{$E$}}
\put(-40,10){$v_M^{min}$}

	\caption{Bifurcation diagram of the inducible operon model \eqref{eq:sysonedel} with $m=4$ and all other parameters as defined in Table~\ref{table: parameter_lac}.
Bifurcations are listed in Table~\ref{table_threess_onepara_MI}.}
	\label{threess_onepara_MI}
\end{figure}

\begin{table}[htp]
	\centering
	\begin{tabular}{| l | l | l | l |}
		\hline
		Bifurcation & Bifurcation parameter value & Unstable eigenvalues& $E^*$ value \\
		\hline \hline
		Hopf & $v_M^{min}=0.19603,$ period = 6.4488 & 0 to 2 & 1.7609 \\
		\hline
		Fold & $v_M^{min}=0.063903$ & 2 to 3 & 1.2317 \\
		\hline
		Hopf & $v_M^{min}=0.12279,$ period = 4.7796 & 3 to 1 & 1.0759 \\
		\hline
		Fold & $v_M^{min}=0.28543$ & 1 to 0 & 0.9249 \\
		\hline \hline
        \rule[-8pt]{0pt}{21pt}\parbox{2.0cm}{Fold of\\ periodic orbits} & $v_M^{min}=0.20812$, period = 5.8831 & 1 to 0 & - \\
        \hline
	\end{tabular}
	\caption{Steady state and Periodic Orbit Bifurcation information
for the Example shown in Figure~\ref{threess_onepara_MI}. }
	\label{table_threess_onepara_MI}
\end{table}

Comparing Figure~\ref{threess_onepara_MI} with the previous example in Figure~\ref{threess_onepara}
we see that changing the value of $m$ from $2$ to $4$ results in two important changes in the bifurcations.
Firstly, in Figure~\ref{threess_onepara_MI} both the fold bifurcations on the branch of steady states now occur for positive values of $v_M^{min}$. Consequently for
$v_M^{min}\in(0.064,0.285)$ there are three co-existing steady states, while for both larger and smaller
values of $v_M^{min}>0$ there is a unique stable steady state.

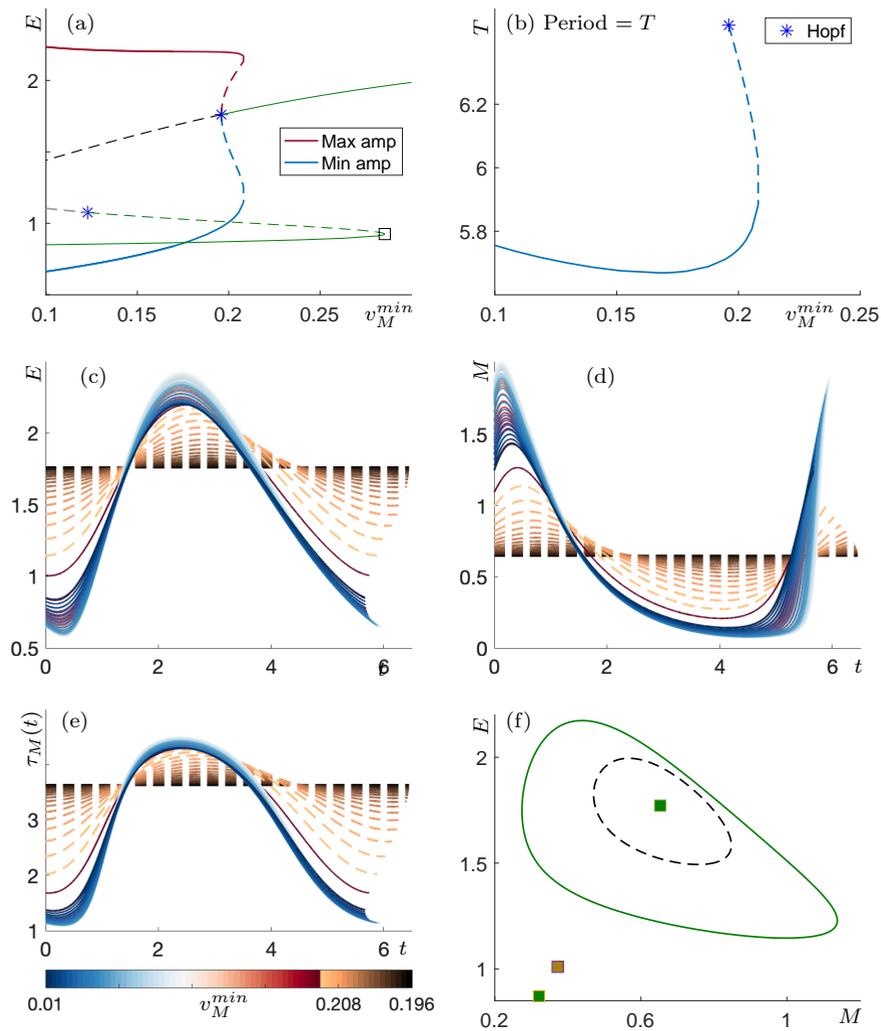
\begin{figure}[htp!]
\hspace*{-1em}\includegraphics[scale=0.42]{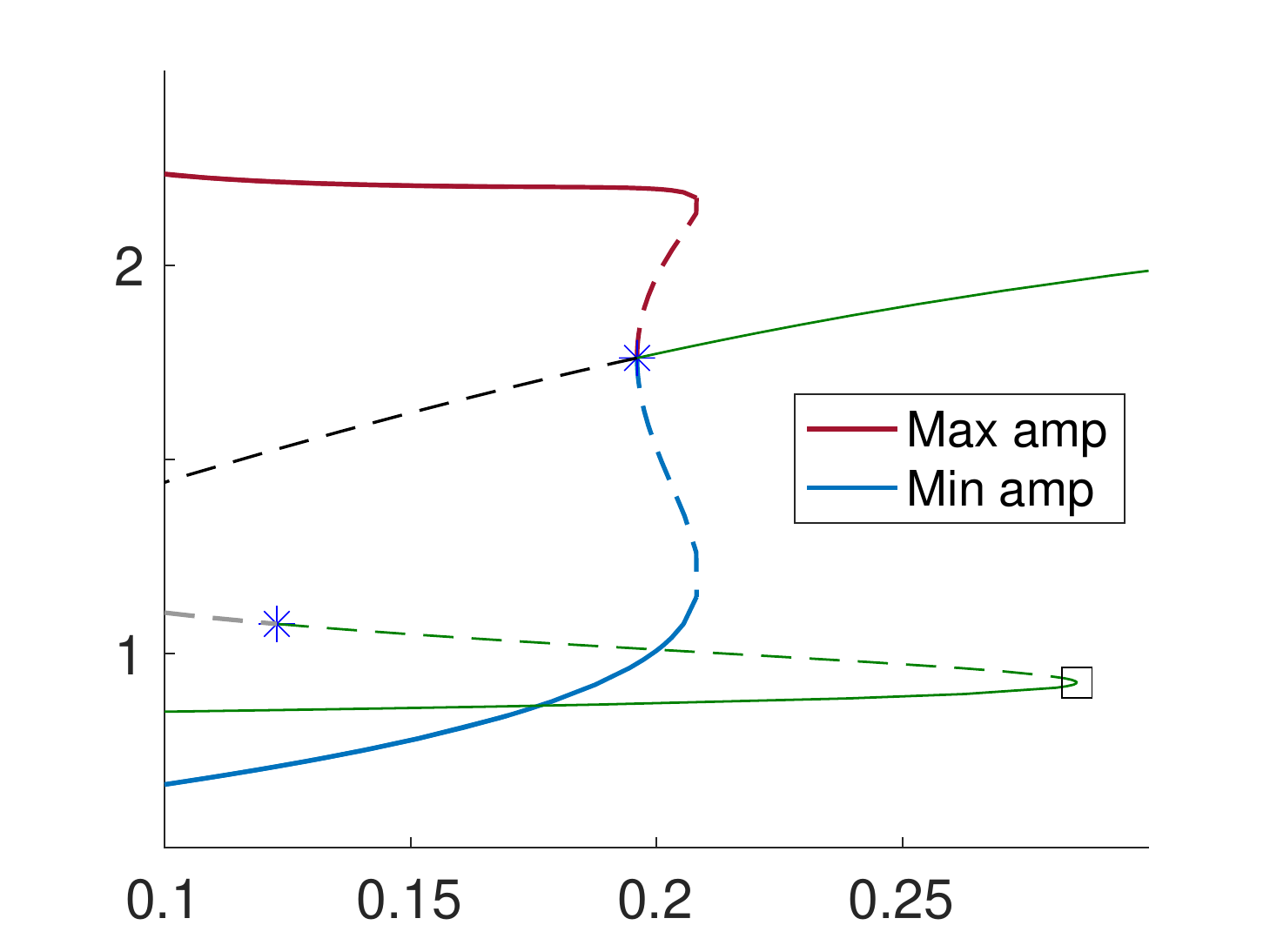}\hspace*{-1em}\includegraphics[scale=0.42]{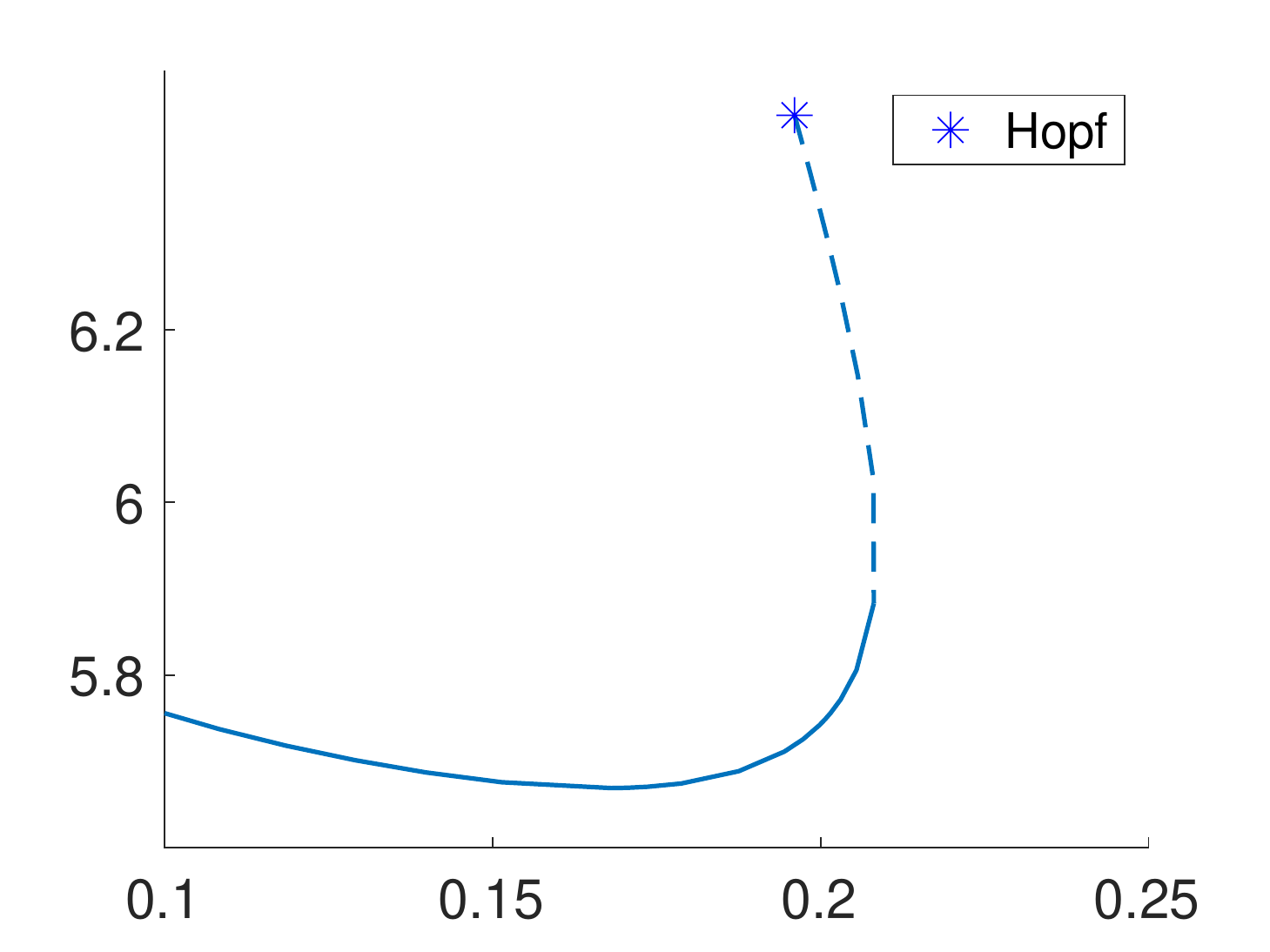}
\put(-315,115){(a)}
\put(-150,115){(b) Period $=T$}
\put(-45,5){$v_M^{min}$}
\put(-203,5){$v_M^{min}$}
\put(-330,114){\rotatebox{90}{$E$}}
\put(-162,113){\rotatebox{90}{$T$}}
	
\hspace*{-1em}\includegraphics[scale=0.42]{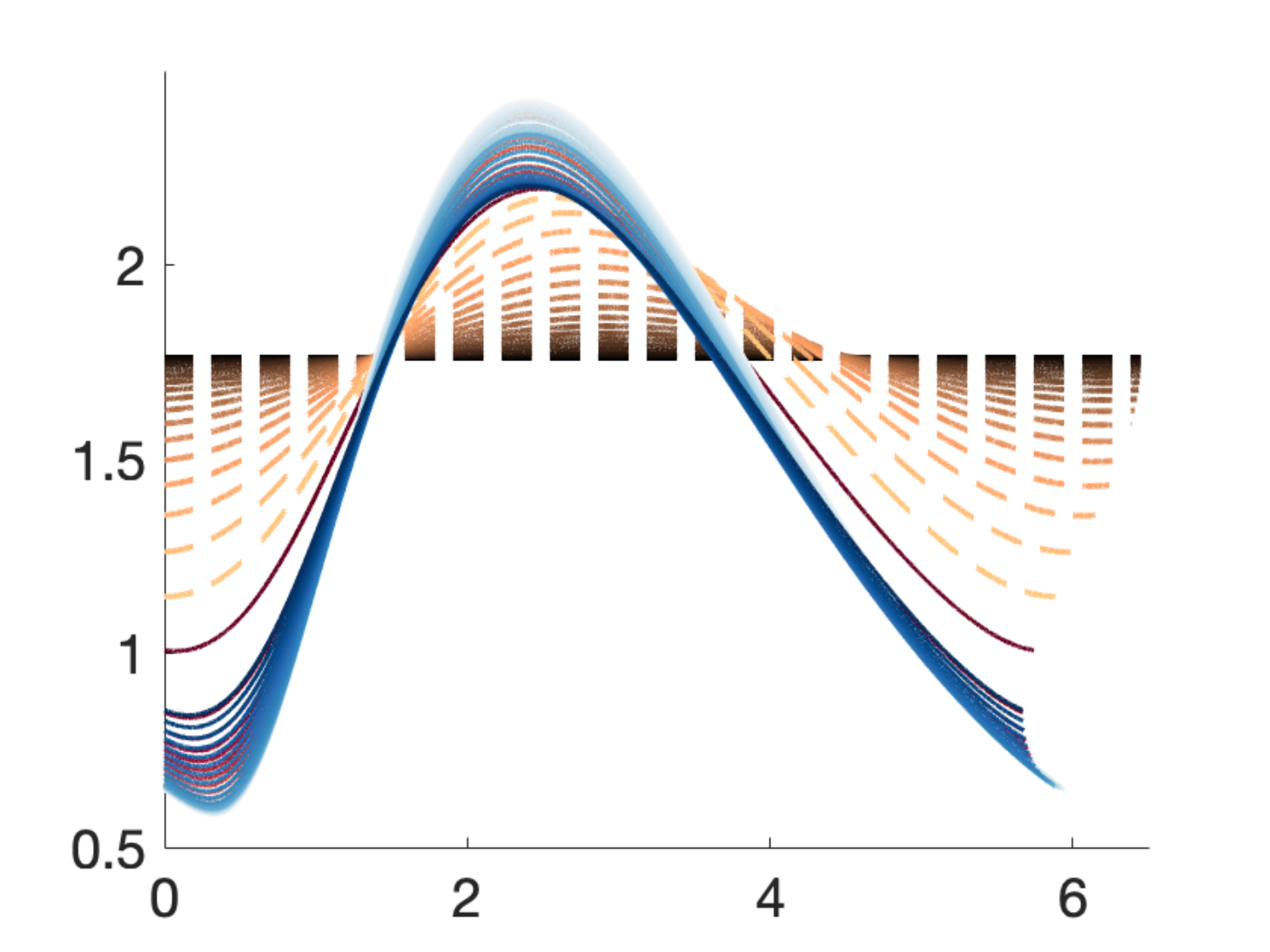}\hspace*{-1em}\includegraphics[scale=0.42]{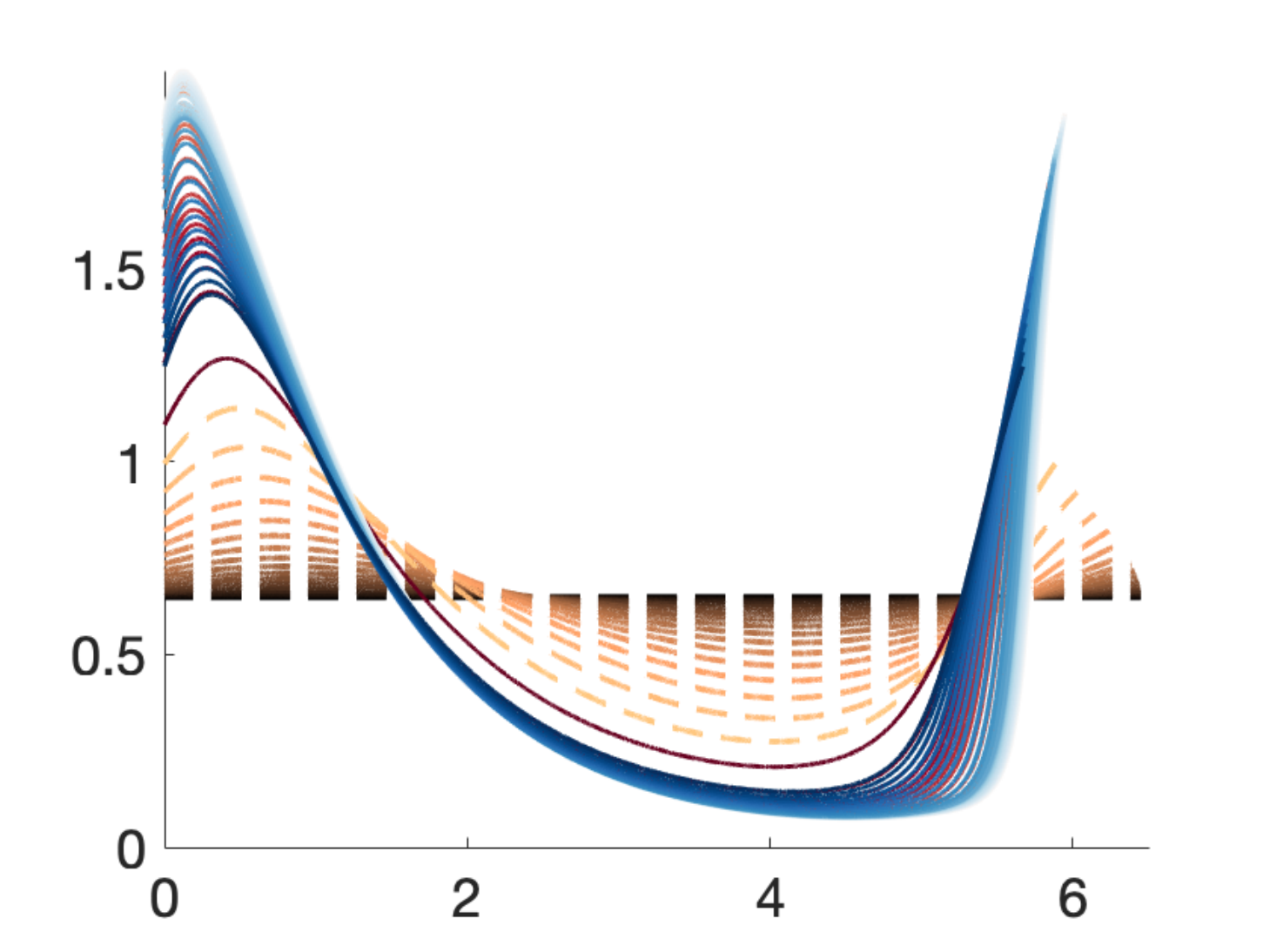}
\put(-310,115){(c)}
\put(-120,115){(d)}
\put(-19,5){$t$}
\put(-198,4){$t$}
\put(-330,115){\rotatebox{90}{$E$}}
\put(-162,115){\rotatebox{90}{$M$}}

\hspace*{-1em}\includegraphics[scale=0.42]{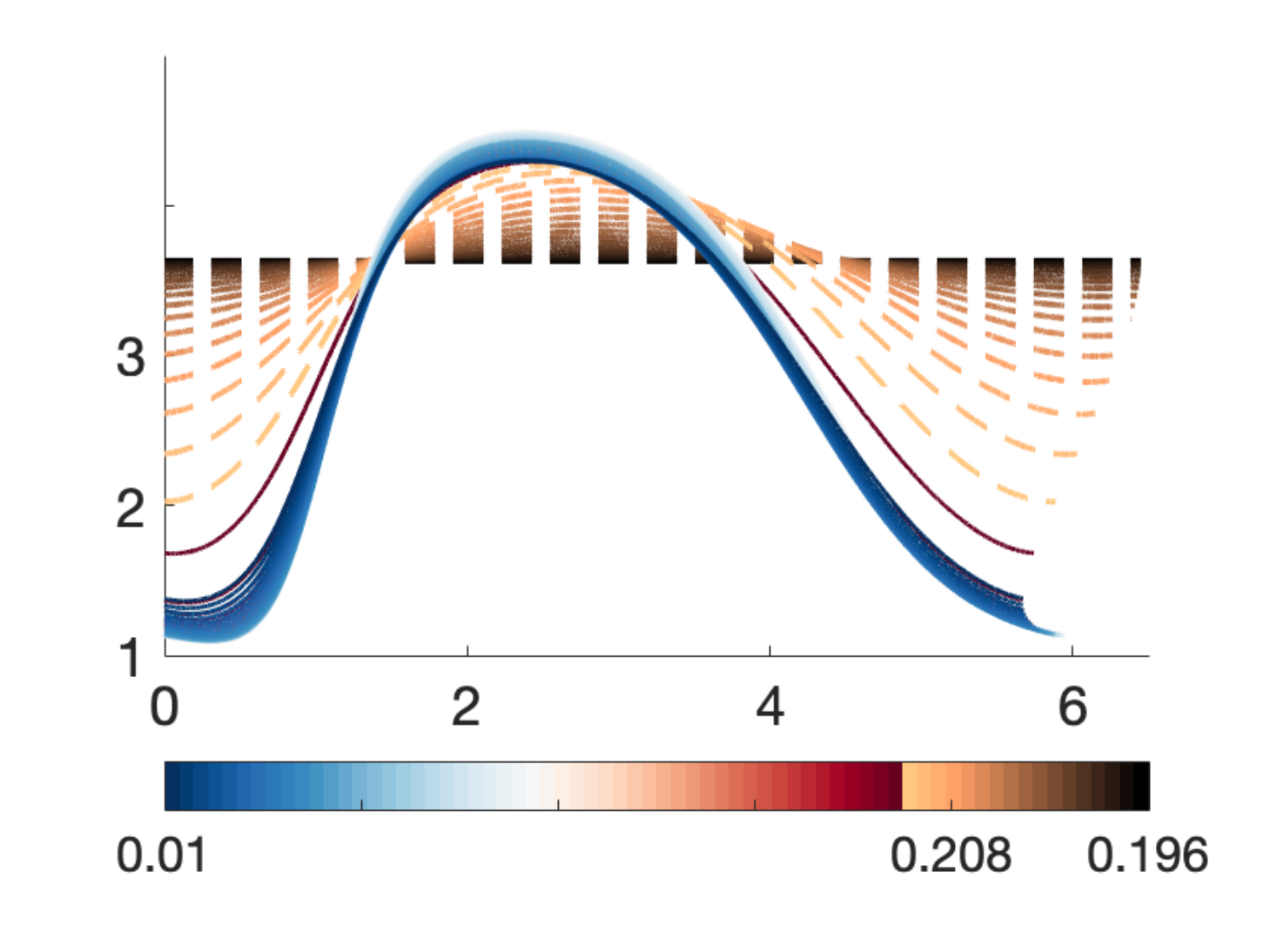}\hspace*{-1em}\includegraphics[scale=0.42]{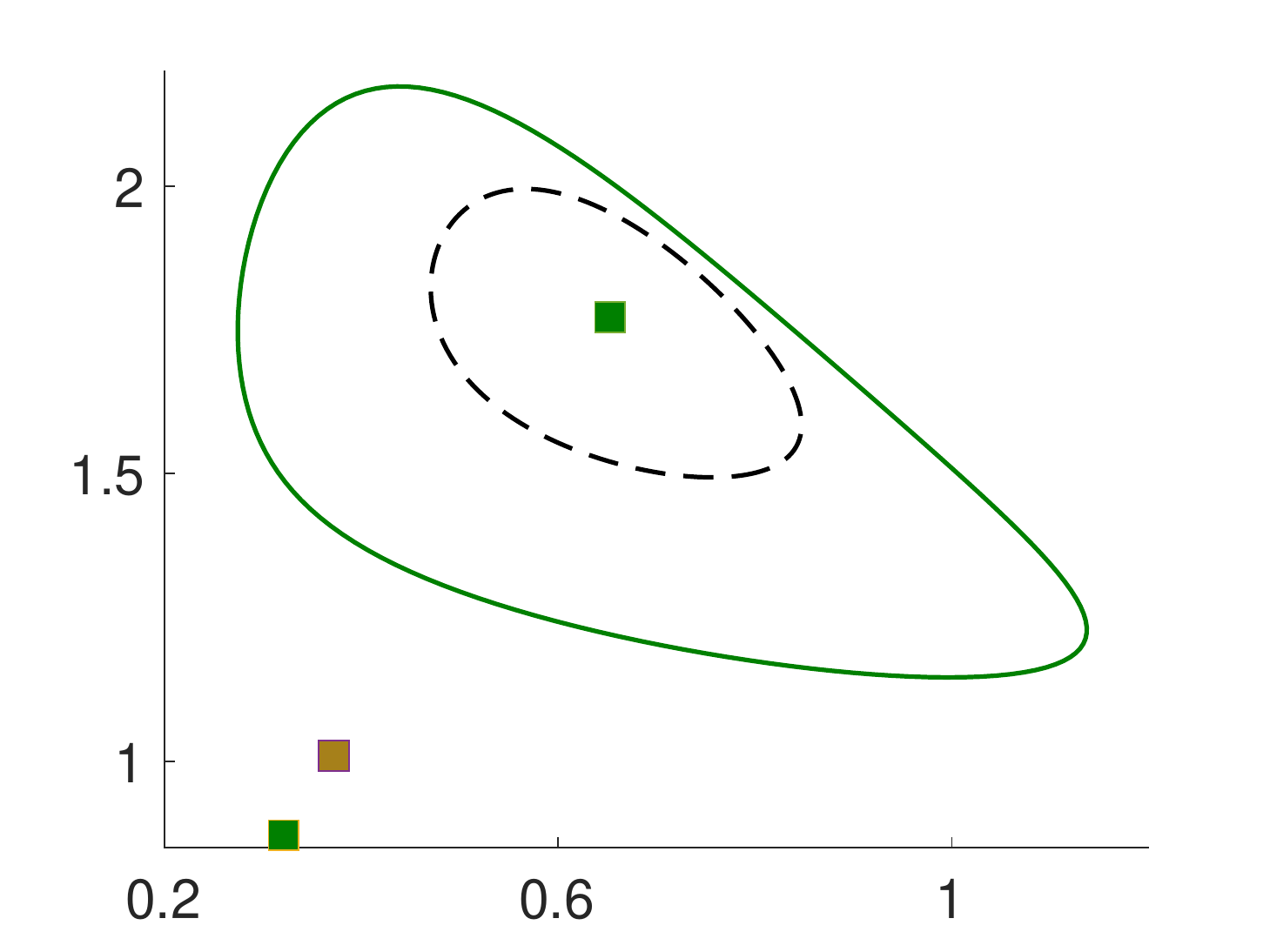}
\put(-150,118){(f)}
\put(-317,118){(e)}
\put(-25,6){$M$}
\put(-162,115){\rotatebox{90}{$E$}}
\put(-190,32){$t$}
\put(-332,103){\rotatebox{90}{$\tau_M(t)$}}
\put(-265,11){$v_M^{min}$}

\caption{Stable and unstable periodic orbits on the branch of periodic orbits emanating from the subcritical Hopf bifurcation in Figure~\ref{threess_onepara_MI}. The colormap indicates
the value of the continuation parameter $v_M^{min}$.
The periodic orbits are shown in (a) amplitude of $E-$component, (b) period, (c) profile in $E$, (d) profile in $M,$ (e) delay $\tau_M$.
Panel (f) shows a projection of phase space onto the $M$-$E$ plane when $v_M^{min} = 0.2$ and tristability occurs.
Periodic orbits are represented by closed curves, and steady states by squares (whose colour indicates the number of unstable eigenvalues as in Figure~\ref{stablehm}).
}
\label{phase_tristability_II}
\end{figure}

The second important difference between the two examples is that the Hopf bifurcation on the upper segment of steady states at $v_M^{min}=0.196$ in Figure~\ref{threess_onepara_MI} is subcritical resulting in a branch of unstable periodic orbits. The change in the criticality of this Hopf bifurcation
between the two examples implies that for some intermediate value $m\in(2,4)$ there is a Bautin bifurcation at which the criticality switches. Bautin bifurcations are well studied \citep{Kuznetsov2004} and
in a two-parameter unfolding generate a branch of fold bifurcations of periodic orbits.

The branch of unstable periodic orbits emanating from the subcritical Hopf bifurcation terminates
in the fold bifurcation of periodic orbits seen in Figure~\ref{threess_onepara_MI} at $v_M^{min}=0.20812$, at which the periodic orbit becomes stable.
As a consequence of the subcritical Hopf bifurcation and fold of periodic orbits there are stable periodic orbits for $v_M^{min}<0.208$ (to the left of the fold bifurcation of periodic orbits) and co-existing stable steady states for $v_M^{min}\in(0.196,285)$ (to the right of the Hopf bifurcation). This creates a
small parameter interval of tristability for $v_M^{min}\in(0.196,0.208)$ between the Hopf bifurcation and fold of periodic orbits bifurcation for which a stable periodic orbit coexists with two stable steady states.
Figure~\ref{phase_tristability_II}(f) shows the dynamics when $v_M^{min}=0.2$ in the tristability region, in a projection of phase space onto the $M$-$E$ plane. The branch of periodic orbits emanating from the Hopf bifurcation at $v_M^{min}=0.19603$ crosses $v_M^{min}=0.2$ twice and both the stable and unstable periodic orbit are shown in the phase portrait.

The other panels of Figure~\ref{phase_tristability_II} show the evolution of the periodic orbit from the Hopf bifurcation on this branch with separate colour maps for the stable and unstable legs of the branch.

For $v_M^{min} \in (0.20812, 0.28543)$ there is bistability between two steady states,
and for $v_M^{min} < 0.19603$ there is bistability between a periodic orbit and a steady state.
There is also a second Hopf bifurcation at $v_M^{min}=0.1228$ which generates small amplitude unstable periodic orbits shown on the bifurcation diagram in Figure~\ref{threess_onepara_MI}.

\subsubsection*{Fold Bifurcation of Periodic Orbits}

For our final example of inducible operon dynamics
we return to the example from Section~\ref{ssec:equilibria} and consider
the one state-dependent delay system \eqref{eq:sysonedel} with the inducible parameter set defined in Table~\ref{table:multsspars}.

\begin{figure}[htp!]
\centering
\includegraphics[scale=0.6]{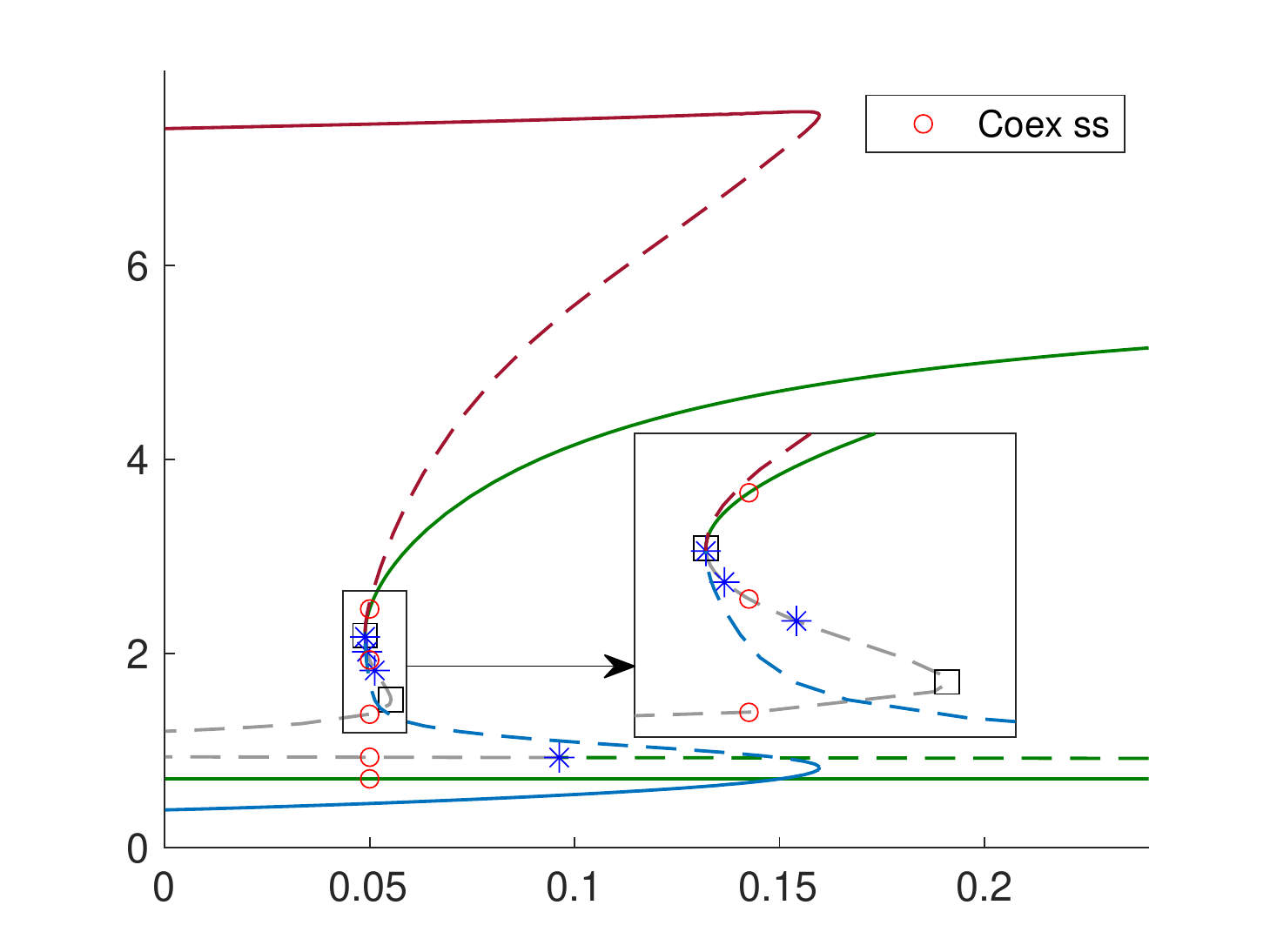}
\put(-230,165){\rotatebox{90}{$E$}}
\put(-40,10){$v_M^{min}$}
\caption{Bifurcation diagram of the model \eqref{eq:sysonedel}
for an inducible system with parameters defined in Table~\ref{table:multsspars}
except $v_M^{min}$ which is taken as the bifurcation parameter.
Red circles denote the five co-existing steady states at $v_M^{min}=0.05$. All other lines and
symbols are defined as in Figures~\ref{fig:threess_onepara_trp} and~\ref{fig:threess_onepara_trp_II}.
}
\label{fivess_onepara}
\end{figure}

\begin{table}[htp!]
	\centering
	\begin{tabular}{| l | l | l | l |}
		\hline
		Bifurcation & Bifurcation parameter value & \rule[-8pt]{0pt}{21pt}\parbox{2.0cm}{Unstable\\ eigenvalues} & $E^*$ value \\
		\hline
		Hopf & $v_M^{min}=0.09624$, period = 4.1286 & 1 to 3 & 0.9270\\
		\hline \hline
		Fold & $v_M^{min}=0.055205$ & 6 to 7 & 1.5229 \\
		\hline
		Hopf & $v_M^{min}=0.051247$, period = 7.1786 & 7 to 5 & 1.8245 \\
		\hline
		Hopf & $v_M^{min}=0.049356$, period = 11.4903 & 5 to 3 & 2.0162 \\
		\hline
		Hopf & $v_M^{min}=0.048868$, period = 23.539 & 3 to 1 & 2.1697 \\
		\hline
		Fold & $v_M^{min}=0.048865$ & 1 to 0 & 2.1830 \\
		\hline \hline
        \rule[-8pt]{0pt}{21pt}\parbox{2.0cm}{Fold of\\ periodic orbits} & $v_M^{min}=0.1597$, period = 7.6675 & 1 to 0 & - \\
        \hline
	\end{tabular}
	\caption{Bifurcation information associated with Figure~\ref{fivess_onepara}.}
    \label{tab:fivess_onepara}
\end{table}

The bifurcation diagram in Figure~\ref{fivess_onepara} extends the diagram previously shown in Figure~\ref{fig:gcont}(b) to
show steady state solutions and periodic orbits along with their stability, as well as Hopf and fold bifurcations. The bifurcations are listed in Table~\ref{tab:fivess_onepara}.

We already saw in Section~\ref{ssec:equilibria} that when $v_M^{min}=v_M^{max}$ and thus both delays are constant, there are three co-existing steady states. Two of these are stable and the intermediate steady state is unstable.
When $v_M^{min}$ is reduced the delay $\tau_M$ becomes state-dependent and a number of bifurcations may occur. The lower stable steady state remains stable  and does not undergo any bifurcations. The intermediate
steady state remains unstable for all $v_M^{min}>0$ but does undergo a Hopf bifurcation. The upper branch of steady states loses stability in a fold bifurcation. There is also another fold bifurcation and several Hopf bifurcations on this branch. Considering all the branches together
there may be up to five co-existing steady states, but as the bifurcation diagram shows, there are only ever one or two co-existing stable steady states.

\begin{figure}[ht!]
\centering
\hspace*{-1em}\includegraphics[scale=0.42]{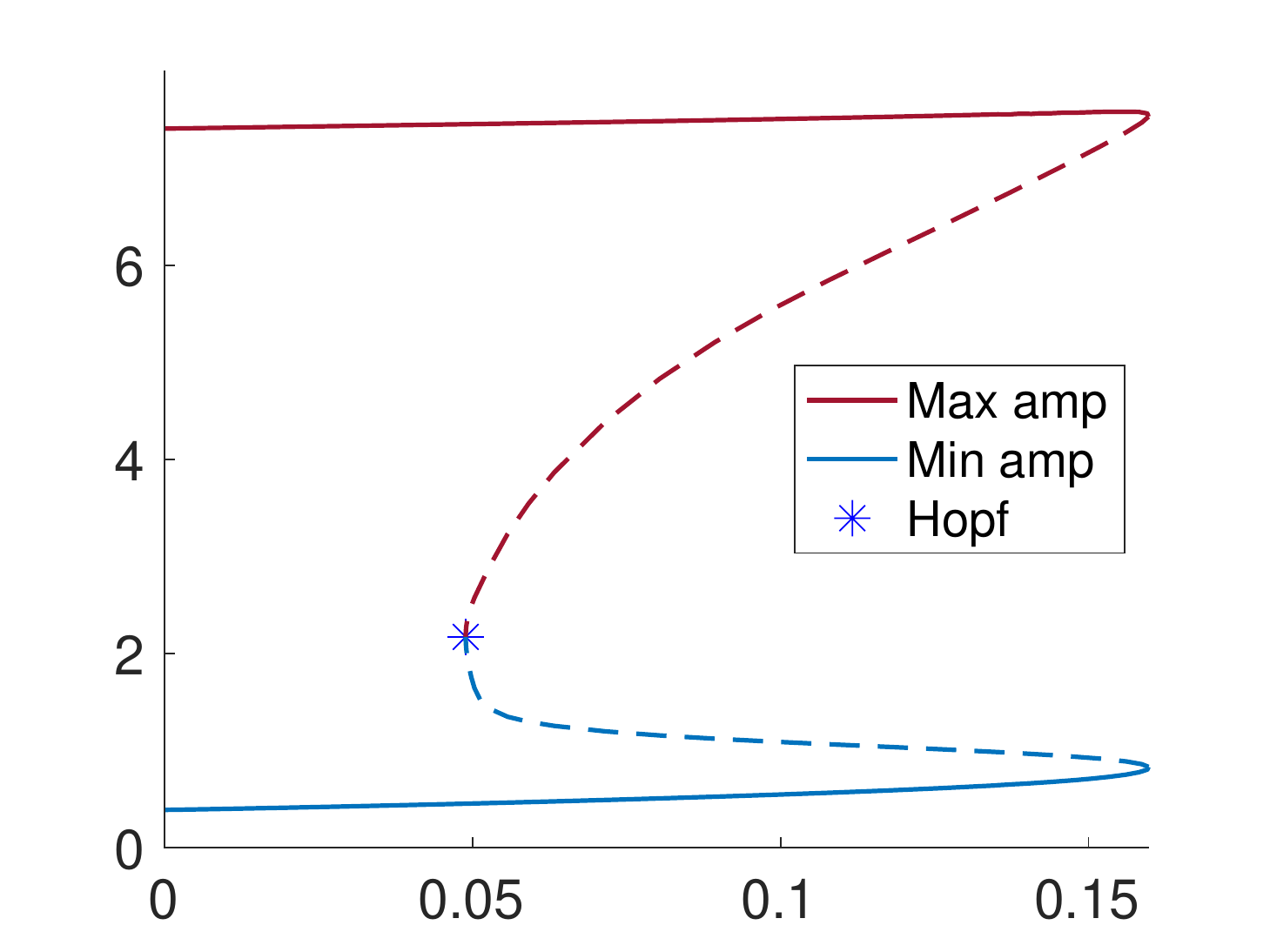}\hspace*{-1em}\includegraphics[scale=0.42]{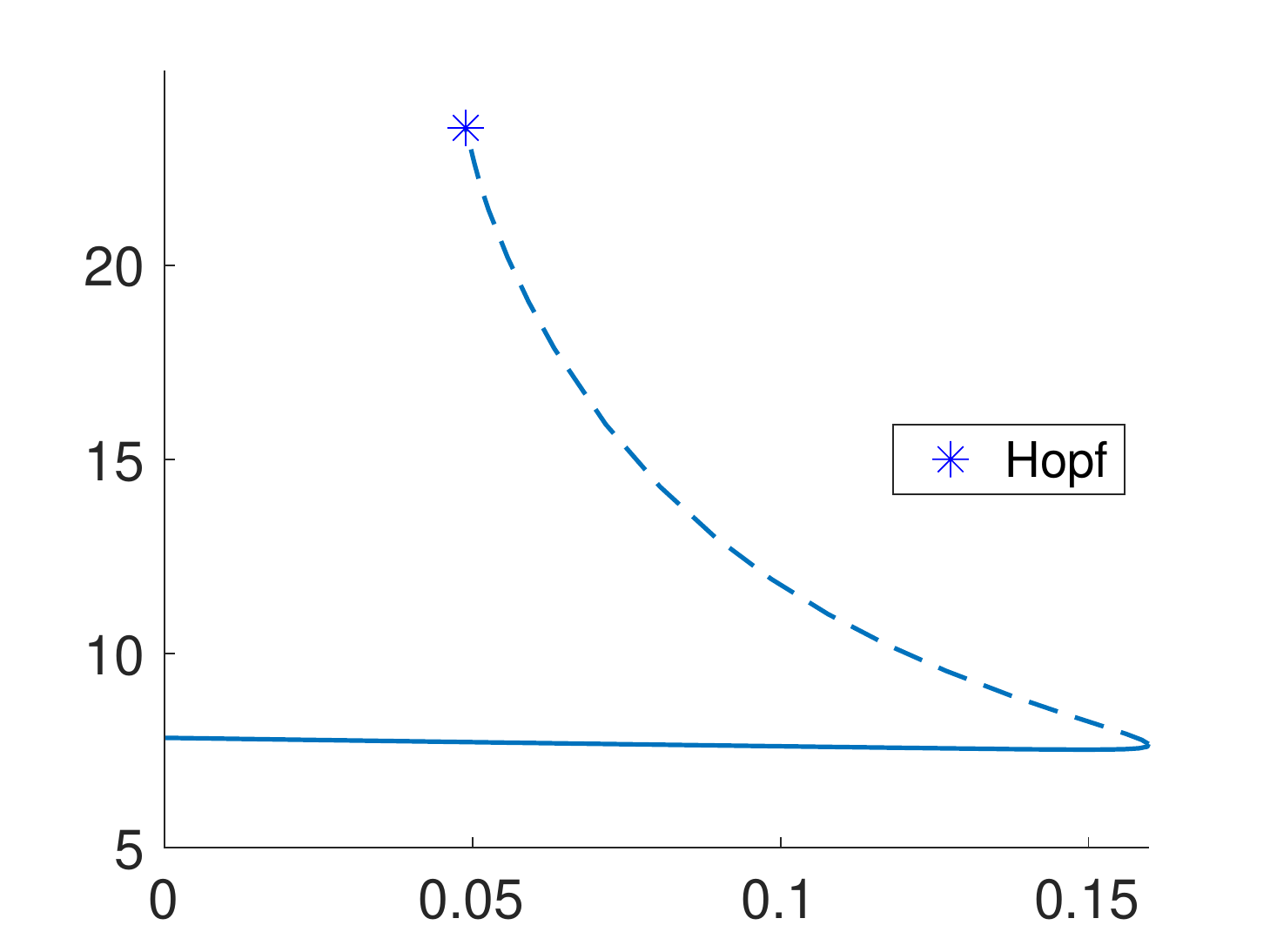}
\put(-319,120){(a)}
\put(-90,115){(b) Period $=T$}
\put(-57,5){$v_M^{min}$}
\put(-225,5){$v_M^{min}$}
\put(-330,114){\rotatebox{90}{$E$}}
\put(-162,113){\rotatebox{90}{$T$}}

\hspace*{-1em}\includegraphics[scale=0.42]{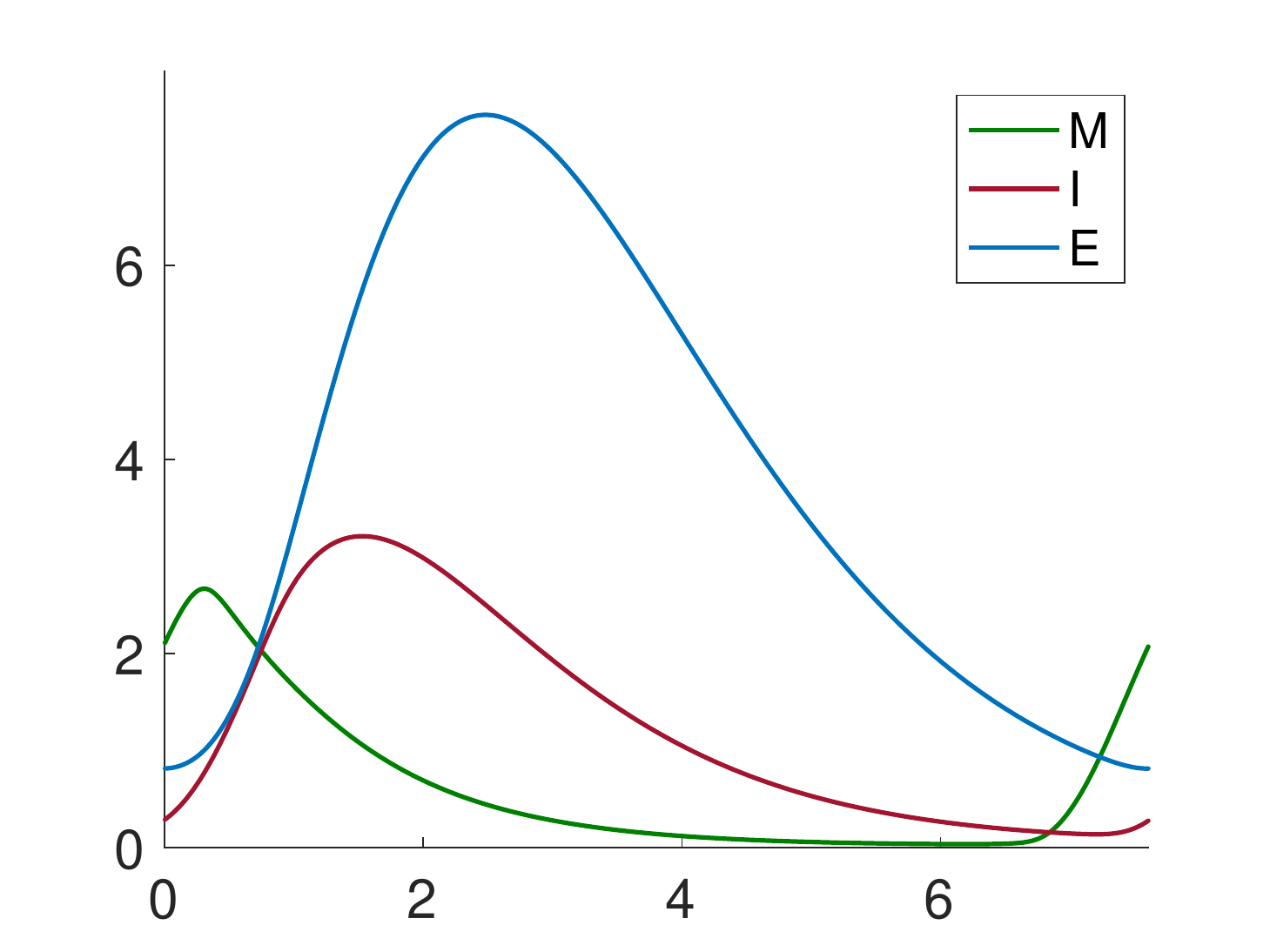}\hspace*{-1em}\includegraphics[scale=0.42]{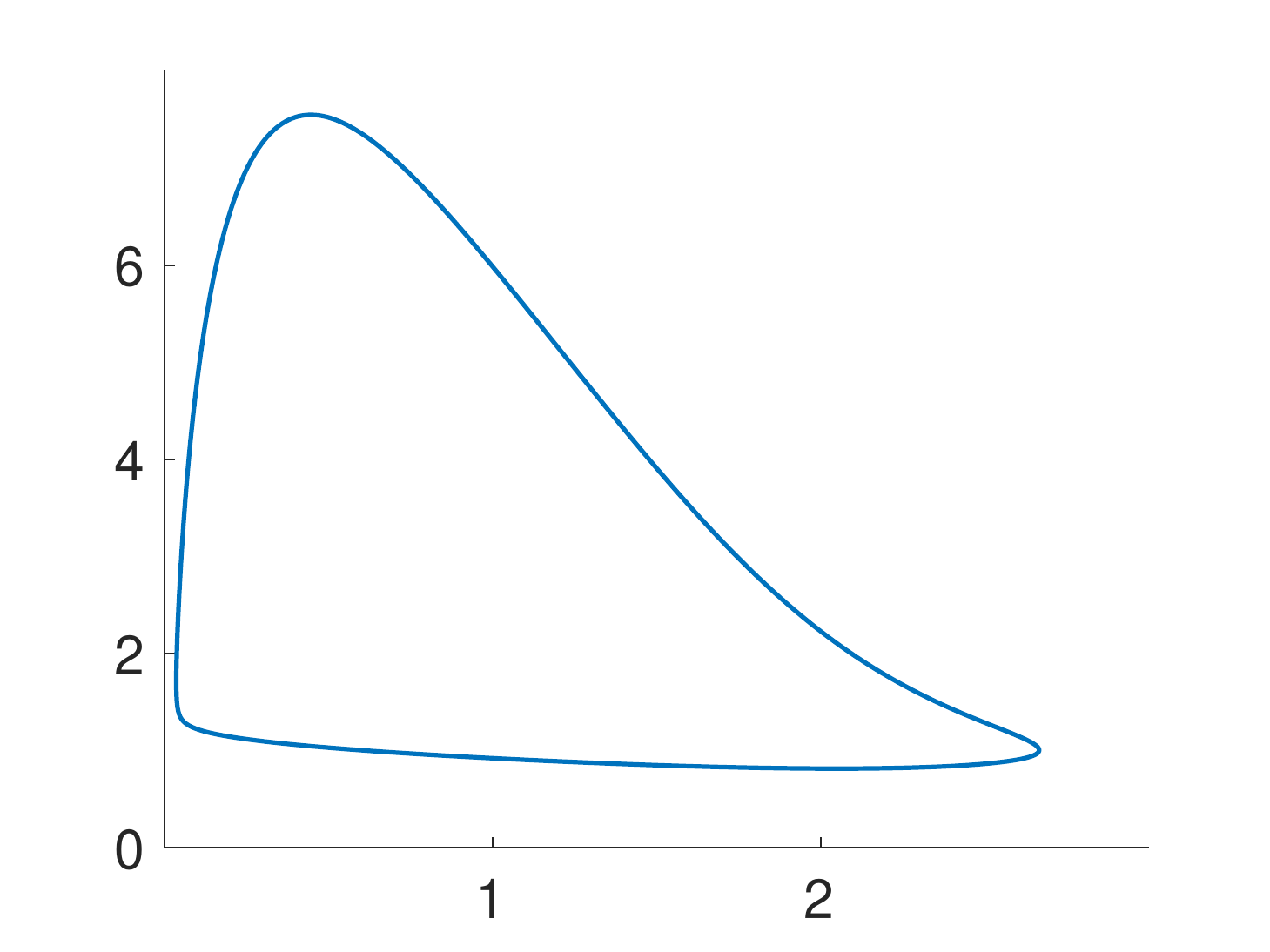}
\put(-320,114){(c)}
\put(-100,113){(d)}
\put(-26,6){$M$}
\put(-162,113){\rotatebox{90}{$E$}}
\put(-190,4){$t$}

\hspace*{-1em}\includegraphics[scale=0.42]{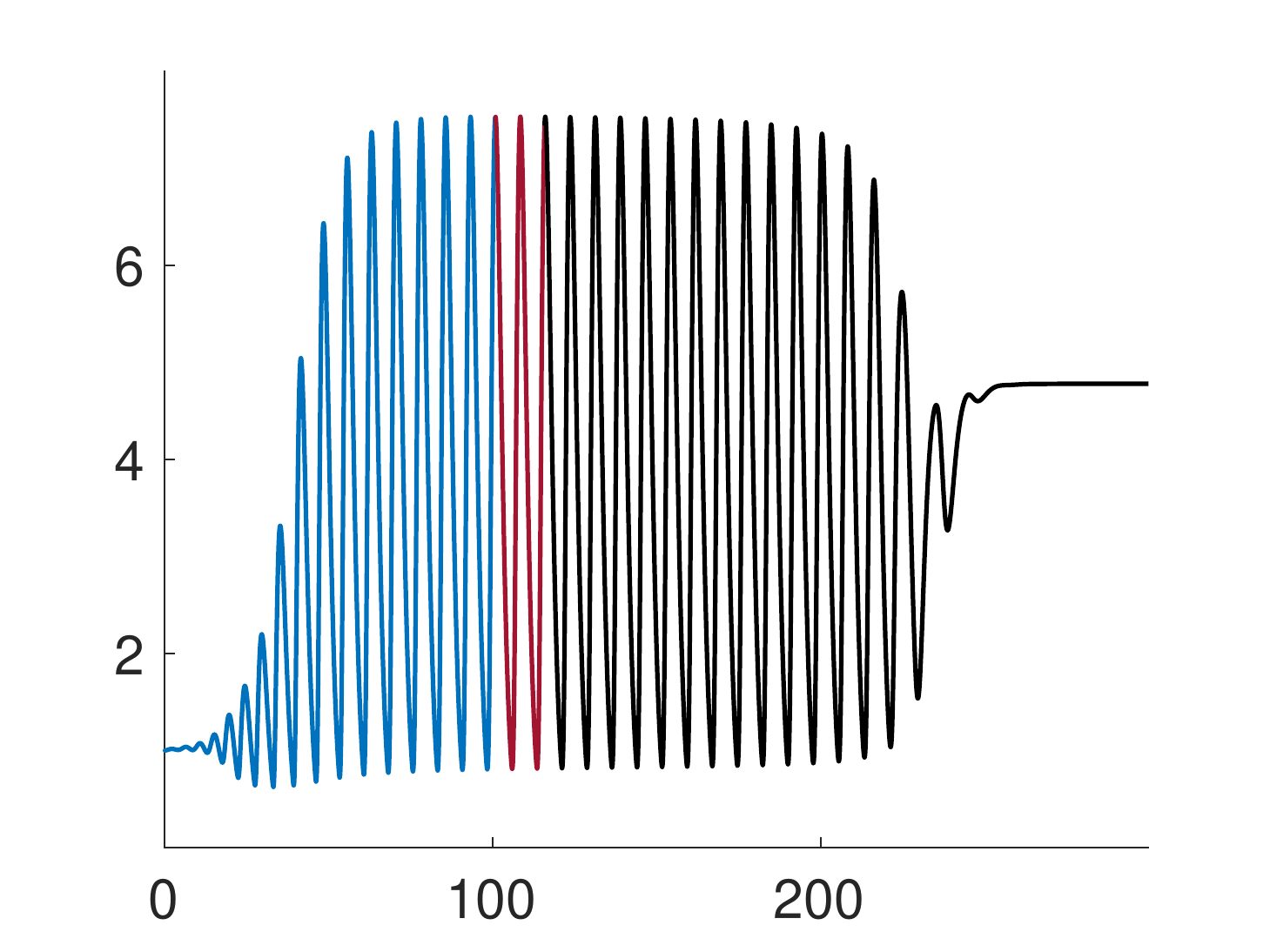}\hspace*{-1em}\includegraphics[scale=0.42]{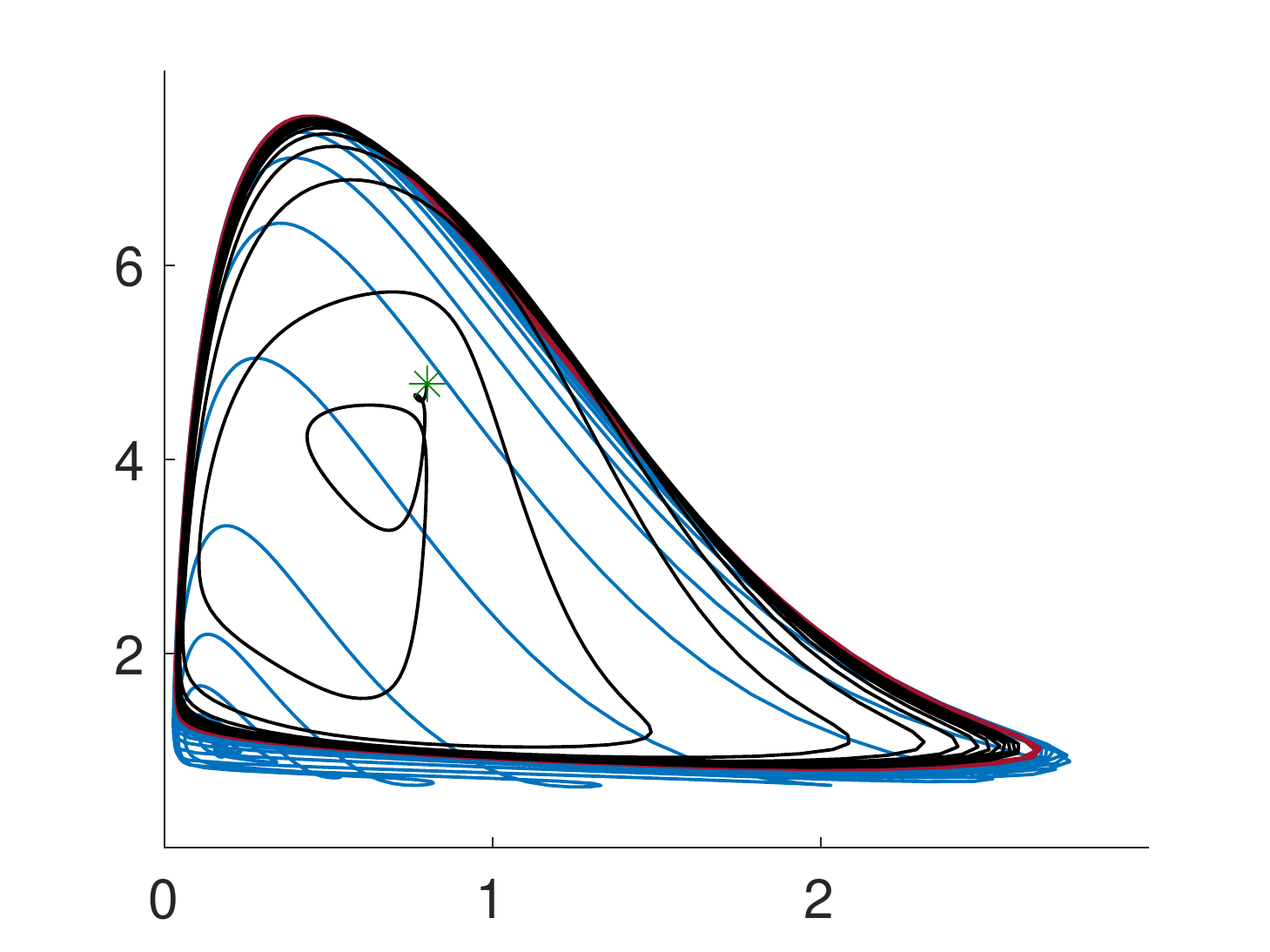}
\put(-318,117){(e)}
\put(-151,117){(f)}
\put(-330,113){\rotatebox{90}{$E$}}
\put(-190,4){$t$}
\put(-162,113){\rotatebox{90}{$E$}}
\put(-26,5){$M$}

\caption{(a) Amplitude and (b) Period for the branch of periodic orbits from
Figure~\ref{tab:fivess_onepara}. (c) The solution components and (d) the projection into phase space
of the periodic orbit at the fold bifurcation $v_M^{min}=0.1597$.
(e) $E$ component and (f) phase space projection for a
 simulation with $v_M^{min}=0.1605$ and initial function close to the intermediate steady state.}

\label{fig:unstabschopf}
\end{figure}
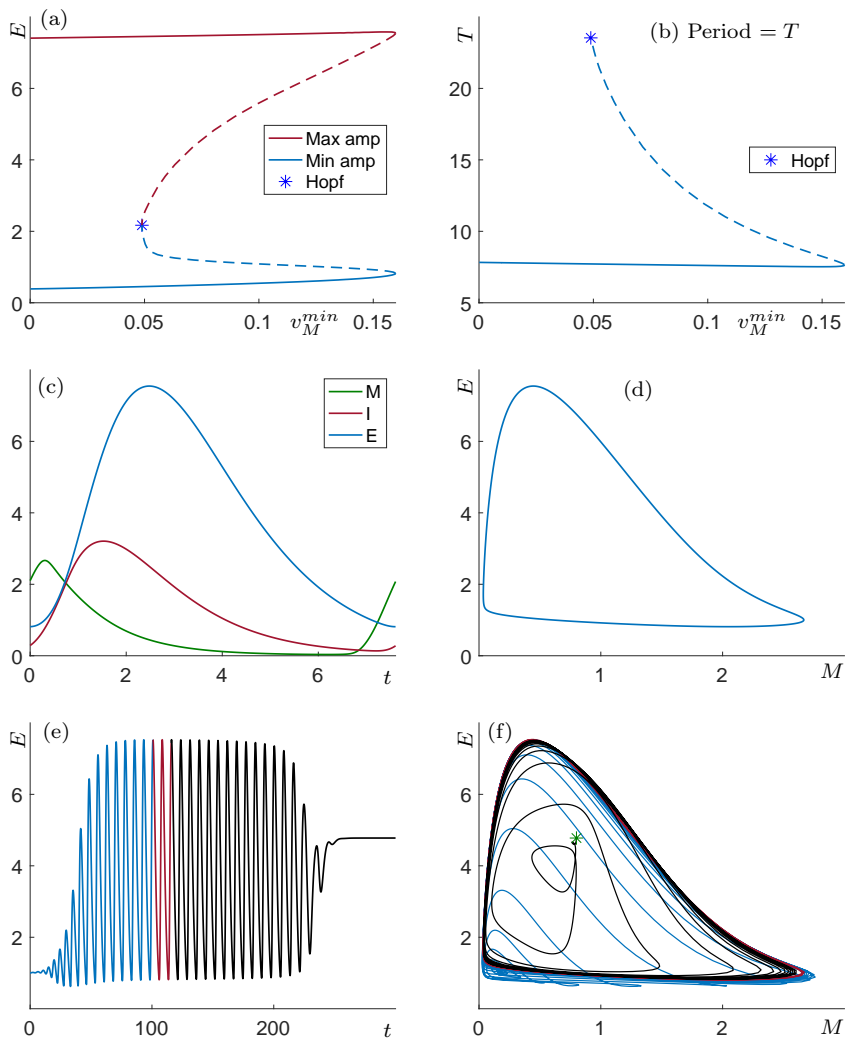

The fold bifurcation at which the steady state loses stability (at $v_M^{min}=0.048865$)
is immediately followed by a Hopf bifurcation (at $v_M^{min}=0.048868$), indicating that this inducible operon is close to a zero-Hopf bifurcation (in Section~\ref{sec:repexamples} we  inferred existence of a zero-Hopf bifurcation for a repressible operon).

The branch of periodic orbits emanating from the Hopf bifurcation is shown in
Figure~\ref{fig:unstabschopf}.
The bifurcation is a supercritical Hopf bifurcation from an unstable steady state, which gives rise to a branch of unstable periodic orbits bifurcating to the right.
The amplitude and period of these orbits are shown in Figure~\ref{fig:unstabschopf}(a) and (b). Interestingly, moving along the branch away from the Hopf bifurcation the period decreases as the amplitude increases until there is a fold bifurcation of periodic orbits at $v_M^{min}=0.1597$
creating a segment of stable periodic orbits on the branch. The periodic orbit
at the fold bifurcation is shown in Figure~\ref{fig:unstabschopf}(c) and (d).
For $v_M^{min}>0.1597$ there is no longer a periodic orbit, but it is still possible to have transient oscillatory dynamics. Figure~\ref{fig:unstabschopf}(e) and (f) show an example of this for
$v_M^{min}>0.1605$, where an initial function close to the unstable intermediate steady state generates a solution with large oscillations for 200 time units before the solution converges to the stable steady state. When the phase space projection of this solution in Figure~\ref{fig:unstabschopf}(f) is compared to the periodic orbit at the fold of periodic orbits (in Figure~\ref{fig:unstabschopf}(d)) it is clear that we are seeing a ghost of the periodic orbit.

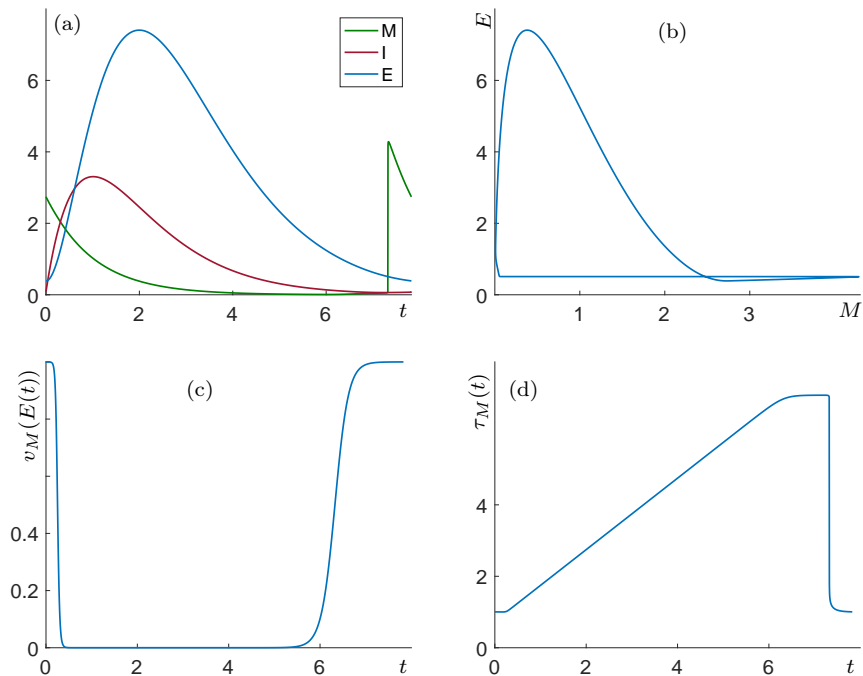
\begin{figure}[htp!]
\hspace*{-1em}\includegraphics[scale=0.42]{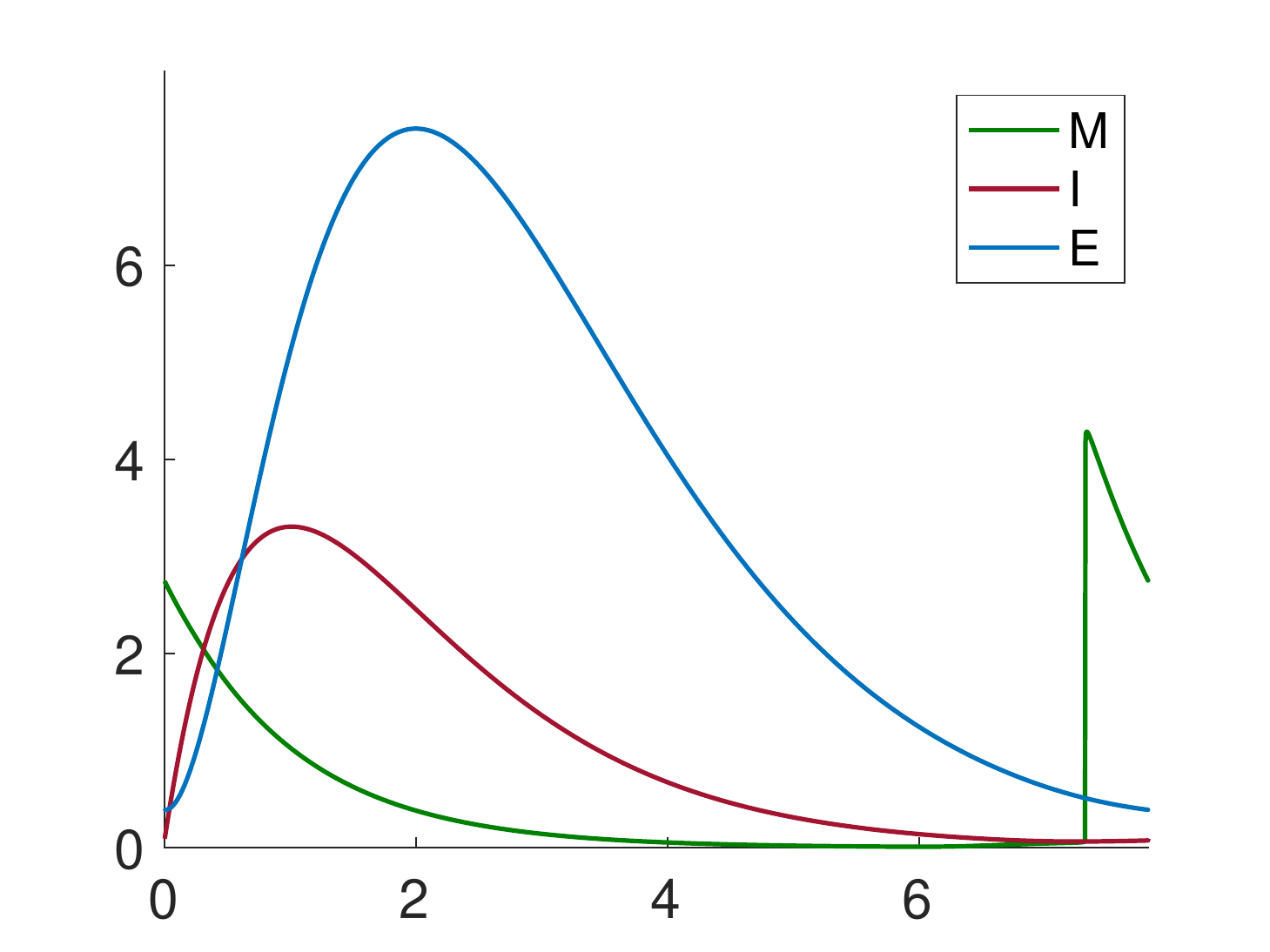}\hspace*{-1em}\includegraphics[scale=0.42]{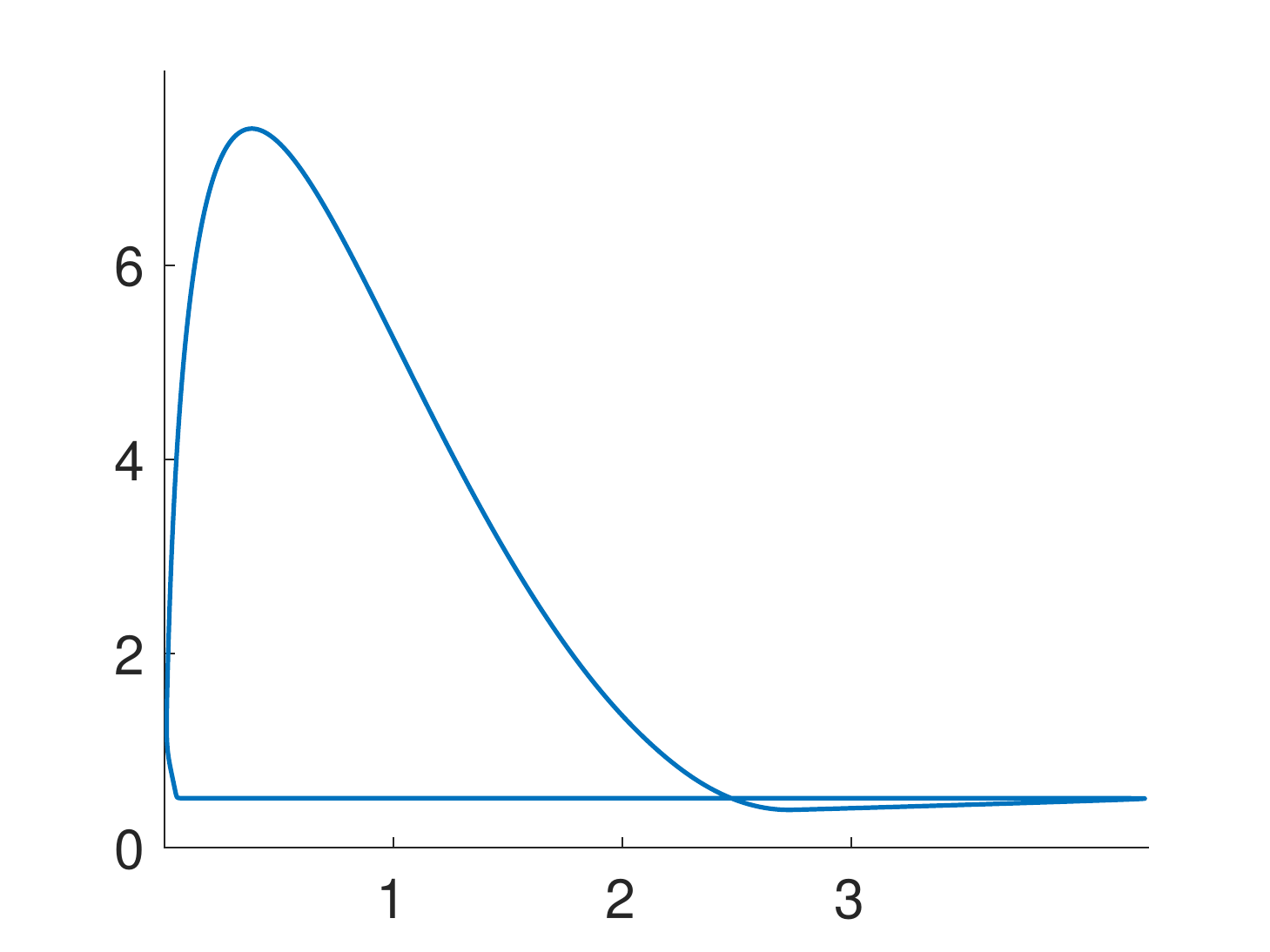}
\put(-320,115){(a)}
\put(-93,112){(b)}
\put(-190,6){$t$}
\put(-25,6){$M$}
\put(-161,115){\rotatebox{90}{$E$}}

\hspace*{-1em}\includegraphics[scale=0.42]{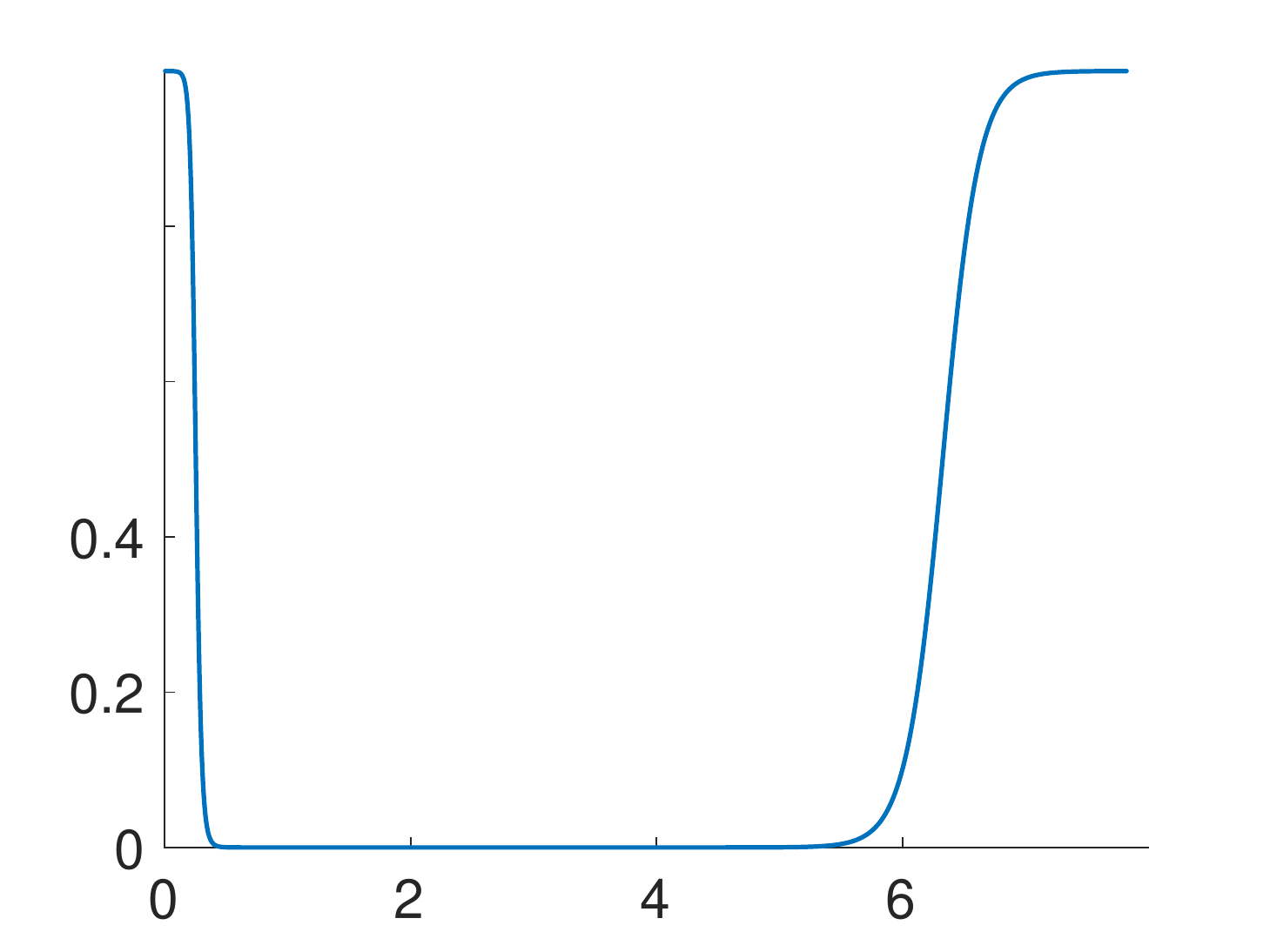}\hspace*{-1em}\includegraphics[scale=0.42]{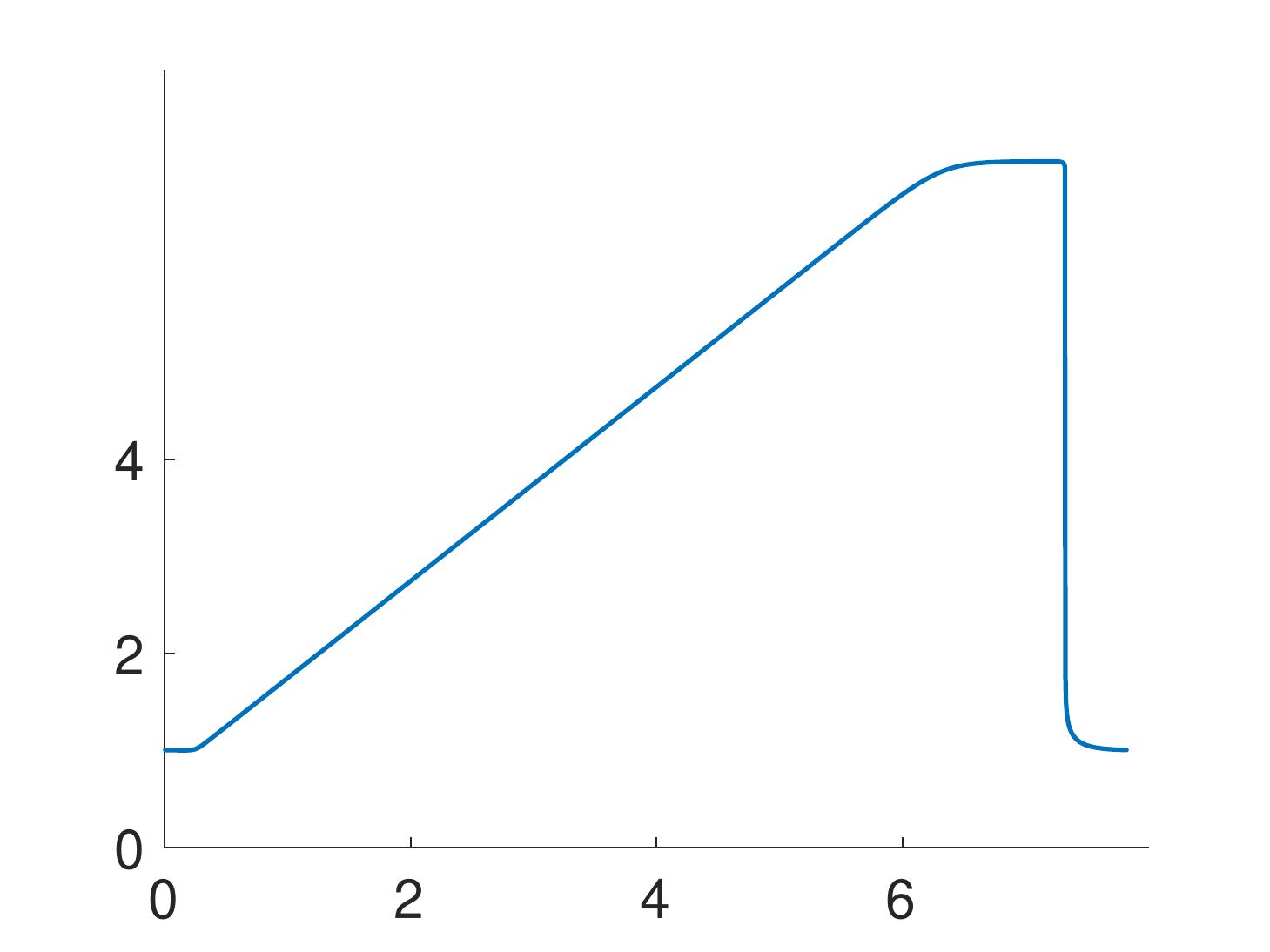}
\put(-270,110){(c)}
\put(-149,110){(d)}
\put(-22,5){$t$}
\put(-190,5){$t$}
\put(-164,97){\rotatebox{90}{$\tau_M(t)$}}
\put(-333,85){\rotatebox{90}{$v_M(E(t))$}}

\caption{
Stable periodic orbit at $v_M^{min}=0.0001$ seen in Figure \ref{fivess_onepara}. The periodic orbit is shown in
(a) solution in all three components over one period,
(b) projection of phase space into $M$-$E$ plane,
(c) velocity $v_M(E(t))$ and
(d) delay $\tau_M(t)$.
}
\label{sim_vmin00001}
\end{figure}

For $v_M^{min}\in(0.048865,0.1597)$ the branch of stable periodic orbits coexists with two stable steady states, creating another example of tristability of solutions. The stable periodic orbit at the left end of the branch with $v_M^{min}=0.0001$ is shown in Figure~\ref{sim_vmin00001}.


The transcription velocity is essentially zero for nearly all of the period, with just a short burst of transcription when $E$ is close to its minimum.
This sudden release of mRNA gives the $M$ component of the solution the characteristic form of a relaxation oscillator, even though the other components of the solution are smooth.

The variation of the delay as a function of time seen in Figure~\ref{sim_vmin00001}(d) shows that the delay  is very far  from being  constant.  The delay is  increasing linearly on the segment of the orbit for which the transcription velocity is zero, and so no transcripts are being completed. During this time the effector $E$  concentration is high and thus the transcription initiation rate $f(E)$ is high; at the same time though, the delay $\tau(E)$ is also increasing. Only when the concentration of the effector $E$ drops sufficiently, does transcription proceed during the last quarter of the period.




Finally, we remark that it is highly delicate to numerically compute a branch of periodic orbits emanating from a Hopf bifurcation close to a co-dimension two zero-Hopf bifurcation in the system \eqref{eq:sysonedel} with
the threshold integral discretized as in Section~\ref{subsec:discretization}. For this reason  we were not able to compute this branch starting from the Hopf bifurcation. Instead, noting that for small values of $v_M^{min}$ there are three steady states, but only the lower one is stable, we performed a numerical simulation  of dynamics as described in Section~\ref{ssec:ddesd} starting close to the upper unstable steady state. This simulation converged to the stable periodic orbit. This periodic orbit was then continued in DDE-BIFTOOL to find the fold bifurcation of periodic orbits and follow the branch of unstable periodic orbits back to the Hopf bifurcation.

\subsection{Two State-Dependent Delays}

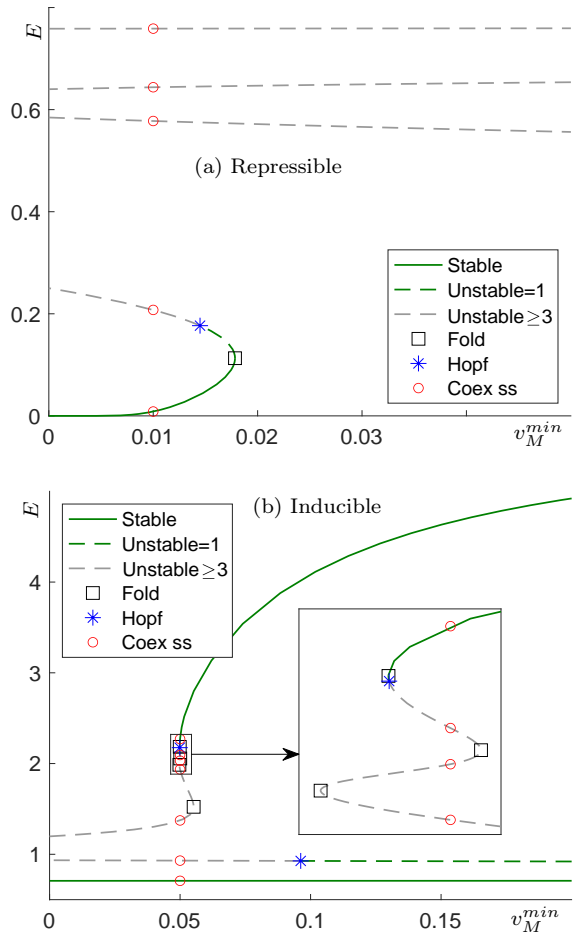
\begin{figure}[thp!]
\centering
\vspace*{-2ex}
\includegraphics[scale=0.6]{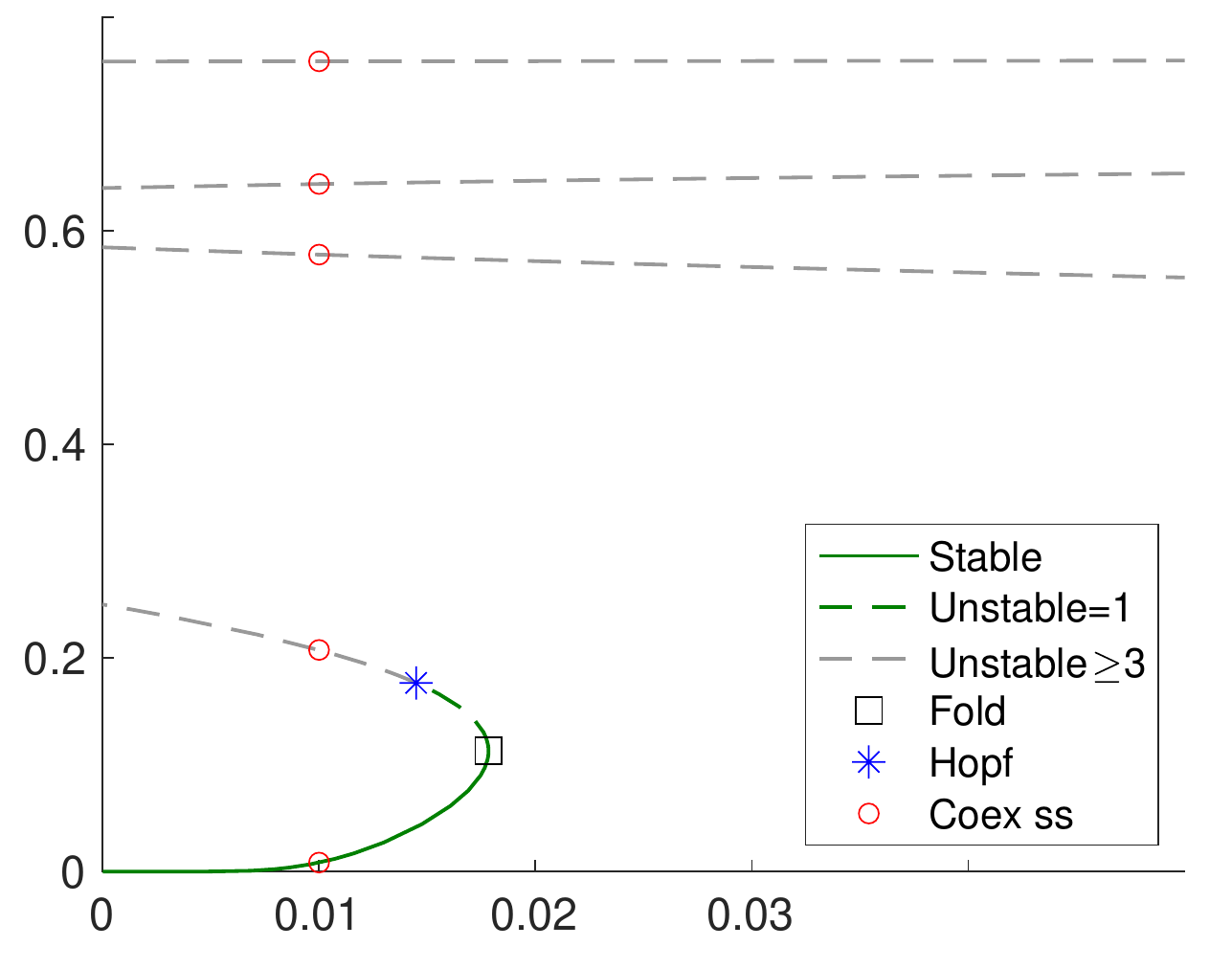}
\put(-213,160){\rotatebox{90}{$E$}}
\put(-30,9){$v_M^{min}$}
\put(-150,110){(a) Repressible}
\vspace*{2ex}
\includegraphics[scale=0.6]{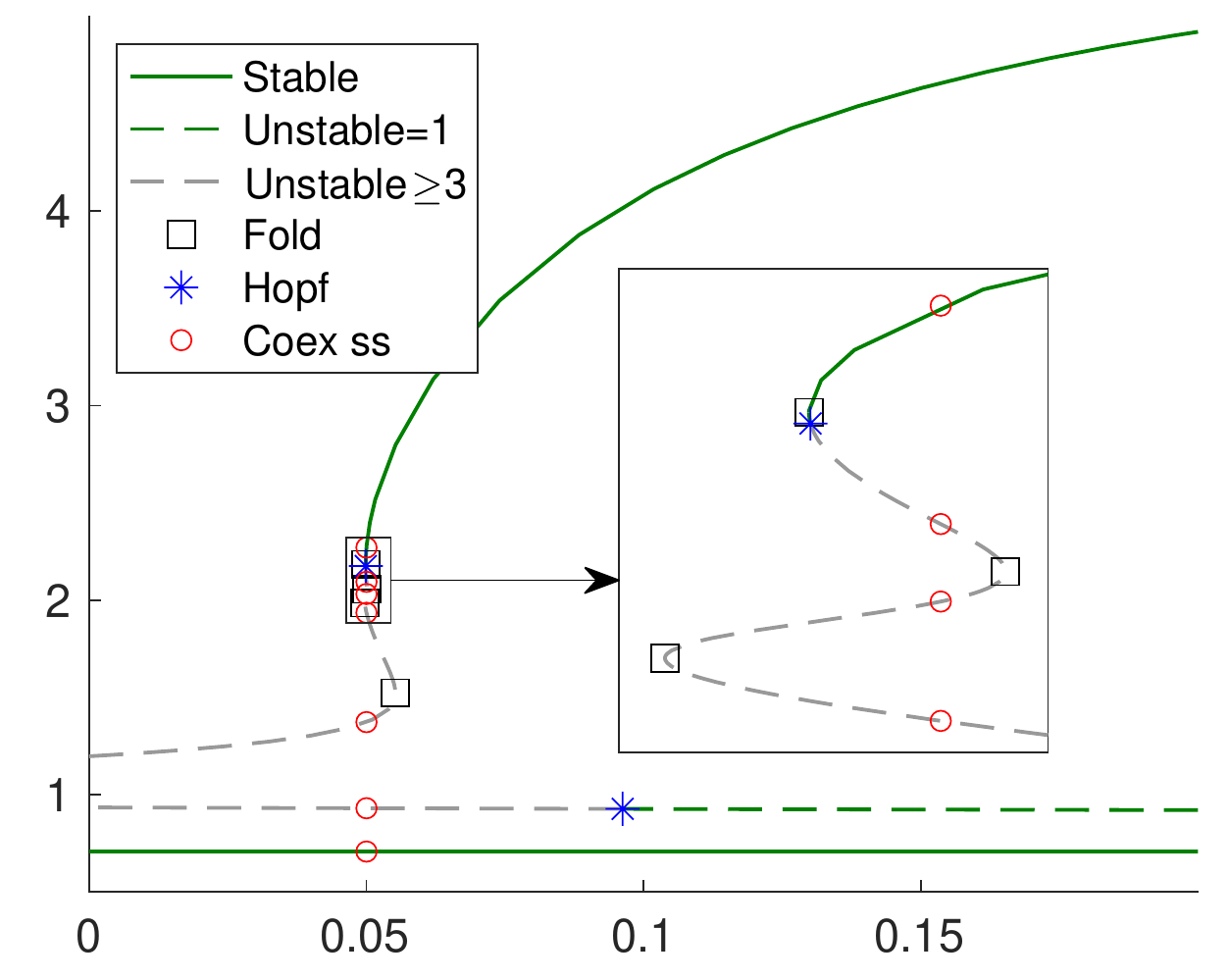}
\put(-28,6){$v_M^{min}$}
\put(-125,162){(b) Inducible}
\put(-211,160){\rotatebox{90}{$E$}}
\caption{Bifurcation diagram of the model \eqref{eq:mrna-delay-var}-\eqref{eq:delay by stateI}
with $\tau_M$ and $\tau_I$ both state-dependent.
(a) Repressible case with same parameters as in Figure~\ref{fig:gcont}(c).
(b) Inducible case with same parameters as in Figure~\ref{fig:gcont}(d).
Symbols and lines are as defined in Figure~\ref{fig:threess_onepara_trp}, except for the red circles which denote the co-existing steady states shown in Figure~\ref{fig:gmultsols}(c) and (d).
}
\label{fig:2delstab}
\end{figure}

\begin{table}[htp!]
	\centering
	\begin{tabular}{| l | l | l | l |}
		\hline
		Bifurcation & Bifurcation parameter value & Unstable eigenvalues& $E^*$ value \\
		\hline
		Fold & $v_M^{min}=0.017832$ & 1 to 0 & 0.1130 \\
		\hline
		Hopf & $v_M^{min}=0.014483,$ period = 34.8874 & 3 to 1 & 0.1768 \\
		\hline
	\end{tabular}
	\caption{Bifurcation information associated with Figure~\ref{fig:2delstab}(a). }
\label{tab:bifs2delrep}
\end{table}

\begin{table}[htp!]
	\centering
	\begin{tabular}{| l | l | l | l |}
		\hline
		Bifurcation & Bifurcation parameter value & Unstable eigenvalues& $E^*$ value \\
		\hline
		Hopf & $v_M^{min}=0.096224,$ period = 4.1287  & 1 to 3 & 0.9270 \\
		\hline \hline
		Fold & $v_M^{min}=0.055205$ & 6 to 7 & 1.5230 \\
		\hline
		Fold & $v_M^{min}=0.049743$ & 15 to 14 & 1.9858 \\
		\hline
		Fold & $v_M^{min}=0.05006$ & 12 to 13 & 2.0558 \\
		\hline
		Hopf & $v_M^{min}=0.049879,$ period = 23.5795 & 3 to 1 & 2.1751 \\
		\hline
		Fold & $v_M^{min}=0.049877$ & 1 to 0 & 2.1841 \\
		\hline
	\end{tabular}
	\caption{Bifurcation information associated with Figure~\ref{fig:2delstab}(b).
}
\label{tab:bifs2delind}
\end{table}

We briefly return to the model \eqref{eq:mrna-delay-var}-\eqref{eq:delay by stateI} with two state-dependent delays. Figure~\ref{fig:2delstab} shows the stability of the steady states and also the steady state bifurcations for the two examples first considered in
Figure~\ref{fig:gcont}(c) and (d) in Section~\ref{ssec:equilibria}. The principal bifurcations are listed in Tables~\ref{tab:bifs2delrep} and~\ref{tab:bifs2delind}.
For the inducible case our numerical code found many pairs of complex conjugate characteristic values crossing the imaginary axis indicating the possibility of many Hopf bifurcations. In Table~\ref{tab:bifs2delind} we only list with Hopf bifurcations that generate steady states with three or fewer unstable eigenvalues, as Hopf bifurcations with more unstable directions will never change the stability of the stability of the steady state and will only generate periodic orbits with multiple unstable Floquet multipliers. 

Compared to the examples from the preceding sections we see that allowing the second delay to also be state-dependent can result in additional co-existing steady states (consistent with \eqref{eq:numss}), however these extra steady states are unstable and there do not seem to be additional stable
invariant objects.
In Figure~\ref{fig:2delstab}(a) there are four unstable and one stable equilibrium, suggesting existence of one or more stable periodic orbits. In Figure~\ref{fig:2delstab}(b) for low $v_{min}$ there are two unstable and one stable equilibrium, again suggesting existence of a stable periodic orbit.

%% file: discussion.tex
\section{Discussion and summary}\label{sec:disc}

This paper studies the Goodwin model of operon dynamics in the presence of state dependent delays in the processes of transcription and translation. The dependence of delays on the state of the system was considered previously \citep{monk03,verdugo07,verriest17,wang21} and justified by the existence of transportation delays  of mRNA export through the nuclear membrane.  We argue that the availability of building blocks for mRNA and protein synthesis, as well as traffic jams of transcribing polymerases and translating ribosomes, affect the velocity of these processes and depend on the state of the cell. In contrast with membrane transportation delays, these effects may influence also  prokaryotic operons.

The focus of the paper is on exploring potential operon dynamics in the presence of state dependent delays both in transcription and translation and contrasting its richness with the dynamics of constant delay systems. We consider two different situations of repressible and inducible operon.

In the repressible case with state dependent transcriptional delay $\tau_M$ and constant translational delay $\tau_I$ we find bistability either between  two steady states or a steady state and a stable periodic orbit (Figures~ \ref{fig:gcont}, \ref{fig:threess_onepara_trp}   and Figures~\ref{fig:threess_onepara_trp_III}).
This periodic orbit has the characteristic of a relaxation oscillator, where the velocity of transcription is very low for a majority of the period, only to produce a brief 'spike' of transcription which results in subsequent spikes of the translated protein and the effector protein, Figure~\ref{fig:vMmin02_trp_II}.

We also found compelling evidence for the existence of complicated dynamics in the repressible case. In particular,  by tracing the trajectory of the stable periodic orbit just before it disappears, we found that it passes  very close to a saddle point which it approaches along the dominant linearized stable direction and leaves along the
one-dimensional unstable manifold (Figure~\ref{stablehm}). An unstable periodic orbit behaviours very similarly at a near identical paranmeter value, except it
leaves along the
one-dimensional unstable manifold in the opposite direction(Figure~\ref{unstablehm}).
This suggests that these periodic orbits disappear in a homoclinic bifurcation.  Furthermore,  linearization at the saddle point at nearby parameters shows that the dominant stable real eigenvalue  becomes a dominant complex pair (Figure~\ref{3DLtrans})) whose magnitude satisfies the assumptions guaranteeing  existence of Shilnikov type chaotic set  \citep{Shilnikov65,Kuznetsov2004}.
We found two parameter values near each other such that at one the stable periodic orbit loses stability through a Hopf bifurcation and at the other at a fold bifurcation and inferred that at some intermediate parameter value those two bifurcations will coincide in a codimension-two zero-Hopf bifurcation.

When both $\tau_M$ and $\tau_I$ are variable we show that the system can admit 5 steady states, Figures~\ref{fig:gcont}(c) and~\ref{fig:2delstab}(a). Further  exploration of the range of  dynamics in this case is left for future studies.

In the inducible case when the delays are constant we can have either a single steady state or  bistability between two steady states. The dynamics are richer in presence of variable delays.

We start with two situations when   $\tau_M$ and $\tau_I$ are constant: either (1) there is a unique steady state or (2) there are 3 steady states, two of which are stable. We then allow $\tau_M$ to vary while   $\tau_I$ remains constant. In case (1) we find coexistence of the steady state
with a  stable periodic orbit, or another steady state, Figure~\ref{threess_onepara};  and in case (2) we find there are an additional 2 equilibria for a total of 5. This results in tristability between two equilibria and a period orbit (Figures~\ref{threess_onepara_MI} and ~\ref{fivess_onepara}).
The stable periodic orbit in (1) has features of a relaxation oscillator; the velocity of transcription remains close to zero for the majority of the period, only to show a rapid increase in a short burst Figures~\ref{sim_vmin00001}. 

When both $\tau_M$ and $\tau_I$ are variable the system can have up to 7 steady states,
Figures~\ref{fig:gcont}(d) and~\ref{fig:2delstab}(b).
Further  exploration of the dynamics in this case is again left for future studies.

The presence of relaxation type oscillations in both inducible and repressible operons provides an intriguing source of pulse generation on a subcellular level. Among several periodic behaviors that have been experimentally observed we focus on transcriptional bursting~\citep{chong2014,Lenstra2016,Chubb2020}. The production of mRNA from some genes does not produce a steady stream of mRNAs, but rather proceeds in bursts of production interspersed by periods of quiescence. The most popular model that describes this phenomena is the \textit{ telegraph model}~\citep{Peccoud1995} where, during the periods when the transcription factor (TF) binds to promoter RNA, polymerases repeatedly initiate transcription, while when transcription factor is off the promoter, the initiation stops. While the data supports temporal coupling between TF binding and initiation of transcriptional bursts, the durations of the binding times and bursts are not equal. For instance, in yeast, an average TF (GAL4) binding time of 34 seconds initiates a mean burst duration of around 2.5 minutes~\citep{Chubb2020}.

We propose that the bursting periodic solutions that we observed in this paper may be one of the mechanisms supporting or enhancing transcriptional bursting. This may be in addition to other proposed mechanisms related to DNA supercoiling~\citep{chong2014}, chromatin opening, scaffold presence at initiation site or  pulses of nuclear localization~\citep{Lenstra2016}.

%% file: app_semiflow.tex
\section{The semiflow of differentiable solution operators generated by the system \eqref{eq:mrna-delay-var}-\eqref{eq:effector-delay-var} and \eqref{eq:delay by stateM}-\eqref{eq:delay by stateI}
}\label{app:semiflow}

For delay differential equations a familiar state space is given by continuous functions on a compact interval, see e.g. \citet{HVL,DvGVLW}.
In case of variable, state-dependent delays, however, there is a specific lack of smoothness which means that in general the initial value problem is not well-posed for only continuous initial data, not to speak of, say, smoothness of solutions with respect to initial data and linearization \citet{W1,HKWW}.

In this section we reformulate the system \eqref{eq:mrna-delay-var}-\eqref{eq:delay by stateI}
as a delay differential equation
\be
x'(t)=G(x_t)
\label{eq:form F}
\ee
with a vector-valued functional $G:C^1_3\to\R^3$ on the Banach space $C^1_3=C^1([-r,0],\R^3)$ of continuously differentiable maps $[-r,0]\to\R^3$, for some $r>0$ which is to be determined. The norm on $C^1_3$ is given by
$$
|\phi|_1=\max_{-r\le t\le 0}|\phi(t)|+\max_{-r\le t\le 0}|\phi'(t)|,
$$
with a chosen norm on $\R^3$. The argument $x_t\in C^1_3$ in \eqref{eq:form F} is defined for $t\in\R$ with $[t-r,t]$ in the domain of the solution $x$  and given by $x_t(s)=x(t+s)$ for $-r\le s\le0$. In other words, $x_t$ is the restriction of $x$ to  $[t-r,t]$ shifted to the interval $[-r,0]$.

In the sequel shifted segments $y_t$ of maps $y$ from an interval $I\subset\R$ into a set $M$ are defined accordingly: For $r>0$ and $t\in\R$ with $[t-r,t]\subset I$ the map $y_t:[-r,0]\to M$ is given by $y_t(s)=y(t+s)$, $-r\le s\le0$.

In addition to the space $C^1_3$ we also need the Banach space $C_3\supset C^1_3$ of continuous maps $[-r,0]\to\R^3$, with the norm given by
$$
|\phi|=\max_{-r\le t\le 0}|\phi(t)|,
$$
and the Banach spaces $C^1$ and $C$  of scalar functions $[-r,0]\to\R$  which are defined analogously to $C^1_3$ and $C_3$.

We shall verify the hypotheses from \citet{W1,W2,HKWW} which guarantee existence, uniqueness, and differentiability with respect to initial data of solutions to an initial value problem which is associated with \eqref{eq:form F} in a submanifold of the space $C^1_3$.

First we rewrite the system \eqref{eq:mrna-delay-var}-\eqref{eq:delay by stateI}
in a form which is more convenient for our purpose. Let positive constants
\begin{gather*}
c_k\quad\text{and}\quad b_k\quad\text{for}\quad k=1,2,3,\\
a_1\quad\text{and}\quad a_2,\\
\mu
\end{gather*}
and a continuously differentiable function $g:\R\to\R$ be given, either with
$$
g(0)=1\quad\text{and g decreasing on}\;[0,\infty)\;\text{with}\quad\lim_{y\to\infty}g(y)>0,
$$
or
$$
0<g(0)<1\quad\text{and g increasing on}\;[0,\infty)\;\text{with}\quad\lim_{y\to\infty}g(y)=1.
$$
Let also continuously differentiable functions
$$
v_k:\R\to[v_{k,min},\infty)
$$
with $0<v_{k,min}$ be given, for $k=1,2$.

In place of \eqref{eq:mrna-delay-var}-\eqref{eq:delay by stateI},  we consider the system
\begin{align}
x_1'(t) & =  c_1\frac{v_1(x_3(t))}{v_1(x_3(t-\tau_1(t))}e^{-\mu\tau_1(t)}g(x_3(t-\tau_1(t)))-b_1x_1(t)\label{eq:F-1},\\
x_2'(t) & =  c_2\frac{v_2(x_1(t))}{v_2(x_1(t-\tau_2(t))}e^{-\mu\tau_2(t)}x_1(t-\tau_2(t))-b_2x_2(t)\label{eq:F-2},\\
x_3'(t) & =  c_3x_2(t)-b_3x_3(t)\label{eq:F-3},\\
a_1 & =  \int^t_{t-\tau_1(t)}v_1(x_3(s))ds=\int^0_{-\tau_1(t)}v_1(x_3(t+s))ds\label{eq:delay-1},\\
a_2 & =  \int^t_{t-\tau_2(t)}v_2(x_1(s))ds=\int^0_{-\tau_2(t)}v_2(x_1(t+s))ds\label{eq:delay-2}.
\end{align}
Assume
$$
r>\frac{a_k}{v_{k,min}}\quad\text{for}\quad k=1,2.
$$
Notice that $r$ is an a-priori bound for both delays $\tau_k(t)$, $k\in\{1,2\}$. Using segment notation  the equations \eqref{eq:delay-1}-\eqref{eq:delay-2} become
$$
a_1=\int^0_{-\tau_1(t)}v_1(x_{3,t}(s))ds\quad\text{and}\quad a_2=\int^0_{-\tau_2(t)}v_2(x_{1,t}(s))ds
$$
with segments $x_{k,t}\in C$, $k\in\{1,2\}$. More generally, we consider the equations
\begin{equation}
a_k=\int^0_{-u}v_k(\phi(s))ds\label{eq:delay-k}
\end{equation}
for $u\in[0,\infty)$ and $\phi\in C$, $k\in\{1,2\}$. Using positivity of the functions $v_k$ and the Intermediate Value Theorem we infer  that for every $k\in\{1,2\}$ and for every $\phi\in C$ there is a uniquely determined solution $u=\delta_k(\phi)\in(0,r)$ of Eq. \eqref{eq:delay-k}.
This yields maps $\delta_k:C\to(0,r)$, $k\in\{1,2\}$. The next proposition guarantees that these maps are continuously differentiable and provides formulae for the derivatives.

\begin{prop} \label{prop:1}
Let $a>0$ and $v_{min}>0$ be given with $\frac{a}{v_{min}}<r$, and let $v:\R\to[v_{min},\infty)$ be continuously differentiable. Then the map $\delta:C\to(0,r)$ given by $\delta(\phi)=u$ with
\be
a=\int^0_{-u}v(\phi(s))ds\label{eq:define delay}
\ee
is continuously differentiable with
$$
D\delta(\phi)\chi=-\frac{\int^0_{-\delta(\phi)}v'(\phi(s))\chi(s)ds}{v(\phi(-\delta(\phi)))}.
$$
In case $\phi(s)=\xi$ for all $s\in[-r,0]$,
$$
\delta(\phi)=\frac{a}{v(\xi)},
\qquad\text{and}\qquad
D\delta(\phi)\chi=-\frac{v'(\xi)}{v(\xi)}\int^0_{-a/v(\xi)}\chi(s)ds.
$$
\end{prop}

Before giving the proof recall from  \citet[p. 47]{W1} or \citet[p. 466]{HKWW}
the elementary facts that the evaluation map
$$
ev_C:C\times[-r,0]\ni(\chi,u)\mapsto\chi(u)\in\R.
$$
is continuous (but not
locally Lipschitz continuous, let alone differentiable),
and that the restricted
evaluation map
$$
ev:C^1\times(-r,0)\ni(\phi,u)\mapsto\phi(u)\in\R
$$
is continuously differentiable with
\be
D\,ev(\phi,u)(\hat{\phi},\hat{u})=D_1ev(\phi,u)\hat{\phi}+D_2ev(\phi,u)\hat{u}=\hat{\phi}(u)+\phi'(u)\hat{u}
\label{eq:eval deriv}
\ee
where $D_1$ and $D_2$ denote partial derivatives with respect to the argument in $C^1$ and in $(-r,0)$, respectively.
For $v:\R\to\R$ continuously differentiable the substitution operator
$$
V:C\ni\phi\mapsto v\circ\phi\in C
$$
is continuously differentiable with
$$
(DV(\phi)\hat{\phi})(s)=v'(\phi(s))\hat{\phi}(s)\quad\text{for all}\quad\hat{\phi}\in C,\quad s\in[-r,0],
$$
see for example \citet[Appendix IV, Lemma 1.5]{DvGVLW}.

\begin{proof}[of Proposition~\ref{prop:1}]
For every $\phi\in C$ the value  $u=\delta(\phi)$ is the unique solution of the equation
$$
h(u,\phi)=0
$$
where $h:(0,r)\times C\to\R$ is given by
$$
h(u,\phi)=a-\int^0_{-u}v(\phi(s))ds=a - ev(I(V(\phi))),-u)
$$
with the continuous linear integration operator
$$
I:C\to C^1,\quad (I\psi)(t)=\int^0_t\psi(s)ds.
$$
The map $h$ is continuously differentiable with
$$
D_1h(u,\phi)1=-v(\phi(-u))<0
$$
and
\begin{align*}
D_2h(u,\phi)\chi & =  -D_1ev(I(V(\phi)),-u)DI(V(\phi))DV(\phi)\chi
 =  -(DI(V(\phi))DV(\phi)\chi)(-u)\\
& =  -(I(DV(\phi)\chi))(-u)=-\int^0_{-u}v'(\phi(s))\chi(s)ds.
\end{align*}
The Implicit Function Theorem applies at every $(\delta(\phi),\phi)\in(-r,0)\times C$ and yields that locally,
$\delta$ is given by a continuously differentiable map. Differentiation of  the equation $h(\delta(\phi),\phi)=0$ gives
\begin{align}
D\delta(\phi)\chi & =  -\frac{D_2h(u,\phi)\chi}{D_1h(u,\phi)1}
\quad\text{(for}\quad u=\delta(\phi))\nonumber\\
& =  -\frac{\int^0_{-u}v'(\phi(s))\chi(s)ds}{v(\phi(-u))}.\nonumber
\end{align}
In case $\phi$ is constant with value $\xi\in\R$ Eq.
\eqref{eq:define delay}
gives $a=\delta(\phi)v(\xi)$, and by the previous formula,
$$
D\delta(\phi)\chi=-\frac{v'(\xi)}{v(\xi)}\int^0_{-a/v(\xi)}\chi(s)ds. $$

\vspace*{-4ex} \mbox{} 
\end{proof}

Using the delay functionals $\delta_k:C\to(0,r)$, $k\in\{1,2\}$  the  system
\eqref{eq:F-1}-\eqref {eq:delay-2}
is reduced to the system
\begin{align}
x_1'(t) & =  c_1\frac{v_1(x_3(t))}{v_1(x_3(t-\delta_1(x_{3,t})))}e^{-\mu\delta_1(x_{3,t})}g(x_3(t-\delta_1(x_{3,t})))-b_1x_1(t),\label{eq:explicit-form-1}\\
x_2'(t) & =  c_2\frac{v_2(x_1(t))}{v_2(x_1(t-\delta_2(x_{1,t})))}e^{-\mu\delta_2(x_{1,t})}x_1(t-\delta_2(x_{1,t}))-b_2x_2(t),\label{eq:explicit-form-2}\\
x_3'(t) & =  c_3x_2(t)-b_3x_3(t).\label{eq:explicit-form-3}
\end{align}
with segments $x_{3,t}\in C$ and $x_{1,t}\in C$. The right hand side of these equations is of the form $G_C(x_t)\in\R^3$ with the map $G_C:C^3\to\R^3$ given by
\begin{align}
G_{C,1}(\phi) & =  c_1\frac{v_1(\phi_3(0))}{v_1(ev_C(\phi_3,-\delta_1(\phi_3)))}e^{-\mu\delta_1(\phi_3)}g(ev_C(\phi_3,-\delta_1(\phi_3)))-b_1\phi_1(0),\label{eq:G-C-1}\\
G_{C,2}(\phi) & =  c_2\frac{v_2(\phi_1(0))}{v_2(ev_C(\phi_1,-\delta_2(\phi_1)))}e^{-\mu\delta_2(\phi_1)}ev_C(\phi_1,-\delta_2(\phi_1))-b_1\phi_2(0),\label{eq:G-C-2}\\
G_{C,3}(\phi)& =  c_3\phi_2(0)-b_3\phi_3(0).\label{eq:G-C-3}
\end{align}
We observe that $G_C$ is continuous. The restriction $G$ of $G_C$ to $C^1_3$ is continuously differentiable because, for $\phi\in C^1_3\subset C_3$, we have
$ev_C(\phi_3,-\delta_1(\phi_3))=ev(\phi_3,-\delta_1(\phi_3))$ and $ev_C(\phi_1,-\delta_2(\phi_1))=ev(\phi_1,-\delta_2(\phi_1))$.  Here
$ev:C^1\times(0,r)\to\R$ is continuously differentiable, the maps $\delta_1$ and $\delta_2$ are continuously differentiable, and the projections $C^1_3\to C^1$ to components as well as the evaluation map
$ev_0:C\ni\chi\mapsto\chi(0)\in\R$ are linear and continuous.

In order to simplify calculations below we now introduce the
continuous maps
$$
E_{C,31}:C_3\to\R,\quad E_{C,31}(\phi)=ev_C(\phi_3,-\delta_1(\phi_3)),
$$
and
$$
E_{C,12}:C_3\to\R,\quad E_{C,12}(\phi)=ev_C(\phi_1,-\delta_2(\phi_1)),
$$
and the continuously differentiable maps
$$
E_{31}:C^1_3\to\R,\quad E_{31}(\phi)=ev(\phi_3,-\delta_1(\phi_3)),
$$
and
$$
E_{12}:C^1_3\to\R,\quad E_{12}(\phi)=ev(\phi_1,-\delta_2(\phi_1)),
$$
with derivatives at $\phi\in C^1$ given by
\begin{align*}
	DE_{31}(\phi)\hat{\phi}& =  D_1ev(\phi_3,-\delta_1(\phi_3))\hat{\phi}_3+D_2ev(\ldots)D(-\delta_1)(\phi_3)\hat{\phi}_3\\	
	& =  \hat{\phi}_3(-\delta_1(\phi_3))-\phi'_3(-\delta_1(\phi_3))D\delta_1(\phi_3)\hat{\phi}_3,\\	
	DE_{12}(\phi)\hat{\phi}& =  D_1ev(\phi_1,-\delta_2(\phi_1))\hat{\phi}_1+D_2ev(\ldots)D(-\delta_2)(\phi_1)\hat{\phi}_1\\	
	& =  \hat{\phi}_1(-\delta_2(\phi_1))-\phi'_1(-\delta_2(\phi_1))D\delta_2(\phi_1)\hat{\phi}_1	
\end{align*}
for all $\hat{\phi}\in C^1_3$. Notice that the right hand sides of these equations make sense also for arguments $\chi\in C_3$ instead of $\hat{\phi}\in C^1_3$. Thus they define linear extensions $D_eE_{31}(\phi):C_3\to\R$ of $DE_{31}(\phi):C^1_3\to\R$ and
$D_eE_{12}(\phi):C_3\to\R$ of $DE_{12}(\phi):C^1_3\to\R$. Using the continuity of the map $ev_C$, and the fact that differentiation $C^1\ni\phi\mapsto\phi'\in C$ is linear and continuous we obtain the next result.

\begin{prop} \label{prop:2}
The maps
$C^1_3\times C_3\ni(\phi,\chi)\mapsto D_eE_{31}(\phi)\chi\in\R$
and $C^1_3\times C_3\ni(\phi,\chi)\mapsto D_eE_{12}(\phi)\chi\in\R$	
are continuous.	\qed
\end{prop}

Incidentally, in the case $\phi\in C^1_3$ is constant with value $(\xi_1,\xi_2,\xi_3)\in\R^3$ we have
$$
E_{31}(\phi)=\xi_3,\quad DE_{31}(\phi)\hat{\phi}=\hat{\phi}_3(-a_1/v_1(\xi_3))
$$
and
$$
E_{12}(\phi)=\xi_1,\quad DE_{12}(\phi)\hat{\phi}=\hat{\phi}_1(-a_2/v_2(\xi_1))
$$

With the linear continuous evaluation map $ev_0:C\ni\phi\mapsto\phi(0)\in\R$ we obtain for the restriction $G$ of $G_C$
\begin{align}
G_1(\phi) & =  c_1\frac{v_1(ev_0\phi_3)}{v_1(E_{31}(\phi))}e^{-\mu\delta_1(\phi_3)}g(E_{31}(\phi))
-b_1ev_0\phi_1\label{eq:form F-1},\\
G_2(\phi) & =  c_2\frac{v_2(ev_0\phi_1)}{v_2(E_{12}(\phi))}e^{-\mu\delta_2(\phi_1)}E_{12}(\phi)-b_2ev_0\phi_2\label{eq:form F-2},\\
G_3(\phi) & =  c_3ev_0\phi_2-b_3ev_0\phi_3.\label{eq:form F-3}
\end{align}
For $G_C$ we have analogous formulas, with $E_{C,31}$ and $E_{C,12}$ instead of $E_{31}$ and $E_{12}$, respectively. In the sequel we will show that the initial value problem
\begin{equation}
x'(t)=G(x_t)\quad\text{for}\quad t>0,\quad x_0=\phi\label{eq:ivp}
\end{equation}
is well-posed on the set
$$
X_G=\{\phi\in C^1_3:\phi'(0)=G(\phi)\}
$$
which is a continuously differentiable submanifold of codimension $3$ in the space $C^1_3$. This result  follows from results in  \cite{W1,W2,HKWW}
provided that  the following two assertions are verified:
$$
X_G\neq\emptyset
$$
and

\begin{itemize}
\item[(e)] {\it each derivative $DG(\phi):C^1_r\to\R^3$, $\phi\in C^1_3$, has a linear extension $D_eG(\phi):C_3\to\R^3$, and the map
$$
C^1_3\times C_3\ni(\phi,\hat{\phi})\mapsto D_eG(\phi)\hat{\phi}\in\R^3
$$
is continuous.}
\end{itemize}

Property (e) is a version of being {\it almost Fr\'echet differentiable}
from \cite{M-PNP}.

We now proceed to the verification of $X_G\neq\emptyset$.
Choose $\phi_3\in C^1$ to be a constant with value $1$. Then $a_1\delta_1(\phi_3)=v_1(1)$. Choose $\xi_1>0$ with
$0=c_1e^{-\mu v_1\!(1)/a_1}g(1)-\xi_1$
and define $\phi_1\in C^1$ to be constant with value $\xi$. Define $\xi_2>0$ by $0=c_3\xi_2-b_3$. There exists $\phi_2\in C^1$ such that $\phi_2(0)=\xi_2$ and $\phi_2'(0)=
c_3\xi_2-b_3$. Then $\phi\in C^1_3$ with the components $\phi_1,\phi_2,\phi_3$ satisfies the equation $\phi'(0)=G(\phi)$ defining the set $X_G$, so $X_G\neq\emptyset$.

To prepare for the proof of extension property (e)
we compute the derivatives $DG(\phi)$, $\phi\in C^1_3$. For $\hat{\phi}\in C^1_3$ we have
$$
DG(\phi)\hat{\phi}=(DG_1(\phi)\hat{\phi},DG_2(\phi)\hat{\phi},DG_3(\phi)\hat{\phi})\in\R^3.
$$
For $\delta_1:C\to\R$ and $\delta_2:C\to\R$ we use the fact that restrictions of differentiable maps $m:C\to\R$ to $C^1$ remain differentiable, with derivatives $D(m|C^1)(\phi):C^1\to\R$ being restrictions of the derivatives $Dm(\phi):C\to\R$, $\phi\in C^1\subset C$, and obtain
\begin{align}
& DG_1(\phi)\hat{\phi} = -b_1\hat{\phi}_1(0)
 +c_1\bigg\{\frac{1}{[v_1(E_{31}(\phi))]^2} \notag \\
& \times \Big[v_1'(\phi_3(0))\hat{\phi}_3(0)\cdot v_1(E_{31}(\phi))
 -v_1(\phi_3(0))v_1'(E_{31}(\phi))DE_{31}(\phi)\hat{\phi}\Big]\cdot e^{-\mu\delta_1(\phi_3)}g(E_{31}(\phi))\notag\\
& +\frac{v_1(\phi_3(0))}{v_1(E_{31}(\phi))}\cdot\Big[
-\mu e^{-\mu\delta_1(\phi_3)}D\delta_1(\phi_3)\hat{\phi}_3\cdot g(E_{31}(\phi))
+e^{-\mu\delta_1(\phi_3)}Dg(E_{31}(\phi))DE_{31}(\phi)\hat{\phi}\Big]\bigg\}
\end{align}
and
\begin{align*}
& DG_2(\phi)\hat{\phi}  = -b_2\hat{\phi}_2(0)+c_2\bigg\{\frac{1}{[v_2(E_{12}(\phi))]^2} \notag \\
& \times \Big[v_2'(\phi_1(0))\hat{\phi}_1(0)\cdot v_2(E_{12}(\phi))
-v_2(\phi_1(0))v_2'(E_{12}(\phi))DE_{12}(\phi)\hat{\phi}Big]\cdot e^{-\mu\delta_2(\phi_1)}E_{12}(\phi)\notag\\
& +\frac{v_2(\phi_1(0))}{v_2(E_{12}(\phi))}\Big[-\mu\,e^{-\mu\delta_2(\phi_1)}D\delta_2(\phi_1)\hat{\phi}_1\cdot E_{12}(\phi)
+e^{-\mu\delta_2(\phi_1)}DE_{12}(\phi)\hat{\phi}\Big]\bigg\}\nonumber
\end{align*}
and
$$
DG_3(\phi)\hat{\phi}=-b_3\hat{\phi}_3(0)+c_3\hat{\phi}_2(0).
$$

Now we are ready to verify property (e).
In the formula for $DG(\phi)\hat{\phi}$, $\phi$ and $\hat{\phi}$ in $C^1_3$, just obtained replace the real numbers $DE_{31}(\phi)\hat{\phi},DE_{12}(\phi)\hat{\phi}$  by $D_eE_{31}(\phi)\chi,D_eE_{12}(\phi)\chi$, respectively, with $\chi\in C_3$, and replace the functions $\hat{\phi}_3,\hat{\phi}_1$ by $\chi_3$ and $\chi_1$, respectively. This defines $D_eG(\phi)\chi\in\R^3$ for $\phi\in C^1_3$ and $\chi\in C_3$ so that the maps $D_eG(\phi):C_3\to\R^3$, $\phi\in C^1_3$, are linear. Using the continuous differentiability of $\delta_1:C\to\R$ and $\delta_2:C\to\R$, and Proposition~\ref{prop:2}, one shows that the map
$C^1_3\times C_3\ni(\phi,\chi)\mapsto D_eG(\phi)\chi\in\R^3$
is continuous.
This finishes the proof of property (e).

With $X_G\neq\emptyset$ and property (e) verified,
 results from  \citet{W1,W2,HKWW} apply and yield the following. The set $X_G$ is a continuously differentiable submanifold of the Banach space $C^1_3$, with codimension 3.  Each $\phi\in X_G$ uniquely determines a maximal continuously differentiable solution $x:[-r,t_x)\to\R^3$, $0<t_x\le\infty$, of the initial value problem
\eqref{eq:ivp}.  That is, $x$ is continuously differentiable and satisfies $x_0=\phi$ and $x'(t)=G(x_t)$ for all $t\in(0,t_x)$, and any other continuously differentiable function $y:[-r,t_y)\to\R^3$, $0<t_y\le \infty$,  which satisfies $y_0=\phi$ and $y'(t)=G(y_t)$ for all $t\in(0,t_y)$ is a restriction of $x$. All segments $x_t$, $0\le t<\infty$, belong to $X_G$ (because of the differential equation). Write $x^{\phi}=x$ and $t_{\phi}=t_x$.
Let
$$
\Omega_G=\{(t,\phi)\in[0,\infty)\times X_G:0\le t<t_{\phi}\}.
$$
The equation
$$
S_G(t,\phi)=x_t^{\phi},
$$
defines a continuous semiflow  $S_G:\Omega_G\to X_G$. For each $t\ge0$ the set
$$
\Omega_{G,t}=\{\phi\in X_G:t<t_{\phi}\}
$$
is an open subset of $X_G$ (possibly empty), $\Omega_{G,0}=X_G$, and each map
$$
S_{G,t}:\Omega_{G,t}\ni\phi\mapsto S_G(t,\phi)\in X_G,\quad t\ge0,
$$
on a non-empty domain is continuously differentiable.

Moreover, the restriction of the semiflow $S_G$ to the open subset $\{(t,\phi)\in\Omega_G:r<t\}$ of the manifold $\R\times X_G$ is continuously differentiable \citet{W2}.

We end this section with remarks about linearization at {\it stationary points} (equilibria) of the semiflow, for which $t_{\phi}=\infty$ and $S_G(t,\phi)=\phi$ for all $t\ge0$. Such $\phi$ are constant since for every $t\ge0$, $x^{\phi}(t)=x^{\phi}_t(0)=S_G(t,\phi)(0)=\phi(0)$, hence
$$
\phi(s)=S_G(r,\phi)(s)=x^{\phi}_r(s)=x^{\phi}(r+s)=\phi(0)\quad \mbox{  for each } s\in[-r,0].
$$
So assume $\phi\in C^1_3$ is constant with value $(\xi_1,\xi_2,\xi_3)\in\R^3$, and $\hat{\phi}\in C^1_3$. We compute
$DG(\phi)\hat{\phi}$, using Proposition~\ref{prop:1} for the values and for the derivatives of the maps
$\delta_1,\delta_2$ in case of constant arguments, and using the calculations right after Proposition~\ref{prop:2} for $DE_{31}(\phi)\hat{\phi},DE_{12}(\phi)\hat{\phi}$ in case of constant arguments. The result is
\begin{align}
DG_1(\phi)\hat{\phi} & = -b_1\hat{\phi}_1(0)+c_1\Bigg\{\frac{1}{[v_1(\xi_3)]^2}\Big[v_1'(\xi_3)\hat{\phi}_3(0)\cdot v_1(\xi_3) \notag\\
& \qquad\qquad\qquad\qquad -v_1(\xi_3)v_1'(\xi_3)\hat{\phi}_3(-a_1/v_1(\xi_3))\Big]\cdot e^{-\mu a_1/v_1(\xi_3)}g(\xi_3)\nonumber\\
& \qquad +\frac{v_1(\xi_3)}{v_1(\xi_3)}\Bigg[-\mu e^{-\mu a_1/v_1(\xi_3)}
\left(-\frac{v_1'(\xi_3)}{v_1(\xi_3)}\int^0_{-a_1/v_1(\xi_3)}\hspace*{-1em}\hat{\phi}_1(s)ds\right)\cdot g(\xi_3) \notag\\
& \qquad\qquad\qquad\qquad +  e^{-\mu a_1/v_1(\xi_3)}g'(\xi_3)\hat{\phi}_3(-a_1/v_1(\xi_3))\Bigg]\Bigg\}\nonumber\\
& =  -b_1\hat{\phi}_{1\!}(0)+A_1\hat{\phi}_3(0)
+\mu A_1 \int^0_{-a_1/v_1(\xi_3)}\hat{\phi}_3(s)ds \notag\\
& \qquad\quad +\left(c_1e^{-\mu a_1^{\!}/^{\!}v_1(\xi_3)}g'^{\!}(\xi_3)-A_1\!\right)\hat{\phi}_3(-a_1^{\!}/^{\!}v_1(\xi_3))
\label{eq:DF-1-equi}
\end{align}
with
\be \label{eq:A1}
A_1=c_1\frac{v_1'(\xi_3)}{v_1(\xi_3)}e^{-\mu a_1/v_1(\xi_3)}g(\xi_3)
\ee
and
\begin{align}
DG_2(\phi)\hat{\phi} & =  -b_2\hat{\phi}_2(0)+c_2\Bigg\{\frac{1}{[v_2(\xi_1)]^2}\Big[v_2'(\xi_1)\hat{\phi}_1(0)\cdot v_2(\xi_1) \notag \\
& \hspace*{13em} - v_2(\xi_1)v_2'(\xi_1)\hat{\phi}_1(-a_2/v_2(\xi_1))\Big]\cdot e^{-\mu a_2/v_2(\xi_1)}\xi_1\nonumber\\
& \qquad +\frac{v_2(\xi_1)}{v_2(\xi_1)}\Bigg[ -\mu\,e^{-\mu a_2/v_2(\xi_1)}\left(-\frac{v_2'(\xi_1)}{v_2(\xi_1)}\int^0_{-a_2/v_2(\xi_1)}\hspace*{-1em}\hat{ \phi}_1(s)ds\right)\cdot\xi_1 \nonumber\\
& \hspace*{15em} + e^{-\mu a_2/v_2(\xi_1)}\hat{\phi}_1(-a_2/v_2(\xi_1))\Bigg]\Bigg\}\nonumber\\
& =   -b_2\hat{\phi}_2(0)+A_2\hat{\phi}_1(0)+\left(c_2e^{-\mu a_2/v_2(\xi_1)}-A_2\!\right)\hat{\phi}(-a_2^{\!}/^{\!}v_2(\xi_1))\notag\\
& \qquad +\mu A_2 \int^0_{-a_2/v_2(\xi_1)}\hat{ \phi}_1(s)ds
\label{eq:DF-2-equi}
\end{align}
with
\be \label{eq:A2}
A_2=c_2\frac{v_2'(\xi_1)}{v_2(\xi_1)}e^{-\mu a_2/v_2(\xi_1)}\xi_1
\ee
and
\be
DG_3(\phi)\hat{\phi}=-b_3\hat{\phi}_3(0)+c_3\hat{\phi}_2(0).\label{eq:DF-3-equi}
\ee

%% file: app_attractor.tex
\section{Positivity, dissipativity, global attractor}\label{app:dissipativity}

Recall that the function $g$ in Eqs. \eqref{eq:explicit-form-1} and \eqref{eq:form F-1} has a positive infimum $g_{min}>0$, and
$$
g_{min}\le g(\xi)\le 1\quad\text{for all}\quad\xi\in\R.
$$
In addition to the hypotheses made in the preceding
Section \ref{app:semiflow} we assume in the present section that the functions $v_k$, $k\in\{1,2\}$, are also bounded from above by real numbers $v_{k,max}\ge0$. Using Eq. \eqref{eq:delay-k} we infer
$$
\frac{a_k}{v_{k,max}}\le\delta_k(\phi)\le\frac{a_k}{v_{k,min}}
$$
for $k=1,2$ and for all $\phi\in C$.

\begin{prop}\label{prop 4}
	For every $c>0$ there exists $c'>0$ so that for all $\phi\in X_G$ with $|\phi(t)|\le c$ on $[-r,0]$ and for all $t\in[-r,t_{\phi})$ we have $|x^{\phi}(t)|\le c'$.
\end{prop}

\begin{proof}
Let $\phi\in X_G$ with $|\phi(t)|\le c$ on $[-r,0]$  be given, set $x=x^{\phi}$. The first term on the right hand side of Eq. \eqref{eq:explicit-form-1} is positive and bounded by the constant
$$d_1^>=c_1\frac{v_{1,max}}{v_{1,min}}.$$
The variation-of-constants formula yields
\begin{align*}
|x_1(t)| & = \bigg|x_1(0)e^{-b_1t}\\
& \qquad + \int_0^t\! e^{-b_1(t-s)}\! \left[\frac{c_1 v_1(x_3(s))}{v_1(x_3(s-\delta_1(x_{3,s})))}
e^{-\mu\delta_1(x_{3,s})}g(x_3(s-\delta_1(x_{3,s})))\right]^{\!}ds\bigg|\\
	& \le c+e^{-b_1t}\frac{d_1^>}{b_1}\left(e^{b_1t}-1\right)\quad\text{for}\quad 0\le t<t_{\phi}\\
	& \le c+\frac{d_1^>}{b_1}\quad\text{for}\quad 0\le t<t_{\phi}.
\end{align*}
	Set $d_1=c+\frac{d_1^>}{b_1}$. With $|x_1(t)|=|\phi_1(t)|\le|\phi|\le c$ on $[-r,0]$ we get $|x_1(t)|\le d_1$ for all $t\in[-r,t_{\phi}]$. Next, the first term on the right hand side of Eq. \eqref{eq:explicit-form-2} is bounded by
$$
d_2^>=c_2\frac{v_{2,max}}{v_{2,min}}d_1.
$$
Using the variation-of-constants formula as above, and $|x_2(t)|=|\phi_2(t)|\le|\phi|\le c$ on $[-r,0]$ we obtain that $x_2$ is bounded by
$d_2=c+\frac{d_2^>}{b_2}$ for all $t\in[-r,t_{\phi}]$. We turn to Eq. \eqref{eq:explicit-form-3} where the first term on the right hand side is bounded by
$$
d_3^>=c_3d_2.
$$
Using the variation-of-constants formula once more, and $|x_3(t)|=|\phi_3(t)|\le|\phi|\le c$ on $[-r,0]$ we obtain that $x_3$ is bounded by $d_3=c+\frac{d_3^>}{b_3}$ on $[-r,t_{\phi}]$.
\end{proof}

We continue with an observation. For any (continuously differentiable) solution $x:[-r,t_x)\to\R^3$ of the system  \eqref{eq:explicit-form-1}-\eqref{eq:explicit-form-3}  and for any $t\in[0,t_x)$ the first term on the right hand side of Eq. \eqref {eq:explicit-form-1} belongs to the interval
$$
[d_1^<,d_1^>]=\left[c_1\frac{v_{1,min}}{v_{1,max}}e^{-\mu a_1/v_{1,min}}g_{min},c_1\frac{v_{1,max}}{v_{1,min}}\right]
$$
and in case
$$
x_1(t)>\frac{d^>_1}{b_1}\quad\text{we have}\quad x_1'(t)\le d_1^>-b_1x_1(t)<0
$$
while in case
$$
x_1(t)<\frac{d^<_1}{b_1}\quad\text{we have}\quad x_1'(t)\ge d^<_1-b_1x_1(t)>0.
$$
Set
\begin{gather*}
	d_2^<=c_2\frac{v_{2,min}}{v_{2,max}}e^{-\mu\,a_2/v_{2,min}}\frac{d_1^<}{b_1} \quad\text{and}\quad
	d_2^>=c_2\frac{v_{2,max}}{v_{2,min}}e^{-\mu\,a_2/v_{2,max}}\frac{d_1^>}{b_1},\\
	d_3^<=c_3\frac{d_2^<}{b_2}  \quad\text{and}\quad d_3^>=c_3\frac{d_2^>}{b_2},
\end{gather*}
and
$$
Q=\left[\frac{d_1^<}{b_1},\frac{d_1^>}{b_1}\right]\times\left[\frac{d_2^<}{b_2},\frac{d_2^>}{b_2}\right]\times\left[\frac{d_3^<}{b_3},\frac{d_3^>}{b_3}\right]\subset\R^3
$$
and
$$
R=\{\phi\in C^1_3:\phi([-r,0])\subset Q\}.
$$

\begin{prop}\label{prop 5}(Global existence, absorption and positive invariance, positivity)
\begin{enumerate}	
\item[(i)] 
For all $\phi\in X_G$, $t_{\phi}=\infty$.
\item[(ii)] 
For every neighbourhood $N$ of $Q$ in $\R^3$ and for each $\phi\in X_G$ there exists $t(\phi,N)\in[0,\infty)$ with $x^{\phi}(t)\in N$ for all $t\ge t(\phi,N)$.
\item[(iii)] 
If $\phi\in X_G\cap R$ then $x^{\phi}(t)\in Q$ for all $t\ge0$.
\item[(iv)] 
If all 3 components of $\phi\in X_G$ are strictly positive then $x^{\phi}_k(t)>0$ for all $t\ge-r$, $k\in\{1,2,3\}$.
\end{enumerate}
\end{prop}

\begin{proof} 1. On assertion (i). Let $\phi\in X_G$ be given. Due to Proposition~\ref{prop 4} the solution $x=x^{\phi}$ is bounded.
Using this and the system \eqref{eq:explicit-form-1}-\eqref{eq:explicit-form-3} we infer that $x'$ is bounded. It follows that $x$ is Lipschitz continuous. Assume now $t_{\phi}<\infty$. Then Lipschitz continuity yields that $x$ has a limit $\xi\in\R^3$ at $t=t_{\phi}$.  $x$ extends to a continuous map $\hat{x}:[-r,t_{\phi}]\to\R^3$. From uniform continuity on the compact interval $[-r,t_{\phi}]$ it follows that the curve $[0,t_{\phi}]\ni t\mapsto \hat{x}_t\in C_3$ is continuous. Using this and the equation
$$
x'(t)=G(x_t)=G_C(\hat{x}_t)\quad\text{for}\quad0\le t<t_{\phi}
$$
with the continuous map $G_C:C^3\to\R^3$ we also conclude that $x'$ has a limit $\eta\in\R^3$ at $t=t_{\phi}$. It follows that $\hat{x}$ is continuously differentiable (with $\hat{x}'(t_{\phi})=\eta$), and $\hat{x}'(t_{\phi})=G_C(\hat{x}_{t_{\phi}})=G(\hat{x}_{t_{\phi}})$. In particular,
$\psi=\hat{x}_{t_{\phi}}$ belongs to $X_G$, and defines a maximal solution $x^{\psi}:[0,t_{\psi})\to\R^3$ of Eq. \eqref{eq:form F}, with $0<t_{\psi}\le\infty$. By means of the semiflow properties it follows that in case $t_{\psi}=\infty$ we have $t_{\phi}=\infty$, in contradiction to the assumption above, while in case $t_{\psi}<\infty$ we get $t_{\phi}\ge t_{\phi}+t_{\psi}$, which contradicts $t_{\psi}>0$.
	
2. On assertion (ii). Let a neighbourhood $N$ of $Q$ in $\R^3$ and $\phi\in X_G$ be given. Set $x=x^{\phi}$. There exists  $\epsilon>0$ so that
for
$$\begin{array}{lcr}
d_{1,-\epsilon} = \frac{d_1^<}{b_1}-\epsilon, \qquad  &
d^{\ast}_{2,-\epsilon}  = c_2\frac{v_{2,min}}{v_{2,max}}e^{-\mu\,a_2/v_{2,min}}\cdot d_{1,-\epsilon}, &
\qquad  d_{2,-\epsilon}  =  \frac{d^{\ast}_{2,-\epsilon}}{b_2}-\epsilon,\\
d_{1,+\epsilon}  =  \frac{d_1^>}{b_1}+\epsilon, &
d^{\ast}_{2,+\epsilon} =  c_2\frac{v_{2,max}}{v_{2,min}}e^{-\mu\,a_2/v_{2,max}}\cdot d_{1,+\epsilon}, &
d_{2,+\epsilon}  =  \frac{d^{\ast}_{2,+\epsilon}}{b_2}+\epsilon
\end{array}$$
we have
$$
d_{1,-\epsilon}>0\quad\text{and}\quad d_{2,-\epsilon}>0\quad\text{and}\quad\frac{c_3\cdot d_{2,-\epsilon}}{b_3}-\epsilon>0,
$$
and
$$
N\supset[d_{1,-\epsilon},d_{1,+\epsilon}]\times[d_{2,-\epsilon},d_{2,+\epsilon}]
\times\left[\frac{c_3\cdot d_{2,-\epsilon}}{b_3}-\epsilon,\frac{c_3\cdot d_{2,+\epsilon}}{b_3}+\epsilon\right].
$$
	
2.1.  Proof that in case $x_1(t)\le d_{1,+\epsilon}$ for some $t\ge0$ we have
$$
x_1(s)\le d_{1,+\epsilon}\quad\text{for all}\quad s\ge t.
$$
Otherwise $x_1(s)>d_{1,+\epsilon}=\frac{d_1^>}{b_1}+\epsilon$ for some $s>t$. For the smallest $u\in[t,s]$ with $x_1(u)=x_1(s)$ we have $0\le x'(u)$ and, on the other hand,
$$
x_1'(u)<d_1^>-b_1x_1(u)=d_1^>-b_1x_1(s)<-\epsilon\,b_1<0.
$$
	
2.2. Proof that in case $x_1(s)>d_{1,+\epsilon}$ for some $s\ge0$ there exists $t>s$ with
$$
x_1(t)\le d_{1,+\epsilon}.
$$
Otherwise $x_1(t)>d_{1,+\epsilon}=\frac{d_1^>}{b_1}+\epsilon$ on $[s,\infty)$. Hence
$$
x_1'(t)<d_1^>-b_1x_1(t)<-\epsilon\,b_1<0\quad\text{on}\quad[s,\infty),
$$
and consequently $x_1(t)\to-\infty$ as $t\to\infty$, in contradiction to the assumption.
	
2.3. It follows that there exists $t_1^{\ast}\ge0$ with
$$
x_1(t)\le d_{1,+\epsilon}\quad\text{for all}\quad t\ge t_1^{\ast}.
$$
Similarly one finds $t_1\ge t_1^{\ast}$ with
$$
x_1(t)\ge d_{1,-\epsilon}\quad\text{for all}\quad t\ge t_1.
$$
Hence
$$
x_1(t)\in[d_{1,-\epsilon},d_{1,+\epsilon}]\quad\text{for all}\quad t\ge t_1.
$$
	
2.4. We proceed to $x_2$. The result of Part 2.3 yields that for $t\ge t_1+r$ the first term on the right hand side of Eq. \eqref{eq:explicit-form-2} is contained in the interval
$$
\left[c_2\frac{v_{2,min}}{v_{2,max}}e^{-\mu a_2^{\!}/^{\!}v_{2,min}}\!\cdot d_{1,-\epsilon},c_2\frac{v_{2,max}}{v_{2,min}}e^{-\mu a_2^{\!}/^{\!}v_{2,max}}\! \cdot  d_{1,+\epsilon}\right]
^{\!} =[d^{\ast}_{2,-\epsilon},d^{\ast}_{2,+\epsilon}].
$$
Arguing as in Parts 2.1-2.3 we find $t_2\ge t_1+r$ so that for all $t\ge t_2$ we have
$$
x_2(t)\in[d_{2,-\epsilon},d_{2,+\epsilon}].
$$
	
2.5. Consider $x_3$. For $t\ge t_2$ the first term on the right hand side of Eq. \eqref{eq:explicit-form-3} is contained in the interval
$$
[c_3\cdot d_{2,-\epsilon},c_3\cdot d_{2,+\epsilon}].
$$
Arguing as in Parts 2.1-2.3 we find $t_3\ge t_2$ so that for all $t\ge t_3$ we have
$$
x_3(t)\in\left[\frac{c_3\cdot d_{2,-\epsilon}}{b_3}-\epsilon,\frac{c_3\cdot d_{2,+\epsilon}}{b_3}+\epsilon\right].
$$
	
2.6. For all $t\ge t_3$,
$$
x(t)\in[d_{1,-\epsilon},d_{1,+\epsilon}]\times[d_{2,-\epsilon},d_{2,+\epsilon}]\times\left[\frac{c_3\cdot d_{2,-\epsilon}}{b_3}-\epsilon,\frac{c_3\cdot d_{2,+\epsilon}}{b_3}+\epsilon\right]\subset N.
$$
	
3. The proof of assertion (iii) begins with the assumption that for a given $\phi\in X_G\cap R$ there exists $s>0$ with $x^{\phi}(s)>\frac{d_1^>}{b_1}$ and is then accomplished by a simplified version of arguments as in Part 2.

4. On assertion (iv).  Let $\phi\in X_G$ be given with strictly positive components. Set $x=x^{\phi}$. The assumption $x_1(t)\le0$ for some $t>0$ leads to a smallest $t>0$ with $x_1(t)=0$. Necessarily, $x'_1(t)\le0$ while Eq. \eqref{eq:explicit-form-1} yields $x'_1(t)>0$. It follows that $x_1(t)>0$ for all $t\ge -r$. Using this and Eq. \eqref{eq:explicit-form-2} and strict positivity of $\phi_2$ one shows in the same way that also $x_2$ is strictly positive. From this one deduces that $x_3$ is strictly positive.
\end{proof}

The solution manifold $X_G$ is a closed subset of the space $C^1_3$, and thereby a complete metric space with respect to the metric given by the norm on $C^1_3$.   The next result implies that the semiflow $S_G$ on the complete metric space $X_G$ is point dissipative as defined in \cite{hale88}, which means that there exists a bounded set $B$ such that for every $\epsilon>0$ and for each $\phi\in X_G$ there exists $t_{B,\epsilon,\phi}\ge0$ with
$S_G(t,\phi)$ contained in the $\epsilon$-neighbourhood $U_{\epsilon}(B)=\cup_{b\in B}U_{\epsilon}(b)$ of $B$.

\begin{cor}\label{cor 1}
There is a bounded open subset $B_G$  of the submanifold $X_G\subset C^1_3$, with
$$
\phi_k(t)>0\quad\text{for all}\quad \phi\in B_G,\quad t\in[-r,0],\quad k\in\{1,2,3\},
$$
such that for every $\phi\in X_G$ there exists $t(\phi)\ge0$ with
$$
S_G(t,\phi)\in B_G\quad\text{for all}\quad t\ge t(\phi).
$$
\end{cor}

\begin{proof}
 Choose $c>0$ so that $N=(0,c)^3$ is a neighbourhood of $Q$. Set $\tilde{R}=\{\phi\in C^1_3:\phi([-r,0])\subset N\}$.
Le $\phi\in X_G$ be given. Choose $t(\phi,N)$ according to Proposition~\ref{prop 4}(ii).
From Eqs. \eqref{eq:explicit-form-1}-\eqref{eq:explicit-form-3}we see that  the map $G$
sends the set $\tilde{R}$ (which is not a bounded subset of $C^1_3$) into a bounded subset of $\R^3$, say, into
$\{x\in\R^3:|x|<b\}$ for some $b>0$. It follows that for all $t\ge t(\phi,N)+r$ we have $|(x^{\phi})'(t)|<b$. For $t\ge t(\phi,N)+2r$ we obtain $x^{\phi}_t\in\{\phi\in X_G\cap\tilde{R}:|\phi'|<b\}=B_G$. The set $B_G$ is an open and bounded
subset of $X_G$, with $0<\phi_k(t)$ for all $\phi\in B_G,\,k\in \{1,2,3\},\,t\in[-r,0]$.
\end{proof}

Recall from \citet{hale88} the definition of a global attractor of a semiflow, which in case of our semiflow $S_G$ is equivalent to saying that a subset $A_G\subset X_G$ is a global attractor if it is
compact, and
invariant in the sense that for every $\phi\in A_G$ there exists a complete flowline\footnote{A complete flowline is a curve $\xi:\R\to X_G$ with
$\xi(t+s)=S_G(t,\xi(s))$ for all $t\ge0$ and $s\in\R$} with $\xi(0)=\phi$ and $\xi(\R)\subset A_G$,
and if
$A_G$ attracts every bounded set $B\subset X_G$ in the sense that given an open neighbourhood  $U\supset A_G$ of $A_G$ in $X_G$ there exists $t_{B,U}\ge0$ such that
$$
S_G([t_{B,U},\infty)\times B)\subset U.
$$

\citet[Theorem 3.4.8]{hale88} guarantees the existence of such a global attractor provided the semiflow is point-dissipative and there exists $t_1\ge0$ so that the semiflow $S_G$ is completely continuous for $t\ge t_1$. The property of being completely continuous (for $t\ge t_1$) is explained after \citet[Lemma 3.2.1]{hale88}. It means that\\
(1) $S_G$ is {\it conditionally completely continuous for $t\ge t_1$} in the sense that for every $t\ge t_1$ and for each bounded set $B\subset X_G$ for which $S_G([0,t]\times B)$ is bounded the set $S_G(t,B)$ is precompact (has compact closure),\\
and that\\
(2) for each bounded set $B\subset X_G$ and for all $t\ge0$ the set $S_G([0,t]\times B)$ is bounded.

\begin{thm} \label{thm 4-1}
The semiflow $S_G$ has a global attractor $A_G\subset X_G$, with $\phi_k(t)>0$ for all $\phi\in A_G,\,t\in[-r,0],\,k\in\{1,2,3\}$.
\end{thm}

\begin{proof}
1.  We first show that for every bounded subset $B\subset X_G$ there exists $c_B>0$ with
$$
|x^{\phi}(t)|\le c_B\quad\text{and}\quad |(x^{\phi})'(t)|\le c_B\quad\text{for all}\quad t\ge-r.
$$
Let $B\subset X_G$ be a  bounded subset which is bounded with respect to the norm of the space $C^1_3$. Proposition~\ref{prop 4} guarantees a constant $c_{B,0}$ with $|x^{\phi}(t)|\le c_{B,0}$ for all $\phi\in B$, $t\ge-r$. Then the formulae \eqref{eq:form F-1}--\eqref{eq:form F-3} show that the set
$\{G(x^{\phi}_t)\in\R^3:\phi\in B,t\ge0\}$ is bounded, and Eq. \eqref{eq:form F} gives that the set
$\{(x^{\phi})'(t)\in\R^3:\phi\in B,t\ge0\}$ is bounded. Also the set $\{\phi'(t)\in\R^3:\phi\in B,-r\le t\le0\}$ is bounded.

2. Claim: For every bounded subset $B\subset X_G$ the set $S_G(r,B)\subset X_G$ has compact closure in $C^1_3$.

Proof: (a)  Let $B\subset X_G$ be a bounded subset of $C^1_3$. Due to Part 1 the sets $\{x^{\phi}(t)\in R^3:\phi\in B,-r\le t\le r\}$ and
$\{(x^{\phi})'(t)\in\R^3:\phi\in B, -r\le t\le r\}$ are bounded. Using the Mean Value Theorem we see that in particular the set $S_G(r,B)$ is equicontinuous. As it also is bounded in $C_3$ the Ascoli-Arz\`{e}la Theorem implies that its closure in $C_3$ is compact.
		
(b) We turn to the set $\{S_G(r,\phi)'\in C_3:\phi\in B\}$ of derivatives, which is bounded in $C_3$, and proceed to show that it is also equicontinuous. As in Part (a) one sees that the closure $K$ of the set
$$
\{S_G(t,\phi)\in C_3:\phi\in B,0\le t\le r\}
$$
in the space $C_3$ is compact. The map $G:C^1_3\to\R^3$ is the restriction of the continuous map $G_C:C_3\to\R^3$
which is uniformly continuous on the compact set $K\subset C_3$. Using the boundedness of the set $\{(x^{\phi})'(t)\in\R^3:\phi\in B, -r\le t\le r\}$ and the Mean Value Theorem one finds that the curves
$$
[0,r]\ni t\mapsto S_G(t,\phi)\in C_3,\quad\phi\in B,
$$
are uniformly Lipschitz continuous, hence equicontinuous. Now
let $t_0\in 0,r]$ and $\epsilon>0$ be given. There exists $\delta>0$ with
$$
|G_C(\phi)-G_C(\psi)|\le\epsilon\quad\text{for all}\quad\phi,\psi\quad\text{in}\quad S_G([0,r]\times B)\quad\text{with}\quad|\phi-\psi|\le\delta,
$$
due to uniform continuity of $G_C$ on $K$. Due to equicontinuity there exists $\eta>0$ with
$$
|S_G(t,\phi)-S_G(t_0,\phi)|\le\delta\quad\text{for all}\;\phi\in B\;\text{and}\; t\in[0,r]\;\text{with}\;|t-t_0|<\eta.
$$
Hence
\begin{align*}
|(x ^{\phi})'(t)-(x^{\phi})'(t_0)| & =|G(S_G(t,\phi))-G(S_G(t_0,\phi))| \\
& =|G_C(S_G(t,\phi))-G_C(S_G(t_0,\phi))|<\epsilon
\end{align*}
for all $\phi\in B$ and $t\in[0,r]$ with $|t-t_0|<\eta$.
	
(c) The Ascoli-Arz\`{e}la Theorem implies that the closure of $\{S_G(r,\phi)'\in C_3:\phi\in B\}$  in $C_3$ is compact.
For the closure of $S_G(r,B)$ in $C^1_3$ to be compact it is sufficient to show that every sequence of points
$\phi_j\in S_G(r,B)$, $j\in\N$, has a subsequence which converges in $C^1_3$. Let a sequence $(\phi_j)_1^{\infty}$ in $S_G(r,B)$be given. Part a) yields that there is a subsequence $(\phi_{j_k})_1^{\infty}$ which converges in $C_3$ to some $\phi\in C_3$. Part b) yields that there is a further subsequence $(\phi_{j_k})_1^{\infty}$ so that the derivatives $(\phi_{j_k})'$ converge in $C_3$ to some $\psi\in C_3$. It follows that $\phi\in C^1_3$ with $\phi'=\psi$, which in turn yields $\phi_{j_k}\to\phi$ in $C^1_3$ as $k\to\infty$.
	
3.  We now show  that for every bounded subset $B\subset X_G$ and for every $t\ge r$ the set $S_G(t,B)\subset X_G$ has compact closure in $C^1_3$.  Let $B\subset X_G$ be bounded and let $t>r$. Then
$$
S_G(t,B)=S_G(t-r,S_G(r,B)).
$$
Use that the closure of $S_G(r,B)$ in $C^1_3$ is compact and belongs to $X_G$ (sinces $X_G$ is a closed subset of $C^1_3$), and that the map $S_G(t-r,\cdot)$ is continuous, and conclude that the closure of $S_G(t,B)$ in $C^1_3$ is contained in a compact subset of $X_G\subset C^1_3$.
	
4. According to the remark preceding Corollary \ref{cor 1} the semiflow $S_G$ is point dissipative. The results of Parts 1 and 3 combined yield that $S_G$ is {\it completely continuous for $t>r$} as stated before \ref{thm 4-1}. Using \cite[Theorem 3.4.8]{hale88} we infer that $S_G$ has a global attractor.
		
5. We now show  that for all $\phi\in A_G$, $t\in[-r,0]$, $k\in\{1,2,3\}$ we have $\phi_k(t)>0$. Let $\phi\in A_G$ be given. There exists a solution $x:\R\to\R^3$ of the system \eqref{eq:explicit-form-1} - \eqref{eq:explicit-form-3} with $x_0=\phi$ and all segments $x_s$, $s\in\R$, in $A_G$. It follows that $x$ is bounded.
	
(a) We now  show  that  $x_1(t)>0$ for all $t\in\R$. Assume $x_1(t)\le0$ for some $t\in\R$. In case $x_1(t)=0$ Eq. \eqref{eq:explicit-form-1} yields $x_1'(t)>0-b_1x_1(t)=0$. It follows that $x_1(u)<0$ for some $u\in(-\infty,t)$. In case $x_1(t)<0$ set $u=t$. For every $s\le u$ the variation-of-constants formula  yields
$$
x_1(u)\ge x_1(s)e^{-b_1(u-s)}+0,
$$
hence
$$
x_1(s)\le x_1(u)e^{b_1(u-s)}\qquad(\to-\infty\quad\text{as}\quad s\to-\infty),
$$
and we arrive at a contradiction to the boundedness of $x_1$.

(b) Using Part (a) and Eq. \eqref{eq:explicit-form-2} and arguing as in Part (a) one finds $x_2(t)>0$ for all $t\in\R$. Using this and Eq. \eqref{eq:explicit-form-3} and arguing once more as in Part (a) we get $x_3(t)>0$ for all $t\in\R$. In particular,
$\phi_k(t)>0$ for all $t\in[-r,0]$ and every $k\in\{1,2,3\}$.
\end{proof}


%% file: app_linearization.tex
\section{Linearization}\label{app:linearization}

We turn to linearization. The subsequent description is based on results proved in \citet[Sections 3.2 and 3.4]{HKWW} and \citet{W1}, \citet{W2}. At a point $\phi\in X_G$ the tangent space of the manifold $X_G$ is given by
$$
T_{\phi}X_G=\{\chi\in C^1_3:\chi'(0)=DG(\phi)\chi\}.
$$
For $\phi\in\Omega_{G,t}$ the derivative
$$
DS_{G,t}(\phi):T_{\phi}X_G\to T_{S_{G,t}(\phi)}X_G
$$
is given by
$$
DS_{G,t}(\phi)\chi=w^{\phi,\chi}_t
$$
where $w^{\phi,\chi}=w$ is the unique continuously differentiable solution $[-r,t_{\phi})\to\R^3$ of the IVP
\bea
w'(t) & = & DG(S_G(t,\phi))w_t\quad\text{for}\quad t>0,
\label{eq:varsyst1}\\
w_0 & = & \chi\in T_{\phi}X_G.
\label{eq:varsyst2}
\eea
Equation \eqref{eq:varsyst1}  is called the linear variational equation along the solution $x^{\phi}$ or  along the flowline
$$
S_G(\cdot,\phi):[0,t_{\phi})\ni t\mapsto S_G(t,\phi)\in X_G.
$$

Suppose now that $\phi\in X_G$ is a stationary point of the semiflow $S_G$. Then $\phi$ is constant  with value, say,  $\xi=(\xi_1,\xi_2,\xi_3)\in\R^3$. The variational equation \eqref{eq:varsyst1}
along the constant solution $x^{\phi}:[-r,\infty)\ni t\mapsto\xi\in\R^3$  becomes
\begin{align}
w_1'(t) & =  -b_1w_1(t)+A_1w_3(t)
+ \left(c_1e^{-\mu a_1/v_1(\xi_3)}g'(\xi_3)-A_1\!\right)w_3(t-a_1/v_1(\xi_3))\nonumber\\
& \qquad +\mu A_1\int^0_{-a_1/v_1(\xi_3)}w_3(s)ds
\label{eq:varsyst1stat}\\
w_2'(t) & =  -b_2w_2(t)+A_2w_1(t)
+\left(c_2e^{-\mu a_2/v_2(\xi_1)}-A_2\right)w_1(t-a_2/v_2(\xi_1))\nonumber\\
& \qquad +\mu A_2\int^0_{-a_2/v_2(\xi_1)}w_1(s)ds \label{eq:varsyst2stat}\\
w_3'(t) & =  -b_3w_3(t)+c_3w_2(t) \label{eq:varsyst3stat}
\end{align}

The derivatives $T_{G,t}=DS_{G,t}(\phi)$, $t\ge0$, form a strongly continuous semigroup on the closed subspace $T_{\phi}X_G$ of the space $C^1_3$. This semigroup is given by $T_{G,t}\chi=T_{G,e,t}\chi$ where $T_{G,e,t}:C_3\to C_3$ is the solution operator associated with the {\it classical} initial value problem
\bea
w'(t) & = & D_eG(\phi)w_t\quad\text{for}\quad t>0,
\label{eq:ivponC1}\\
w_0 & = & \chi\in C_3,
\label{eq:ivponC2}
\eea
with a continuous linear vector-valued functional $L:C_3\to\R^3$, $L=D_eG(\phi)$, as in the monographs, e. g., \cite{HVL,DvGVLW}. Recall  that by definition the solution $w:[-r,\infty)\to\R^3$ of the initial value problem
\eqref{eq:ivponC1}-\eqref{eq:ivponC2} is only continuous, with the restriction
$w|_{[0,\infty)}$ continuously differentiable and satisfying
\eqref{eq:ivponC1}.

The extended derivative $D_eG(\phi):C_3\to\R^3$ in the case just considered where $\phi\in X_G$ is a stationary point with value $\xi$ is given by \eqref{eq:DF-1-equi}-\eqref{eq:DF-3-equi}, now for $\hat{\phi}\in C_3$. Therefore the equation \eqref{eq:ivponC1}  coincides with the system \eqref{eq:varsyst1stat}-\eqref{eq:varsyst3stat}, considered for continuous maps $[-r,\infty)\to\R^3$ whose restrictions to $[0,\infty)$ are differentiable and satisfy \eqref{eq:varsyst1stat}-\eqref{eq:varsyst3stat} for all $t\ge0$.

%
%

The stability of the zero solution of the linear variational equation
\eqref{eq:varsyst1stat}-\eqref{eq:varsyst3stat} is determined by the spectrum $\sigma\subset\C$ of the generator of the semigroup $(T_{G,t})_{t\ge0}$ on $T_{\phi}X_G\subset C^1_3$, which coincides with the spectrum $\sigma_e\subset\C$ of the generator of the semigroup on $C_3$.

The spectrum $\sigma_e$ consists of the solutions $\lambda\in\C$ of the {\it characteristic equation}, which is obtained from the Ansatz $\R\ni t\mapsto e^{\lambda\,t}z\in\C^3$, $z\in\C^3\setminus\{0\}$, for a complex-valued  solution of the system
\eqref{eq:varsyst1stat}-\eqref{eq:varsyst3stat} as follows. We write down the system for $w:t\mapsto e^{\lambda t}z$, multiply by $e^{-\lambda t}$,  and obtain a linear equation of the form
\be \label{eq:M33eq}
(M(\lambda)-\lambda\cdot I_3)z=0
\ee
with a $3\times 3$-matrix $M(\lambda)$ (complex coefficients) and with $I_3$ denoting the $3\times 3$-unit matrix.
We then set
\begin{equation}
\Delta(\lambda)=\det(M(\lambda)-\lambda\cdot I)=0.\label{eq:char}
\end{equation}
Following this recipe we get the linear system
\begin{align*}
\lambda\,z_1 & =  -b_1z_1+A_1z_3 + \left(c_1e^{-\mu\,a_1/v_1(\xi_3)}g'(\xi_3)-A_1\right)z_3e^{-\lambda\,a_1/v_1(\xi_3)}\nonumber\\
& \qquad +\mu\,A_1\int^0_{-a_1/v_1(\xi_3)}e^{\lambda s}ds\,z_3 \\
\lambda\,z_2 & =  -b_2z_2+A_2z_1+\left(c_2e^{-\mu\,a_2/v_2(\xi_1)}-A_2\right)e^{-\lambda\, a_2/v_2(\xi_1)}z_1\\
& \qquad +\mu\,A_2\int^0_{-a_2/v_2(\xi_1)}e^{\lambda\,s}ds\,z_1 \\
\lambda\,z_3 & = -b_3z_3+c_3z_2,
\end{align*}
and with the functions $k_1:\C\to\C$ and $k_2:\C\to\C$ given by
\begin{align*}
 k_1(\lambda) & =   A_1+\left(c_1e^{-\mu\,a_1/v_1(\xi_3)}g'(\xi_3)-A_1\!\right)e^{-\lambda\,a_1/v_1(\xi_3)}+\mu A_1\! \int^0_{-a_1/v_1(\xi_3)}\hspace*{-1em}e^{\lambda s}ds,\\
  k_2(\lambda) & =   A_2+\left(c_1e^{-\mu\,a_2/v_2(\xi_1)}-A_2\right)e^{-\lambda\,a_2/v_2(\xi_1)}+\mu A_2\! \int^0_{-a_2/v_2(\xi_1)}\hspace*{-1em}e^{\lambda s}ds,
\end{align*}
respectively, we obtain
\begin{equation} \label{eq:MlambdaHO}
M(\lambda)=\left(
\begin{matrix}
-b_1 & 0 & k_1(\lambda)\\
k_2(\lambda) & -b_2
& 0 \\
0 & c_3 & -b_3
\end{matrix}
\right)
\end{equation}
for $\lambda\in\C$, which yields
$$
\Delta(\lambda)=\det(M(\lambda)-\lambda I_3) = (b_1+\lambda)(b_2+\lambda)(b_3+\lambda)+k_1(\lambda)k_2(\lambda)c_3.
$$
So the characteristic equation associated with the stationary point
$\phi$ with value $\xi$ is
\begin{equation} \label{eq:char-eq-at-phi}
\Delta(\lambda)=(b_1+\lambda)(b_2+\lambda)(b_3+\lambda)+k_1(\lambda)k_2(\lambda)c_3=0.
\end{equation}

Next we rewrite \eqref{eq:char-eq-at-phi} in the original notation of \eqref{eq:mrna-delay-var}--\eqref{eq:delay by stateI}.
In terms of the parameters and functions in the system \eqref{eq:mrna-delay-var}--\eqref{eq:delay by stateI}
the parameters and functions in the system \eqref{eq:explicit-form-1}-\eqref{eq:explicit-form-3} are $\mu$ and
$$
\begin{array}{l}
\xi_1  =  M^*\\
\xi_2  =  I^*\\
\xi_3  =  E^*
\end{array} \qquad
\begin{array}{l}
c_1  =  \beta_M\\
c_2  =  \beta_I\\
c_3  =  \beta_E
\end{array} \qquad
\begin{array}{l}
b_1  = \bgamma_M\\
b_2  = \bgamma_I\\
b_3    \bgamma_E
\end{array} \qquad
\begin{array}{l}
\phantom{_1}g  =  f\\
v_1  =  v_M\\
v_2  =  v_I\\
a_1  =  a_M\\
a_2  =  a_I.
\end{array}
$$

The equilibria of the systems \eqref{eq:mrna-delay-var}--\eqref{eq:delay by stateI}
and \eqref{eq:explicit-form-1}-\eqref{eq:explicit-form-3}
are given by the solutions
$(M^*,I^*,E^*)=(\xi_1,\xi_2,\xi_3)\in\R^3$ of the systems
$$\begin{array}{l}
0  =  \beta_Me^{-\mu a_M/v_M(I^*)}_{\!}f(E^*)-\bgamma_{M\!} M^*
\nonumber\\
0  =  \beta_Ie^{-\mu a_I/v_I(M^*)}M^*-\bgamma_I I^*\\
0  =  \beta_E I^*-\bgamma_E E^*.
\end{array} \Longleftrightarrow \quad\!
\begin{array}{l}
0  =  c_1e^{-\mu a_1/v_1(\xi_3)}_{\!}g(\xi_3)-b_1\xi_1\\
0  =  c_2e^{-\mu a_2/v_2(\xi_1)}\xi_1-b_2\xi_2\\
0  =  c_3\xi_2-b_3\xi_3.
\end{array}$$

From \eqref{eq:A1} and \eqref{eq:A2} we can rewrite $A_1$ and $A_2$ in the notation of
\eqref{eq:mrna-delay-var}--\eqref{eq:delay by stateI} as
\begin{align*}
A_1 = \beta_M\frac{v_M'(E^*)}{v_M(E^*)}e^{-\mu a_M/v_M(E^*)}f(E^*), \quad
A_2  = \beta_I\frac{v_I'(M^*)}{v_I(M^*)}e^{-\mu a_I/v_I(M^*)}M^*.
\end{align*}
and hence
\begin{align*}
k_1(\lambda) & =   A_1+\left(\beta_Me^{-\mu a_{M}^{\!}/^{\!}v_{M}(E^*)}f'^{\!}(E^*)-A_1\!\right)
\!e^{-\lambda a_{\!M}^{\!}/^{\!}v_M(E^*)}
+\mu A_1\hspace*{-1.5em}\int\limits^0_{-a_M^{\!}/^{\!}v_M(E^*)}\hspace*{-1.5em}e^{\lambda s}ds,\\
k_2(\lambda) & =  A_2+\left(\beta_Ie^{-\mu a_I/v_I(M^*)}-A_2\right)e^{-\lambda a_I/v_I(M^*)}+\mu A_2 \!\int^0_{-a_I/v_I(M^*)}\hspace*{-1em}e^{\lambda s}ds.
\end{align*}
For $\lambda\neq0$  we evaluate the integrals and obtain
\begin{align*}
k_1(\lambda) & = A_1+\left(\beta_Me^{-\mu a_{M\!}/^{\!}v_{M\!}(E^*)}f'^{\!}(E^*)-A_1\right)e^{-\lambda a_{M\!}/^{\!}v_{M\!}(E^*)}
+\frac{\mu A_1}{\lambda} \left( 1-  e^{-\lambda a_{M\!}/^{\!}v_{M\!}(E^*)} \!\right) \\
& = A_1 ( 1-e^{-\lambda a_M^{\!}/^{\!}v_M(E^*)})\Bigl(1+\frac\mu\lambda\Bigr)
+ \beta_Me^{-\mu a_M^{\!}/^{\!}v_M(E^*)}f'(E^*)e^{-\lambda a_M^{\!}/^{\!}v_M(E^*)}\\
& = \beta_M e^{-\mu\tau_M(E^*)}\!\left(
\frac{v_M'(E^*)}{v_M(E^*)}f(E^*)( 1-e^{-\lambda\tau_M(E^*)})\Bigl(1+\frac\mu\lambda\Bigr)
+ f'^{\!}(E^*)e^{-\lambda\tau_M(E^*)}\! \right)
\end{align*}
and
\begin{align*}
k_2(\lambda) & =  A_2+\left(\beta_Ie^{-\mu a_I^{\!}/^{\!}v_I(M^*)}-A_2\!\right)e^{-\lambda a_I^{\!}/^{\!}v_I(M^*)}+\frac\mu\lambda A_2
(1- e^{-\lambda a_I/v_I(M^*)}) \\
&=A_2(1- e^{-\lambda a_I/v_I(M^*)})\Bigl(1+\frac\mu\lambda\Bigr)
+ \beta_Ie^{-\mu a_I/v_I(M^*)}e^{-\lambda a_I/v_I(M^*)}\\
&=\beta_Ie^{-\mu\tau_{I\!}(M^*)}\!\left(
\frac{v_I'(M^*)}{v_I(M^*)}M^*(1- e^{-\lambda\tau_I(M^*)})\Bigl(1+\frac\mu\lambda\Bigr)
+e^{-\lambda\tau_I(M^*)}
\right).
\end{align*}
In terms of the parameters and functions of the system \eqref{eq:mrna-delay-var}--\eqref{eq:delay by stateI}
the characteristic equation \eqref{eq:char-eq-at-phi} becomes $\Delta(\lambda)=0$ for
\begin{align} \label{eq:char-eq-at-MIE}
\Delta(\lambda) & =  (\bgamma_M+\lambda)(\bgamma_I+\lambda)(\bgamma_E+\lambda)
+\beta_E k_1(\lambda)k_2(\lambda),  \end{align}   
which is identical to \eqref{eq:Delta-lambda-MIE},\eqref{eq:klambda1}.
